\documentclass[12pt,reqno]{amsart}
\usepackage{amssymb}
\usepackage{amsfonts}
\usepackage{hyperref}
\usepackage{graphicx}
\usepackage{color}
\usepackage{amsmath}
\usepackage[font=scriptsize,labelfont=bf]{caption}
\usepackage{wrapfig}

\graphicspath{{./CurtisPictureFile/}}
\definecolor{myurlcolor}{rgb}{0,0.6,0.2}
\hypersetup{colorlinks, linkcolor=myurlcolor, citecolor=myurlcolor, urlcolor=myurlcolor}
\def\equationautorefname~#1\null{(#1)\null}

\theoremstyle{plain}
\newtheorem*{acknowledgment}{Acknowledgment}

\newtheorem{algorithm}{Algorithm}[section]
\newtheorem{thm}{Thm}
\newtheorem{covlem}{Covering Lemma}
\newtheorem*{perelmancovlem}{Perelman Covering Lemma}

\newtheorem{corollary}[algorithm]{Corollary}

\newtheorem{definition}[algorithm]{Definition}
\newtheorem*{Remark-delta}{Remark on all things $\protect\delta$}

\newtheorem{lemma}[algorithm]{Lemma}
\newtheorem{layPerelCovlemma}[algorithm]{Parameterized Perelman Covering Lemma}

\newtheorem{limitlemmalet}[thm]{Limit Lemma}

\newtheorem{theorem} [algorithm] {Theorem}
\newtheorem{theoremlet}[thm]{Theorem}

\newtheorem*{theoremnonum}{Theorem}

\newtheorem{definitionlet}[thm]{Definition}
\newtheorem*{remarknonum}{Remark}

\newtheorem{Induction Statement}[algorithm]{Induction Statement}
\newtheorem*{definitionnonum}{Definition}

\newtheorem{proposition}[algorithm]{Proposition}

\numberwithin{equation}{algorithm}

\begin{document}
\title{Stability, Finiteness and Dimension Four }
\author{Curtis Pro}
\address{{Department of Mathematics, California State University, Stanislaus}%
}
\email{cpro@csustan.edu}
\author{Frederick Wilhelm}
\thanks{This work was supported by a grant from the Simons Foundation
(\#358068, Frederick Wilhelm)}
\address{{Department of Mathematics, University of California, Riverside}}
\email{fred@math.ucr.edu}
\urladdr{https://sites.google.com/site/frederickhwilhelmjr/home}

\begin{abstract}
We prove that for any $k\in \mathbb{R},$ $v>0,$ and $D>0$ there are only
finitely many diffeomorphism types of closed Riemannian $4$--manifolds with
sectional curvature $\geq k,$ volume $\geq v,$ and diameter $\leq D.$
\end{abstract}

\maketitle

\pdfbookmark[1]{Introduction}{Introduction}Let $\mathcal{M}%
_{k,v,d}^{K,V,D}\left( n\right) $ denote the class of closed Riemannian $n$%
--manifolds $M$ with 
\begin{equation*}
\begin{array}{cclccc}
k & \leq & \sec \,M & \leq & K, &  \\ 
v & \leq & \mathrm{vol}\,M & \leq & V, & \mathrm{and} \\ 
d & \leq & \mathrm{diam}\,M & \leq & D, & 
\end{array}%
\end{equation*}%
where $\sec $ is sectional curvature, $\mathrm{vol}$ is volume, and $\mathrm{%
diam}$ is diameter.

Cheeger's finiteness theorem says that $\mathcal{M}_{k,v,0}^{K,\infty
,D}\left( n\right) $ contains only finitely many diffeomorphism types (\cite%
{Cheeg1}, \cite{Cheeg2}, \cite{KirSie}, \cite{Pete}). Removing the
hypothesis of an upper curvature bound, Grove, Petersen, and Wu showed that
for $n\neq 3$, $\mathcal{M}_{k,v,0}^{\infty ,\infty ,D}\left( n\right) $ has
only finitely many homeomorphism types (\cite{GrovPetWu}). Perelman then
proved his stability theorem (\cite{Perel1},\cite{Kap2}) which, together
with Gromov's precompactness theorem, extends the Grove, Petersen, Wu result
to give finiteness of homeomorphism types in $\mathcal{M}_{k,v,0}^{\infty
,\infty ,D}\left( n\right) $ in all dimensions. Using \cite{KirSie}, the
Grove-Petersen-Wu/Perelman finiteness theorem says for $n\neq 4$, $\mathcal{M%
}_{k,v,0}^{\infty ,\infty ,D}\left( n\right) $ contains only finitely many
diffeomorphism types.

In contrast to other dimensions, there are compact topological $4$%
--manifolds that admit infinitely many distinct smooth structures (\cite%
{FriMorg}, cf. also \cite{FintStern}, \cite{FintStr2}). So it is natural
that Problem 1.1 in \cite{AIM} asks whether $\mathcal{M}_{k,v,0}^{\infty
,\infty ,D}\left( 4\right) $ also contains only finitely many diffeomorphism
types. We resolve this here and thus get the following generalization of
Cheeger's finiteness theorem.

\begin{theoremlet}
\label{finiteness} For any $k\in \mathbb{R},$ $v>0,$ $D>0,$ and $n\in 
\mathbb{N}$ the class $\mathcal{M}_{k,v,0}^{\infty ,\infty ,D}\left(
n\right) $ contains only finitely many diffeomorphism types.
\end{theoremlet}

Theorem \ref{finiteness} is optimal in the sense that the conclusion is
false if any one of the hypotheses is removed (\cite{GrovPet1}). Perelman
showed that the conclusion is also false if the lower bound on sectional
curvature is weakened to a lower bound on Ricci curvature (\cite{Perel3}).
On the other hand, Cheeger and Naber proved that the conclusion of Theorem %
\ref{finiteness} holds in dimension 4, if the lower bound on sectional
curvature is replaced by a two sided bound on Ricci curvature (\cite%
{CheegNab}). To do this, they prove the limits of manifolds in their class
can have only isolated singularities, in the appropriate sense. In contrast,
the limit spaces that correspond to Theorem \ref{finiteness} can have
bifurcating singular sets of all dimensions $\leq n-2$ (see e.g. Figure \ref%
{bifurc}).

A crucial distinction between $\mathcal{M}_{k,v,0}^{K,\infty ,D}\left(
n\right) $ and $\mathcal{M}_{k,v,0}^{\infty ,\infty ,D}\left( n\right) $ is
that there is a uniform lower bound for the injectivity radii of all $M\in 
\mathcal{M}_{k,v,0}^{K,\infty ,D}\left( n\right) $, whereas the spaces of $%
\mathcal{M}_{k,v,0}^{\infty ,\infty ,D}\left( n\right) $ do not even have a
uniform lower bound for their contractibility radii (\cite{Cheeg1}, \cite%
{Cheeg2}, \cite{GrovPet1}). This leads to a further contrast: if $X$ is in
the Gromov--Hausdorff closure of $\mathcal{M}_{k,v,0}^{K,\infty ,D}\left(
n\right) $ then all points of $X$ have neighborhoods that are almost
isometric to small euclidean balls. Simple examples show that this is false
for $\mathcal{M}_{k,v,0}^{\infty ,\infty ,D}\left( n\right) .$

In the 1980's, Gromov's Compactness Theorem provided a new perspective on
Cheeger's Finiteness Theorem. In particular, we learned that if $M_{i}\in 
\mathcal{M}_{k,v,0}^{K,\infty ,D}\left( n\right) $ converge to $X$, then all
but finitely many of the $M_{i}$ must be diffeomorphic to $X.$ (See, for
example, \cite{Cheeg3}). As Gromov also showed that $\mathcal{M}%
_{k,0,0}^{\infty ,\infty ,D}\left( n\right) $ is precompact, Cheeger's
finiteness theorem becomes a corollary. Combined with Perelman's Stability
Theorem, this leads naturally to the

\bigskip

\noindent \textbf{Diffeomorphism Stability Question}. \emph{(Problem 32 in 
\cite{Petersen1}) Given }$k\in \mathbb{R},$\emph{\ }$v,D>0,$\emph{\ and }$%
n\in N,$\emph{\ let }$\left\{ M_{\alpha }\right\} _{\alpha =1}^{\infty
}\subset M_{k,v,0}^{\infty ,\infty ,D}\left( n\right) $\emph{\ be a
Gromov--Hausdorff convergent sequence. Are all but finitely many of the }$%
M_{\alpha }$\emph{s diffeomorphic to each other? }\bigskip

An affirmative answer to this question would have far reaching consequences.
For example, Grove and the second author showed

\begin{theoremnonum}
(\cite{GrovWilh2}) If the answer to the Diffeomorphism Stability Question is
\textquotedblleft yes\textquotedblright , then every Riemannian $n$%
--manifold $M$ with $\sec \,M\geq 1$ and $\mathrm{diam}\,M>\frac{\pi }{2}$
is diffeomorphic to $\mathbb{S}^{n}.$
\end{theoremnonum}

Theorem \ref{finiteness} follows by combining Gromov's Precompactness
Theorem with the following result.

\begin{theoremlet}
\label{dif stab- dim 4}Diffeomorphism Stability holds in dimension $4$.
\end{theoremlet}

The remainder of the paper is a proof of Theorem \ref{dif stab- dim 4}. To
start, we fix a sequence%
\begin{equation*}
\left\{ M_{\alpha }\right\} _{\alpha =1}^{\infty }\subset \mathcal{M}%
_{k,v,0}^{\infty ,\infty ,D}\left( 4\right) \text{ that converges to }X.
\end{equation*}%
Our goal is to prove that all but finitely many of the $M_{\alpha }$ are
diffeomorphic to each other. By work of Kuwae, Machigashira, and Shioya, the
desired diffeomorphism is known to exist away from the singularities of $X$
(see \cite{KMS} and Theorem \ref{KMS diffeo thm}, below). We call this the
KMS\ diffeomorphism. Our strategy is to extend a version of the KMS
diffeomorphism to progressively more singular regions. To do this we show in
Limit Lemma \ref{cover of X lemma} (below) that any $M_{\alpha }$
sufficiently close to $X$ admits a very special type of handle decomposition
that is related to the singular structure of $X.$ The handles of Limit Lemma %
\ref{cover of X lemma} will allow us to use theorems of Cerf and Smale to
carry out the required extension. We start with a review of the Cerf and
Smale theorems. 
\begin{wrapfigure}[6]{r}{0.33\textwidth}\centering
\vspace{-8pt}
\includegraphics[scale=0.1]{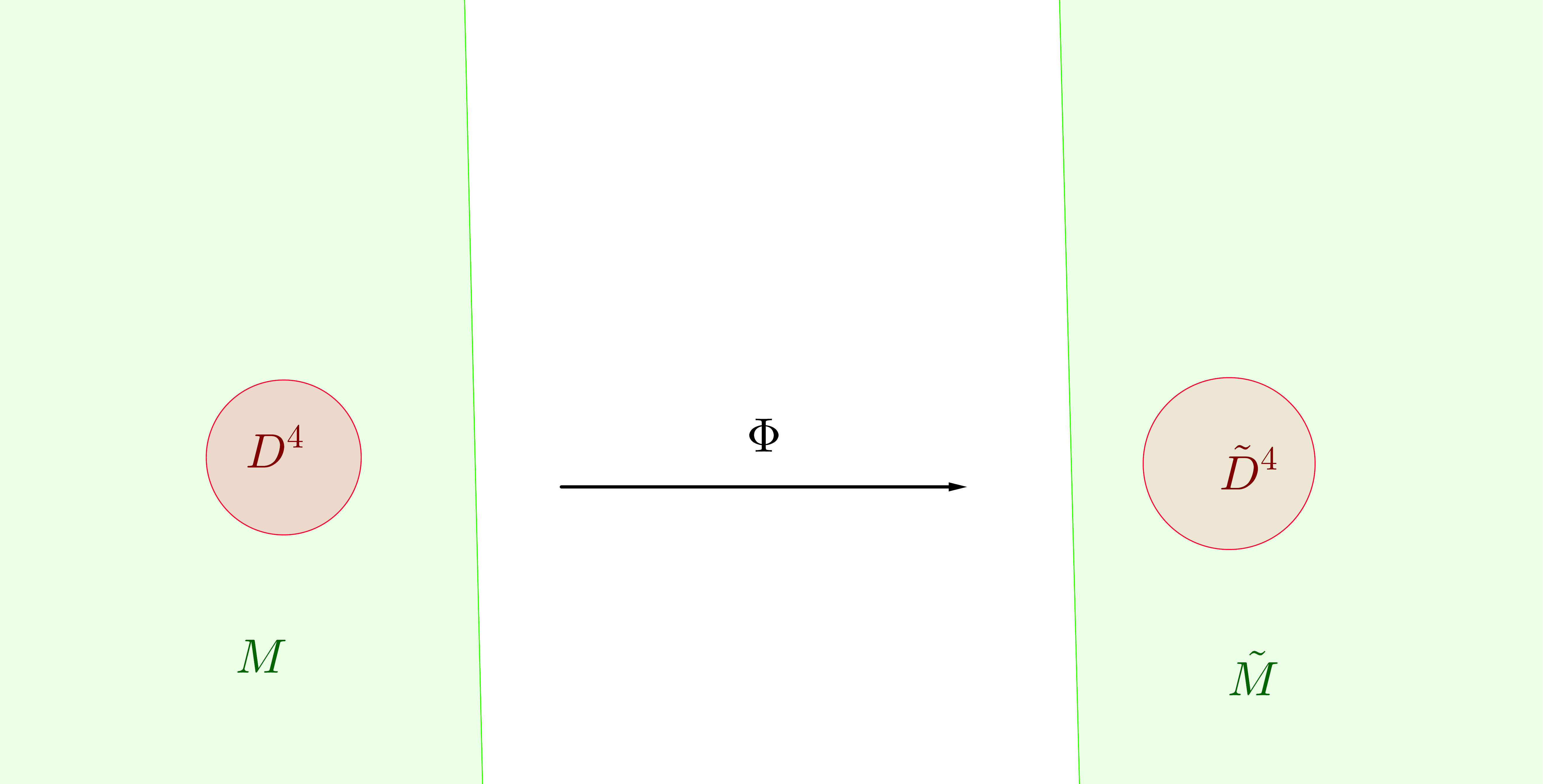}
\caption{{\protect\tiny {Cerf's Extension Theorem} }}
\label{cerf}
\end{wrapfigure}
Cerf showed that any diffeomorphism of $\mathbb{S}^{3}$ extends to $\mathbb{D%
}^{4}$ (see \cite{Cerf} cf also \cite{BalmKlin} and \cite{Hat}). The
following is a consequence

\vspace{6pt}

\noindent \textbf{Cerf's Extension Theorem. } \emph{Let }$M$\emph{\ and }$%
\tilde{M}$\emph{\ be two smooth }$4$--\emph{manifolds with smoothly embedded
closed }$4$\emph{--disks, }$D^{4}\subset M$\emph{\ and }$\tilde{D}%
^{4}\subset \tilde{M}.$\emph{\ Let }$\Phi :\overline{M\setminus D^{4}}
\longrightarrow \overline{\tilde{M}\setminus \tilde{D}^{4}} $ \emph{be a
diffeomorphism so that } $\Phi \left( \partial D^{4}\right) =\partial \tilde{%
D}^{4}.$\emph{\ Then }$\Phi $\emph{\ extends to a diffeomorphism }$%
M\longrightarrow \tilde{M}.$ \bigskip

\begin{wrapfigure}[16]{r}{0.33\textwidth}\centering
\includegraphics[scale=.19]{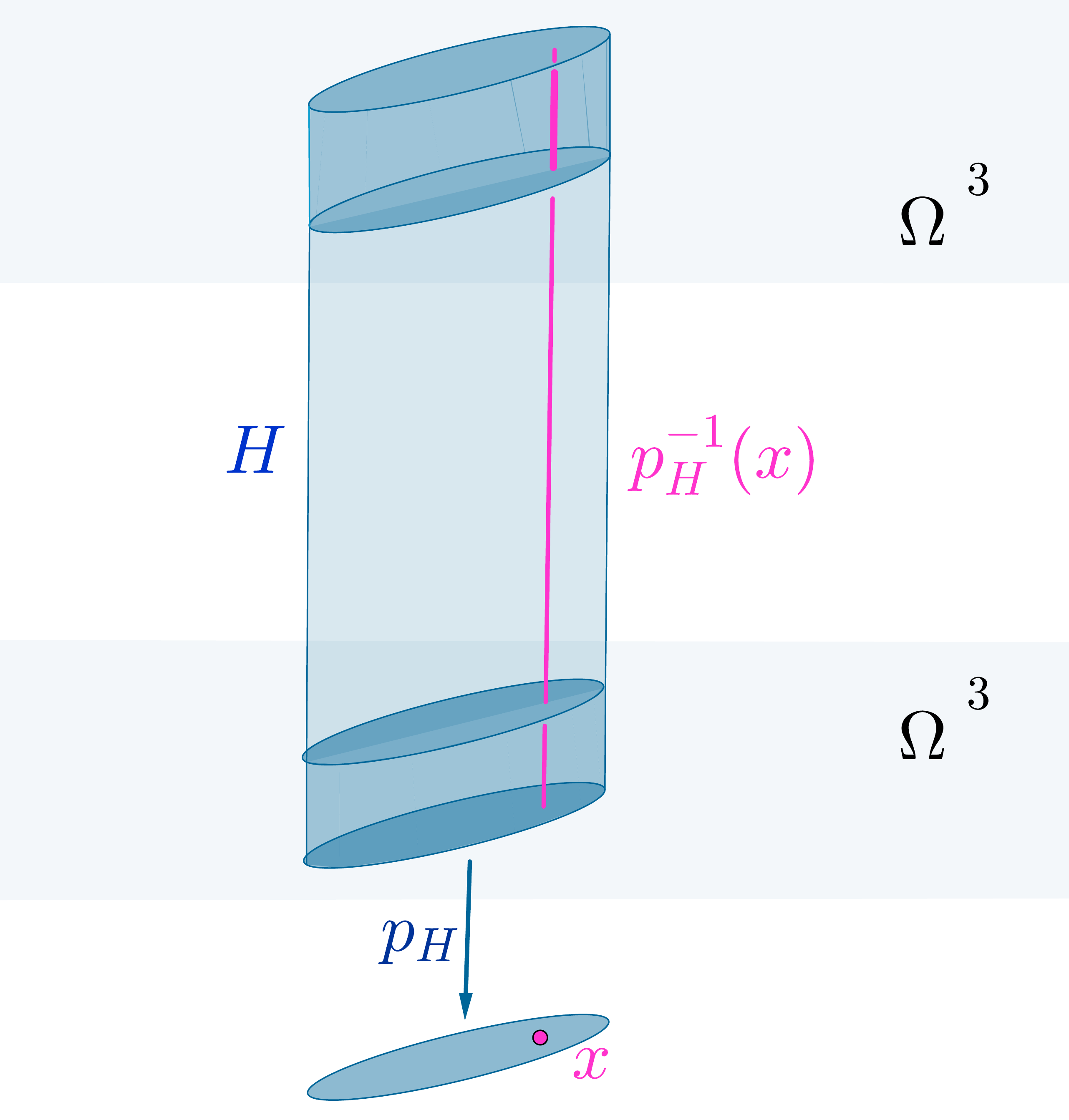}
\vspace{3pt}
\caption{{\protect\tiny {\ This is a $1$--handle $H$ attached to a $3$%
--manifold $\Omega^3$. The projection $p_H$ is topologically from ${\ 
\mathbb{R}}^3 $ to ${\ \mathbb{R}}^2 $. The dark blue depicts the
intersection of $\Omega^3$ and $H $. }}}
\label{handle}
\end{wrapfigure}
We also use an extension theorem that is essentially due to Smale. It
follows from the fact that, for $n=1$ or $2,$ the diffeomorphism group of
the $n$--sphere deformation retracts to the orthogonal group (\cite{Sm}). To
state Smale's Theorem, we recall that a smooth $n$--manifold $M$ is obtained
from a smooth $n$--manifold $\Omega $ by \textbf{attaching an }$\left(
n-j\right) $\textbf{--handle}, $H,$ provided:

\begin{enumerate}
\item $M=\Omega \cup H.$

\item $\Omega $ and $H$ are open in $M.$

\item $H\ $is the total space of a trivial vector bundle $%
p_{H}:H\longrightarrow \mathbb{B}^{j}$ with fiber $\mathbb{R}^{n-j}.$ Here $%
\mathbb{B}^{j}$ is an open ball in $\mathbb{R}^{j}.$

\item The restriction of $p_{H}$ to $H\setminus \Omega $ is a trivial bundle
over $\mathbb{B}^{j}$ whose fibers are diffeomorphic to the closed ball in $%
\mathbb{R}^{n-j}$. In particular, the restriction of $p_{H}$ to $\partial
\left( H\setminus \Omega \right) $ is a trivial sphere bundle over $\mathbb{B%
}^{j}.$
\end{enumerate}

\noindent \textbf{Smale's Extension Theorem. }\emph{Let }$M$\emph{, }$\tilde{%
M},$\emph{\ }$\Omega ,$\emph{\ and }$\tilde{\Omega}$\emph{\ be smooth }$4$--%
\emph{manifolds. For }$k\in \left\{ 1,2\right\} ,$ \emph{suppose that }$M$%
\emph{\ is obtained from }$\Omega $\emph{\ by attaching a }$\left(
4-k\right) $\emph{--handle, }$\left( H,p_{H}\right) ,$\emph{\ and that }$%
\tilde{M}$\emph{\ is obtained from }$\tilde{\Omega}$\emph{\ by attaching an }%
$\left( 4-k\right) $\emph{--handle, }$( \tilde{H},p_{\tilde{H}}) . $\emph{\ }

\emph{Let}%
\begin{equation*}
\Phi :\mathrm{closure}( \Omega ) \longrightarrow \mathrm{closure}( \tilde{%
\Omega})
\end{equation*}%
\emph{be a diffeomorphism so that}%
\begin{equation}
\left. p_{\tilde{H}}\circ \Phi \right\vert _{H\cap \Omega }=\left.
p_{H}\right\vert _{H\cap \Omega }.  \label{smale respt}
\end{equation}%
\emph{Then }$\Phi $\emph{\ extends to a diffeomorphism }$M\longrightarrow 
\tilde{M}$\emph{\ so that }$p_{\tilde{H}}\circ \Phi =p_{H}.$

\begin{figure}[h]
\includegraphics[scale=0.18]{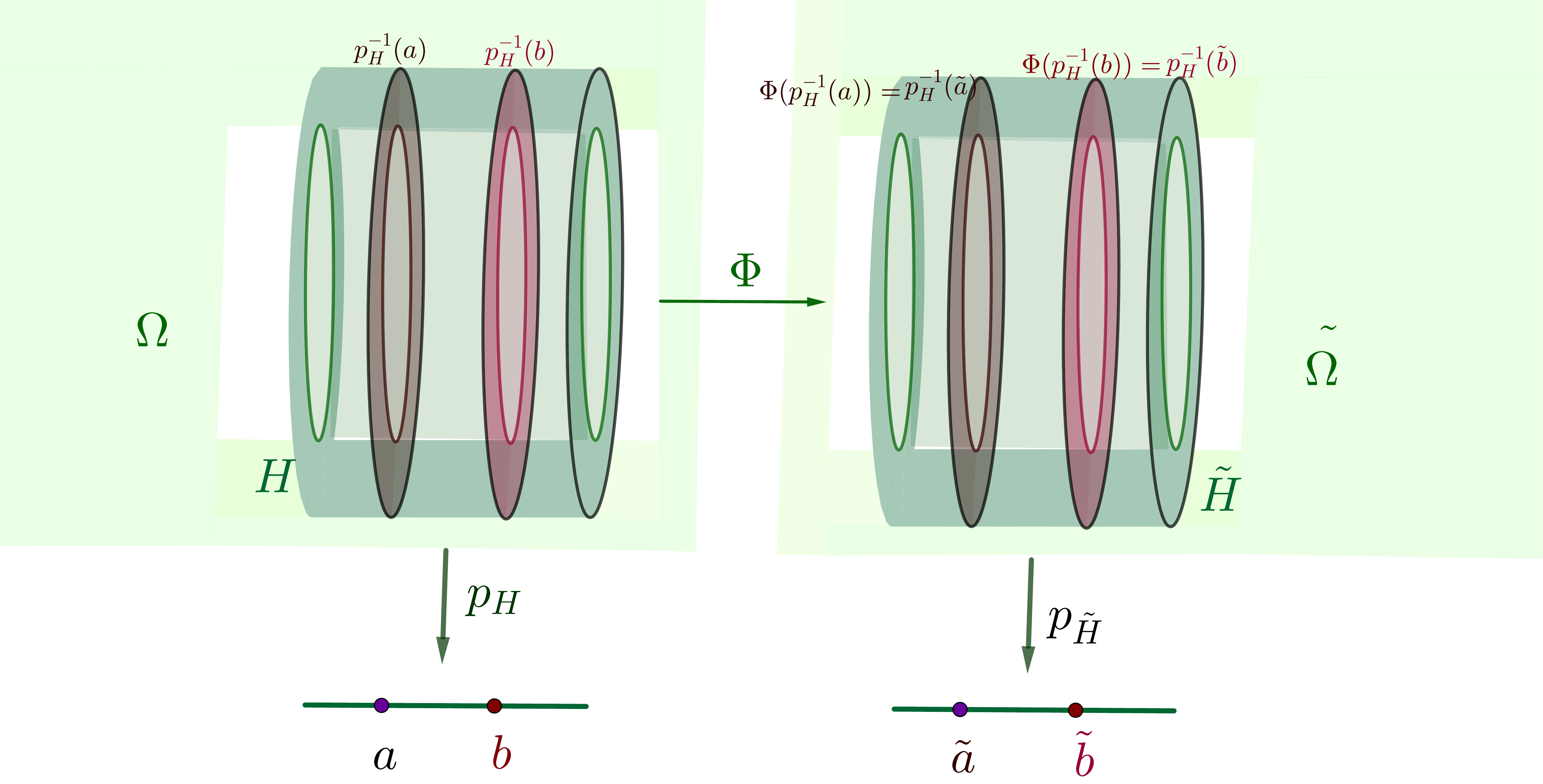}
\caption{{\protect\tiny {\ Pictured are $2$-handles $H$, $\tilde{H}$
attached to $3$--manifolds $\Omega$, $\tilde{\Omega}$. The dark green
depicts the intersections. Originally the diffeomorphism $\Phi$ takes $%
\Omega $ to $\tilde{\Omega}$; however, because it also takes fibers of $p_H$
to fibers of $p_{\tilde{H}}$, and because $\mathrm{Diff}(S^1)$ deformation
retracts to $O(2)$, $\Phi$ extends to a diffeomorphism from $\Omega \cup H
\longrightarrow \tilde{\Omega} \cup \tilde{H}$. Since $\mathrm{Diff}(S^2)$
deformation retracts to $O(3)$, the same statement applies if $H$ and $%
\tilde{H}$ are $3$--handles attached to $4 $--manifolds. }}}
\label{smale}
\end{figure}

Together with Invariance of Domain, the hypotheses of Smale's Extension
Theorem imply that $\Phi $ takes $\partial \Omega $ to $\partial \tilde{%
\Omega}.$ Each of these is a trivial sphere bundle that bounds the trivial
disk bundles $H\setminus \Omega $ and $\tilde{H}\setminus \tilde{\Omega},$
respectively. This, together with (\ref{smale respt}) and Smale's
deformation retraction of \textrm{Diff}$\left( S^{4-k-1}\right) $ to $%
O\left( 4-k\right) ,$ gives the method to extend $\Phi $ to a diffeomorphism
from $M$ to $\tilde{M}.$

The strategy to prove Theorem \ref{dif stab- dim 4}, which we outline in the
next section, is to show any two manifolds $M_{\alpha }$ and $M_{\beta }$
that are sufficiently close to $X$ admit a special type of handle
decomposition that is related to the singular structure of $X$. This handle
decomposition, together with Cerf and Smale's Extension Theorems then allows
us to extend a version of the KMS diffeomorphism to a diffeomorphism from $%
M_{\alpha }$ to $M_{\beta }$.

\section*{Outline of the Proof}

Central to the proof of Theorem \ref{dif stab- dim 4} is the approximation
of geometric structures in $X$ with corresponding structures in all but
finitely many of the $M_{\alpha }$. We call these approximations in the $%
M_{\alpha }$ \emph{stable} structures. Our objective is to show that the
singular structure of $X$ leads to a stable handle decomposition in the $%
M_{\alpha }$. We formulate this precisely in Limit Lemma \ref{cover of X
lemma}, below. First we establish necessary terminology.

We will use the $0$, $1$, and $2$-strained subsets of $X$ to construct $4$, $%
3$, and $2$-handles in the $M_{\alpha }$. So, for the remainder of the paper
we assume that the reader has a working knowledge of the concept of an $%
\left( l,\delta \right) $--strained point in an Alexandrov space. These were
originally defined in \cite{BGP}, where they were called burst points (see
also, e.g., \cite{ProWilh1}).

The starting point in the proof of Theorem \ref{dif stab- dim 4} is to
construct a diffeomorphism between open, almost-dense subsets of any two of
the $M_{\alpha }$ that are sufficiently close to the following subset of $X$.

\begin{definitionnonum}
Given $\delta >0,$ we let $X_{\delta }^{4}$ be the points of $X$ that are $%
\left( 4,\delta \right) $--strained.
\end{definitionnonum}

\begin{definitionnonum}
Let $U\subset X$ be an open subset and $p:U\rightarrow \mathbb{R}^{k}$ a
map. Given a sequence of sets $\{U_{\alpha }\}_{\alpha =1}^{\infty }$ and
maps $\{p_{\alpha }:U_{\alpha }\rightarrow \mathbb{R}^{k}\}_{\alpha
=1}^{\infty }$ where $U_{\alpha }\subset M_{\alpha }$ is open, we call $%
\{\left( U_{\alpha },p_{\alpha }\right) \}_{\alpha =1}^{\infty }$ an \textbf{%
approximation} of $\left( U,p\right) $ provided $U_{\alpha }\rightarrow U$
and $p_{\alpha }\rightarrow p$ in the Gromov-Hausdorff topology.
\end{definitionnonum}

Kuwae, Machigashira, and Shioya have already constructed diffeomorphisms
between approximations of the top stratum. More precisely

\begin{theoremlet}[Kuwae-Machigashira-Shioya \protect\cite{KMS}]
\label{KMS diffeo thm}There is a $\delta >0$ so that every open $\Omega
\subset X_{\delta }^{4}$ with $\mathrm{closure}(\Omega )\subset X_{\delta
}^{4}$ has approximations $\Omega _{\alpha }\subset M_{\alpha }$ so that for
all but finitely many $\alpha $ and $\beta $ there is a diffeomorphism 
\begin{equation*}
\Phi _{\beta ,\alpha }:\Omega _{\alpha }\rightarrow \Omega _{\beta }.
\end{equation*}
\end{theoremlet}

To prove Theorem \ref{dif stab- dim 4} we show that there is a choice of $%
\Omega \subset X_{\delta }^{4}$ so that all but finitely many of the $%
M_{\alpha }$ are obtained from the corresponding $\Omega _{\alpha }$ by
first attaching a finite number of $2$-handles, then a finite number of $3$%
-handles, and lastly a finite number of $4$--handles. We also show, using
Cerf and Smale's Extension Theorems, that at each stage of attachment, $\Phi
_{\beta ,\alpha }$ can be extended over the attached handle. The end result
is a diffeomorphism between $M_{\alpha }$ and $M_{\beta }$. To accomplish
this, we will need to establish several compatibility conditions between $%
\Phi _{\beta ,\alpha }$ and the handle decompositions of $M_{\alpha }$ and $%
M_{\beta }$. Next we give several definitions related to these compatibility
conditions which will allow us to state Limit Lemma \ref{cover of X lemma}.

\begin{definitionnonum}
\label{k framed dfn}A subset $H\subset X$ together with a map $%
p_{H}:H\rightarrow \mathbb{R}^{k}$ is called $(k,\varkappa )$\textbf{--framed%
} provided there are approximations $H_{\alpha }\subset M_{\alpha }$ of $H$
and $p_{\alpha }:H_{\alpha }\rightarrow \mathbb{R}^{k}$ of $p_{H}$ with the
following properties:

\begin{enumerate}
\item For all but finitely many $\alpha $, the maps $p_{\alpha }:H_{\alpha
}\rightarrow \mathbb{R}^{k}$ are $C^{1}$, $\varkappa $-almost Riemannian
submersions for some $\varkappa >0$ whose images are diffeomorphic to an
open ball in $\mathbb{R}^{k}$.

\item For all but finitely many $\alpha $, the pair $(H_{\alpha },p_{\alpha
})$ is a trivial vector bundle.
\end{enumerate}
\end{definitionnonum}

\begin{definitionnonum}
Let $H\subset X$ be $k$--framed. A function $r_{H}:H\longrightarrow \mathbb{R%
}$ is called a \textbf{fiber exhaustion function} for $H$ if and only if for
all but finitely many $\alpha $ there are approximations $H_{\alpha
},p_{\alpha }$ and $r_{\alpha }$ of $H,$ $p_{H},$ and $r_{H}$ so that the
restriction of $r_{\alpha }$ to each fiber of $p_{\alpha }$ has exactly one
critical point which is nondegenerate and maximal.
\end{definitionnonum}

The fiber exhaustion function $r_{H_{\alpha }}$ plays the role of the \emph{%
negative} of the radial function for each fiber of $\left( H_{\alpha
}^{k},r_{H_{\alpha }}\right) ,$ and hence has a maximum rather than a
minimum. This is an artifact of our construction, which gives $r_{H}$ as a
concave down function.

Officially, a framed set is a triple, $\left( H,p,r\right) ,$ but to
simplify the exposition we will take the liberty of writing it either as a
pair, $\left( H,p\right) ,$ or a singleton, $H$, depending on the context.
When $(H,p_{H})$ is $(k,\varkappa )$-framed and neither $p_{H}$ or $%
\varkappa $ are needed for discussion, we will just say $H$ is $k$-framed.
We use the superscript $^{k},$ as in $H^{k}$, to indicate that $H^{k}$ is $k$%
-framed. We adopt the convention that any expression involving the symbol $%
dp_{H}$ implicitly asserts that the expression holds with the symbol $dp_{H}$
replaced by the differential $dp_{\alpha }$ for all but finitely many $%
\alpha $.

\begin{definitionnonum}
Given $\delta ,\varkappa >0$, a cover $\mathcal{H=H}^{0}\cup \mathcal{H}%
^{1}\cup \mathcal{H}^{2}$ of $X\setminus X_{\delta }^{4}$ is called $%
\varkappa $-\textbf{framed} if the $\mathcal{H}^{k}$ consist of a finite
number of $(k,\varkappa )$-framed sets $(H,p_{H},r_{H})$. In the event that $%
\mathcal{H}$ does not cover $X\setminus X_{\delta }^{4},$ we will call it a $%
\varkappa $-\textbf{framed collection}.
\end{definitionnonum}

Although the notion of a framed cover depends on the parameters $\varkappa $
and $\delta $, in the interest of simplicity, we will usually not mention
this explicitly, and we adopt the convention that all statements about
framed covers are only valid for sufficiently large $\alpha $ and
sufficiently small $\varkappa $ and $\delta .$

The proof of Theorem \ref{dif stab- dim 4} is broken down into the following
steps. First, we construct a framed cover $\mathcal{H}$ of $X\setminus
X_{\delta }^{4}$, opens sets 
\begin{equation*}
\Omega _{\alpha }\subset M_{\alpha },\Omega _{\beta }\subset M_{\beta }\text{%
,}
\end{equation*}%
as in Theorem \ref{KMS diffeo thm}, and a diffeomorphism $\Phi _{\beta
,\alpha }:\Omega _{\alpha }\rightarrow \Omega _{\beta }$ so that for $\gamma
=\alpha $ or $\beta $,

\begin{itemize}
\item $M_{\gamma }=\Omega _{\gamma }\bigcup \cup \mathcal{H}_{\gamma }.$

\item All of the $2$-framed sets $H^{2}\in \mathcal{H}^{2}$ have
approximations $H_{\gamma }^{2}\in \mathcal{H}_{\gamma }^{2}$ that are
disjoint $2$-handles attached to $\Omega _{\alpha }$, and

\item $\Phi _{\beta ,\alpha }$ satisfies the hypotheses of Smale's Extension
Theorem and hence extends over each of these disjoint $2$-handles.
\end{itemize}

This produces a diffeomorphism $\Phi _{\beta ,\alpha }^{2}:\Omega _{\alpha
}^{2}\rightarrow \Omega _{\beta }^{2}$ where $\Omega _{\gamma }^{2}=\Omega
_{\gamma }\cup \left( \cup \mathcal{H}_{\gamma }^{2}\right) $ for $\gamma
=\alpha $ or $\beta .$ For the next step, we require that

\begin{itemize}
\item all of the $1$-framed sets $H^{1}\in \mathcal{H}^{1}$ have
approximations $H_{\alpha }^{1}\in \mathcal{H}_{\alpha }^{1}$ that are
disjoint $3$-handles attached to $\Omega _{\alpha }^{2}$, and

\item $\Phi _{\beta ,\alpha }^{2}$ satisfies the hypotheses of Smale's
Extension Theorem and hence extends over each of these disjoint $3$-handles.
\end{itemize}

This produces a diffeomorphism $\Phi _{\beta ,\alpha }^{3}:\Omega _{\alpha
}^{3}\rightarrow \Omega _{\beta }^{3}$ where $\Omega _{\gamma }^{3}=\Omega
_{\gamma }^{2}\cup \left( \cup \mathcal{H}_{\gamma }^{1}\right) $ for $%
\gamma =\alpha $ or $\beta .$ The last step requires that

\begin{itemize}
\item all $0$-framed sets $H^{0}\in \mathcal{H}^{0}$ have approximations $%
H_{\alpha }^{0}\in \mathcal{H}_{\alpha }^{0}$ that are disjoint $4$-handles
attached to $\Omega _{\alpha }^{3}$, and

\item $\Phi _{\beta ,\alpha }^{3}$ satisfies the hypotheses of Cerf's
Extension Theorem and hence extends over each of these disjoint $4$-handles.
\end{itemize}

The final product is the desired diffeomorphism $\Phi :M_{\alpha
}\rightarrow M_{\beta }$. This only works because our framed cover $\mathcal{%
H}$ satisfies certain other constraints. For example, the Smooth Schoenflies
Problem is open in dimension $4,$ so to perform the final extension, we will
identify $M_{\beta }\setminus \Phi _{\beta ,\alpha }^{3}\left( \Omega
_{\alpha }^{3}\right) $ with a disjoint union of $4$--disks. This is the
main role of the fiber exhaustion functions $r_{H}$ and the conditions that
we impose on them in Definitions \ref{top stratum resp dfn} and \ref{dfn
resp framed collect} below.

The next definition provides the condition on the framed sets that allows us
to identify them as handles attached in the manner that our extension
process requires.

\begin{definitionlet}
\label{dfn proper}Let $\mathcal{H=H}^{0}\cup \mathcal{H}^{1}\cup \mathcal{H}%
^{2}$ be a framed cover and let $\Omega \subset $ $X_{\delta }^{4}$ be an
open set. We say that $\mathcal{H}$ is $\Omega $--\textbf{proper} provided

\begin{enumerate}
\item For each $(H_{2},p_{2})\in \mathcal{H}^{2}$ and each $x\in H_{2},$ 
\begin{equation}
\partial \left( \mathrm{closure}\left( p_{2}^{-1}\left( p_{2}\left( x\right)
\right) \right) \right) \subset \Omega .  \label{proper 2 framed}
\end{equation}

\item For each $(H_{1},p_{1})\in \mathcal{H}^{1}$ and each $x\in H_{1},$ 
\begin{equation}
\partial \left( \mathrm{closure}\left( p_{1}^{-1}\left( p_{1}\left( x\right)
\right) \right) \right) \subset \Omega \bigcup \cup \mathcal{H}^{2}.
\label{proper 1 framed}
\end{equation}

\item For each $(H_{0},p_{0})\in \mathcal{H}^{0}$ and each $x\in H^{0},$ 
\begin{equation*}
\partial \left( \mathrm{closure}\left( p_{0}^{-1}\left( p_{0}\left( x\right)
\right) \right) \right) \subset \Omega \bigcup \cup \mathcal{H}^{2}\bigcup
\cup \mathcal{H}^{1}.
\end{equation*}
\end{enumerate}
\end{definitionlet}

To ensure that our various diffeomorphisms can be extended over handles, we
require two further conditions. The first imposes a constraint on how the
KMS--diffeomorphism $\Phi _{\beta ,\alpha }:\Omega _{\alpha }\rightarrow
\Omega _{\beta }$ interacts with our framed cover.

\begin{definitionlet}
\label{top stratum resp dfn}Let $\mathcal{H=H}^{0}\cup \mathcal{H}^{1}\cup 
\mathcal{H}^{2}$ be a framed cover, and let $\Omega \subset $ $X_{\delta
}^{4}$ be an open set with $\mathrm{closure}\left( \Omega \right) \subset
X_{\delta }^{4}.$ We say that $\mathcal{H}$ \textbf{respects }$\Omega $%
\textbf{\ }if and only if 
\begin{equation*}
X=\cup \mathcal{H}\cup \Omega ,
\end{equation*}%
and the following hold:

\begin{enumerate}
\item There are approximations $\Omega _{\alpha }$ of $\Omega $ and, for all
but finitely many $\alpha ,\beta $, diffeomorphisms 
\begin{equation*}
\Phi _{\beta ,\alpha }:\Omega _{\alpha }\longrightarrow \Omega _{\beta }.
\end{equation*}

\item For each $k\in \{0,1,2\}$ and $k$-framed set $(H^{k},p,r)\in \mathcal{H%
}^{k}$, there are approximations $(H_{\alpha }^{k},p_{\alpha },r_{\alpha })$
and $(H_{\beta }^{k},p_{\beta },r_{\beta })$ so that 
\begin{equation*}
\Phi _{\beta ,\alpha }\left( \Omega _{\alpha }\cap H_{\alpha }\right)
=\Omega _{\beta }\cap H_{\beta },
\end{equation*}%
\begin{equation*}
p_{\alpha }\circ \Phi _{\alpha ,\beta }=p_{\beta },\text{ and }r_{\alpha
}\circ \Phi _{\alpha ,\beta }=r_{\beta }.
\end{equation*}
\end{enumerate}
\end{definitionlet}

\begin{wrapfigure}[12]{r}{0.55\textwidth}\centering
\vspace{-15pt}
\includegraphics[scale=0.215]{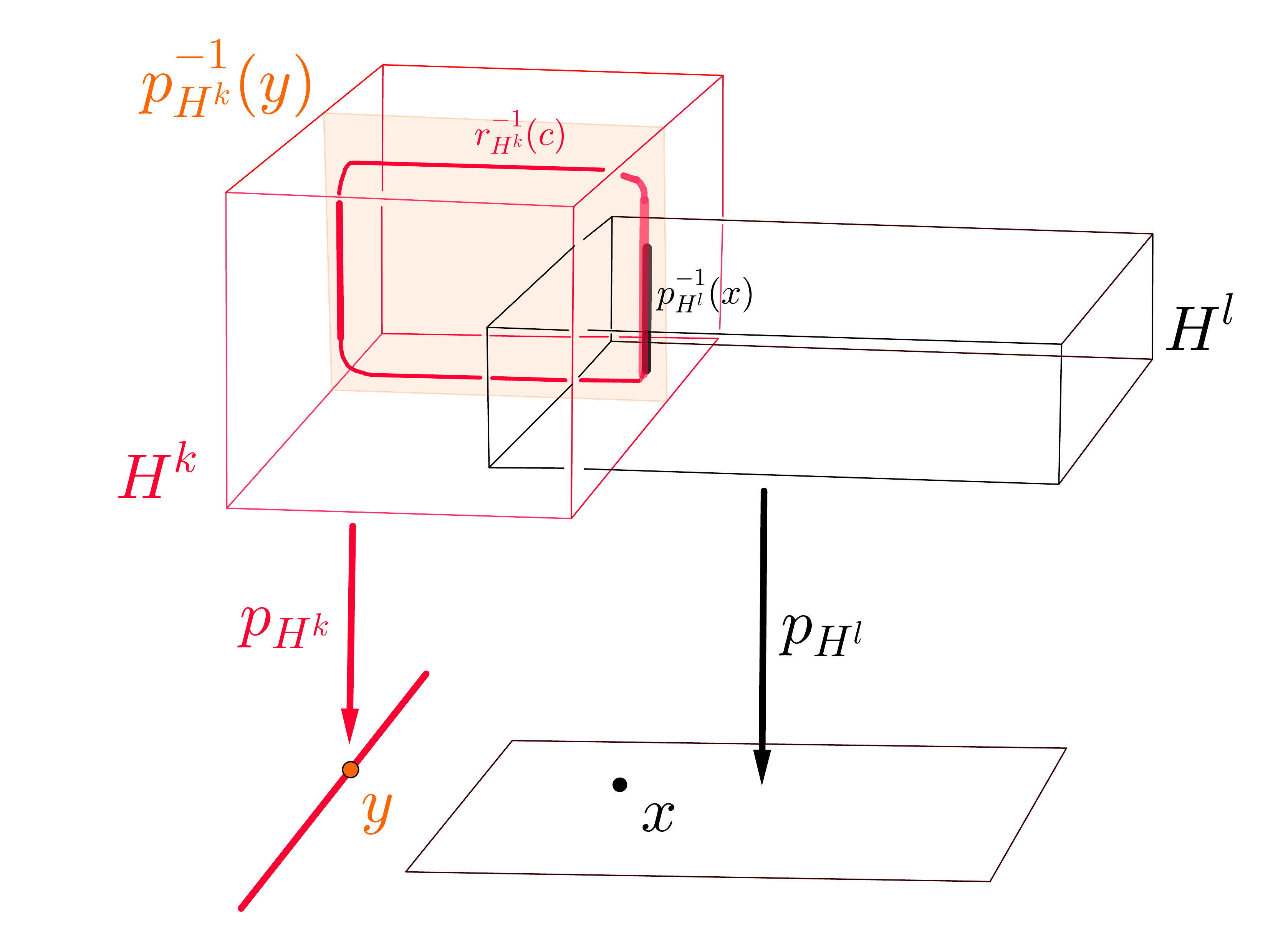}
\vspace{-11pt}
\caption{{\protect\tiny {The framed sets $H^k$ and $H^l$ intersect
respectfully, because the fibers of $p_{H^l}$ that intersect $H^k$ are
contained in both the fibers of $p_{H^k}$ and level sets of the fiber
exhaustion function $r_{H^k}$ of $H^k$. }}}
\label{respect figure}
\end{wrapfigure}

Our collection of $2$--framed sets $\mathcal{H}^{2}$ is disjoint, so if the
previous definition is satisfied, we can use Smale's Theorem to extend the
diffeomorphism $\Phi _{\beta ,\alpha }$ to the elements of $\mathcal{H}^{2}.$
To allow the extension process to continue to the $1$-- and $0$--framed
sets, the intersections of our framed sets will satisfy the following
additional constraint.

\begin{definitionlet}
\label{dfn resp framed collect}(Figure \ref{respect figure}) For $k<l,$ let $%
\left( H^{k},p_{H^{k}},r_{H^{k}}\right) \ $be $k$--framed and let $\left(
H^{l},p_{H^{l}},r_{H^{l}}\right) $ be $l$--framed. We say that $H^{k}$ and $%
H^{l}$ \textbf{intersect respectfully} if and only if there is a
neighborhood $N$ of $\partial \left( \mathrm{closure}\left( p_{H^{l}}\left(
H^{l}\right) \right) \right) $ so that for all $x\in p_{H^{l}}^{-1}\left(
N\right) \cap H^{k},$ 
\begin{eqnarray}
p_{H^{l}}^{-1}\left( p_{H^{l}}\left( x\right) \right) &\subset
&p_{H^{k}}^{-1}\left( p_{H^{k}}\left( x\right) \right) \text{ and\label%
{fiber contain}} \\
p_{H^{l}}^{-1}\left( p_{H^{l}}\left( x\right) \right) &\subset
&r_{H^{k}}^{-1}\left( r_{H^{k}}\left( x\right) \right) .
\label{schenfl cont}
\end{eqnarray}%
A framed collection is \textbf{respectful }provided all intersections of
sets with different framing dimensions are respectful.
\end{definitionlet}

In Section \ref{sect kasrten Keuker}, we combine the preceding ideas and
show that Theorem \ref{dif stab- dim 4} follows from

\begin{limitlemmalet}
\label{cover of X lemma}If $\delta $ is sufficiently small, then there is a
respectful framed cover $\mathcal{H=H}^{0}\cup \mathcal{H}^{1}\cup \mathcal{H%
}^{2},$ of $X\setminus X_{\delta }^{4}$ with the following properties:

\begin{enumerate}
\item For $k=0,$ $1,$ or $2,$ the closures of the elements of $\mathcal{H}%
^{k}$ are pairwise disjoint.

\item There is an open $\Omega \subset X_{\delta }^{4}$ so that $\mathcal{H}$
respects $\Omega $ and is $\Omega $--proper.
\end{enumerate}
\end{limitlemmalet}

Section \ref{not and conve} establishes notations and conventions. The
remainder and vast majority of the paper is divided into three parts and an
appendix, which together are a proof of Limit Lemma \ref{cover of X lemma}.
The parts are divided into subordinate sections. The parts are numbered $%
0,1, $ and $2$ because the main goal of the respective parts is to make the
collections of $0$--, $1$--$,$ and $2$--framed sets disjoint. The appendix
completes the proof of Limit Lemma \ref{cover of X lemma} by showing that $%
X_{\delta }^{4}$ respects $\mathcal{H}.$

To prove Limit Lemma \ref{cover of X lemma} we establish four other,
successively better covering lemmas. The first is essentially due to
Perelman and is proven in Section \ref{perlem framed cover sect} (cf \cite%
{GrovWilh2}, \cite{Kap1}, and \cite{ProWilh1}).

\begin{perelmancovlem}
Given $\varkappa ,\delta >0,$ there is a $\varkappa $-framed cover $\mathcal{%
P}=\mathcal{P}^{0}\cup \mathcal{P}^{1}\cup \mathcal{P}^{2}$ of $X\setminus
X_{\delta }^{4}$.
\end{perelmancovlem}

Each of Parts \ref{Part 0}, \ref{Part 1} and \ref{Part 2} begins with a
detailed summary. We refer the reader to these summaries for an outline of
the rest of the paper.

\begin{remarknonum}
The idea of using the Cerf and Smale extension theorems to extend
diffeomorphisms is also used in e.g. \cite{GrovWilh1}, \cite{GrovWilh2}, 
\cite{ProSilWilh}, and \cite{ProWilh1}.
\end{remarknonum}

\begin{acknowledgment}
Theorem \ref{finiteness} answers Problem 1.1 of the American Institute of
Mathematics Problem list for manifolds of nonnegative curvature, which was
discussed at the AIM workshop with the same name in September, 2007. The
second author is grateful to AIM for supporting his attendance at this
workshop and to the participants for the many insights that they shared. We
thank Vitali Kapovitch for showing us the example in Figure \ref{bifurc},
for discussions related to Corollary \ref{real grad pres}, and, in general,
for several extensive conversations about this problem over the years. We
are grateful to Michael Sill for multiple preliminary discussions about this
work, to Stefano Vidussi for conversations about the topological issues, and
to Karsten Grove and Catherine Searle for listening to extended outlines of
the proof and providing valuable feedback and advice.
\end{acknowledgment}

\section{Why Theorem \protect\ref{dif stab- dim 4} Follows From Limit Lemma 
\protect\ref{cover of X lemma}\label{sect kasrten Keuker}}

In this section, we prove Theorem \ref{dif stab- dim 4} assuming Limit Lemma %
\ref{cover of X lemma}. The key ingredients are the Cerf and Smale Extension
Theorems stated above.

\begin{proof}[Proof of Theorem \protect\ref{dif stab- dim 4} assuming Limit
Lemma \protect\ref{cover of X lemma}]
Let $\Omega $ be an open subset of $X_{\delta }^{4}$ that satisfies Part 2
of Limit Lemma \ref{cover of X lemma}. Then there are open sets $\Omega
_{\beta }\subset M_{B}$ and $\Omega _{\alpha }\subset M_{\beta }$ and a
diffeomorphism $\Phi _{\beta ,\alpha }^{4}:\Omega _{\alpha }\rightarrow
\Omega _{\beta }$ as in Definition \ref{top stratum resp dfn}. In
particular, for $k\in \left\{ 0,1,2\right\} $ and all $\left( H,p,r\right)
\in \mathcal{H}^{k}$, there are approximations $(H_{\alpha },p_{\alpha
},r_{\alpha })\in \mathcal{H}_{\alpha }^{k}$ and $(H_{\beta },p_{\beta
},r_{\beta })\in \mathcal{H}_{\beta }^{k}$ so that 
\begin{equation}
\Phi _{\beta ,\alpha }^{4}\left( \Omega _{\alpha }\cap H_{\alpha }\right)
=\Omega _{\beta }\cap H_{\beta },\text{ }  \label{sammmme im}
\end{equation}%
\begin{equation}
\text{ }p_{\beta }\circ \Phi _{\beta ,\alpha }=p_{\alpha },\text{ and }%
r_{\beta }\circ \Phi _{\beta ,\alpha }=r_{\alpha }.  \label{bbbase equiv}
\end{equation}

For $\gamma =\alpha $ or $\beta $, set $\Omega _{\gamma }^{2}:=\Omega
_{\gamma }\cup (\cup \mathcal{H}_{\gamma }^{2})$ and $\Omega _{\gamma
}^{1}:=\Omega _{\gamma }^{2}\cup (\cup \mathcal{H}_{\gamma }^{1}).$ It
follows that for all sufficiently large $\gamma ,$ 
\begin{equation*}
M_{\gamma }=\Omega _{\gamma }^{1}\cup (\cup \mathcal{H}_{\gamma }^{0}).
\end{equation*}%
Since the $2$-framed sets of $\mathcal{H}_{\gamma }^{2}$ are disjoint,
Equations (\ref{sammmme im}) and (\ref{bbbase equiv}), together with the
Smale Extension Theorem and the fact that $\mathcal{H}$ is $\Omega $%
--proper, allow us to successively extend $\Phi _{\beta ,\alpha }:\Omega
_{\alpha }\rightarrow \Omega _{\beta }$ to each $2$--framed set to get a
diffeomorphism 
\begin{equation*}
\Phi _{\beta ,\alpha }^{2}:\Omega _{\alpha }^{2}\rightarrow \Omega _{\beta
}^{2}
\end{equation*}%
so that for all $\left( H^{2},p\right) \in \mathcal{H}^{2},$ 
\begin{equation}
\Phi _{\beta ,\alpha }^{2}\left( H_{\alpha }^{2}\right) =H_{\beta }^{2},%
\text{\label{same image II}}
\end{equation}%
and%
\begin{equation}
\left. p_{\beta }\circ \Phi _{\beta ,\alpha }^{2}\right\vert _{H_{\alpha
}^{2}}=\left. p_{\alpha }\right\vert _{H_{\alpha }^{2}}.
\label{bassse equiv}
\end{equation}

Since $\mathcal{H}$ is respectful (see (\ref{fiber contain}),(\ref{schenfl
cont})), Equations (\ref{sammmme im}), (\ref{bbbase equiv}), (\ref{same
image II}), and (\ref{bassse equiv}) give us that for every $\left(
H,p,r\right) \in \mathcal{H}^{0}\cup \mathcal{H}^{1},$

\begin{equation}
\Phi _{\beta ,\alpha }^{2}\left( \Omega _{\alpha }^{2}\cap H_{\alpha
}\right) =\Omega _{\beta }^{2}\cap H_{\beta },  \label{sm im 2 2}
\end{equation}
\begin{equation}
\left. p_{\beta }\circ \Phi _{\beta ,\alpha }^{2}\right\vert _{H_{\alpha
}\cap \Omega _{\alpha }^{2}}=\left. p_{\alpha }\right\vert _{H_{\alpha }\cap
\Omega _{\alpha }^{2}},\text{ and}\left. r_{\beta }\circ \Phi _{\beta
,\alpha }^{2}\right\vert _{H_{\alpha }^{2}}=\left. r_{\alpha }\right\vert
_{H_{\alpha }}.  \label{resp 86}
\end{equation}

Equations (\ref{sm im 2 2}) and (\ref{resp 86}), Smale's Extension Theorem,
and the fact that $\mathcal{H}$ is $\Omega $--proper allow us to
successively extend $\Phi _{\beta ,\alpha }^{2}:\Omega _{\alpha
}^{2}\rightarrow \Omega _{\beta }^{2}$ to each $1$--framed set to get a
diffeomorphism 
\begin{equation*}
\Phi _{\beta ,\alpha }^{1}:\Omega _{\alpha }^{1}\longrightarrow \Omega
_{\beta }^{1}
\end{equation*}%
so that for all $\left( H,p,r\right) \in \mathcal{H}^{1},$%
\begin{equation}
\Phi _{\beta ,\alpha }^{1}\left( H_{\alpha }\right) =H_{\beta },\text{ and}
\label{sammm
immm}
\end{equation}%
\begin{equation}
\left. p_{\beta }\circ \Phi _{\beta ,\alpha }^{1}\right\vert _{H_{\alpha
}\cap \Omega _{\alpha }^{3}}=\left. p_{\alpha }\right\vert _{H_{\alpha }\cap
\Omega _{\alpha }^{3}}.  \label{resp 97}
\end{equation}

Equations (\ref{sammmme im}) through (\ref{resp 97}), together with our
assumption that $\mathcal{H}$ is respectful, give us that for all $\left(
H^{0},r\right) \in \mathcal{H}^{0},$ 
\begin{equation}
\left. r_{\beta }\circ \Phi _{\beta ,\alpha }^{1}\right\vert _{H_{\alpha
}^{0}}=\left. r_{\alpha }\right\vert _{H_{\alpha }^{0}}.  \label{cerfssss}
\end{equation}

The final step is to apply the Cerf Extension Theorem to get the desired
diffeomorphism, $\Phi :M_{\alpha }\longrightarrow M_{\beta }.$
\end{proof}

\section{Notations and Conventions\label{not and conve}}

For the remainder of the paper, we fix a sequence $\left\{ M_{\alpha
}\right\} _{\alpha =1}^{\infty }\subset \mathcal{M}_{k,v,0}^{\infty ,\infty
,D}\left( 4\right) $ that Gromov-Hausdorff converges to $X$, that is, 
\begin{equation*}
\lim_{\alpha \rightarrow \infty }M_{\alpha }=X.
\end{equation*}%
Our goal is to prove that all but finitely many of the $M_{\alpha }$s are
diffeomorphic to each other. We have seen that to do this, it suffices to
prove Limit Lemma \ref{cover of X lemma}. We assume throughout that all
metric spaces are complete, and the reader has a basic familiarity with
Alexandrov spaces, including, but not limited to \cite{BGP}. We adopt all of
the Notations and Conventions from \cite{ProWilh1}, and in addition,

\begin{enumerate}
\item We write $\mathbb{S}^{n}$ for the unit sphere in $\mathbb{R}^{n+1},$ $%
\mathbb{B}^{n}$ for an open unit $n$--ball, and $\mathbb{B}^{n}\left(
r\right) $ for an open $r$--ball in $\mathbb{R}^{n}.$

\item For a metric space $Y$, we let $\lambda Y$ be the metric space
obtained from $Y$ by multiplying all distances by $\lambda .$

\item Given a constant $C\in \mathbb{R},$ we use the symbols $C^{-}$ and $%
C^{+}$ to denote constants that are slightly smaller or slightly larger than 
$C.$ We adopt the convention that any statement involving $C^{-}$ or $C^{+}$
implicitly asserts that there is a choice of $C^{-}\in \left( \frac{C}{2}%
,C\right) $ or a choice of $C^{+}$ in $\left( C,2C\right) $ for which the
statement is true.

\item Specific constants appear throughout the paper. While we have strived
to make the dependencies among these constants clear, they are not chosen to
be optimal, with the exception of $\frac{\pi }{6}$\ in Theorem \ref%
{Gromov-Del Prop}. Rather we have endeavored to choose them in a way that
maximizes the readability of the paper. Even with optimal constants, our
methods do not lead to a practical estimate for the number of diffeomorphism
classes for most cases of Theorem \ref{finiteness}.

\item If $\Omega \left( r\right) $ is a parameterized family of sets with $%
\Omega =\Omega \left( r_{0}\right) ,$ we write $\Omega ^{+}$ or $\Omega ^{-}$
for $\Omega \left( r_{0}^{+}\right) $ or $\Omega \left( r_{0}^{-}\right) $
respectively. We adopt the convention that any statement involving $\Omega
^{+}$ or $\Omega ^{-}$ implicitly asserts that there is a choice of $\Omega
^{+}$ or a choice of $\Omega ^{-}$ for which the statement is true.\label%
{plus minus conve}

\item For $\lambda \geq 1,$ we say that an embedding $\iota $ is $\lambda $%
--bilipshitz if both $\iota $ and $\iota ^{-1}$ are $\lambda $--Lipshitz.

\item Let $\mathcal{C}$ and $\mathcal{D}$ be collections of subsets of $X.$
We say that the Hausdorff distance between $\mathcal{C}$ and $\mathcal{D}$
is $<\varepsilon $ if and only if there is a bijection $\Phi :\mathcal{C}%
\longrightarrow \mathcal{D}$ so that $\mathrm{dist}_{\mathrm{Haus}}\left(
C,\Phi \left( C\right) \right) <\varepsilon $ for all $C\in \mathcal{C}.$

\item For a simplicial complex $\mathcal{T}$, we adopt the convention that $%
\mathcal{T}_{i}$ is the collection of $i$--simplices of $\mathcal{T}$ and
that $\left\vert \mathcal{T}\right\vert $ is the underlying polyhedron.

\item For $A\subset X,$ we write $\bar{A}$ for the closure of $A$ in $X.$
\end{enumerate}

\begin{definition}
\label{stable subm dfn}A subset $U\subset X$ together with a map $%
p_{U}:U\rightarrow \mathbb{R}^{k}$ is called a stable $\varkappa $-almost
Riemannian submersion provided:

\begin{enumerate}
\item There are approximations $U_{\alpha }\subset M_{\alpha }$ of $U$ and $%
p_{\alpha }:U_{\alpha }\rightarrow \mathbb{R}^{k}$ of $p_{U}$.

\item For all but finitely many $\alpha $ and some $\varkappa >0,$ the maps $%
p_{\alpha }:U_{\alpha }\rightarrow \mathbb{R}^{k}$ are $C^{1}$, $\varkappa $%
-almost Riemannian submersions whose images are diffeomorphic to an open
ball in $\mathbb{R}^{k}$.
\end{enumerate}
\end{definition}

\begin{remarknonum}
When there is no possibility of confusion, we will omit the word \emph{stable%
} and call such a map a Riemannian submersion.
\end{remarknonum}

\begin{definition}
\label{C1 close S-subms dfn}Given $\varkappa >0$ and two maps $\pi
,p:U\longrightarrow \mathbb{R}^{k},$ we say that 
\begin{equation*}
\left\vert \pi -p\right\vert _{C_{\mathrm{stab}}^{1}}<\varkappa
\end{equation*}%
if and only if for all but finitely many $\alpha ,$ there are approximations 
$\left\{ \pi _{\alpha }\right\} $ and $\left\{ p_{\alpha }\right\} $ of $\pi 
$ and $p$ so that 
\begin{equation*}
\left\vert \pi _{\alpha }-p_{\alpha }\right\vert _{C^{1}}<\varkappa .
\end{equation*}%
Similarly, we say that 
\begin{equation*}
\left\vert d\pi -dp\right\vert _{\mathrm{stab}}<\varkappa
\end{equation*}%
if and only if for all but finitely many $\alpha ,$ there are approximations 
$\left\{ \pi _{\alpha }\right\} $ and $\left\{ p_{\alpha }\right\} $ of $\pi 
$ and $p$ so that 
\begin{equation*}
\left\vert d\pi _{\alpha }-dp_{\alpha }\right\vert <\varkappa .
\end{equation*}
\end{definition}

\setcounter{part}{-1}

\part{The Perelman-Plaut Covering Lemma\label{Part 0}}

In the first section of this part, we prove the Perelman Covering Lemma
(Section \ref{perlem framed cover sect} of the paper). In the rest of the
paper, we re-do the Perelman cover multiple times, obtaining at each stage
covers that are progressively closer to satisfying the conclusions of Limit
Lemma \ref{cover of X lemma}. The end result is a respectful framed cover so
that for $k=0,1,$ or $2,$ the collection of $k$--framed sets is disjoint.

Before attempting to make the $1$-- and $2$-- framed sets disjoint, we
constrain the way that they intersect in Section \ref{perl plt sect}. The
argument uses a sphere theorem of Plaut's from \cite{Pla},%
\begin{wrapfigure}[11]{r}{0.4\textwidth}\centering
\vspace{-12pt}
\includegraphics[scale=0.3]{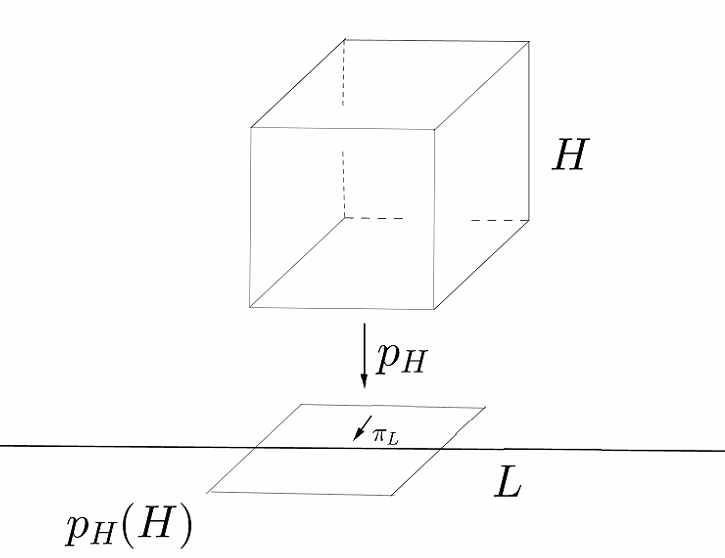}
\vspace{-8.5pt}
\caption{{\protect\tiny {\ $(H,p_{H})$ is 2-framed. $(H,\protect\pi %
_{L}\circ p_{H})$ is artificially 1-framed.}}}
\label{art fr basic}
\end{wrapfigure} so we call the resulting cover the Perelman-Plaut Framed
Cover (see Lemma \ref{Ms. Pacman limit lemma} and Figures \ref{lined up}, %
\ref{lined up, but not swallowing}, \ref{notlinedup}, and \ref{vert gap fig}%
).

The focus of Parts \ref{Part 1} and \ref{Part 2} of the paper is to replace
the collections of $1$-- and $2$-- framed sets of the Perelman-Plaut Framed
Cover by collections that are pairwise disjoint. To separate the $1$--framed
sets, we add additional $0$--framed sets, and to separate the $2$--framed
sets, we add additional $0$-- and $1$-- framed sets. The additional $0$--
and $1$-- framed sets are subsets of existing $1$-- and $2$-- framed sets.
The main idea is that \newpage \noindent 
\begin{wrapfigure}[9]{r}{0.33\textwidth}\centering
\vspace{-16pt}
\includegraphics[scale=0.13]{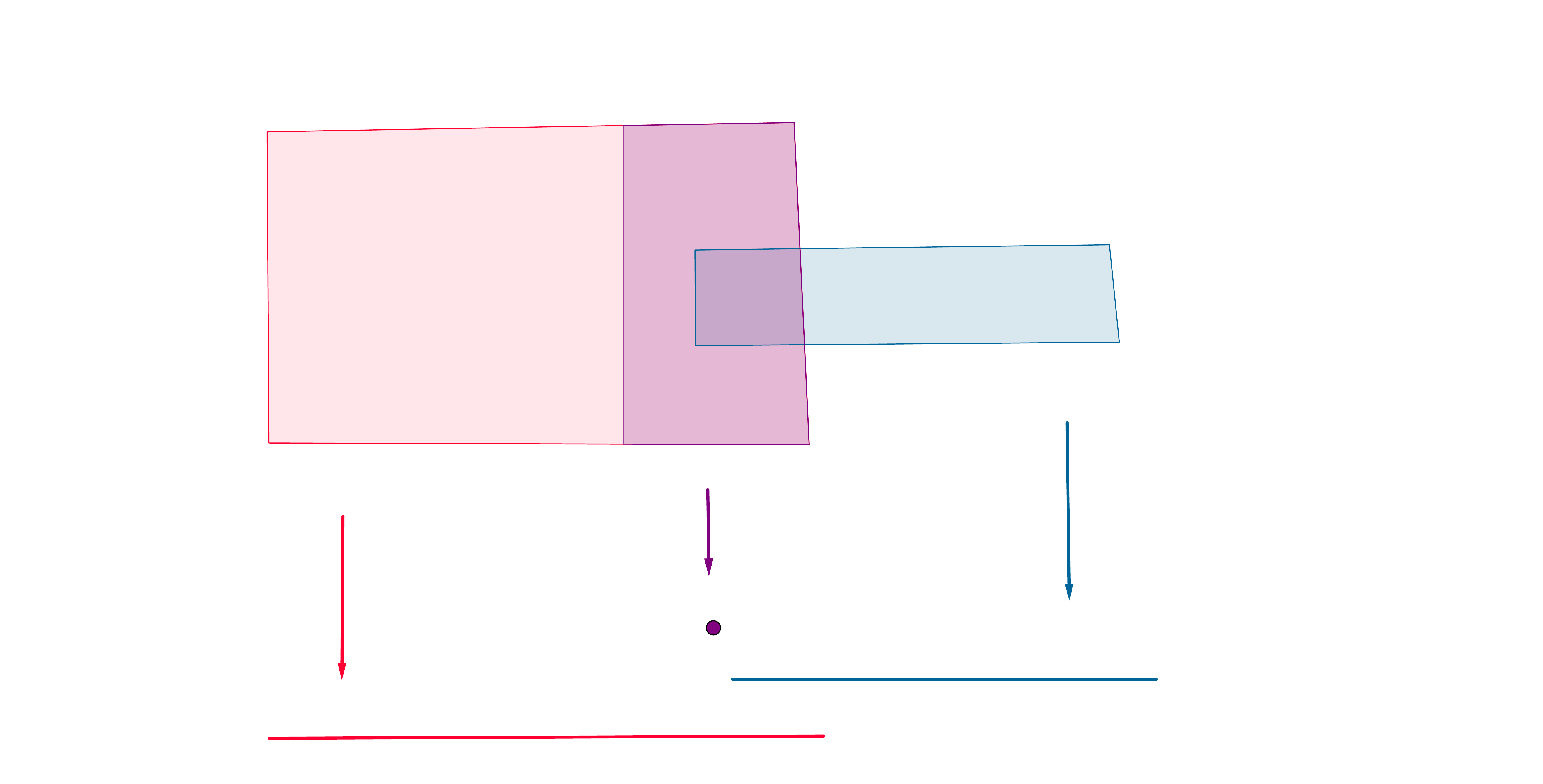}
\vspace{-22pt}
\caption{{\protect\tiny {The purple artificially $0$--framed set is
separating the red and blue $1$--framed sets.}}}
\label{art fr sep}
\end{wrapfigure} if $(H,p_{H})$ is $l$--framed, then $H$ is often $k$%
--framed for all $k\in \left\{ 0,1,\ldots ,l\right\} .$ Indeed if $L\subset 
\mathbb{R}^{l}$ is an affine $k$--dimensional subspace, then by composing $%
p_{H}$ with orthogonal projection $\pi _{L}:\mathbb{R}^{l}\longrightarrow L,$
we get a submersion $\pi _{L}\circ p_{H}$ to the lower-dimensional Euclidean
space, $L.$ With this and some other mild hypotheses, $H$ becomes $k$%
--framed. When we take this perspective on $H$, we say that $H$ is \textbf{%
artificially }$k$\textbf{-framed } (see Figures \ref{art fr basic} and \ref%
{art fr sep}). The concept of artificially framed sets is formally
introduced in Subsection \ref{atrf framed subset}.

In the remainder of Part \ref{Part 0}, we develop tools to deform framed
covers to respectful covers and to verify the disjointness of framed sets.
These are the topics of Subsections \ref{creating resp} and \ref{quant
openness}.

\section{The Perelman Framed Cover\label{perlem framed cover sect}}

In this section, we prove the Perelman Covering Lemma.

\begin{lemma}[Perelman Covering Lemma]
\label{perelmanframedcover} Given $\varkappa ,\delta >0$, there is a finite $%
\varkappa $-framed cover $\mathcal{P}=\mathcal{P}^{0}\cup \mathcal{P}%
^{1}\cup \mathcal{P}^{2}$ of $X\setminus X_{\delta }^{4}$.
\end{lemma}

We will show this follows from next result (cf \cite{GrovWilh2, Kap1,
ProWilh1}).

\begin{lemma}[Framed Lemma]
\label{framedsetprop} Given$\ $a $(k,\delta )$-strained point $q\in X$ and $%
d>0$, there is a $(k,\tau (\delta ))$-framed set $(\digamma ,p_{\digamma
},r_{\digamma })$ with $q\in \digamma $ so that 
\begin{equation}
\mathrm{diam}(\digamma )\leq d\mathrm{.\label{diam small}}
\end{equation}
\end{lemma}

With just a little extra effort, we will also take the following major step
toward getting a proper cover.

\begin{lemma}[Proper Lemma]
\label{highstrainbdry}Let $q\in X$ be as in the Framed Lemma. We can choose $%
\digamma $ to satisfy the conclusion of the Framed Lemma and to also have
the property that for every $y\in p_{\digamma }(\digamma )$, every point in
the boundary of the closure of the fiber $p_{\digamma }^{-1}(y)$ is $%
(k+1,\tau (\delta ))$-strained$.$
\end{lemma}

\begin{proof}[Proof of the Perelman Covering Lemma assuming Lemma \protect
\ref{framedsetprop}]
Let $\tau _{F}$ and $\tau _{P}$ be the $\tau s$ of the Framed and Proper
Lemmas, respectively. Set 
\begin{equation*}
\tau :=\max \left\{ \tau _{P},\tau _{F}\right\} .
\end{equation*}%
Given $\varkappa ,\delta >0$, choose $\delta _{1},\delta _{2}$ so that $%
0<\delta _{1}<\delta _{2}<\delta $, 
\begin{equation}
\text{ }\tau \left( \delta _{2}\right) <\min \{\delta ,\varkappa \},\text{
and }\tau \left( \delta _{1}\right) <\min \{\delta _{2},\delta ,\varkappa \}.
\label{delta
choices}
\end{equation}

Let $\mathfrak{s}^{0}$ be the points of $X$ that are not $\left( 1,\delta
_{1}\right) $--strained. Since $\mathfrak{s}^{0}$ is compact, we cover $%
\mathfrak{s}^{0}$ by $0$--framed sets and let $\mathcal{P}^{0}$ be a finite
subcollection. We show in Proposition \ref{0 str isol prop} (below) that $%
s^{0}$ is isolated. We can thus choose the diameters of the $\digamma $
small enough so that the closures of the $\digamma $ are disjoint.

By the definition of $\mathfrak{s}^{0},$ all points in $X\setminus (\cup 
\mathcal{P}^{0})^{-}$ are at least $(1,\delta _{1})$-strained. Let $%
\mathfrak{s}^{1}$ be all points in $X\setminus (\cup \mathcal{P}^{0})^{-}$
that are not $\left( 2,\delta _{2}\right) $-strained. Use compactness and
the Framed Lemma to choose a finite collection $\mathcal{P}^{1}$ of $(1,\tau
_{F}\left( \delta _{1}\right) )$-framed sets that cover $\mathfrak{s}^{1}.$

All points in $X\setminus ((\cup \mathcal{P}^{0})^{-}\cup (\cup \mathcal{P}%
^{1})^{-})$ are then at least $(2,\delta _{2})$-strained in $X$. Let $%
\mathfrak{s}^{2}$ be all points in $X\setminus ((\cup \mathcal{P}%
^{0})^{-}\cup (\cup \mathcal{P}^{1})^{-})$ that are not $(4,\delta )$%
-strained. Again using compactness, construct a finite cover $\mathcal{P}%
^{2} $ of $\mathfrak{s}^{2}$ by $(2,\tau _{F}\left( \delta _{2}\right) )$%
-framed sets. It follows that $X\setminus X_{\delta }^{4}\subset (\cup 
\mathcal{P}^{0})^{-}\cup (\cup \mathcal{P}^{1})^{-}(\cup \mathcal{P}%
^{2})^{-},$ as desired.
\end{proof}

\addtocounter{algorithm}{1}

\subsection{Lemmas \protect\ref{framedsetprop} and \protect\ref%
{highstrainbdry}}

Except for the statement about the diameter, the Framed Lemma is a
consequence of results in, e.g., \cite{GrovWilh2, Kap1, ProWilh1}. Since we
will also use several explicit facts from the proof, we give a
self-contained exposition below.

The main tool is a reformulation of Perelman's construction of a strictly
concave function on small neighborhoods of any point $q$ in an Alexandrov
space $X$ (see Proposition 3.1 of \cite{ProWilh1}, cf pages 389-390 of \cite%
{Grmv}).

\begin{proposition}[Perelman's Concave Function]
\label{Grisha prop}For $q\in X$, let $p$\ be any point in the interior\ of a
segment with endpoint $q$. Given $\varepsilon >0$, there is a $\rho >0$\ and
a $C^{1}$-function 
\begin{equation*}
f_{p}:B(q,\rho )\longrightarrow \mathbb{R}
\end{equation*}%
that is strictly $(-1)$-concave and satisfies 
\begin{equation}
|f_{p}-\mathrm{dist}_{p}|_{C_{\mathrm{stab}}^{1}}<\varepsilon .
\label{C1
stab approx}
\end{equation}%
\emph{\ }
\end{proposition}

\begin{proof}
Apart from Inequality (\ref{C1 stab approx}), this is Proposition 3.1 of 
\cite{ProWilh1}. To establish Inequality (\ref{C1 stab approx}), note that
in \cite{ProWilh1} we saw that 
\begin{equation}
f_{p}\left( \cdot \right) :=\frac{1}{\mathrm{vol}\left( B\left( p,\eta
\right) \right) }\int_{z\in B\left( p,\eta \right) }\psi \circ \mathrm{dist}%
_{z}\left( \cdot \right) ,  \label{formuala for f}
\end{equation}%
where $\psi :\mathbb{R}$ $\longrightarrow \mathbb{R}$ is an appropriately
chosen, smooth, strictly increasing function. Inequality (\ref{C1 stab
approx}) follows from (\ref{formuala for f}) and the fact that we have
chosen $p$ to be in the interior of a segment with endpoint\emph{\ }$q.$
\end{proof}

Since Perelman's functions are determined by distance functions, they are
stable under Gromov-Hausdorff convergence. Our proof of the Framed and
Proper Lemmas combines this fact together with the following
\textquotedblleft ideal\textquotedblright\ model for local Alexandrov
geometry.

\begin{proposition}
\label{local}Suppose $q$ is $(k,\delta ,r)$-strained in $X.$ If $\lambda $
is sufficiently small, there is a $\tau (\delta ,\lambda )$-embedding%
\begin{equation}
\Phi _{\lambda }:\frac{1}{\lambda }B(q,\lambda )\longrightarrow \mathbb{R}%
^{k}\times C(E)  \label{ideal embedding}
\end{equation}%
so that $\Phi _{\lambda }\left( q\right) =\left( 0,\ast \right) ,$ where $%
C(E)$ is the Euclidean cone on some Alexandrov space $E$ with $\mathrm{curv}%
\geq 1$, and $\ast $ is the cone point in $C(E).$
\end{proposition}

Motivated by this result, we call $\mathbb{R}^{k}\times C(E)$ the \textbf{%
ideal tangent cone} and denote it by $T_{q}^{\mathrm{ideal}}X$. In Corollary %
\ref{alm join} (below), we show that on compact sets, $T_{q}^{\mathrm{ideal}%
}X$ is almost isometric to $T_{q}X$ and that the Euclidean factor of $T_{q}^{%
\mathrm{ideal}}X$ corresponds to an almost Euclidean factor of $T_{q}X.$
Thus $T_{q}^{\mathrm{ideal}}X$ is the idealized local model we would have if
the $\delta $\ of $(k,\delta ,r)$ were $0.$ Before proving Proposition \ref%
{local}, we exploit it to prove the Framed and Proper Lemmas.

\begin{proof}[Proof of the Framed and Proper Lemmas assuming Proposition 
\protect\ref{local}]
First we use Perelman's functions to construct a $k$-framed neighborhood of $%
(0,\ast )$ in $\mathbb{R}^{k}\times C(E)$. Combining this with the $\tau
(\delta ,\lambda )$-embedding (\ref{ideal embedding}), we will get our $k$%
-framed neighborhood of $q.$

Motivated by euclidean geometry, for $v,w\in T_{q}^{\mathrm{ideal}}X$ with $%
\mathrm{dist}\left( \left( 0,\ast \right) ,v\right) =\mathrm{dist}\left(
\left( 0,\ast \right) ,w\right) =1,$ we let $\sphericalangle (v,w)$ denote
the intrinsic distance between $v$ and $w$ in 
\begin{equation*}
S\left( \left( 0,\ast \right) ,1\right) :=\left\{ \left. u\in T_{q}^{\mathrm{%
ideal}}X\text{ }\right\vert \text{ }\mathrm{dist}\left( \left( 0,\ast
\right) ,u\right) =1\right\} .
\end{equation*}

To construct $k$-framed sets about $(0,\ast )$ in $T_{q}^{\mathrm{ideal}}X$,
we select points $\{v_{1},\ldots ,v_{k}\}$ in $T_{q}^{\mathrm{ideal}}X$ that
are close to an orthonormal basis for $\mathbb{R}^{k}$ and satisfy $%
\sphericalangle (v_{i},v_{j})>\pi /2$ for all $i\neq j$. Next, we select a
finite set of points $\{e_{i}\}\subset S\left( \left( 0,\ast \right)
,1\right) $ that are close to $E$ and satisfy

\begin{enumerate}
\item $\sphericalangle (e_{i},v_{j})>\pi /2$ for all $i$ and all $j$.

\item $\{e_{i}\}$ is no more than $\pi /100$ from a subset of $E$ that is $%
\frac{\pi }{8}$--dense in $E.$
\end{enumerate}

We then take the Perelman functions $f_{v_{i}},f_{-v_{i}},$ and $f_{e_{j}}$
from Proposition \ref{Grisha prop} and adjust them by adding constants so
that they all agree at $(0,\ast ).$ By Proposition \ref{Grisha prop}, there
is a $\tilde{\beta}_{\max }$ and a $\tilde{\nu}_{\max }$ so that the $%
f_{v_{i}}$ and the $f_{e_{i}}$ are $\left( -1\right) $--concave on 
\begin{equation*}
\digamma :=\mathbb{B}^{l}\left( \tilde{\beta}_{\max }\right) \times B\left(
\ast ,\tilde{\nu}_{\max }\right) .
\end{equation*}%
The density of the configuration of the $\{v_{i}\}$ and $\{e_{j}\}$ then
implies that the strictly concave down function 
\begin{equation*}
h=\min_{i,j}\{f_{v_{i}},f_{-v_{i}},f_{e_{j}}\}
\end{equation*}%
has a unique maximum at $(0,\ast )$. Moreover, if $\tilde{\beta}_{\max }$
and $\tilde{\nu}_{\max }$ are sufficiently small, then \textrm{diam}$\left(
\digamma \right) \leq d.$ Since the $\{v_{i}\}$ are close to an orthonormal
basis, the map $p_{\digamma }:\digamma \rightarrow \mathbb{R}^{k}$ defined
as 
\begin{equation}
p_{\digamma }=(f_{v_{1}},\ldots ,f_{v_{k}})  \label{subm dfn}
\end{equation}%
is a $\tau (\delta )$-almost Riemannian Submersion, if $\tilde{\beta}_{\max
} $ and $\tilde{\nu}_{\max }$ are sufficiently small. Since $\sphericalangle
(v_{i},v_{j})>\pi /2,$ $\sphericalangle (e_{i},v_{j})>\pi /2$, and the $%
f_{e_{i}}$ are $(-1)$-- concave, the function $r_{\digamma }:\digamma
\rightarrow \mathbb{R}$ defined as%
\begin{equation*}
r_{\digamma }=\min_{i}\{f_{e_{i}}\}
\end{equation*}%
is $(-1)$--concave when restricted to any fiber of $p_{\digamma }$ (see
e.g., pages 663--664 in \cite{Wilh}). By construction, the restriction of $%
r_{\digamma }$ to the fiber $p_{\digamma }^{-1}(p_{\digamma }(0,\ast ))$ is
maximal at $(0,\ast ),$ and by concavity, this maximum is unique. By
continuity, for any $b\in \mathbb{B}^{l}\left( \tilde{\beta}_{\max }\right) $%
, the unique maximum of $r_{\digamma }|_{p_{\digamma }^{-1}\left( b\right)
\cap \digamma }$ will be no more than $\frac{\tilde{\nu}_{\max }}{100}$ from 
$\mathbb{R}^{k}\times \{\ast \}$ provided $\tilde{\beta}_{\max }$ is
sufficiently small. In particular, $\left( \digamma ,p_{\digamma }\right) $
is a fiber bundle whose fibers are disks. Since the base $\mathbb{B}%
^{l}\left( \tilde{\beta}_{\max }\right) $ of $\left( \digamma ,p_{\digamma
}\right) $ is contractible, $\left( \digamma ,p_{\digamma }\right) $ is a
trivial vector bundle.

Since every point in $\mathbb{R}^{k}\times \left( C(E)\setminus \left\{ \ast
\right\} \right) $ is $\left( k+1\right) $--strained and the maxima of $%
r_{\digamma }$ along the fibers of $p_{\digamma }$ are close to $\mathbb{R}%
^{k}\times \{\ast \}$, the Proper Lemma \ref{highstrainbdry} holds for the
framed set $(\digamma ,p_{\digamma },r_{\digamma })\subset T_{q}^{\mathrm{%
ideal}}X$. Since being $(k+1)$-strained is stable under Gromov-Hausdorff
approximation$,$ the Proper Lemma \ref{highstrainbdry} follows from
Proposition \ref{local}.
\end{proof}

For the proof of Proposition \ref{local}, we use the Join Lemma from \cite%
{GrovWilh1} and a Theorem of Plaut from \cite{Pla}. Plaut calls a set of $2n$
points $x_{1},y_{1},\ldots ,x_{n},y_{n}$ in a metric space $Y$ \textbf{%
spherical} if $\mathrm{dist}(x_{i},y_{i})=\pi $ for all $i$ and $\det [\cos 
\mathrm{dist}(x_{i},x_{j})]>0.$

\begin{remarknonum}
If $x_{1},\ldots ,x_{n}$ are points in $\mathbb{S}^{n+k}\subset \mathbb{R}%
^{n+k+1},$ then $\sqrt{\mathrm{det}[\cos \mathrm{dist}_{\mathbb{S}%
^{n+k}}(x_{i},x_{j})]}$ is the $n$--dimensional volume of the parallelepiped
spanned by $\left\{ x_{1},\ldots ,x_{n}\right\} .$ So Plaut's condition
should be viewed as a quantification of linear independence.
\end{remarknonum}

\begin{theorem}[Plaut \protect\cite{Pla}]
\label{plauttheorem} If $X$ has curvature $\geq 1$ and contains a spherical
set $\sigma $ of $2(n+1)$ points, then there is a subset $S\subset X$ that
contains $\sigma $ and is isometric to $\mathbb{S}^{n}$.
\end{theorem}

\begin{definition}
As in \cite{ProWilh1} (cf also \cite{BGP, OSY}), we say an Alexandrov space $%
\Sigma $ with $\mathrm{curv\,}\Sigma \geq 1$ is \textbf{globally }$(m,\delta
)$\textbf{-strained} by pairs of subsets $\{A_{i},B_{i}\}_{i=1}^{m}$
provided 
\begin{equation*}
\begin{array}{ll}
|\mathrm{dist}(a_{i},b_{j})-\frac{\pi }{2}|<\delta , & \mathrm{dist}%
(a_{i},b_{i})>\pi -\delta , \\ 
|\mathrm{dist}(a_{i},a_{j})-\frac{\pi }{2}|<\delta , & |\mathrm{dist}%
(b_{i},b_{j})-\frac{\pi }{2}|<\delta%
\end{array}%
\end{equation*}%
for all $a_{i}\in A_{i}$, $b_{i}\in B_{i}$ and $i\neq j$.
\end{definition}

\begin{lemma}
\label{alm join lemma}Given natural numbers $k\leq n+1$ and $v>0,$ there is
a $\delta >0$ with the following property: If $\Sigma $ is a globally $%
(k,\delta )$-strained $n$--dimensional Alexandrov space with $\mathrm{curv}%
\geq 1$ and $\mathrm{vol}>v,$ then there is an Alexandrov space $E$ of $%
\mathrm{curv}\geq 1$ and a $\tau \left( \delta |v\right) $--homeomorphism 
\begin{equation*}
h:\Sigma \longrightarrow \mathbb{S}^{k-1}\ast E,
\end{equation*}%
where $\mathbb{S}^{k-1}\ast E$ has the spherical join metric.
\end{lemma}

\begin{proof}
If not, then there is a sequence of $n$--dimensional Alexandrov spaces $%
\left\{ \Sigma _{i}\right\} $ with $\mathrm{curv}\geq 1$ and $\mathrm{vol}>v$
that are globally $(k,\delta _{i})$-strained with $\delta _{i}\rightarrow 0,$
and that are not $\tau \left( \delta _{i}|v\right) $--homeomorphic to a
space of the form $\mathbb{S}^{k-1}\ast E.$

After passing to a convergent subsequence, we have $\lim_{i\rightarrow
\infty }\Sigma _{i}=\Sigma ,$ where $\Sigma $ is globally $(k,0)$-strained.
This global strainer for $\Sigma $ is a spherical set of $2k$ points. By
Plaut's theorem (above), $\Sigma $ contains a metrically embedded copy of $%
\mathbb{S}^{k-1}.$ By the Join Lemma, there is an Alexandrov space $E$ of $%
\mathrm{curv}\geq 1$ so that $\Sigma $ is isometric to the spherical join $%
\mathbb{S}^{k-1}\ast E.$ By Perelman's Stability Theorem, all but finitely
many of the $\Sigma _{i}$s are $\tau \left( \delta _{i}|v\right) $%
--homeomorphic to $\Sigma ,$ a contradiction.
\end{proof}

A lower volume bound and an upper diameter bound for an Alexandrov space $X$
give a lower bound for the volume of all spaces of directions of $X$. So
applying the previous result with $\Sigma $ being a space of directions of $%
X $ gives us

\begin{corollary}
\label{alm join}Let $X$ be an $n$--dimensional Alexandrov space with $%
\mathrm{curv}\geq -1.$ Given $k\in \left\{ 0,1,\ldots ,n\right\} ,$ there is
a $\delta >0$ with the following property:

If $q\in X$ is $\left( k,\delta \right) $--strained, then there is a space $%
E $ of $\mathrm{curv}\geq 1$ and a pointed homeomorphism 
\begin{equation*}
\Psi :\left( T_{q}X,\ast \right) \longrightarrow \left( T^{\mathrm{ideal}%
}_{q}X,\left( 0,\ast \right) \right)
\end{equation*}%
whose restriction to any $r$--ball around $\ast $ $\in T_{q}X$ is a $\tau
\left( \delta \right) r$--embedding.
\end{corollary}

This, together with the Parameterized Stability Theorem (\cite{Kap2},
Theorem 7.8) and the fact $\left( \lambda X,q\right) \rightarrow \left(
T_{q}X,\ast \right) $ as $\lambda \rightarrow \infty ,$ gives us Proposition %
\ref{local}.

We close this subsection with a proof of the fact that $0$-strained points
in $X$ are isolated.

\begin{proposition}
\label{0 str isol prop}Given $\delta >0,$ the set of points in $X$ that are
not $\left( 1,\delta \right) $--strained is isolated.
\end{proposition}

\begin{proof}
Suppose $q\in X$ is not $(1,\delta )$-strained and $x\in X$ satisfies $%
\mathrm{dist}\left( q,x\right) =\frac{\lambda }{2}.$ Let 
\begin{equation*}
\Omega _{\lambda }:\frac{1}{\lambda }B\left( q,\lambda \right)
\longrightarrow B\left( \ast ,1\right) \subset T_{q}X
\end{equation*}%
be a $\tau \left( \lambda \right) $--Gromov-Hausdorff approximation with 
\begin{equation}
\mathrm{dist}\left( \Omega _{\lambda }\left( x\right) ,\frac{1}{2}\uparrow
_{q}^{x}\right) <\tau \left( \lambda \right) .  \label{and of blow up}
\end{equation}%
In $T_{p}X,$ we have that $\frac{1}{2}\uparrow _{q}^{x}$ is $\left( k,\delta
,r\right) =\left( 1,0,\frac{1}{4}\right) $-strained by $\ast $ and the
direction $\uparrow _{q}^{x}.$ For $\delta >0,$ the property of being $%
\left( k,\delta ,r\right) $--strained is Gromov-Hausdorff stable, so (\ref%
{and of blow up}) gives us that $x$ is $\left( 1,\tau \left( \delta ,\lambda
\right) ,\frac{1}{4}\right) $--strained in $\frac{1}{\lambda }X,$ and hence $%
x$ is $\left( 1,\tau \left( \delta ,\lambda \right) ,\frac{\lambda }{4}%
\right) $--strained in $X.$ Since $x$ was an arbitrary point with $\mathrm{%
dist}\left( q,x\right) =\frac{\lambda }{2},$ the result follows.
\end{proof}

\addtocounter{algorithm}{1}

\subsection{The Parameterized Perelman Framed Cover}

In this section, we record further properties of the Perelman cover that are
consequences of its construction.

\begin{lemma}[Parameterized Framed Lemma]
\label{cor two parame exh} Let $q\in X$ be $(k,\delta ,r)$-strained. There
are numbers $\beta _{\max }$ and $\nu _{\max }$ and a $(k,\tau \left( \delta
\right) )$-framed neighborhood $(F(\beta _{\max },\nu _{\max }),p,r)$ of $q$
so that the following hold:

\begin{enumerate}
\item For all $(\beta ,\nu )\in \lbrack \beta _{\max }/1000,\beta _{\max
}]\times \lbrack \nu _{\max }/1000,\nu _{\max }]$, there is a set $\digamma
(\beta ,\nu )\subset F(\beta _{\max },\nu _{\max })$ so that $(\digamma
(\beta ,\nu ),p,r)$ is $(k,\tau \left( \delta \right) )$-framed and $%
\digamma (\beta ,\nu )$ satisfies the conclusion of the Proper Lemma.

\item If $\beta _{1}<\beta _{2}$ and $\nu _{1}<\nu _{2}$, then 
\begin{equation*}
\digamma (\beta _{1},\nu _{1})\Subset \digamma (\beta _{2},\nu _{2}).
\end{equation*}

\item The Gromov-Hausdorff distance 
\begin{equation*}
\mathrm{dist}(\digamma (\beta ,\nu ),\mathbb{B}^{k}(\beta )\times B(\ast
,\nu ))<\tau (\delta )(\beta _{\max }+\nu _{\max }),
\end{equation*}%
where $\mathbb{B}^{k}(\beta )\times B(\ast ,\nu )$ is a product of metric
balls in $T_{q}^{\mathrm{ideal}}X=\mathbb{R}^{k}\times C(E)$.
\end{enumerate}
\end{lemma}

\begin{proof}
We get Parts 1 and 2 by combining the construction of $\digamma $ from the
proof of the Framed Lemma with the Proper Lemma. To prove Part 3, we first
observe that since the statement is Gromov--Hausdorff stable, as in the
proof of Lemma \ref{framedsetprop}, it suffices to prove the result with $X$
replaced by $T_{q}^{\mathrm{ideal}}X$ and $q$ replaced by $\left( 0,\ast
\right) .$ The rest of the proof is very similar to the the proof of Lemma %
\ref{framedsetprop}, the difference being we replace $\{v_{1},\ldots
,v_{k}\} $ with a $\xi $-dense subset of $\left\{ s_{i}\right\} \subset
S\left( \left( 0,\ast \right) ,1\right) \cap \mathbb{R}^{k},$ and we require
that $\{e_{i}\}$ be a set that is no more than $\xi $ from a $\xi $-dense
subset of $E.$ We let 
\begin{equation*}
h_{B}=\min \{f_{s_{i}}\}\text{~~~ and ~~~}r=\min \{f_{e_{j}}\}.
\end{equation*}%
We again add constants to $f_{s_{i}}$ and $f_{e_{j}}$ so that they all agree
at $(0,\ast ).$ This allows us to assume that $h_{B}|_{\mathbb{R}^{k}\times
\ast }$ and $r|_{0\times C\left( E\right) }$ have their unique maximums at $%
(0,\ast )$. We define $p_{\digamma }$ via (\ref{subm dfn}), as before, and
let $\digamma (\beta _{\max },\nu _{\max })$ be the intersection of
superlevels of $h_{B}$ and $r$.

For each of the $f_{s_{i}}$ and the $f_{e_{j}},$ let $\varepsilon $ be as in
(\ref{C1 stab approx}). By choosing both $\varepsilon $ and $\xi $ to be
small, we force the directions of maximal decrease of $h_{B}$ and $r$ at a
point $(u,v)\in \digamma (\beta _{\max },\nu _{\max })$ to be arbitrarily
close to the radial fields of $\mathbb{R}^{k}\times \{v\}$ and $\{u\}\times
C(E),$ uniformly on compact subsets of $\left\{ \mathbb{B}^{k}(\beta _{\max
})\setminus \left\{ 0\right\} \right\} \times \left\{ B(\ast ,\nu _{\max
})\setminus \left\{ 0\right\} \right\} .$ Part 3 follows by combining this
with the Fundamental Theorem of Calculus.
\end{proof}

\begin{wrapfigure}[8]{r}{0.5\textwidth}\centering
\vspace{-15pt}
\includegraphics[scale=0.14]{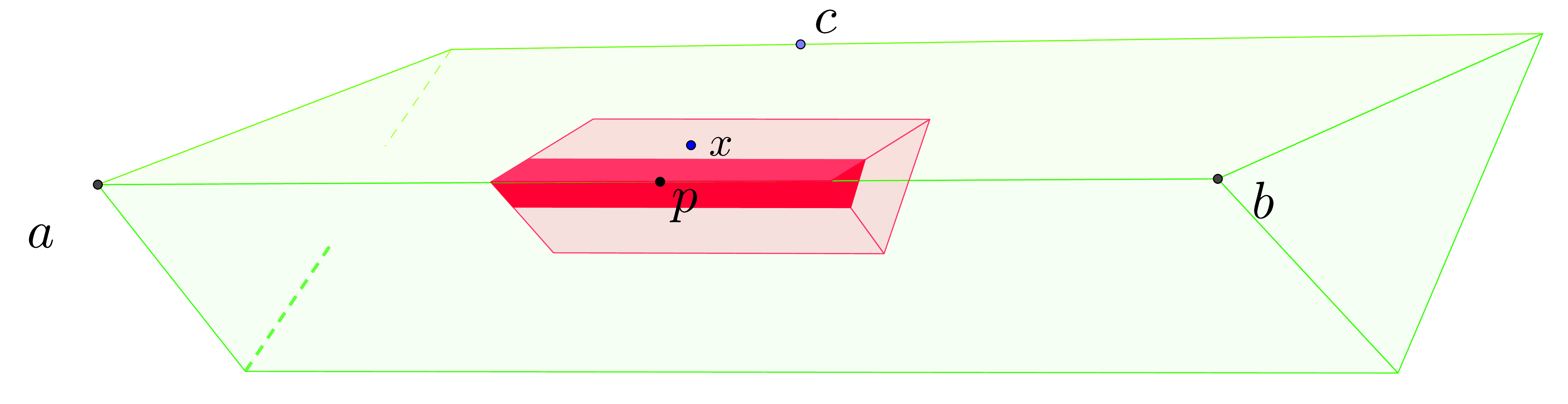}
\caption{{\protect\tiny {Let $X$ be the double of the convex subset of $%
\mathbb{R}^3$ in the picture. The double of the red set is $1$--strained in $%
X$ by $(a,b)$. Lemma \protect\ref{real grad pres} guarantees that the light
red points are also $2 $--strained. For example, $x$ is $2$ strained by $\{
(a,b) , (p,c)\}$.}}}
\label{vitali}
\end{wrapfigure}

\noindent\textbf{Remark.} \emph{In many cases, the sets from Lemma \ref{cor
two parame exh} that we consider will be those for which the two parameters
are equal. When this occurs, we will write $\digamma \left( \beta \right) $
for $\digamma \left( \beta ,\beta \right) .$ When the precise values of the
parameters are inconsequential, we will simply write $\digamma .$} \vspace{%
.1 in}

Part 1 of Lemma \ref{cor two parame exh} gives us

\begin{corollary}[Quantitatively Proper]
\label{real grad pres} The sets $(\digamma (\beta _{\max },\nu _{\min
}),p,r) $ from Lemma \ref{cor two parame exh} have the property that every
point of 
\begin{equation*}
\digamma (\beta _{\max },\nu _{\max })\setminus \digamma (\beta _{\max },\nu
_{\max }/1000)
\end{equation*}%
is $(k+1,\tau (\delta ))$-strained.
\end{corollary}

In fact, our choices of $\varepsilon $ and $\xi $ in the proof of the
Parameterized Framed Lemma give us that $\left( r_{\digamma },p_{\digamma
}\right) $ is an explicit $\left( k+1\right) $--framing for the points in $%
\digamma (\beta _{\max },\nu _{\max })\setminus \digamma (\beta _{\max },\nu
_{\max }/1000).$

\begin{corollary}
\label{cor expli k+1 framing}For each $x\in \digamma \left( \beta _{\mathrm{%
\max }},\nu _{\mathrm{\max }}\right) \setminus \digamma \left( \beta _{%
\mathrm{\max }},\frac{\nu _{\max }}{1,000}\right) ,$ there is a neighborhood 
$N_{x}$ of $x$ that is the total space of a stable trivial vector bundle%
\begin{eqnarray}
p_{N_{x}} &:&N_{x}\longrightarrow \mathbb{R}\times \mathbb{R}^{k}  \notag \\
p_{N_{x}}\left( x\right) &=&\left( r_{\digamma }\left( x\right) ,p_{\digamma
}\left( x\right) \right)  \label{k+1 submer dfn}
\end{eqnarray}%
whose projection is a $\tau \left( \delta \right) $--almost Riemannian
submersion. In particular, $\left\{ \digamma ,N_{x}\right\} $ is a
respectful framed collection.
\end{corollary}

\begin{figure}[tbp]
\includegraphics[scale=0.145]{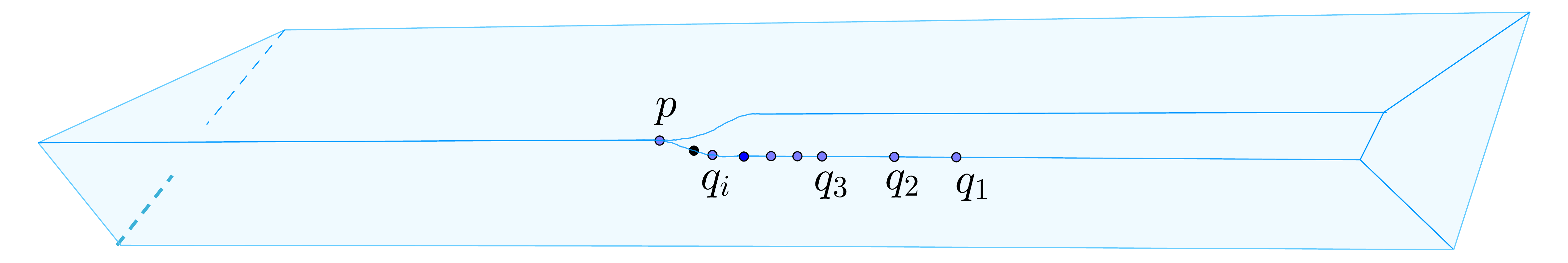}
\caption{{\protect\tiny {In Figure \protect\ref{vitali}, the points that
cannot be $2$--strained are precisely those on the segment $ab$. However,
the double of the convex set in this picture explains why it is difficult to
improve on the statement in Corollary \protect\ref{real grad pres}, apart
from the comment that $\frac{1}{1,000}$ can be any a priori constant.
Consider applying Lemma \protect\ref{cor two parame exh} to each of the
points $q_{i}$. There are two ways to avoid the upper branch of $1$%
--strained points: the quantity $\protect\beta _{max}$ can converge to $0$
as the $q_{i}\rightarrow p$, or the branch can be \textquotedblleft
swallowed\textquotedblright\ by $\digamma (\protect\beta _{max},\frac{%
\protect\nu _{max}}{1,000})$. These considerations together with the fact
that the bifurcation in this figure can be repeated any finite number of
times leads us to think that a practical improvement of Corollary \protect
\ref{real grad pres} is not possible.}}}
\label{bifurc}
\end{figure}

\begin{remarknonum}
The constant $\frac{1}{1,000}$ in (\ref{real grad pres}) can of course be
replaced by any positive constant we like, but as far as we know, this
constant must be chosen before $\digamma $ is constructed. We believe this
is necessary because of the possibly of bifurcating strata. We are grateful
to Vitali Kapovitch for pointing this out to us (cf Figure \ref{bifurc}).
\end{remarknonum}

Together the results of this section give us a parameterized version of the
Perelman Covering Lemma.

\begin{definition}
\label{Perelman cover}A framed cover $\mathcal{H}=\mathcal{H}^{0}\cup 
\mathcal{H}^{1}\cup \mathcal{H}^{2}$ of $X\setminus X_{\delta }^{4}$ is 
\textbf{parameterized }if there is a $b>1$ so that for $t\in \left[ \frac{1}{%
2},b\right] ,$ 
\begin{equation*}
\mathcal{H}^{0}\left( t\right) \cup \mathcal{H}^{1}\left( t\right) \cup 
\mathcal{H}^{2}\left( t\right) :=\left\{ \digamma _{i}^{0}\left( t\beta
_{i}^0\right) \right\} _{\digamma _{i}^{0}\in \mathcal{H}^{0}}\cup \left\{
\digamma _{i}^{1}\left( t\beta _{i}^{1}\right) \right\} _{\digamma
_{i}^{1}\in \mathcal{H}^{1}}\cup \left\{ \digamma _{i}^{2}\left( t\beta
_{i}^{2}\right) \right\} _{\digamma _{i}^{2}\in \mathcal{H}^{2}}
\end{equation*}%
is a framed cover of $X\setminus X_{\delta }^{4},$ and for each $k=0,1,2,$
each member of the $1$--parameter family $\left\{ \digamma _{i}^{k}\left(
t\beta _{i}^{k}\right) \right\} _{t\in \left[ \frac{1}{2},b\right] }$ of $k$%
--framed sets satisfies the conclusions of Lemma \ref{cor two parame exh}
and Corollary \ref{real grad pres}. We say that $\mathcal{H}^{k}$\textbf{\
is parameterized on }$\left[ \frac{1}{2},b\right] .$
\end{definition}

By combining the proof of the Perelman Covering Lemma with Lemma \ref{cor
two parame exh} and Corollary \ref{real grad pres} we get the

\begin{layPerelCovlemma}
\label{layered Perel cov lemma}For all $\delta >0,$ there is a parameterized%
\textbf{\ }framed cover $\mathcal{P}=\mathcal{P}^{0}\cup \mathcal{P}^{1}\cup 
\mathcal{P}^{2}$ of $X\setminus X_{\delta }^{4}$ so that the closures of
elements of $\mathcal{P}^{0}$ are pairwise disjoint, and for each $k\in
\left\{ 0,1,2\right\} ,$ $\mathcal{P}^{k}$\textbf{\ }is parameterized on $%
\left[ \frac{1}{2},25\right] .$
\end{layPerelCovlemma}

\noindent \textbf{Convention: }From this point all of our framed covers will
be parameterized. To keep the nomenclature simple, we will use the term
framed cover for parameterized framed cover. All of our assertions about
framed covers will also be valid for slightly larger and slightly smaller
versions of the covers, but in the interest of simplicity, typically, we
will not mention this.

\section{The Perelman-Plaut Framed Cover\label{perl plt sect}}

\begin{wrapfigure}[7]{r}{0.4\textwidth}\centering
\vspace{-16pt}
\includegraphics[scale=.14]{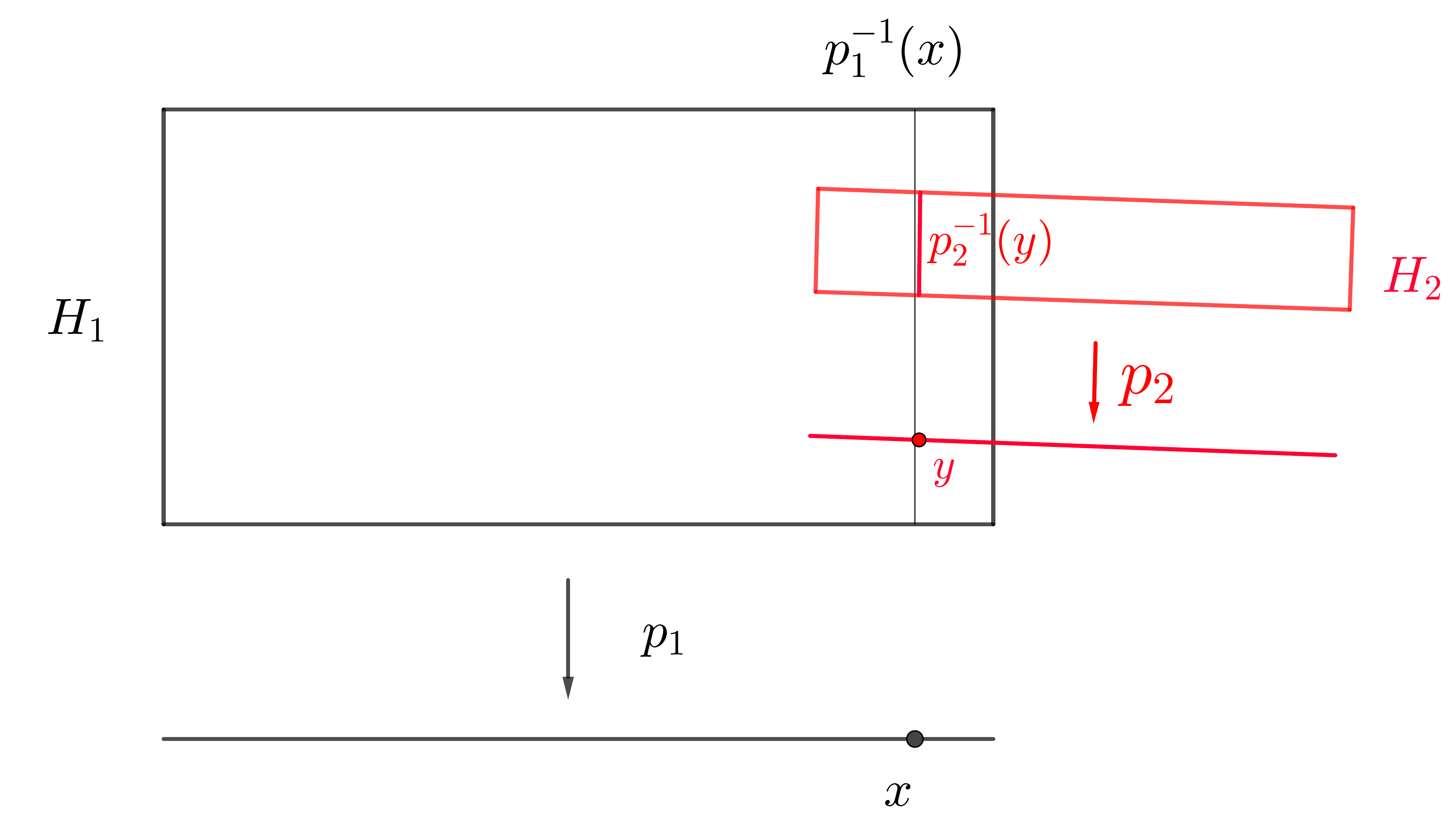}
\vspace{-6pt}
\caption{Fiber swallowing and $\protect\kappa $-lined up}
\label{lined up}
\end{wrapfigure}

\begin{wrapfigure}[8]{l}{0.37\textwidth}\centering
\vspace{-0.27in}
\includegraphics[scale=.15]{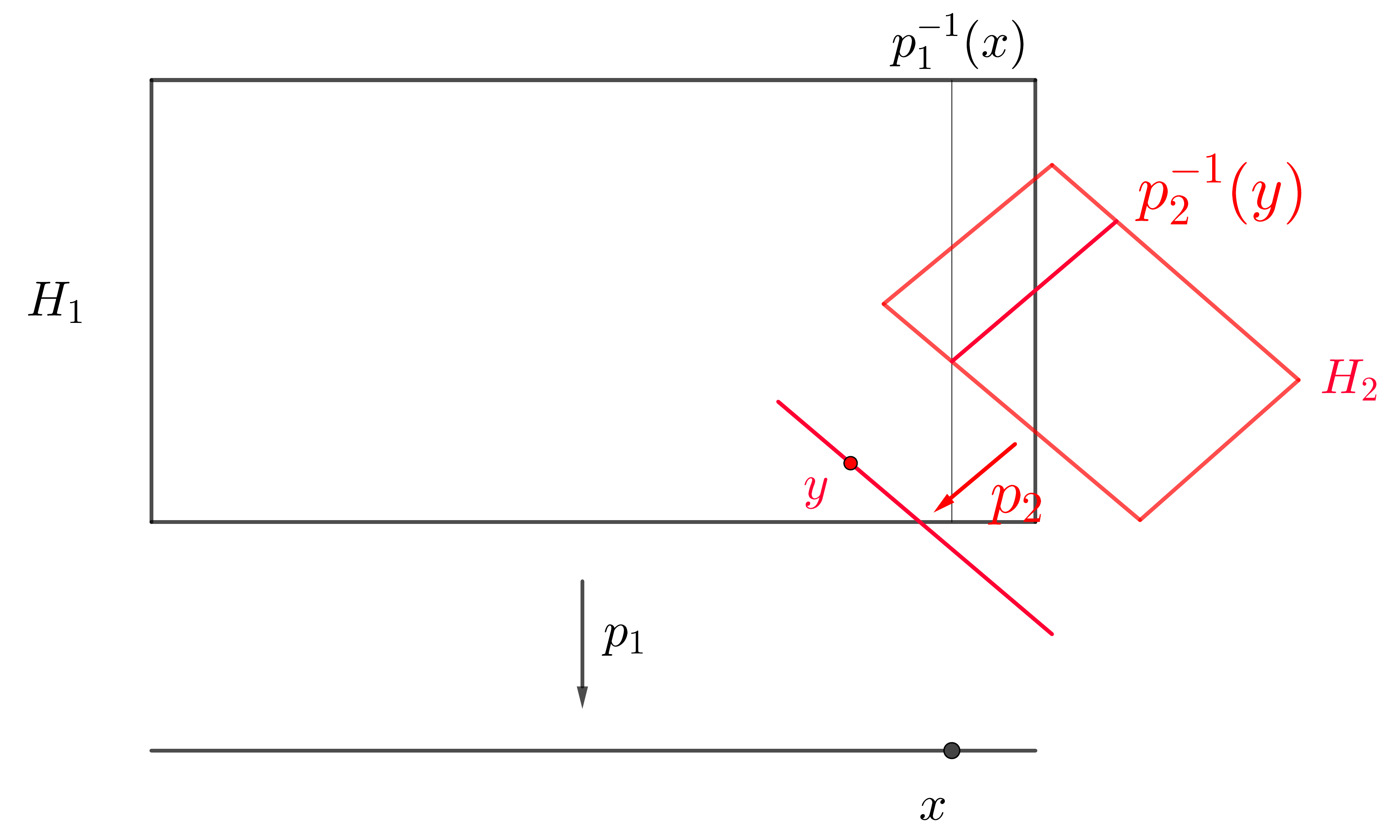}
\vspace{-0.12in}
\caption{Not $\protect\kappa $-lined up}
\label{notlinedup}
\end{wrapfigure}

\begin{wrapfigure}[4]{r}{0.28\textwidth}\centering
\vspace{-16pt}
\includegraphics[scale=0.09]{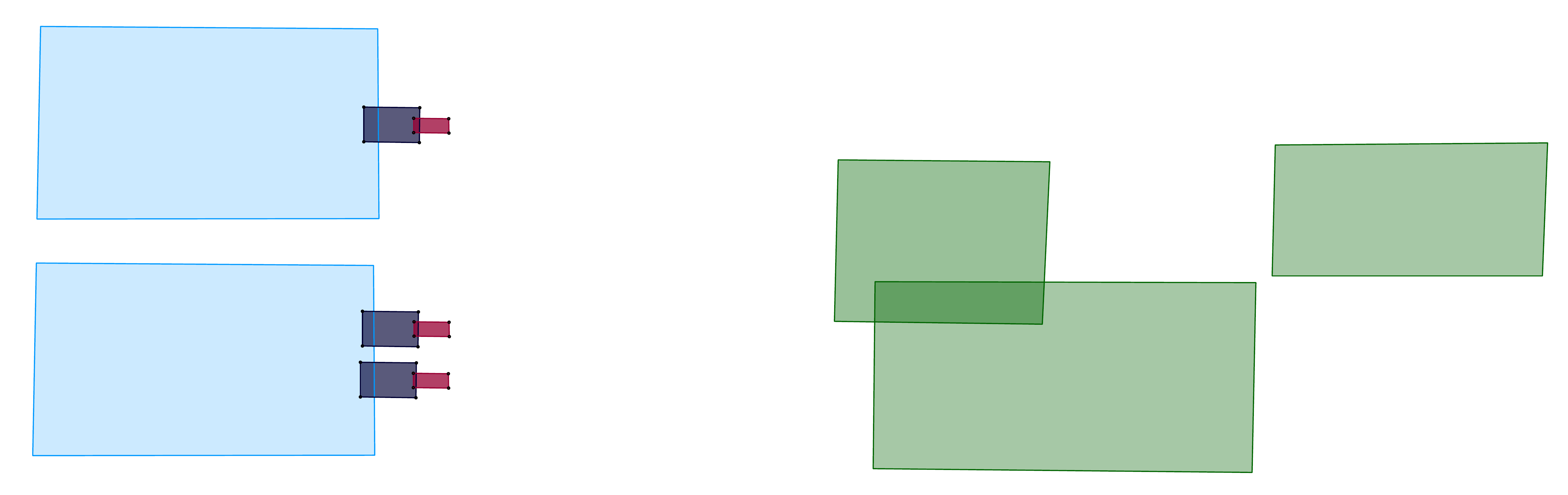}
\vspace{-11pt}
\caption{{\protect\tiny {The framed sets on the left are vertically
separated. Those on the right are not vertically separated. }}}
\label{vert gap fig}
\end{wrapfigure}

In this section, we improve the Perelman Framed Cover to one for which the
intersections of the $1$- and $2$-framed sets are more tightly constrained.
We call these constraints\textbf{\ }$\kappa $\textbf{--lined up, vertically
separated,} and \textbf{fiber swallowing }(see Figures \ref{lined up}, \ref%
{lined up, but not swallowing}, \ref{notlinedup}, and \ref{vert gap fig}).
Before defining these new terms, we state the result.

\begin{lemma}[Perelman-Plaut Covering Lemma]
\label{Ms. Pacman limit lemma} Given $\delta ,\kappa ,\mu >0,$ there is a $%
\varkappa _{0}>0$ with the following property: For all $\varkappa \in \left(
0,\varkappa _{0}\right) $, there is a parameterized, $\varkappa $-framed
cover $\mathcal{Z}=\mathcal{Z}^{0}\cup \mathcal{Z}^{1}\cup \mathcal{Z}^{2}$
of $X\setminus X_{\delta }^{4}$ so that

\begin{enumerate}
\item The closures of the elements of $\mathcal{Z}^{0}$ are pairwise
disjoint.

\item $\mathcal{Z}^{1}$ and $\mathcal{Z}^{2}$ are $\kappa $--lined up,
vertically separated, and fiber swallowing.

\item Given $k\in \left\{ 0,1,2\right\} ,$ all $k$-framed sets $\digamma
_{i}^{k}\left( \beta _{i}^{k}\right) $ in $\mathcal{Z}^{k}$ satisfy $\beta
_{i}^{k}\leq \mu . 
$
\end{enumerate}
\end{lemma}

\noindent We start with the definition of $\kappa $--lined up.

\begin{definition}[See Figures \protect\ref{lined up}, \protect\ref%
{notlinedup}, and \protect\ref{lined up, but not swallowing} ]
\label{Dfn eps line up} Given $\kappa >0$, we say two almost Riemannian
submersions $q,p:U\longrightarrow \mathbb{R}^{k}$ are $\kappa $\textbf{%
--lined up} provided 
\begin{equation}
\hspace*{2.2 in} \left\vert dp-I_{p,q}\circ dq\right\vert _{\mathrm{stab}%
}<\kappa  \label{c1 close inequal}
\end{equation}
for some isometry $I_{p,q}:\mathbb{R}^{k}\longrightarrow \mathbb{R}^{k}$. We
say that $p$ and $q$ are \textbf{lined up }if there is a diffeomorphism $%
I_{p,q}:\mathbb{R}^{k}\longrightarrow \mathbb{R}^{k}$ so that $%
p=I_{p,q}\circ q.$
\end{definition}

\begin{definition}
Given $\kappa >0$, let $\mathcal{U}$ be a collection of $k$-framed sets
parameterized on $[\frac{1}{2},25].$ We say that $\mathcal{U}$ is $\kappa $%
\textbf{--lined up} provided for all $(\digamma _{1}(\beta
_{1}),p_{1}),(\digamma _{2}(\beta _{2}),p_{2})\in \mathcal{U}$ with $%
\digamma _{1}(\beta _{1})\cap \digamma _{2}(\beta _{2})\neq \emptyset ,$ the
submersions $p_{1}$ and $p_{2}$ are $\kappa $-lined up on $\digamma
_{m}\left( 5\beta _{m}\right) \cap \digamma _{q}\left( 5\beta _{q}\right) .$
\end{definition}

\mbox{}\vspace*{-0.3in}

\begin{definition}
A collection of $k$-framed sets $U$ parameterized on $[\frac{1}{2},25]$ is
called \textbf{vertically separated} if and only if for all $(\tilde{\digamma%
}(\tilde{\beta}),\tilde{p}),(\digamma (\beta ),p)\in U$\emph{\ with }$\tilde{%
\beta}\geq \beta \ $ and $\tilde{\digamma}(4\tilde{\beta},\tilde{\beta})\cap
\digamma (\beta ,\beta )=\varnothing ,$ 
\begin{equation*}
\hspace*{-1.75in} \mathrm{dist} \left( \tilde{\digamma}(\tilde{\beta}%
),\digamma (\beta )\right) \geq \frac{\tilde{\beta}}{6}.
\end{equation*}
\end{definition}

\begin{definition}
\label{fiber swalloing} A collection of $k$-framed sets $\mathcal{U}$ is
called \textbf{fiber swallowing} if for all $(\digamma (\beta ),p),(\tilde{%
\digamma}(\tilde{\beta}),\tilde{p})\in \mathcal{U}$ and all $x\in \digamma
(\beta )\cap \tilde{\digamma}(\tilde{\beta})$, 
\begin{equation}
p^{-1}\left( p\left( x\right) \right) \cap \digamma \left( \beta ,\frac{1}{%
100}\beta \right) \subset \tilde{\digamma}\left( \tilde{\beta}^{+},\frac{5}{%
100}\tilde{\beta}\right) ,  \label{cores intersect}
\end{equation}%
provided $\tilde{\beta}\geq \beta .$
\end{definition}

Let $\tau _{PF}$ and $\tau _{QP}$ be the $\tau s$ of the Parameterized
Framed Lemma (\ref{cor two parame exh}) and the Quantitatively Proper
Corollary (\ref{real grad pres}), respectively. Set $\tau :=\max \left\{
\tau _{PF},\tau _{QP}\right\} . $ Given $\delta ,\varkappa >0$, choose $%
\delta _{1},\delta _{2}$ so that $0<\delta _{1}<\delta _{2}<\delta $, 
\begin{equation}
\tau \left( \delta _{2}\right) <\min \{\delta ,\varkappa \}\text{ and }\tau
\left( \delta _{1}\right) <\min \{\delta _{2},\varkappa ,\delta \}.
\label{grown up
delta choices}
\end{equation}

By adding a non-redundancy criterion to the construction of $\mathcal{P}$,
we will get a framed cover $\mathcal{Z}=\mathcal{Z}^{0}\cup \mathcal{Z}%
^{1}\cup \mathcal{Z}^{2}$ so that $\mathcal{Z}^{1}$ and $\mathcal{Z}^{2}$
are $\kappa $--lined up, vertically separated, and fiber swallowing. 
\begin{wrapfigure}[9]{l}{0.423\textwidth}\centering
\vspace{-13.5pt}
\includegraphics[scale=0.155]{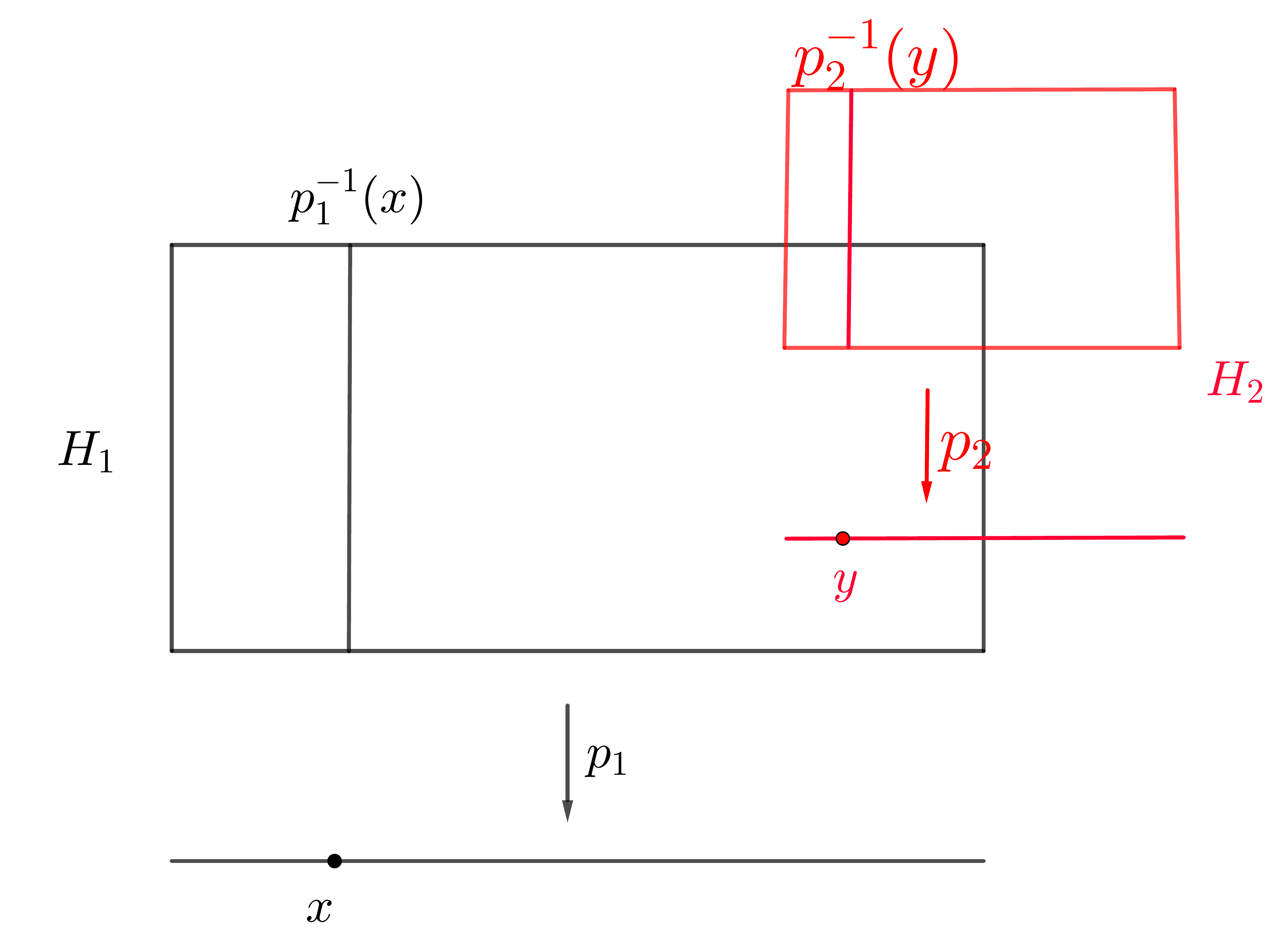}
\vspace{-14pt}
\caption{$\protect\kappa $-lined up but not fiber swallowing}
\label{lined up, but not swallowing}
\end{wrapfigure}
\hspace*{.132in} We set $\mathcal{Z}^{0}=\mathcal{P}^{0},$ and repeat the
construction of the $1$--framed sets in the proof of the Perelman Covering
Lemma (\ref{perelmanframedcover}), with the added constraint that we require
each $Z\in \mathcal{Z}^{1}$ to actually contain a point of $X$ that is not $%
(2,\delta _{2})$-strained. Having obtained $\mathcal{Z}^{1},$ we impose the
analogous requirement on $\mathcal{Z}^{2}$, namely that each of its elements
contain a point of $X$ that is not $(4,\delta )$-strained. $\mathcal{Z}$ is
thus a $\varkappa $-framed cover of $X\setminus X_{\delta }^{4}$, and since $%
\mathcal{Z}^{0}=\mathcal{P}^{0}$, the 
\begin{wrapfigure}[9]{r}{0.4\textwidth}\centering
\vspace{-.13in}
\includegraphics[scale=0.16]{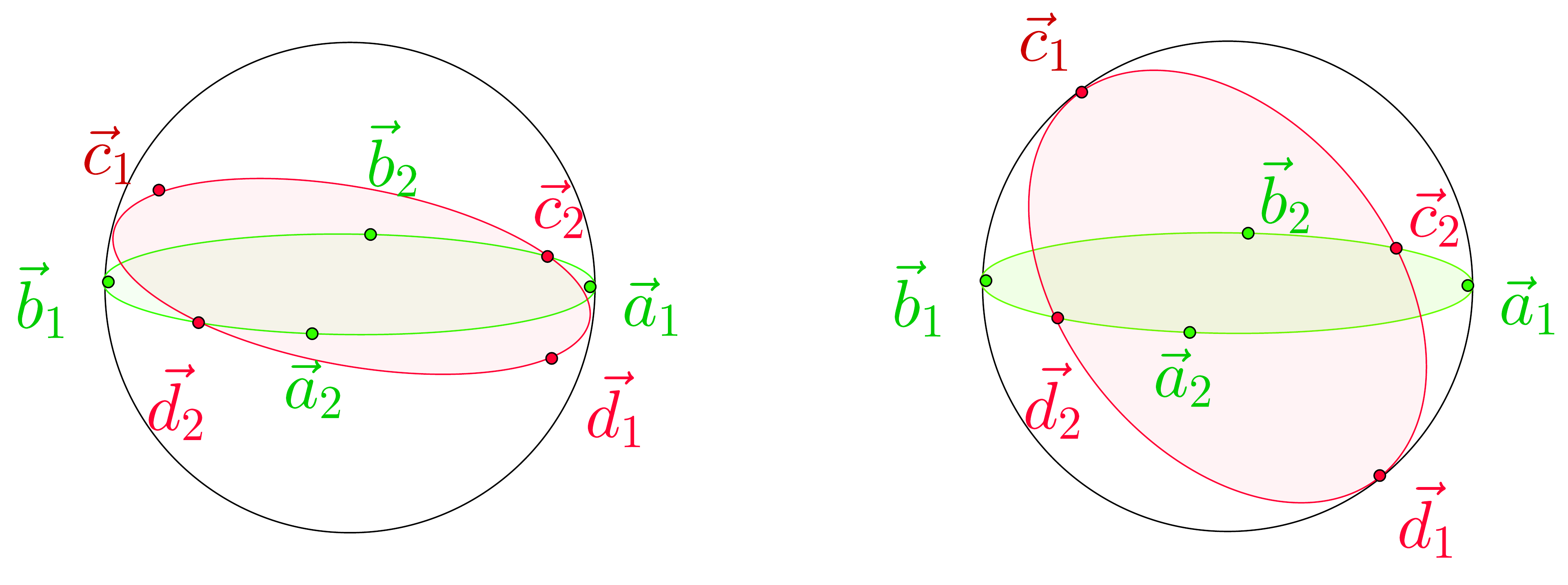}
\vspace{-.07in}
\caption{{\protect\tiny {\ Both pictures are of directions that correspond
to two sets of strainers $\{ (a_1, b_1), (a_2, b_2) \}$ and $\{ (c_1, d_1),
(c_2, d_2) \}$ in an unknown space of directions $\Sigma_p$. In the picture
on the left, the corresponding submersions are $\protect\kappa$--lined up.
In the picture on the right, the submersions are not $\protect\kappa$--lined
up, so $p$ is actually 3--strained. }}}
\label{linedup-or-strained}
\end{wrapfigure}


\noindent closures of the elements of $\mathcal{Z}^{0}$ are disjoint. It
only remains to prove Part 2 of Lemma \ref{Ms. Pacman limit lemma}, that is,
we must show $\mathcal{Z}^{1}$ and $\mathcal{Z}^{2}$ are both $\kappa $%
-lined up, vertically separated, and fiber swallowing. The following
technical result is the key to proving this.

\begin{proposition}[Line up--see Figure \protect\ref{linedup-or-strained}]
\label{line up} Given $\varepsilon _{1},\kappa >0,$ there is a $\varkappa >0$
so that if $(\digamma _{1},p_{\digamma _{1}})$, and $(\digamma
_{2},p_{\digamma _{2}})$ are $(k,\varkappa )$--framed sets from the Lemma %
\ref{cor two parame exh}, then either

\noindent 1. all points in $\digamma _{1}\cap \digamma _{2}$ are $\left(
k+1,\varepsilon _{1}\right) $--strained, or

\noindent 2. $p_{\digamma _{1}}$ and $p_{\digamma _{2}}$ are $\kappa $%
--lined up,

\noindent provided $\max \left\{ \mathrm{diam}\left( \digamma _{1}\right) ,%
\mathrm{diam}\left( \digamma _{2}\right) \right\} $ is small enough compared
with $r$.
\end{proposition}

Before giving the proof we show that this proposition implies that $\mathcal{%
Z}$ is $\kappa $-lined up, vertically separated, and fiber swallowing. We
start with the $\kappa $--lined up assertion.

\begin{proposition}
Assume that Proposition \ref{line up} holds. Given $\delta ,\kappa >0$,
there is a $\varkappa >0$ so that the $\varkappa $-framed cover $\mathcal{Z}$
constructed above is $\kappa $-lined up, provided the diameters of all
elements of $\mathcal{Z}$ are sufficiently small.
\end{proposition}

\begin{proof}
Since the elements of $\mathcal{Z}$ are close to a product of balls, if $%
(\digamma _{1}\left( \beta _{1}\right) ,p_{1}),(\digamma _{2}\left( \beta
_{2}\right) ,p_{2})\in \mathcal{Z}^{k}$, $\beta _{1}\leq \beta _{2}$, and $%
\digamma _{1}(\beta _{1})\cap \digamma _{2}(\beta _{1})\neq \emptyset $,
then 
\begin{equation*}
\digamma _{1}(\beta _{1})\subset \digamma _{2}(5\beta _{2}).
\end{equation*}

Set $\delta =\delta _{3}=\delta _{4}$ and for $k=1,2$, let $\varepsilon
_{1}=\delta _{k+1}$. Assume $\varkappa $ and the diameter of each element of 
$\mathcal{Z}$ is small enough so that the line up Proposition \ref{line up}
holds for all pairs $\left\{ \digamma _{1}(5\beta _{1}),\digamma _{2}(5\beta
_{2})\right\} $ with $\digamma _{1},\digamma _{2}\in \mathcal{Z}^{k}$. If $%
\digamma _{1}(\beta _{1})\cap \digamma _{2}(\beta _{2})\neq \emptyset $, we
claim $p_{1}$ and $p_{2}$ are $\kappa $-lined up. If not, by Proposition \ref%
{line up}, all points of $\digamma _{1}(5\beta _{1})\cap \digamma
_{2}(5\beta _{2})$ are $(k+1,\delta _{k+1})$-strained. However, if $\beta
_{1}\leq \beta _{2}$, then $\digamma _{1}(\beta _{1})\subset \digamma
_{2}(5\beta _{2})$, so all points of $\digamma _{1}(\beta _{1})$ are $%
(k+1,\delta _{k+1})$-strained. Since this contradicts the fact that $%
\digamma _{1}\in \mathcal{Z}^{k}$, $\mathcal{Z}^{k}$ is $\kappa $-lined up.
\end{proof}

Next we show that $\mathcal{Z}$ is vertically separated and fiber swallowing.

\begin{proposition}
Given $\delta >0$, there are $\kappa ,\varkappa >0$ so that if $\mathcal{Z}$
is one of the $\kappa $-lined up $\varkappa $--framed covers of $X\setminus
X_{\delta }^{4}$ constructed above, then $\mathcal{Z}$ is vertically
separated and fiber swallowing.
\end{proposition}

\begin{proof}
For $k=1,2$, suppose $\digamma _{1}(\beta _{1}),\digamma _{2}(\beta _{2})\in 
\mathcal{Z}^{k},$ $\beta _{2}\geq \beta _{1},$ $\digamma _{1}(\beta
_{1},\beta _{1})\cap \digamma _{2}(4\beta _{2},\beta _{2})=\varnothing ,$
and $\mathrm{dist}\left( \digamma _{1}(\beta _{1}),\digamma _{2}(\beta
_{2})\right) <\frac{\beta _{2}}{6}.$ Then $\digamma _{1}(\beta _{1})\subset
\digamma _{2}(4\beta _{2},4\beta _{2})\setminus \digamma _{2}(4\beta
_{2},\beta _{2}).$ Combined with the Quantitatively Proper Corollary (\ref%
{real grad pres}), this implies that all points of $\digamma _{1}(\beta
_{1}) $ are $\tau _{QP}(\delta _{k})<\tau (\delta _{k})<\delta _{k+1}$%
-strained, where, as before, $\delta =\delta _{3}=\delta _{4}.$ This
contradicts $\digamma _{1}\in \mathcal{Z}^{k},$ so $\mathcal{Z}^{k}$ is
vertically separated.

For $k=1,2$, suppose $\digamma _{1}(\beta _{1}),\digamma _{2}(\beta _{2})\in 
\mathcal{Z}^{k}$, $\beta _{2}\geq \beta _{1}$ and $x\in \digamma _{1}(\beta
_{1})\cap \digamma _{2}(\beta _{2})$. Set 
\begin{equation*}
B_{1}(x):=p_{{1}}^{-1}\left( p_{{1}}\left( x\right) \right) \cap \digamma
_{1}\left( \beta _{1},\frac{1}{100}\beta _{1}\right) .
\end{equation*}%
We will show $B_{1}(x)\subset \digamma _{2}(\beta _{2}^{+},\frac{5}{100}%
\beta _{2})$.

By Part 3 of Lemma \ref{cor two parame exh}, if $y,y^{\prime }\in B_{1}(x)$,
then 
\begin{equation*}
\mathrm{dist}(y,y^{\prime })<\frac{2^{+}}{100}\beta _{1}.
\end{equation*}%
If $\kappa $ and $\varkappa $ are small enough and $y\in B_{1}(x)\setminus
\digamma _{2}(\beta _{2}^{+},\frac{5}{100}\beta _{2})$, then 
\begin{equation*}
\mathrm{dist}\left( y,\digamma _{2}\left( \beta _{2},\frac{1}{100}\beta
_{2}\right) \right) >\left( \frac{5}{100}-\frac{1^{+}}{100}\right) \beta
_{2}=\left( \frac{4^{-}}{100}\right) \beta _{2}.
\end{equation*}%
The previous two displays give us, 
\begin{equation*}
B_{1}(x)\cap \digamma _{2}\left( \beta _{2}^{+},\frac{2^{-}}{100}\beta
_{2}\right) =\emptyset .
\end{equation*}%
Together with the fact that $\left\{ \digamma _{1},\digamma _{2}\right\} $
is $\kappa $--lined up, this gives us that 
\begin{equation*}
\digamma _{1}\left( \beta _{1},\frac{1}{100}\beta _{1}\right) \subset
\digamma _{2}\left( 5\beta _{2}\right) \setminus \digamma _{2}\left( 5\beta
_{2},\frac{1}{1,000}\beta _{2}\right) ,
\end{equation*}%
provided $\kappa $ and $\varkappa $ are small enough. Combined with the
Quantitatively Proper Corollary (\ref{real grad pres}), this implies that
all points of $\digamma _{1}$ are $\tau _{QP}(\delta _{k})<\tau (\delta
_{k})<\delta _{k+1}$-strained, where, as before, $\delta =\delta _{3}=\delta
_{4}.$ This contradicts $\digamma _{1}\in \mathcal{Z}^{k}$. Thus $%
B_{1}(x)\subset \digamma _{2}(\beta _{2}^{+},\frac{5}{100}\beta _{2}),$ and $%
\mathcal{Z}^{k}$ is fiber swallowing.
\end{proof}

\addtocounter{algorithm}{1}

\subsection{The Line Up Proposition \protect\ref{line up}}

To complete the proof of Lemma \ref{Ms. Pacman limit lemma}, in this
subsection we prove the Line Up Proposition (\ref{line up}).

We will need the next two results on the stability of straining directions
from \cite{ProWilh1}. Here we assume that $\left\{ M_{\alpha }\right\}
_{\alpha }$ is a sequence of $n$--dimensional Alexandrov spaces that
converge to $X$ and that $B\subset X$ is $\left( k,\delta ,r\right) $%
-strained by $\left\{ \left( a_{i},b_{i}\right) \right\} _{i=1}^{k}.$

\begin{proposition}[Lemma 1.4 in \protect\cite{ProWilh1}]
\label{angle convergence} For any $x\in B$ and $i\neq j,$ if $\theta $ is
any of the angles $\sphericalangle \left( a_{i},x,b_{j}\right)
,\sphericalangle \left( b_{i},x,b_{j}\right) ,\sphericalangle \left(
a_{i},x,a_{j}\right) $, then 
\begin{equation*}
|\pi /2-\theta |<\tau (\delta ).
\end{equation*}
\end{proposition}

To state the next result, we adopt Petrunin's convention (\cite{Pet}) that $%
\Uparrow _{x}^{p}$ denotes the set of directions of segments from $x$ to $p.$
Later we will use $\uparrow _{x}^{p}$ to denote the direction of a single
segment from $x$ to $p$.

\begin{proposition}[Proposition 1.5 in \protect\cite{ProWilh1})]
\label{angggl conv prop} Suppose $y\in B$ and $\left\{ \left( a_{i}^{\alpha
},b_{i}^{\alpha }\right) \right\} _{i=1}^{k}\subset M_{\alpha }\times
M_{\alpha }$ converge to $\left\{ \left( a_{i},b_{i}\right) \right\}
_{i=1}^{l}.$ Let $c^{\alpha }\in M_{\alpha }$ converge to $c\in B\setminus
\left\{ y\right\} .$

Then for $y^{\alpha }\in M^{\alpha }$ with $y^{\alpha }\rightarrow y,$ 
\begin{equation*}
\left\vert \sphericalangle \left( \Uparrow _{y^{\alpha }}^{a_{i}^{\alpha
}},\Uparrow _{y^{\alpha }}^{c^{\alpha }}\right) -\sphericalangle \left(
\Uparrow _{y}^{a_{i}},\Uparrow _{y}^{c}\right) \right\vert <\tau \left(
\delta \right) +\tau \left( \frac{1}{\alpha }|r\right) .
\end{equation*}
\end{proposition}

Let $\left( \digamma _{1},p_{\digamma _{1}}\right) $ and $\left( \digamma
_{2},p_{\digamma _{2}}\right) $ be as in the Line Up Proposition (\ref{line
up}). Suppose that $\digamma _{1}$ and $\digamma _{2}$ are constructed via $%
(k,\delta ,r)$-strainers $\left\{ \left( a_{i},b_{i}\right) \right\}
_{i=1}^{k}$ and $\left\{ \left( c_{i},d_{i}\right) \right\} _{i=1}^{k}$, and 
$\mathrm{diam}\left\{ \digamma _{1}\right\} ,\mathrm{diam}\left\{ \digamma
_{2}\right\} \leq d,$ for some $d>0$. Use a superscript index $\alpha $ to
denote approximations of these structures in $M_{\alpha }$. Applying the
previous result to this situation gives us

\begin{corollary}
\label{constant directions} For $x^{\alpha },y^{\alpha }\in \digamma
_{1}^{\alpha }\cap \digamma _{2}^{\alpha },$ 
\begin{equation*}
\left\vert \sphericalangle \left( \Uparrow _{x^{\alpha }}^{a_{i}^{\alpha
}},\Uparrow _{x^{\alpha }}^{c_{j}^{\alpha }}\right) -\sphericalangle \left(
\Uparrow _{y^{\alpha }}^{a_{i}^{\alpha }},\Uparrow _{y^{\alpha
}}^{c_{j}^{\alpha }}\right) \right\vert <\tau \left( \delta \right) +\tau
\left( \frac{1}{\alpha },d|r\right) .
\end{equation*}
\end{corollary}

In the next proposition, we show $p_{\digamma _{1}^{\alpha }}$ and $%
p_{\digamma _{2}^{\alpha }}$ are $\kappa $-lined up provided the
differential $dp_{\digamma _{2}^{\alpha }}$ at some point $x^{\alpha }\in
\digamma _{1}^{\alpha }(\beta _{1})\cap \digamma _{2}^{\alpha }(\beta _{2})$
is almost an isometry when restricted to the horizontal space of $%
p_{\digamma _{1}^{\alpha }}$.

\begin{proposition}
\label{line up sort of} Let $H_{1}^{\alpha }$ be the horizontal distribution
of $p_{\digamma _{1}^{\alpha }}$. Given $\kappa >0,$ there is a $\lambda
_{0}<1$ so that if for some $x\in \digamma _{1}^{\alpha }\cap \digamma
_{2}^{\alpha }$ and some $\lambda \in \left( \lambda _{0},1\right) ,$ $%
\left( dp_{\digamma _{2}^{\alpha }}\right) _{x}$ is $\frac{1}{\lambda }$%
--bilipschitz on the unit sphere of $H_{1}^{\alpha },$ then $p_{\digamma
_{1}^{\alpha }}$ and $p_{\digamma _{2}^{\alpha }}$ are $\kappa $--lined up,
provided $\tau \left( \delta \right) +\tau \left( \frac{1}{\alpha }%
,d|r\right) <\kappa .$
\end{proposition}

\begin{proof}
Let $\left( p_{\digamma _{1}^{\alpha }}\right) ^{i}$ and $\left( p_{\digamma
_{2}^{\alpha }}\right) ^{i}$ be the $i^{th}$--coordinate functions of $%
p_{\digamma _{1}^{\alpha }}$ and $p_{\digamma _{2}^{\alpha }}.$

Since $p_{\digamma _{2}^{\alpha }}$ is a $\tau \left( \delta \right) $%
-almost Riemannian submersion, the hypothesis that $\left( dp_{\digamma
_{2}^{\alpha }}\right) _{x}|_{H_{1}^{\alpha }}$ is $\frac{1}{\lambda }$%
--bilipschitz tells us that at $x,$ the horizontal distributions 
\begin{equation*}
H_{1}^{\alpha }=\mathrm{span}\left\{ \nabla \left( p_{\digamma _{1}^{\alpha
}}\right) ^{i}\right\} _{i}\text{ and }H_{2}^{\alpha }=\mathrm{span}\left\{
\nabla \left( p_{\digamma _{2}^{\alpha }}\right) ^{j}\right\} _{j}
\end{equation*}%
nearly coincide. On the other hand, Corollary \ref{constant directions}
tells us that the angles between the above gradients do not vary by more
than $\tau \left( \delta \right) +\tau \left( \frac{1}{\alpha },d|r\right) $
on $\digamma _{1}^{\alpha }\cap \digamma _{2}^{\alpha }.$ It follows that $%
H_{1}^{\alpha }$ and $H_{2}^{\alpha }$ nearly coincide throughout $\digamma
_{1}^{\alpha }\cap \digamma _{2}^{\alpha }.$ Since both of $p_{\digamma
_{1}^{\alpha }}$ and $p_{\digamma _{2}^{\alpha }}$ are $\tau \left( \delta
\right) $-almost Riemannian submersions, it follows that that $p_{\digamma
_{1}^{\alpha }}$ and $p_{\digamma _{2}^{\alpha }}$ are $\kappa $--lined up,
provided 
\begin{equation*}
\tau \left( \delta \right) +\tau \left( \frac{1}{\alpha },d|r\right) <\kappa
,
\end{equation*}%
and $\lambda _{0}$ is close enough to $1.$
\end{proof}

\mbox{}\vspace{-.5in} 
\begin{wrapfigure}[11]{r}{0.43\textwidth}\centering
\vspace{-12pt}
\includegraphics[scale=.16]{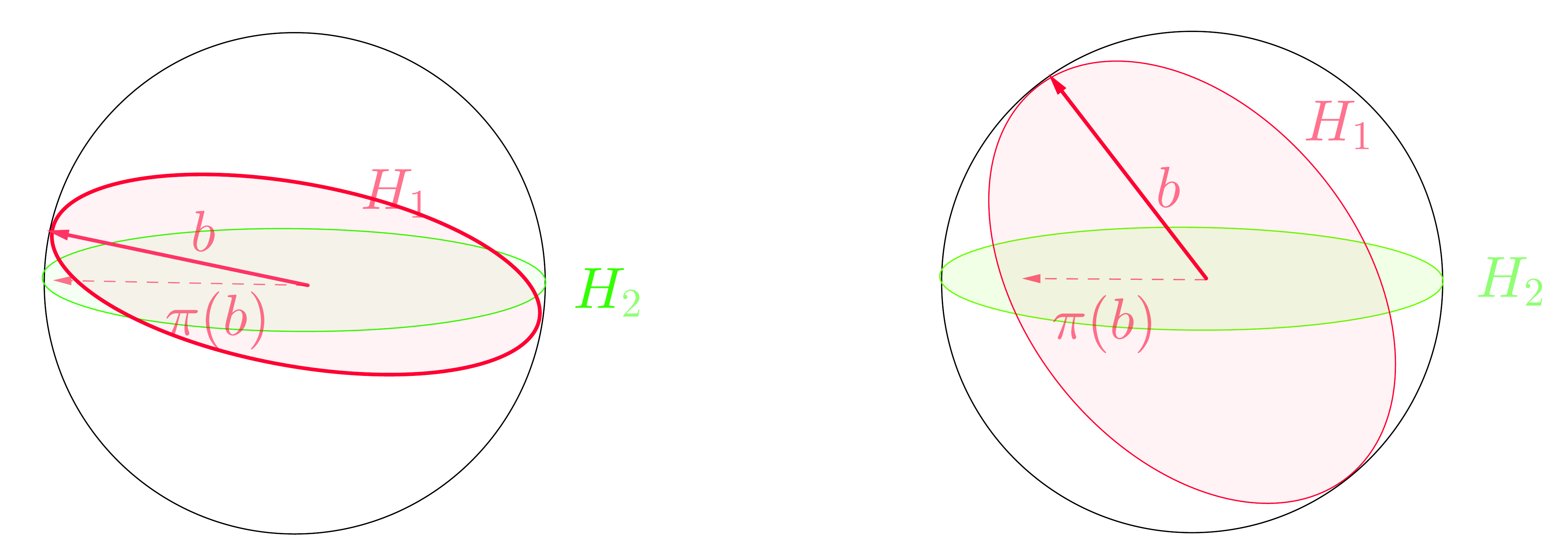}
\vspace{-18pt}
\caption{{\protect\tiny {\ In the figure on the left, the orthogonal
projection of the unit sphere in the subspace $H_1$ to the subspace $H_2$ is
almost an isometry. The fact that this is not the case on the right is
detected by the vector $b$ whose orthogonal projection to $H_2$ has length
significantly less than $1$. }}}
\label{Bilipsorsmall}
\end{wrapfigure}

\begin{proposition}[see Figure \protect\ref{Bilipsorsmall}]
\label{not bilip prop} Let $V$ be an inner product space with $k$%
--dimensional subspaces $H_{1}$ and $H_{2}.$ Let $S_{1}$ be the unit sphere
in $H_{1}.$ Let $B$ be an orthonormal basis for $H_{1},$ and let 
\begin{equation*}
\hspace*{-3in}\pi :V\longrightarrow H_{2}
\end{equation*}%
be orthogonal projection.

For all $\lambda \in \left( \frac{1}{2},1\right) ,$ there is a $\tilde{%
\lambda}\in $ $\left( \lambda ,1\right) $ so that either $\pi |_{S_{1}}$ is $%
\frac{1}{\lambda }$--bilipschitz, or for some $b\in B,$%
\begin{equation*}
\hspace*{-3in}\left\vert \pi \left( b\right) \right\vert <\tilde{\lambda}.
\end{equation*}
\end{proposition}

\begin{proof}
Since $\pi $ is orthogonal projection, it is $1$--Lipshitz on unit vectors,
and $\left\vert \pi \left( b\right) \right\vert \leq 1$ for all $b\in B.$

If the result is false, there is a $\lambda \in \left( \frac{1}{2},1\right) $
and a sequence $\left\{ \tilde{\lambda}_{i}\right\} _{i=1}^{\infty }\subset
\left( \lambda ,1\right) $ with $\tilde{\lambda}_{i}\rightarrow 1$ so that $%
\pi |_{S_{1}}$ is not $\frac{1}{\lambda }$--bilipschitz and%
\begin{equation*}
\tilde{\lambda}_{i}\leq \left\vert \pi \left( b\right) \right\vert \leq 1%
\text{ for all }b\in B.
\end{equation*}%
Sending $i\rightarrow \infty $, it follows that%
\begin{equation}
\left\vert \pi \left( b\right) \right\vert =1\text{ for all }b\in B.
\label{norm 1 on basis}
\end{equation}%
Since $\pi |_{S_{1}}$ is $1$--Lipshitz but not $\frac{1}{\lambda }$%
--bilipshitz, there is a $s\in S_{1}$ so that $\left\vert \pi \left(
s\right) \right\vert <\lambda .$ Combining this with (\ref{norm 1 on basis}%
), we see that for at least two distinct elements of $B,$ say $b_{1}$ and $%
b_{2},$%
\begin{equation}
\left\vert \left\langle \pi \left( b_{1}\right) ,\pi \left( b_{2}\right)
\right\rangle \right\vert >\tau \left( 1-\lambda \right) .
\label{inner products}
\end{equation}%
(\ref{norm 1 on basis}) and (\ref{inner products}) together imply that $\pi
|_{S_{1}}$ is not $1$--Lipschitz. Since $\pi $ is an orthogonal projection,
this is a contradiction.
\end{proof}

The following is a natural deformation of the hypotheses of Plaut's Theorem %
\ref{plauttheorem}.

\begin{definition}
\label{dfn alm sph}A set of $2n$ points $x_{1},y_{1},\ldots ,x_{n},y_{n}$ in
a metric space $Y$ is called $\left( \delta |v\right) $\textbf{--almost
spherical} if $\mathrm{dist}(x_{i},y_{i})>\pi -\delta $ for all $i$ and $%
\mathrm{det}[\cos \mathrm{dist}(x_{i},x_{j})]>v>0$.
\end{definition}

The following result gives a connection between strainers and almost
spherical sets.

\begin{proposition}[Proposition 1.14 in \protect\cite{ProWilh1}]
\label{almost Plauter}Let $X$ have curvature $\geq 1,$ dimension $n,$ and
contain a $\left( \delta |v\right) $--almost spherical set $S$ of $2(m+1)$
points$.$ Then there is an $\left( m+1,\tau \left( \delta |v\right) \right) $%
--global strainer $\left\{ \left( a_{i},b_{i}\right) \right\} _{i=1}^{m+1}$
for $X$.
\end{proposition}

We are ready to prove the Line Up Proposition.

\begin{proof}[Proof of Proposition \protect\ref{line up}]
Let $H_{1}^{\alpha }$ and $H_{2}^{\alpha }$ be the horizontal distributions
of $p_{\digamma _{1}^{\alpha }}$ and $p_{\digamma _{2}^{\alpha }}.$ Let $%
S_{1}$ be the unit sphere in $H_{1}^{\alpha }$. Given $\kappa >0,$ we use
Proposition \ref{line up sort of} to choose $\delta ,d>0,\ \alpha ,$ and $%
\lambda \in \left( \frac{1}{2},1\right) $ so that if $dp_{\digamma
_{2}^{\alpha }}|_{S_{1}}$ is $\frac{1}{\lambda }$--bilipschitz at a point,
then $p_{\digamma _{1}}$ and $p_{\digamma _{2}}$ are $\kappa $--lined up
everywhere. Therefore, if $p_{\digamma _{1}}$ and $p_{\digamma _{2}}$ are
not $\kappa $-lined up, then $dp_{\digamma _{2}^{\alpha }}|_{H_{1}}$ is not $%
\frac{1}{\lambda }$-bilipschitz everywhere on $\digamma _{1}^{\alpha }\cap
\digamma _{2}^{\alpha }.$ By Proposition \ref{not bilip prop}, we then have
that for all $y^{\alpha }\in \digamma _{1}^{\alpha }\cap \digamma
_{2}^{\alpha }$, there is $\tilde{\lambda}\in \left( \lambda ,1\right) $ and
a straining direction, say $\uparrow _{y^{\alpha }}^{a_{1}^{\alpha }}\in
\Sigma _{y^{\alpha }}$, so that 
\begin{equation}
\left\vert dp_{\digamma _{2}^{\alpha }}\left( \uparrow _{y^{\alpha
}}^{a_{1}^{\alpha }}\right) \right\vert <\tilde{\lambda}.  \label{small proj}
\end{equation}%
By further constraining $\delta ,$ $d,$ and $\alpha ,$ and appealing to
Proposition \ref{angggl conv prop}, we see that $\tilde{\lambda}$ can be
chosen independent of $y^{\alpha }\in \digamma _{1}^{\alpha }\cap \digamma
_{2}^{\alpha }$.

It follows from Lemma \ref{angle convergence} that there is a choice of $%
\tau $ that we will call $\tau _{\mathrm{vol}}$ so that 
\begin{equation}
\det \left( \left\{ \cos \sphericalangle \left( \uparrow _{y^{\alpha
}}^{c_{i}^{\alpha }},\uparrow _{y^{\alpha }}^{c_{j}^{\alpha }}\right)
\right\} _{i,j}\right) >1-\tau _{\mathrm{vol}}\left( \delta \right) .
\label{k d vol}
\end{equation}%
Writing a $\left( k+1\right) \times \left( k+1\right) $ matrix in the form 
\begin{equation*}
\left( 
\begin{array}{cc}
k\times k & k\times 1 \\ 
k\times 1 & 1\times 1%
\end{array}%
\right)
\end{equation*}%
and interpreting the following determinant as the square of the $(k+1)$%
-volume of the parallelepiped determined by $\left\{ \uparrow _{y^{\alpha
}}^{c_{i}^{\alpha }}\right\} _{i=1}^{l}\cup $ $\left\{ \uparrow _{y^{\alpha
}}^{a_{1}^{\alpha }}\right\} $ embedded in $\mathbb{R}^{k+1}$, we conclude,
using (\ref{small proj}) and (\ref{k d vol}), that

\begin{equation*}
\det \left( 
\begin{array}{cc}
\left\{ \cos \sphericalangle \left( \uparrow _{y^{\alpha }}^{c_{i}^{\alpha
}},\uparrow _{y^{\alpha }}^{c_{j}^{\alpha }}\right) \right\} _{i,j} & 
\left\{ \cos \sphericalangle \left( \uparrow _{y^{\alpha }}^{a_{1}^{\alpha
}},\uparrow _{y^{\alpha }}^{c_{j}^{\alpha }}\right) \right\} _{1,j} \\ 
\left\{ \cos \sphericalangle \left( \uparrow _{y^{\alpha }}^{c_{i}^{\alpha
}},\uparrow _{y^{\alpha }}^{a_{1}^{\alpha }}\right) \right\} _{i,1} & \cos
\sphericalangle \left( \uparrow _{y^{\alpha }}^{a_{1}^{\alpha }},\uparrow
_{y^{\alpha }}^{a_{1}^{\alpha }}\right)%
\end{array}%
\right) >1-\tau _{\mathrm{vol}}\left( \delta \right) -\tilde{\lambda}^{2}.
\end{equation*}%
Letting $\alpha \rightarrow \infty $ gives us 
\begin{equation*}
\det \left( 
\begin{array}{cc}
\left\{ \cos \sphericalangle \left( \uparrow _{y}^{c_{i}},\uparrow
_{y}^{c_{j}}\right) \right\} _{i,j} & \left\{ \cos \sphericalangle \left(
\uparrow _{y}^{a_{1}},\uparrow _{y}^{c_{j}}\right) \right\} _{1,j} \\ 
\left\{ \cos \sphericalangle \left( \uparrow _{y}^{c_{i}},\uparrow
_{y}^{a_{1}}\right) \right\} _{i,1} & \cos \sphericalangle \left( \uparrow
_{y}^{a_{1}},\uparrow _{y}^{a_{1}}\right)%
\end{array}%
\right) \geq 1-\tau _{\mathrm{vol}}\left( \delta \right) -\tilde{\lambda}%
^{2}.
\end{equation*}%
So for each $y\in \digamma _{1}\cap \digamma _{2}$, $\left\{ \left( \uparrow
_{y}^{c_{i}},\uparrow _{y}^{d_{i}}\right) \right\} _{i=1}^{k}$ together with 
$\left( \uparrow _{y}^{a_{1}},\uparrow _{y}^{b_{1}}\right) $ forms a $\left(
\delta |1-\tau _{\mathrm{vol}}\left( \delta \right) -\tilde{\lambda}%
^{2}\right) $-almost spherical set $S\subset \Sigma _{y}$ of $2(k+1)$ points.

To complete the proof, we further constrain $\delta $ so that 
\begin{equation*}
\tau _{\mathrm{plaut}}\left( \delta |1-\tau _{\mathrm{vol}}\left( \delta
\right) -\tilde{\lambda}^{2}\right) <\varepsilon _{1},
\end{equation*}%
where $\tau _{\mathrm{plaut}}$ is the $\tau $ of Proposition \ref{almost
Plauter}. Thus Proposition \ref{almost Plauter} gives that each point of $%
\digamma _{1}\cap \digamma _{2}$ is $\left( k+1,\varepsilon _{1}\right) $%
--strained, as claimed.
\end{proof}

\subsection{Improving the Perelman-Plaut Cover}

\begin{wrapfigure}[13]{l}{0.3\textwidth}\centering
\vspace{-1pt}
\includegraphics[scale=0.12]{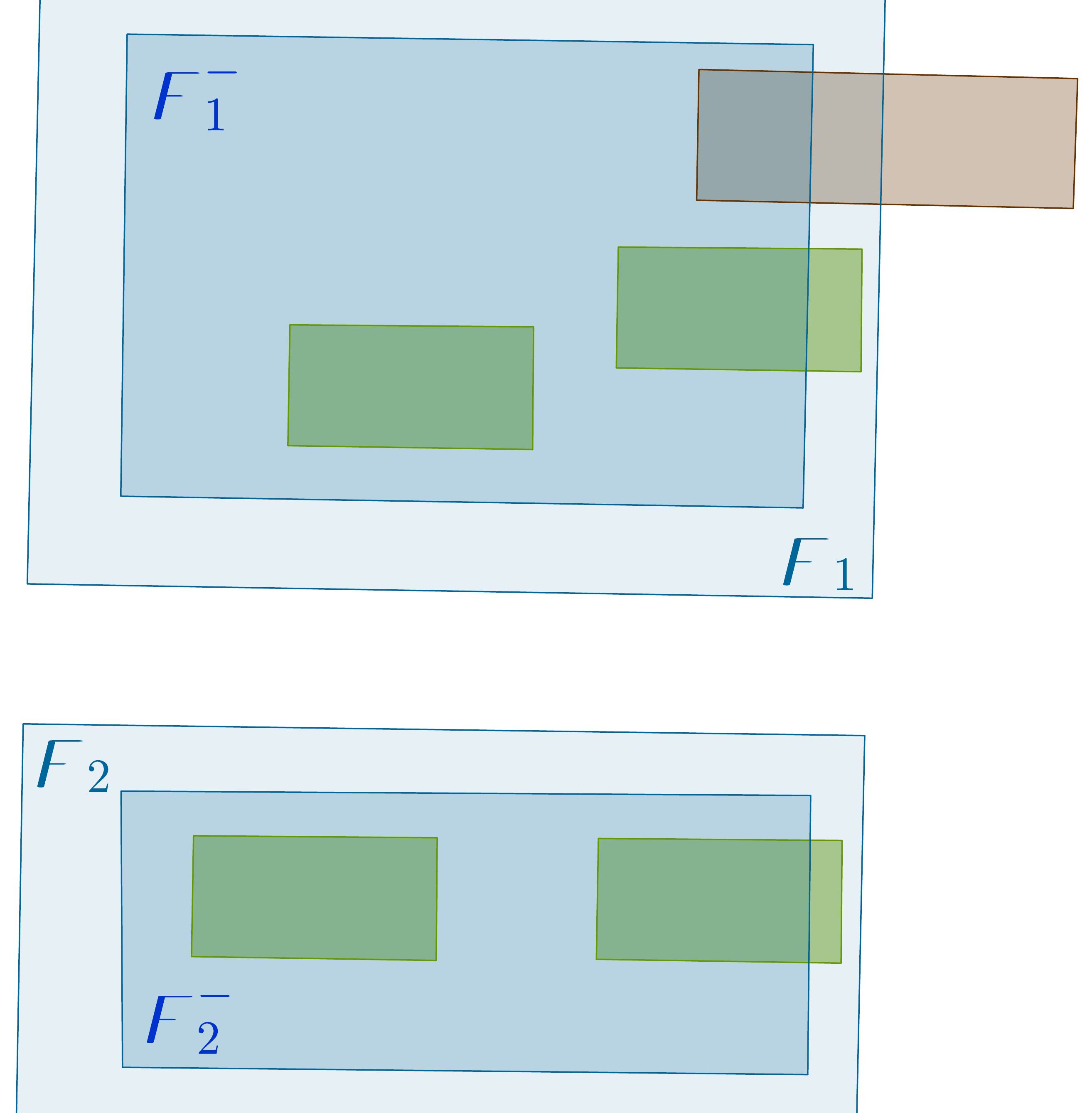}
\vspace{-1pt}
\caption{{\protect\tiny {The collection of green $1$--framed sets thoroughly
respects the collection of $0$--framed sets $\{ \digamma_1, \digamma_2\}$;
however, the brown set does not thoroughly respect $\{ \digamma_1,
\digamma_2\}$. } } }
\label{thor resp}
\end{wrapfigure}
The next step in the proof of Limit Lemma \ref{cover of X lemma} is to
modify the Perelman-Plaut Covering Lemma (\ref{Ms. Pacman limit lemma}) by
constraining the way that the $1$-- and $2$--framed sets intersect the $0$%
--framed sets and imposing further geometric constraints on the $1$--framed
sets. The results are stated in Lemma \ref{art Ms Pacman lemma} (below). The
constraint on intersections is called \textbf{thorough respect}.

\begin{definition}[see Figure \protect\ref{thor resp}]
\label{nat resp 0 str dfn} Let $\mathcal{U}=\mathcal{U}^{0}\cup \mathcal{U}%
^{1}\cup \mathcal{U}^{2}$ be a framed collection. For $0\leq k<l\leq 2$,
suppose that$\ $for each $U^{l}\in \mathcal{U}^{l}$ there is at most one set 
$\digamma \in \mathcal{U}^{k}$ such that 
\begin{equation}
\hspace*{2in} U^{l}\cap \digamma ^{-}\neq \emptyset .  \label{oneinter cond}
\end{equation}%
In this event, suppose further that 
\begin{equation}
\hspace*{1.9in}U^{l}\subset \digamma  \label{conta if intercond}
\end{equation}%
and $U^{l}$ intersects $\digamma $ respectfully. When this occurs for a
single $l>k,$ we say that $\mathcal{U}^{l}$ \textbf{thoroughly respects} $%
\mathcal{U}^{k}.$ When this occurs for all $l>k,$ we say that $\mathcal{U}%
^{k}$ is \textbf{thoroughly respectful}.
\end{definition}

\begin{wrapfigure}[9]{r}{0.3\textwidth}\centering
\vspace{-9pt}
\includegraphics[scale=0.16]{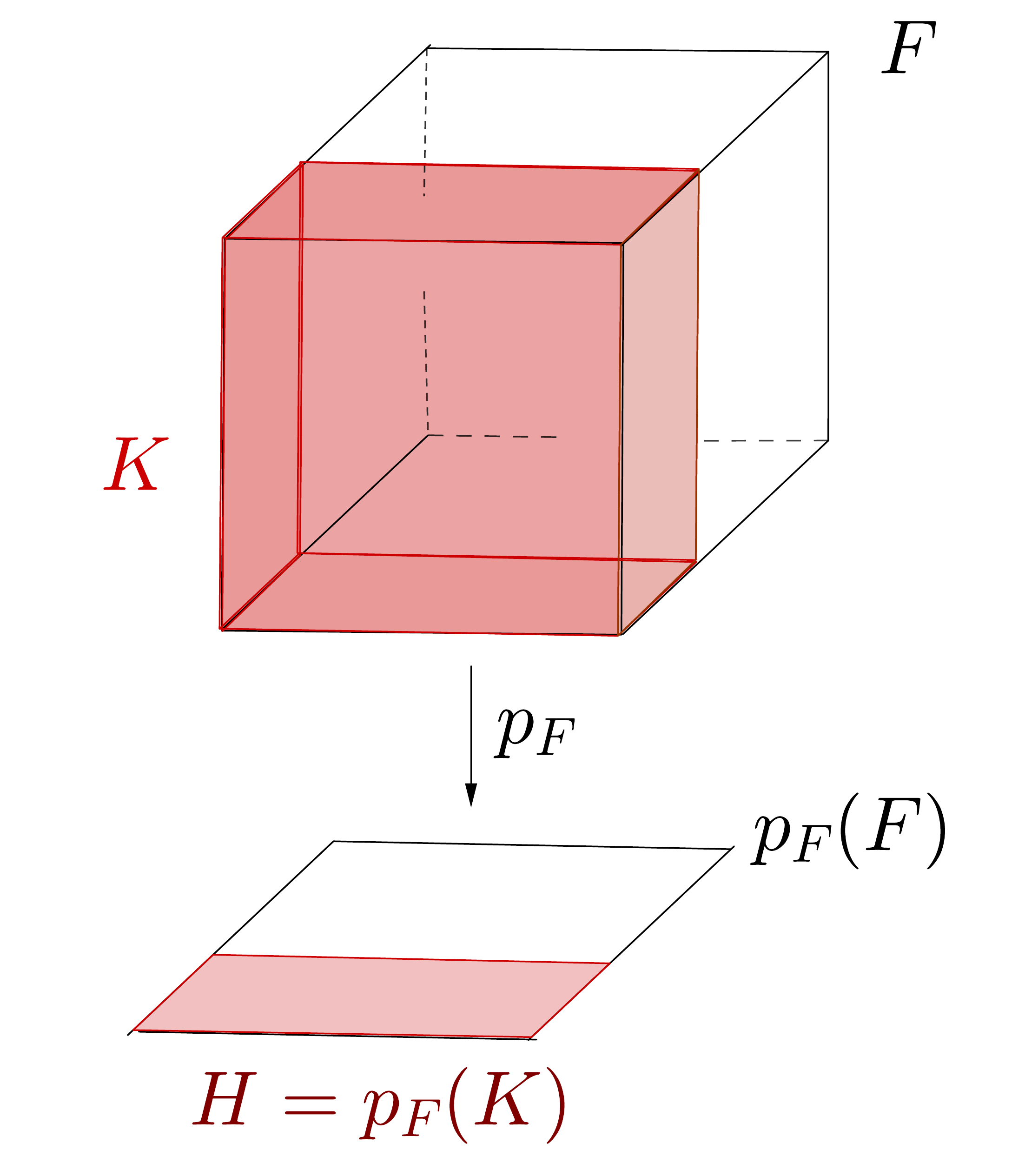}
\vspace{-9pt}
\caption{ {\protect\tiny {A framed half space} }}
\label{half space fig}
\end{wrapfigure}

The main result of this subsection asserts that the cover of Lemma \ref{Ms.
Pacman limit lemma} can be chosen to have a thoroughly respectful $0$%
--stratum. To establish this, we will modify some of the framed sets as
follows.

\begin{definition}[see Figure \protect\ref{half space fig}]
\label{perelman half space}Let $\digamma $ and $p_{\digamma }:\digamma
\longrightarrow \mathbb{R}^{k}$ be as in Lemma \ref{framedsetprop}. Given a
half space $H\subset \mathbb{R}^{k}$ so that $H\cap p_{\digamma }\left(
\digamma \right) $ is diffeomorphic to an open ball, we call 
\begin{equation*}
\hspace*{-2in}K\equiv p_{\digamma }^{-1}\left( H\cap p_{\digamma }\left(
\digamma \right) \right)
\end{equation*}%
a \textbf{framed half space. }
\end{definition}

\mbox{}\vspace{-.1in} 
\begin{wrapfigure}[18]{l}{0.3\textwidth}\centering
\vspace{-5pt}
\includegraphics[scale=0.16]{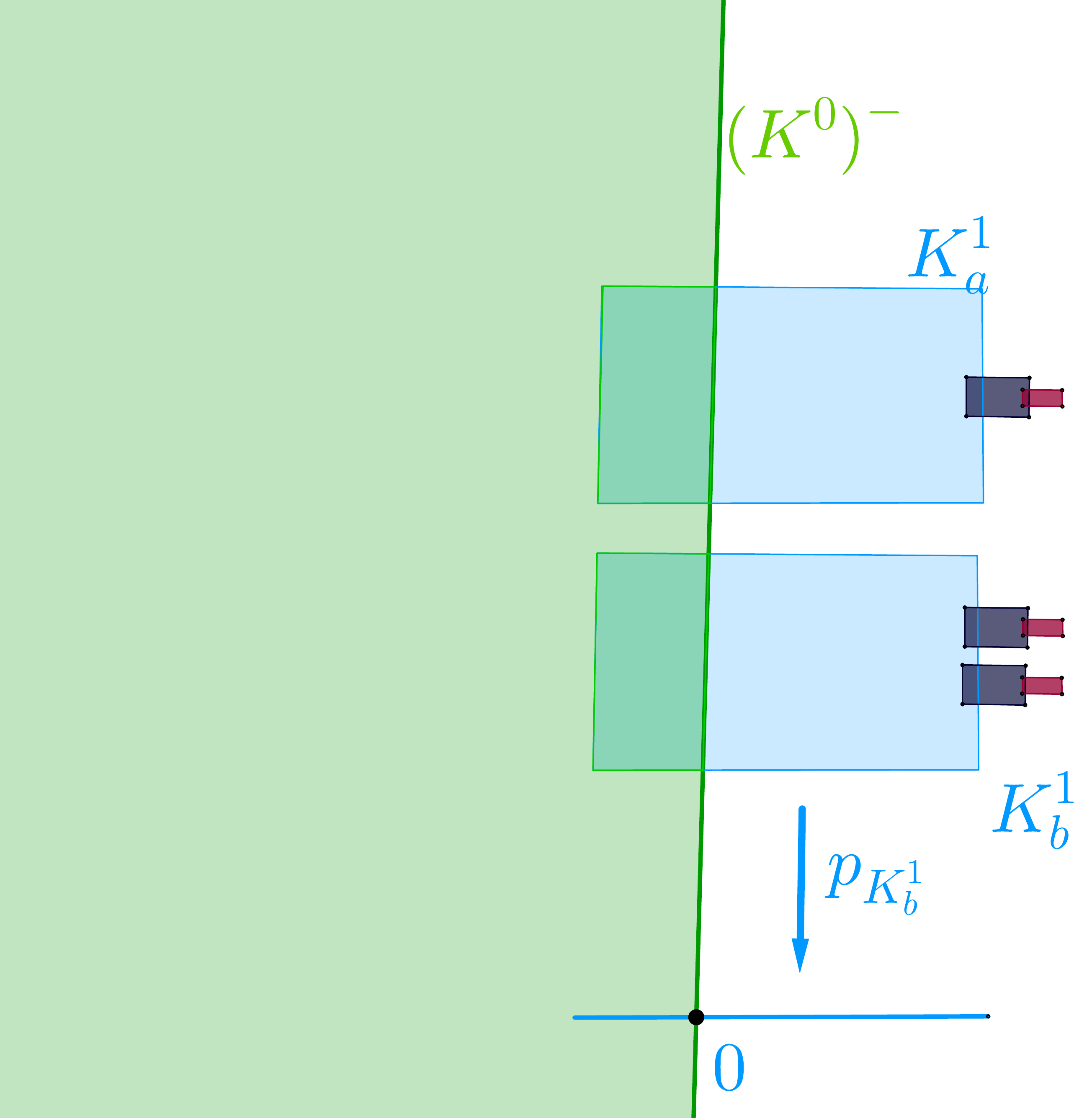}
\vspace{5pt}
\caption{{\protect\tiny {This depicts how elements $K^{1}_a$ and $K_{b}^1$
of $\mathcal{K}^1$ intersect an element $K^0$ of $\mathcal{K}^0$, as
described in Lemma \protect\ref{art Ms Pacman lemma}.  }}}
\label{Covering Lemma 0+}
\end{wrapfigure}
\vspace{-.4in}

\begin{lemma}
\label{art Ms Pacman lemma} There is a framed cover $\mathcal{K}=\mathcal{K}%
^{0}\cup \mathcal{K}^{1}\cup \mathcal{K}^{2}$ that has all of the properties
of the cover $\mathcal{Z}$ of Lemma \ref{Ms. Pacman limit lemma}. In
addition, $\mathcal{K}^{0}$ is thoroughly respectful, and \vspace{4pt}

\noindent 1. For each $\digamma \in \mathcal{K}^{1},$ there is a $\beta
_{\digamma }>0$ so that 
\begin{equation}
\hspace*{2in} p_{\digamma }\left( \bar{\digamma}\left( \beta \right) \right)
=\left[ 0,\beta _{\digamma }\right] .  \label{1 framed image}
\end{equation}

\noindent 2. Set $\Omega _{0}:=\cup \mathcal{K}^{0}.$ If $\digamma \left(
\beta _{\digamma }\right) ,\tilde{\digamma}\left( \tilde{\beta}_{\digamma
}\right) \in \mathcal{K}^{1}$, $\digamma \cap \Omega _{0}^{-}\neq
\varnothing ,$ and $\tilde{\digamma}\cap \Omega _{0}^{-}=\varnothing ,$ then 
\begin{equation}
\hspace*{2in} \beta _{\digamma }\geq 100\tilde{\beta}_{\digamma }\text{,%
\label{much smaller away}}
\end{equation}%
and if $\left( Z,r_{Z}\right) $ is the element of $\mathcal{K}^{0}$ so that $%
\digamma \cap Z\neq \emptyset ,$ then 
\begin{equation}
\hspace*{2in} p_{\digamma }\left( \bar{\digamma}\cap \bar{Z}^{-}\right)
=\left\{ 0\right\} ,\text{ }  \label{0 1 interface}
\end{equation}%
and $p_{\digamma }$ is lined up with $r_{Z}.$

\noindent 3. If $\digamma _{l}\in \mathcal{K}^{1}$ and $\overline{\digamma
_{l}\left( \beta _{l}\right) }\cap \Omega _{0}^{-}\neq \varnothing ,$ then 
\begin{equation}
\digamma _{l}\left( \beta _{l}\right) \cap
\bigcup\limits_{i=1}^{l-1}\digamma _{i}\left( \beta _{i}^{+}\right)
=\varnothing  \label{disj beggining}
\end{equation}%
and 
\begin{equation}
p_{\digamma _{l}}^{-1}\left[ 0,\frac{\beta _{l}}{2}\right] \cap
\bigcup\limits_{i=l+1}^{L}\digamma _{i}\left( \beta _{i}^{+}\right)
=\varnothing ,  \label{disj end}
\end{equation}%
where $\mathcal{K}^{1}=\left\{ \digamma _{i}\left( \beta _{i}^{+}\right)
\right\} _{i=1}^{L}.$
\end{lemma}

\begin{proof}
For $\delta _{1}$ and $\delta _{2}$ as in (\ref{grown up delta choices}), let%
\begin{eqnarray}
X^{0} &=&\{p\in X\mid p\text{ is not }(1,\delta _{1})\text{--strained}\}, 
\notag \\
X^{1} &=&\{p\in X\mid p\text{ is }(1,\delta _{1})\text{--strained, but not }%
(2,\delta _{2})\text{--strained}\},  \notag \\
X^{2} &=&\{p\in X\mid p\text{ is }(2,\delta _{2})\text{--strained, but not }%
(4,\delta )\text{--strained}\},\text{ and }  \notag \\
X^{4} &=&\{p\in X\mid p\text{ is }(4,\delta )\text{--strained}\}.
\label{dfn of Xi}
\end{eqnarray}

Let 
\begin{equation*}
\mathcal{Z}=\mathcal{Z}^{0}\cup \mathcal{Z}^{1}\cup \mathcal{Z}^{2}
\end{equation*}%
be the framed cover of Lemma \ref{Ms. Pacman limit lemma}, and set 
\begin{equation*}
\mathcal{K}^{0}\equiv \mathcal{Z}^{0}.
\end{equation*}%
For $l\in \left\{ 1,2\right\} ,$ when we construct $\mathcal{K}^{l}\equiv
\left\{ \digamma _{i}^{l}\left( \beta _{i}^{l}\right) \right\} ,$ we require
that the $\beta _{i}^{l}$ be small enough so that (\ref{oneinter cond}) and (%
\ref{conta if intercond}) hold with $k=0$. This is possible because
Inequality (\ref{diam small}) allows us to choose our framed sets to be as
small as we please, and we have a fixed collection, $\mathcal{Z}^{0},$ whose
closures are disjoint.

To each $Z\in \mathcal{Z}^{0}$ and each $x\in \partial \bar{Z}^{-}\cap
X^{l}, $ we apply Corollary \ref{cor expli k+1 framing}. This yields a
collection of framed half spaces $\digamma _{\mathrm{near}}^{l}\equiv
\left\{ K_{\mathrm{near,}i}^{l}\left( \beta _{\mathrm{near,}i}^{l}\right)
\right\} $ that cover $X^{l}\cap \partial \Omega _{0}^{-},$ and satisfy (\ref%
{k+1 submer dfn}). It follows that if $\digamma \in \mathcal{K}_{\mathrm{near%
}}^{l}\mathcal{\ }$intersects $Z\in \mathcal{Z}^{0}$, then $r_{Z}$ and $%
p_{\digamma }\ $are lined up.

For each $\digamma \in \mathcal{K}_{\mathrm{near}}^{1},$ we post compose
with an appropriate isometry to arrange that $p_{K}\left( \bar{\digamma}%
\right) =\left[ 0,\beta _{\digamma }\right] .$ In particular, $\mathcal{K}_{%
\mathrm{near}}^{1}$ then satisfies (\ref{1 framed image}). Since our $1$%
--framed sets respect the radial functions of our $0$--framed sets, we can
also arrange for (\ref{0 1 interface}) to hold.

Our collection of half spaces now satisfies Definition \ref{nat resp 0 str
dfn} but is unlikely to cover $\left( X^{1}\cup X^{2}\right) \setminus
\Omega ^{0}.$ To remedy this, we repeat the procedure of the proof of Lemma %
\ref{Ms. Pacman limit lemma} with the extra requirement that if the
additional framed sets are $\digamma _{j}^{l}\left( \beta _{j}^{l}\right) ,$
then 
\begin{equation}
\max_{j}\left\{ \beta _{j}^{l}\right\} \leq \frac{1}{100}\min_{i}\left\{
\beta _{\mathrm{near,}i}^{l}\right\} .  \label{big near small far}
\end{equation}

Call the resulting framed cover $\mathcal{K}=\mathcal{K}^{0}\cup \mathcal{K}%
^{1}\cup \mathcal{K}^{2}.$ Then $\mathcal{K}^{0}$ is thoroughly respectful
and disjoint, and $\mathcal{K}^{1}$ and $\mathcal{K}^{2}$ are fiber
swallowing, vertically separated, and $\kappa $--lined up. Inequality (\ref%
{big near small far}) ensures that (\ref{much smaller away}) holds. Thus
Part 2 holds. To prove Part 1, for each $\digamma \in \mathcal{K}^{1},$ we
post compose $p_{\digamma }$ with an appropriate isometry to arrange that $%
p_{\digamma }\left( \bar{\digamma}\right) =\left[ 0,\beta _{\digamma }\right]
$ for all $\digamma \in \mathcal{K}^{1}.$

To prove (\ref{disj beggining}), order $\mathcal{K}^{1}:=\left\{ \digamma
_{j}\left( \beta _{j}\right) \right\} _{j}$ so that if $i\leq l,$ then $%
\beta _{i}\geq \beta _{l}$. Suppose that $\overline{\digamma _{l}\left(
\beta _{l}\right) }\cap \Omega _{0}^{-}\neq \varnothing $ and $\digamma
_{l}\left( \beta _{l}\right) \cap \digamma _{i}\left( \beta _{i}^{+}\right)
\neq \varnothing $ for some $i\leq l-1.$ Since $\beta _{i}\geq \beta _{l},$
it follows from Part 2 that $\overline{\digamma _{i}\left( \beta _{i}\right) 
}\cap \Omega _{0}^{-}\neq \varnothing .$ Since $\mathcal{K}^{0}$ is
thoroughly respectful, it follows that there are $Z_{l},Z_{i}\subset 
\mathcal{K}^{0}$ so that%
\begin{equation}
\digamma _{l}\left( \beta _{l}\right) \subset Z_{l}\text{ and }\digamma
_{i}\left( \beta _{i}\right) \subset Z_{i}.  \label{DIFFERENT CONTAIN}
\end{equation}%
Taking into account that $\digamma _{l}\left( \beta _{l}\right) \cap
\digamma _{i}\left( \beta _{i}^{+}\right) \neq \varnothing $ and that the
collection $\mathcal{K}^{0}$ has pairwise disjoint closures, it follows,
after a possible further restriction on the $\beta _{j},$ that $%
Z_{l}=Z_{i}=:Z.$ From the thorough respect condition, we conclude 
\begin{equation*}
\digamma _{i}\left( \beta _{i}\right) \cup \digamma _{l}\left( \beta
_{l}\right) \subset Z.
\end{equation*}%
Via Part 2, we have that on their common domains, $r_{Z},$ $p_{\digamma
_{i}},$ and $p_{\digamma _{l}}$ are lined up, and 
\begin{equation*}
p_{\digamma _{i}}\left( \partial \bar{Z}^{-}\right) =p_{\digamma _{l}}\left(
\partial \bar{Z}^{-}\right) =0.
\end{equation*}%
This, together with the fiber swallowing condition, Equation (\ref{1 framed
image}), and the fact that $p_{\digamma _{l}}$ and $p_{\digamma _{i}}$are $%
\kappa $--lined up, $\varkappa $--almost Riemannian submersions gives us 
\begin{equation*}
\digamma _{l}\left( \beta _{l}^{-},\frac{\beta _{l}}{2}\right) \subset
\digamma _{i}\left( \beta _{i}\right) .
\end{equation*}%
Combining this with (\ref{real grad pres}), we see that we may delete such a 
$\digamma _{l}$ from our collection without losing the property that $%
\mathcal{K}^{0}\cup \mathcal{K}^{1}$ covers $X^{1}.$ After performing these
deletions, $\mathcal{K}^{1}$ satisfies (\ref{disj beggining}). Via a similar
argument that also exploits (\ref{much smaller away}), we can remove
redundant sets from $\mathcal{K}^{1}$ until (\ref{disj end}) holds.
\end{proof}

Throughout the remainder of the paper, we set 
\begin{equation*}
\Omega _{0}:=\cup \mathcal{K}^{0}.
\end{equation*}%
We call $\mathcal{K}^{0}$ and $\mathcal{K}^{1}$ the natural $0$ and $1$
strata of $X.$

\section{Tools to Manipulate a Framed Cover\label{tools sect}}

In this section, we develop tools\ used throughout the remainder of the
paper to prove Limit Lemma \ref{cover of X lemma}. In Subsection \ref{atrf
framed subset}, we introduce \textbf{artificially framed sets} which we
later use to construct a cover with the disjointness properties required by
Limit Lemma \ref{cover of X lemma} (see Figures \ref{art fr basic} and \ref%
{art fr sep}). In Subsection \ref{creating resp}, we develop tools to deform
framed covers to respectful framed covers, and in Subsection \ref{quant
openness}, we establish some tools to verify disjointness.

\addtocounter{algorithm}{1}

\subsection{Artificially Framed Sets\label{atrf framed subset}}

In this subsection, we introduce \textbf{artificially framed sets}\textit{. }%
They are the main tool that we will use to make the framed cover of Limit
Lemma \ref{cover of X lemma} have the required disjointness properties.

\begin{definition}
\label{artif str dfn} Let $\left( \digamma \left( \beta ,\nu \right)
,p_{\digamma }\right) $ be an $l$-framed set. Let $q_{S}:\digamma (\beta
,\nu )\rightarrow \mathbb{R}^{l}$ be a submersion that is $\kappa $-lined up
with $p_{\digamma }$. Let $U_{S}\subset q_{S}\left( \digamma (\beta ,\nu
)\right) \subset \mathbb{R}^{l}$ be an open set that is diffeomorphic to an
open ball, and let $\tilde{\nu}$ be a number in $\left( 0,\nu \right] .$ We
say that $(S,q_{S},U_{S})$ is \textbf{artificially $l$-framed subordinate to}
$(\digamma ,p_{\digamma })$ provided 
\begin{equation}
S=q_{S}^{-1}\left( U_{S}\right) \cap \digamma (\beta ,\tilde{\nu}).
\label{dfn of S 1st}
\end{equation}
\end{definition}

A framed half space is an example of an artificially framed set (see
Definition \ref{perelman half space} and Figure \ref{half space fig}). In
the previous definition, the artificially framed set $S$ has the same
framing dimension as $(\digamma ,p_{\digamma }).$ We also need this notion
for smaller framing dimensions. This is accomplished via the following
definition.

\begin{definition}[see Figures \protect\ref{art fr basic} and \protect\ref%
{atrf framed 2}]
\label{artf k framed dfn} Let $\left( \digamma ,p_{\digamma }\right) $ be an 
$l$-framed set and let $(S,q_{S},U_{S})$ be an artificially $l$-framed
subordinate to $(\digamma ,p_{\digamma })$. Let $L_{S}\subset \mathbb{R}^{l}$
be a $k$-dimensional affine subspace and $\mathfrak{s}\subset L_{S}$
diffeomorphic to an open ball. If $U_{S}$ has the form 
\begin{equation*}
U_{S}=\left\{ \left. w+u\text{ }\right\vert \text{ }w\in \mathfrak{s},u\perp 
\mathfrak{s},\text{ and }\left\vert u\right\vert <b_{S}\right\}
\end{equation*}%
for some $b_{S}>0,$ then we say that $S$ (or $(S,q_{S},U_{S},\mathfrak{s})$)
is \textbf{artificially} $k$\textbf{-framed}.
\end{definition}

Let $\pi _{\mathfrak{s}}:\mathbb{R}^{l}\longrightarrow L_{S}=\mathbb{R}^{k}$
be orthogonal projection. To justify the name artificially $k$-framed, we
note that after setting 
\begin{equation}
p_{S}:=\pi _{\mathfrak{s}}\circ q_{S},  \label{pee-ess}
\end{equation}%
$\left( S,p_{S}\right) $ becomes a trivial vector bundle. Thus an
artificially $k$-framed set $(S,q_{S},\mathfrak{s})$ is $k$-framed provided
it has a fiber exhaustion function. This is a consequence of the next result.

\begin{wrapfigure}[13]{r}{0.4\textwidth}\centering
\vspace{-35pt}
\includegraphics[scale=0.2]{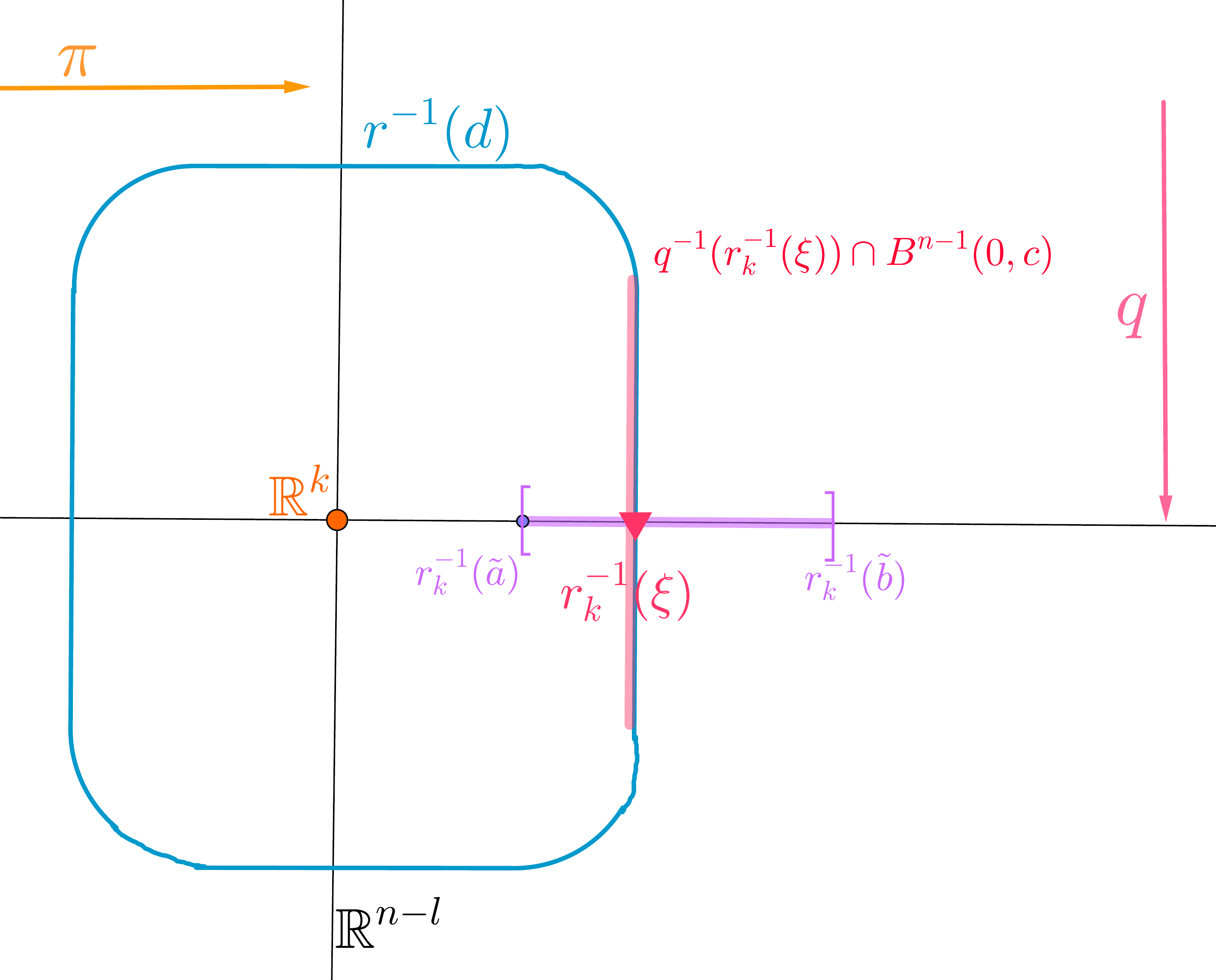}
\vspace{-6pt}
\caption{Depicts Lemma \protect\ref{sucky euclidean lemma} in the case when $%
k=0$.}
\label{fib exh lemma pict}
\end{wrapfigure}

\begin{proposition}[see Figure \protect\ref{fib exh lemma pict} for the
Euclidean analog]
\label{artf framed is framed prop} Let $(\digamma ,p_{\digamma },r_{\digamma
})$ be a $l$--framed set. Suppose that $(S,q_{S},U_{S},\mathfrak{s})$ is an
artificially $k$-framed set that is subordinate to $\digamma $ with $%
S=q_{S}^{-1}\left( U_{S}\right) \cap \digamma (\beta ,\tilde{\nu}),$ and let 
\begin{equation*}
r_{\mathfrak{s}}:U_{S}\longrightarrow \left[ a,b\right) \subset \mathbb{R}
\end{equation*}%
be a fiber exhaustion function for $\pi _{\mathfrak{s}}:\mathbb{R}%
^{l}\longrightarrow L_{S}.$

Given any $\tilde{a},\tilde{b}\in \left( a,b\right) $ and $c\in \left( 0,%
\tilde{\nu}\right) $, there is a fiber exhaustion function $r:S\rightarrow 
\mathbb{R}$ for $p_{S}:=\pi _{\mathfrak{s}}\circ q_{S}$ so that for all $\xi
\in \left[ \tilde{a},\tilde{b}\right] ,$%
\begin{equation*}
q_{S}^{-1}\left( r_{\mathfrak{s}}^{-1}\left( \xi \right) \right) \cap
r_{\digamma }^{-1}\left[ 0,c\right] \subset r^{-1}(d)
\end{equation*}%
for some $d\in \mathbb{R}.$
\end{proposition}

This follows by combining the trivialization for $S$ with the corresponding
result for Euclidean space, which we state here.

\begin{lemma}
\label{sucky euclidean lemma} Let $q:\mathbb{R}^{n}=\mathbb{R}^{l}\times 
\mathbb{R}^{n-l}\longrightarrow \mathbb{R}^{l}\text{ and} $ 
\begin{equation*}
\hspace*{-2.5in} \pi :\mathbb{R}^{l}=\mathbb{R}^{k}\times \mathbb{R}%
^{l-k}\longrightarrow \mathbb{R}^{k}
\end{equation*}%
be orthogonal projections, and let 
\begin{equation*}
\hspace*{-2.5in} r_{k}:\mathbb{R}^{k}\times \mathbb{R}^{l-k}\longrightarrow %
\left[ 0,\infty \right)
\end{equation*}%
be any fiber exhaustion function for $\pi .$ Given any $0<\tilde{a}<\tilde{b}%
\ $and any $c>0$ there is a fiber exhaustion function $r:\mathbb{R}%
^{n}\longrightarrow \mathbb{R} $ for $p:=\pi \circ q:\mathbb{R}^{n}=\mathbb{R%
}^{k}\times \mathbb{R}^{l-k}\times \mathbb{R}^{n-l}\longrightarrow \mathbb{R}%
^{k} $ so that for all $\xi \in \left[ \tilde{a},\tilde{b}\right] ,$ 
\begin{equation*}
q^{-1}\left( r_{k}^{-1}\left( \xi \right) \right) \cap B^{n-l}\left(
0,c\right) \subset r^{-1}\left( d\right)
\end{equation*}%
for some $d\in \mathbb{R}.$
\end{lemma}

\addtocounter{algorithm}{1}

\subsection{Creating Respect\label{creating resp}}

\begin{wrapfigure}[13]{l}{0.4\textwidth}\centering
\vspace{-4pt}
\includegraphics[scale=0.16]{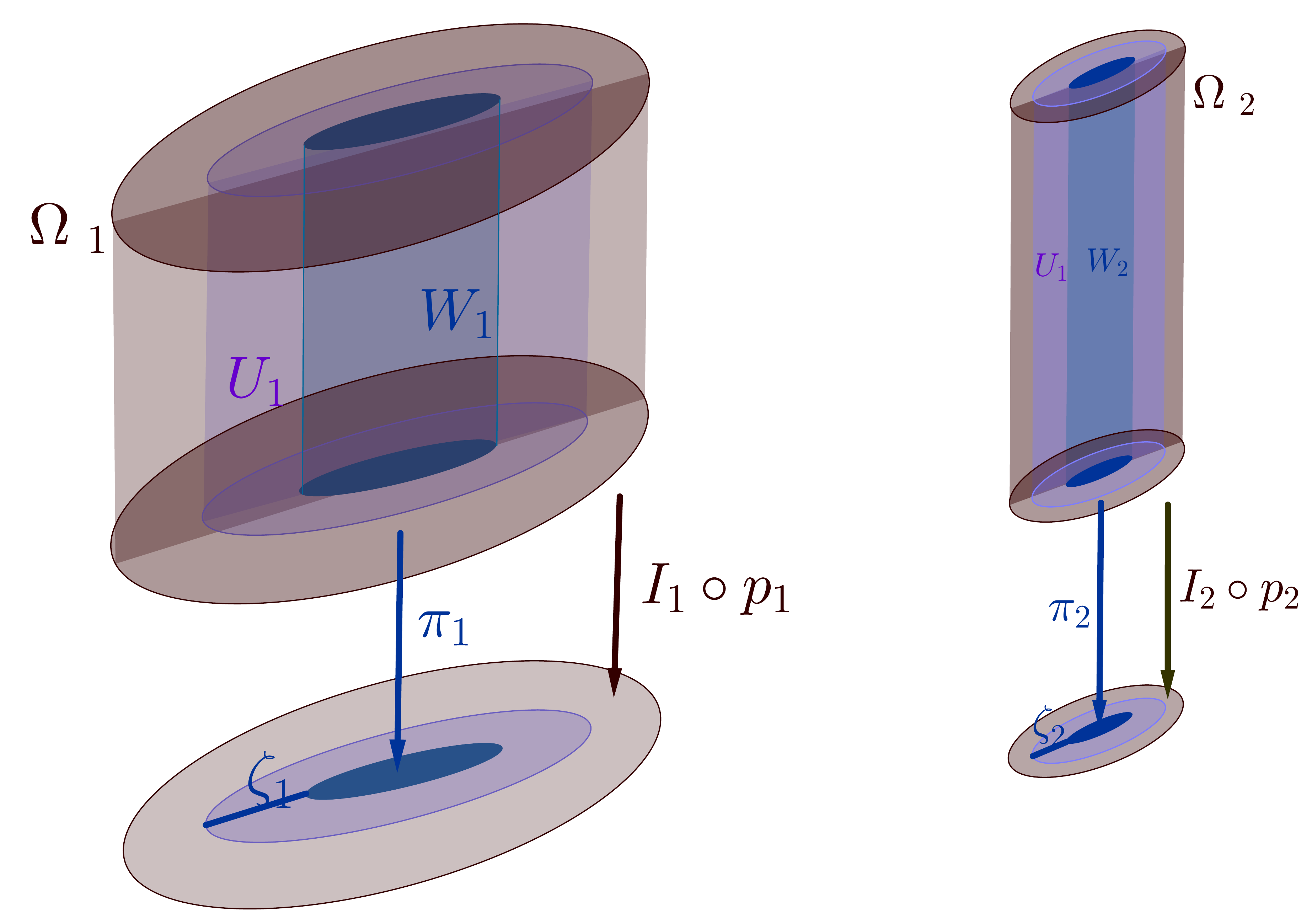}
\caption{{\protect\tiny {\ Suppose that $\protect\pi _{i}$ and $p_{i}$ are $%
\protect\kappa _{i}$--lined up, and we want to glue the submersions as in
the Submersion Deformation Lemma. Since $\frac{\mathrm{Diam}(\Omega _{1})}{%
\protect\zeta _{1}}$ is much smaller than $\frac{\mathrm{Diam}(\Omega _{2})}{%
\protect\zeta _{2}}$, $\protect\kappa _{1}$ is allowed to be much larger
than $\protect\kappa _{2}$. }}}
\label{subm def fig}
\end{wrapfigure}
In this subsection, we outline the tools that will be used to deform a
framed cover to a respectful framed cover. The first step is to glue
together submersions of neighboring framed sets. This is possible because of
the $\kappa $--lined up condition and the following lemma, which is mostly a
consequence of Lemma 7.4 in \cite{ProWilh1}.

\begin{lemma}[Submersion Deformation Lemma]
\label{const rank def} For every $C>1$ and $n\in \mathbb{N},$ there are $%
\kappa _{0},\varkappa _{0}>0$ with the following property. Let $M$ be a
Riemannian $n$--manifold. Let $W\subset U\subset \Omega \subset M$ be three
nonempty, open, precompact sets so that 
\begin{equation}
\hspace*{2.5in}\mathrm{diam}\left( \Omega \right) \leq C\mathrm{dist(}%
W,\Omega \setminus U\mathrm{)},\text{ }  \label{Sub def hyp}
\end{equation}%
and let $p,\pi :\Omega \longrightarrow \mathbb{R}^{k}$ be $\kappa $--lined
up, $\varkappa $--almost Riemannian submersions for some $\kappa \in \left(
0,\kappa _{0}\right) $ and $\varkappa \in \left( 0,\varkappa _{0}\right) $.
Then there is an isometry $I:\mathbb{R}^{k}\longrightarrow \mathbb{R}^{k}$
and a submersion $\psi :\Omega \rightarrow \mathbb{R}^{k}$ with the
following properties.

\begin{enumerate}
\item $\psi|_{W} = \pi$ and $\psi|_{\Omega\setminus U} = I\circ p$.

\item $\left\vert d\psi -d\pi \right\vert <\kappa \left( 1+2C\right) $ and $%
\left\vert d\psi -d\left( I\circ p\right) \right\vert <\kappa \left(
1+2C\right) .$

\item For $\varkappa _{1}:=\varkappa +\kappa \left( 1+2C\right) ,$ $\psi $
is a $\varkappa _{1}$--almost Riemannian submersion.

\item If $F$ is a subset of $\Omega $ with $\pi _{k}\circ I\circ p|_{F}=\pi
_{k}\circ \pi |_{F},$ then $\pi _{k}\circ \psi |_{F}=\pi _{k}\circ I\circ
p|_{F}=\pi _{k}\circ \pi |_{F},$ where $\pi _{k}:\mathbb{R}^{l}\rightarrow 
\mathbb{R}^{k}$ is projection to the first $k$--factors.
\end{enumerate}
\end{lemma}

\begin{proof}
We will show that for 
\begin{equation*}
\zeta :=\mathrm{dist(}W,\Omega \setminus U\mathrm{),}
\end{equation*}%
there is an isometry $I:\mathbb{R}^{k}\longrightarrow \mathbb{R}^{k}$ so
that 
\begin{equation}
\left\vert d\pi -d\left( I\circ p\right) \right\vert +2\frac{\left\vert \pi
-I\circ p\right\vert _{C^{0}}}{\zeta }<\kappa \left( 1+2C\right) .
\label{make it so}
\end{equation}%
Except for Part 3, the result follows from this by applying Lemma 7.4 in 
\cite{ProWilh1}.

To prove (\ref{make it so}), choose a point $w\in W$ and an isometry $I:%
\mathbb{R}^{k}\longrightarrow \mathbb{R}^{k}$ so that 
\begin{equation*}
\left\vert d\pi -d\left( I\circ p\right) \right\vert <\kappa \text{ and }%
I\circ p\left( w\right) =\pi \left( w\right) .\text{ }
\end{equation*}%
Then 
\begin{equation*}
\left\vert \pi -I\circ p\right\vert _{C^{0}}<\kappa \mathrm{diam}\left(
\Omega \right) \leq \kappa C\zeta ,
\end{equation*}%
so 
\begin{eqnarray}
\left\vert d\pi -d\left( I\circ p\right) \right\vert +2\frac{\left\vert \pi
-I\circ p\right\vert _{C^{0}}}{\zeta } &<&\kappa +2\frac{\kappa C\zeta }{%
\zeta }  \notag \\
&=&\kappa \left( 1+2C\right) .  \label{62}
\end{eqnarray}%
With the exception of Part 3, this completes the proof.
\end{proof}

Part 3 of the Submersion Deformation Lemma asserts that $\psi $ is $%
\varkappa _{1}=\varkappa +\kappa \left( 1+2C\right) $-almost Riemannian.
Since $\pi $ is $\varkappa $-almost Riemannian, this follows by combining
the next lemma with the assertion of Part 2 of the Submersion Deformation
Lemma that $\left\vert d\psi -d\pi \right\vert <\kappa \left( 1+2C\right) .$

\begin{lemma}
Let $\pi $ be a $\varkappa $--almost Riemannian submersion. Suppose $\psi $
is a $C^{1}$--map that satisfies%
\begin{equation*}
\left\vert d\pi -d\psi \right\vert <\delta .
\end{equation*}%
Then $\psi $ is a $\left( \varkappa +\delta \right) $--almost Riemannian
submersion, provided $\varkappa $ and $\delta $ are sufficiently small.
\end{lemma}

\begin{proof}
For any unit $w,$ 
\begin{equation}
\left\vert d\psi \left( w\right) \right\vert \leq \left\vert d\pi \left(
w\right) \right\vert +\delta \leq 1+\varkappa +\delta .  \label{norm}
\end{equation}%
If $x$ is $\pi $--horizontal and unit, then%
\begin{equation}
\left\vert \left\vert d\psi \left( x\right) \right\vert -1\right\vert \leq
\left\vert \left\vert d\pi \left( x\right) \right\vert -1\right\vert +\delta
<\varkappa +\delta .  \label{duh}
\end{equation}%
So $\psi $ is a submersion. Let $x^{h}$ be the $\psi $--horizontal component
of $x.$ Then 
\begin{equation*}
\left\vert d\psi \left( x^{h}\right) \right\vert =\left\vert d\psi \left(
x\right) \right\vert ,
\end{equation*}%
so by (\ref{duh}), 
\begin{equation*}
\left\vert \left\vert d\psi \left( x^{h}\right) \right\vert -1\right\vert
\leq \delta +\varkappa .
\end{equation*}%
Since $\left\vert x^{h}\right\vert \leq \left\vert x\right\vert =1$, this
gives us 
\begin{eqnarray*}
\left\vert x^{h}\right\vert \left( 1-\delta -\varkappa \right) &\leq
&1-\delta -\varkappa \\
&\leq &\left\vert d\psi \left( x^{h}\right) \right\vert \\
&\leq &\left\vert x^{h}\right\vert \left( 1+\varkappa +\delta \right) ,\text{
by }(\ref{norm}).
\end{eqnarray*}%
Since the map $x\longmapsto x^{h}$ is an isomorphism, the result follows.
\end{proof}

After using the Submersion Deformation Lemma to line up our submersions, we
will create fiber exhaustion functions so that (\ref{schenfl cont}) is
satisfied. We accomplish this via the next two results.

\begin{proposition}
\label{deform funct}Let $M$ be a Riemannian $n$--manifold. Let $W\subset
U\subset \Omega \subset M$ be three nonempty, open sets. Let $\lambda
:\Omega \longrightarrow \left[ 0,1\right] $ be $C^{\infty }$ so that $%
\lambda |_{W}=0$ and $\lambda |_{\Omega \setminus U}=1$. Let $p,\pi :\Omega
\longrightarrow \mathbb{R}$ be submersions so that $\pi \geq p,$ $%
\sphericalangle \left( \nabla p,\nabla \pi \right) \leq \frac{11}{12}\pi $,
and wherever $\nabla \lambda \neq 0$, $\max \left\{ \sphericalangle \left(
\nabla p,\nabla \lambda \right) ,\sphericalangle \left( \nabla \pi ,\nabla
\lambda \right) \right\} \leq \frac{11}{12}\pi $. Then there is a submersion 
$\psi :\Omega \longrightarrow \mathbb{R}$ such that 
\begin{equation}
\psi |_{W}=p\text{ and }\psi |_{\Omega \setminus U}=\pi .  \label{42}
\end{equation}
\end{proposition}

\begin{proof}
Set 
\begin{equation*}
\psi =\left( 1-\lambda \right) p+\lambda \pi .
\end{equation*}%
Then $\psi $ satisfies (\ref{42}) and 
\begin{equation*}
\nabla \psi =\left( 1-\lambda \right) \nabla p+\lambda \nabla \pi +\left(
\pi -p\right) \nabla \lambda .
\end{equation*}%
Since all three summands on the right lie in a half space, $\nabla \psi $ is
nowhere $0$, and $\psi $ is a submersion. 
\begin{figure}[b]
\includegraphics[scale=0.24]{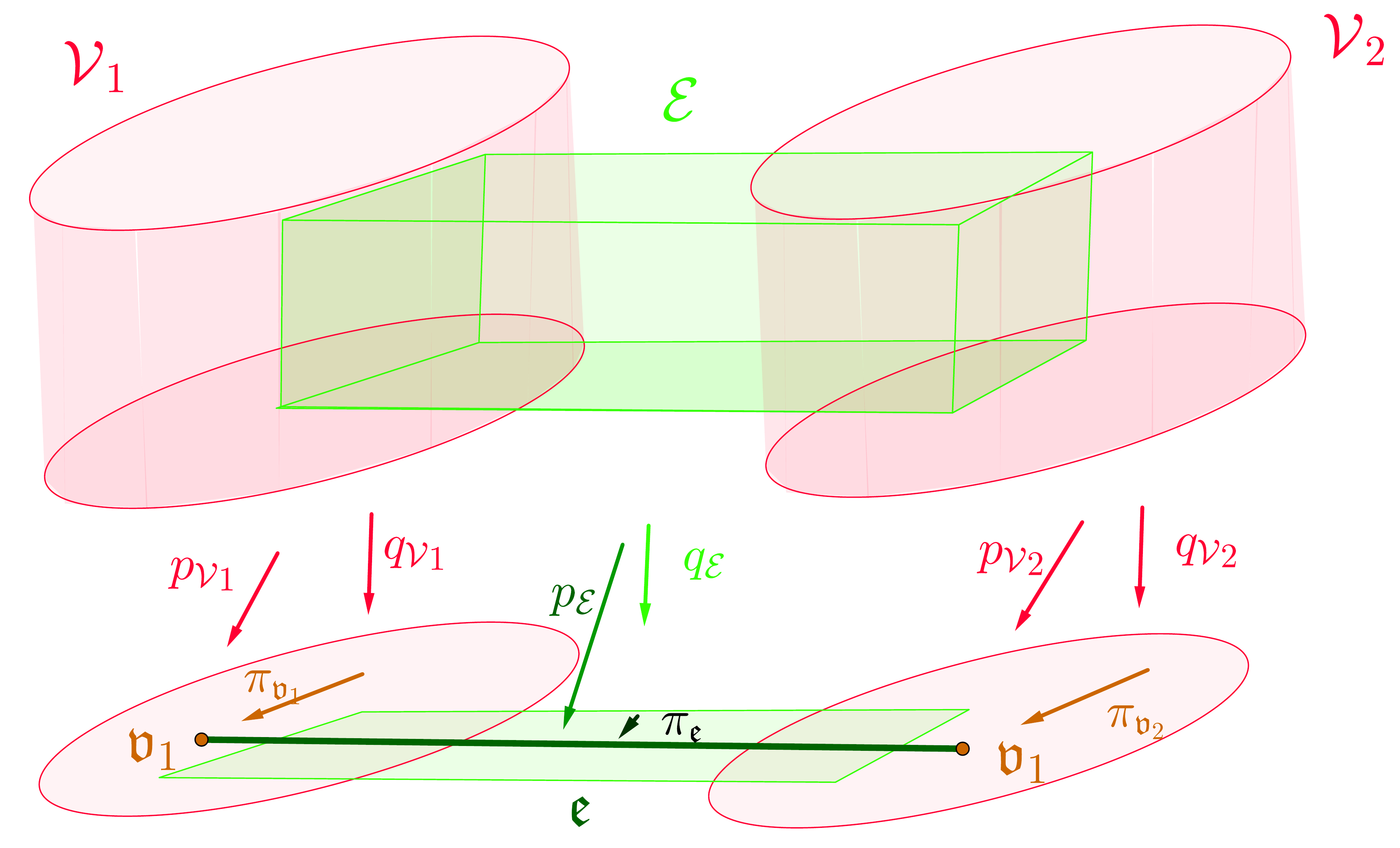}
\caption{{\protect\tiny {The green box $\mathcal{E}$ is artificially $1$%
--framed. The red cylinders, $\mathcal{V}_{1}$ and $\mathcal{V}_{2}$, are
artificially $0$--framed. The red cylinders are taller than the green box,
so (\protect\ref{major swallow}) is satisfied. The vertices $\mathfrak{v}%
_{1} $ and $\mathfrak{v}_{2}$ are disjoint from the green rectangle, so (%
\protect\ref{disjoint core}) is satisfied. }}}
\label{atrf framed 2}
\end{figure}
\end{proof}

\begin{lemma}
\label{line up respct lemma}For $l,h\in \left\{ 0,1,2\right\} $ with $l<h$,
let $\left( L,q_{L},U_{L},\mathfrak{l}\right) $ and $\left( H,q_{H},U_{H},%
\mathfrak{h}\right) $ be artificially $l$-framed and $h$-framed sets,
respectively, which are subordinate to $k$--framed sets, where $h\leq k\leq
2.$

Suppose that $q_{L}$ and $q_{H}$ are lined up on $(H\setminus H^{-})\cap L,$
and for $p_{H}:=\pi _{\mathfrak{h}}\circ q_{H}$, 
\begin{equation}
p_{H}^{-1}\left( p_{H}\left( H\cap L\right) \right) \Subset L
\label{major swallow}
\end{equation}%
and 
\begin{equation}
\mathfrak{l}\cap \bar{U}_{H}=\varnothing \text{ (see Figure \ref{atrf framed
2})}.  \label{disjoint core}
\end{equation}%
Then

\begin{enumerate}
\item The pair $\left\{ \left( L,p_{L}\right) ,\left( H,p_{H}\right)
\right\} $ satisfies the fiber containment condition (\ref{fiber contain}).

\item There is a choice of fiber exhaustion function $r_{L}$ for $L$ so that 
$\left\{ L,H\right\} $ is respectful.

\item The fiber exhaustion function $r_{L}$ for $L$ can be constructed as
follows: Given any fiber exhaustion function $r_{\mathfrak{l}}$ for $\pi _{%
\mathfrak{l}}:U_{L}\longrightarrow \mathrm{span}\left( \mathfrak{l}\right) $
whose 
\begin{equation}
\text{gradient is horizontal for }\pi _{\mathfrak{h}}:U_{H}\longrightarrow 
\mathrm{span}\left( \mathfrak{h}\right)  \label{horiz grad condition}
\end{equation}
on $q_{H}\left( H\cap L\right) ,$ any of the corresponding fiber exhaustion
functions for $L$ constructed from $r_{\mathfrak{l}}$ via Proposition \ref%
{artf framed is framed prop} make $\left\{ L,H\right\} $ respectful.

\item When $h=2,$ (\ref{horiz grad condition}) is satisfied for all fiber
exhaustion functions $r_{\mathfrak{l}}$ for $\pi _{\mathfrak{l}%
}:U_{L}\longrightarrow \mathrm{span}\left( \mathfrak{l}\right) .$
\end{enumerate}
\end{lemma}

\begin{proof}
For $x\in \left( H\setminus H^{-}\right) \cap L,$ Hypothesis (\ref{major
swallow}) together with $p_{H}=\pi _{\mathfrak{h}}\circ q_{H}$ gives us 
\begin{equation*}
q_{H}^{-1}\left( q_{H}\left( x\right) \right) \subset p_{H}^{-1}\left(
p_{H}\left( x\right) \right) \subset L.
\end{equation*}%
Combining this with the fact that $q_{H}$ is lined up with $q_{L}$ gives us 
\begin{equation}
q_{H}^{-1}\left( q_{H}\left( x\right) \right) \subset q_{L}^{-1}\left(
q_{L}\left( x\right) \right) .  \label{fiber containment 1}
\end{equation}

For each of the three cases $\left( l,h\right) \in \left\{ \left( 0,1\right)
,\left( 0,2\right) ,\left( 1,2\right) \right\} ,$ we must show that (\ref%
{fiber contain}) and (\ref{schenfl cont}) hold for $\left\{ L,H\right\} .$
In the two cases when $l=0,$ $p_{L}$ has only one four-dimensional fiber, $%
p_{L}^{-1}\left( \mathfrak{l}\right) =L,$ so the fiber containment condition
(\ref{fiber contain}) follows from Hypothesis (\ref{major swallow}). When $%
\left( l,h\right) =\left( 1,2\right) ,$ we have that $p_{H}=q_{H}.$ Together
with (\ref{fiber containment 1}), this gives us that for all $x\in \left(
H\setminus H^{-}\right) \cap L,$%
\begin{equation}
p_{H}^{-1}\left( p_{H}\left( x\right) \right) =q_{H}^{-1}\left( q_{H}\left(
x\right) \right) \subset q_{L}^{-1}\left( q_{L}\left( x\right) \right)
\subset p_{L}^{-1}\left( p_{L}\left( x\right) \right) .
\label{more fiber containment}
\end{equation}%
So (\ref{fiber contain}) holds in all three cases.

To establish (\ref{schenfl cont}), suppose that $L$ is subordinate to the $2$%
-framed set $(\digamma (\beta _{,}\nu ),p_{\digamma })$ and has the form 
\begin{equation*}
L=q_{L}^{-1}\left( U_{L}\right) \cap \digamma \left( \beta ,\tilde{\nu}%
\right) \text{ for }\tilde{\nu}<\nu .
\end{equation*}%
Together with (\ref{major swallow}), this implies that there is a $c<\tilde{%
\nu}$ so that 
\begin{equation}
q_{H}^{-1}\left( q_{H}\left( H\cap L\right) \right) \subset p_{H}^{-1}\left(
p_{H}\left( H\cap L\right) \right) \subset \digamma \left( \beta ,c\right) .
\label{swallow aTM}
\end{equation}

Suppose, for the moment, that we have a fiber exhaustion function $r_{%
\mathfrak{l}}$ for $\left( U_{L},\pi _{\mathfrak{l}}\right) $ as in Part $3$%
. Use $r_{\mathfrak{l}}$ and (\ref{disjoint core}) to apply Proposition \ref%
{artf framed is framed prop} and conclude that there is a fiber exhaustion
function $r_{L}$ for $\left( L,p_{L}\right) $ so that for all $x\in H\cap L,$
there is a $d$ so that 
\begin{equation}
q_{L}^{-1}\left( r_{\mathfrak{l}}^{-1}\left( r_{\mathfrak{l}}\left(
q_{H}\left( x\right) \right) \right) \right) \cap \digamma \left( \beta
,c\right) \subset r^{-1}(d).  \label{sucky outcome}
\end{equation}%
Since $\nabla r_{\mathfrak{l}}$ is horizontal for $\pi _{\mathfrak{h}%
}:U_{H}\longrightarrow \mathrm{span}\left( \mathfrak{h}\right) $ on $%
q_{H}\left( H\cap L\right) ,$ for all $x\in H\cap L,$ 
\begin{equation*}
\pi _{\mathfrak{h}}^{-1}\left( \pi _{\mathfrak{h}}\left( q_{H}\left(
x\right) \right) \right) \subset r_{\mathfrak{l}}^{-1}\left( r_{\mathfrak{l}%
}\left( q_{H}\left( x\right) \right) \right) .
\end{equation*}%
Combining this with (\ref{fiber containment 1}), (\ref{sucky outcome}), and $%
p_{H}=\pi _{\mathfrak{h}}\circ q_{H},$ we get that for all $x\in \left(
H\setminus H^{-}\right) \cap L,$ 
\begin{equation*}
p_{H}^{-1}\left( p_{H}\left( x\right) \right) =q_{H}^{-1}\left( \pi _{%
\mathfrak{h}}^{-1}\left( \pi _{\mathfrak{h}}q_{H}\left( x\right) \right)
\right) \cap \digamma \left( \beta ,c\right) \subset r^{-1}(d),
\end{equation*}%
and (\ref{schenfl cont}) holds, provided we can show that in each of our
three cases $\left( l,h\right) \in \left\{ \left( 0,1\right) ,\left(
0,2\right) ,\left( 1,2\right) \right\} ,$ there is a fiber exhaustion
function $r_{\mathfrak{l}}$ for $\left( U_{L},\pi _{\mathfrak{l}}\right) $
as in Part $3.$ When $h=2$, $\pi _{\mathfrak{h}}=\mathrm{id}_{\mathbb{R}%
^{2}},$ so any fiber exhaustion function for $\pi _{\mathfrak{l}}$ is
horizontal for $\pi _{\mathfrak{h}}.$ When $h=1,$ we construct $r_{\mathfrak{%
l}}$ by combining (\ref{disjoint core}) and Proposition \ref{deform funct}.
\end{proof}

\addtocounter{algorithm}{1}

\subsection{Quantitative Fiber Disjointness\label{quant openness}}

In this subsection, we establish two tools to verify the disjoint property
of our ultimate framed cover. The first says that just as a Riemannian
submersion has equidistant fibers, an almost Riemannian submersion has
almost equidistant fibers.

\begin{proposition}
\label{alm equi dist fibers}Let 
\begin{equation*}
\pi :M\longrightarrow B
\end{equation*}%
be a $\varkappa $--almost Riemannian submersion. For all $b_{1},b_{2}\in B$
and $x\in \pi ^{-1}\left( b_{1}\right) $, 
\begin{equation}
\left\vert \mathrm{dist}\left( x,\pi ^{-1}\left( b_{2}\right) \right) -%
\mathrm{dist}\left( b_{1},b_{2}\right) \right\vert \leq \varkappa \mathrm{%
dist}\left( b_{1},b_{2}\right) .  \label{Eqn eps open}
\end{equation}
\end{proposition}

\begin{proof}
Let $\gamma $ be a minimal geodesic from $x$ to $\pi ^{-1}\left(
b_{2}\right) .$ Then 
\begin{eqnarray*}
\mathrm{dist}\left( x,\pi ^{-1}\left( b_{2}\right) \right) &=&\mathrm{len}%
\left( \gamma \right) \\
&\geq &\left( 1-\varkappa \right) \mathrm{len}\left( \pi \circ \gamma \right)
\\
&\geq &\left( 1-\varkappa \right) \mathrm{dist}\left( b_{1},b_{2}\right) .
\end{eqnarray*}

To prove the opposite inequality, let $\gamma :\left[ 0,l\right]
\longrightarrow B$ be a minimal geodesic from $b_{1}$ to $b_{2}$ with $%
\left\vert \gamma ^{\prime }\left( t\right) \right\vert \equiv 1.$ Let $%
\tilde{\gamma}$ be a horizontal lift of $\gamma $ starting at $x.$ Then $%
\tilde{\gamma}\left( l\right) \in \pi ^{-1}\left( b_{2}\right) $. Since $%
\tilde{\gamma}$ is everywhere horizontal, 
\begin{equation*}
\left\vert \left\vert \tilde{\gamma}^{\prime }\left( t\right) \right\vert
-1\right\vert \leq \varkappa ,
\end{equation*}%
so%
\begin{equation*}
\left\vert \mathrm{len}\left( \tilde{\gamma}\right) -\mathrm{len}\left(
\gamma \right) \right\vert \leq \varkappa \mathrm{len}\left( \gamma \right)
=\varkappa \mathrm{dist}\left( b_{1},b_{2}\right) .
\end{equation*}%
Thus 
\begin{eqnarray*}
\mathrm{dist}\left( x,\pi ^{-1}\left( b_{2}\right) \right) &\leq &\mathrm{len%
}\left( \tilde{\gamma}\right) \\
&\leq &\mathrm{len}\left( \gamma \right) +\varkappa \mathrm{len}\left(
\gamma \right) \\
&=&\left( 1+\varkappa \right) \mathrm{dist}\left( b_{1},b_{2}\right) ,
\end{eqnarray*}%
as desired.
\end{proof}

\begin{lemma}
\label{Cor Q disj}For $i=1,2,$ let $\left( \digamma _{i}\left( 5\beta
_{i}\right) ,p_{i}\right) $ be $\left( l,\varkappa \right) $--framed sets
with $\beta _{2}\geq \beta _{1}.$ For all $b_{1},b_{2}\in p_{1}\left(
\digamma _{1}\left( \beta _{1}\right) \cap \digamma _{2}\left( \beta
_{2}\right) \right) $ and $x_{1}\in p_{1}^{-1}\left( b_{1}\right) \cap
\digamma _{1}\left( \beta _{1}\right) $,

\begin{enumerate}
\item 
\begin{equation}
\mathrm{dist}\left( x_{1},p_{2}^{-1}\left( b_{2}\right) \cap \digamma
_{2}\left( 5\beta \right) \right) \leq \left( 1+\varkappa \right) \mathrm{%
dist}\left( b_{1},b_{2}\right) .  \label{fibers close}
\end{equation}

\item Suppose that $\left\{ \left( \digamma _{1}\left( 5\beta _{1}\right)
,p_{1}\right) ,\left( \digamma _{2}\left( 5\beta _{2}\right) ,p_{2}\right)
\right\} $ is $\kappa $--lined up and fiber swallowing. If $I:\mathbb{R}%
^{l}\longrightarrow \mathbb{R}^{l}$ is an isometry so that 
\begin{equation}
p_{1}\left( x_{1}\right) =I\circ p_{2}\left( x_{1}\right) \text{ and }%
\left\vert dp_{1}-I\circ dp_{2}\right\vert <\kappa ,  \label{JLP}
\end{equation}%
then for all $\tilde{x}_{1}\in p_{1}^{-1}\left( b_{1}\right) \cap \digamma
_{1}\left( \beta _{1}\right) ,$ 
\begin{equation}
\left\vert \mathrm{dist}\left( b_{1},b_{2}\right) -\mathrm{dist}\left( 
\tilde{x}_{1},\left( I\circ p_{2}\right) ^{-1}\left( b_{2}\right) \cap
\digamma _{2}\left( 5\beta _{2}\right) \right) \right\vert \leq \varkappa 
\mathrm{dist}\left( b_{1},b_{2}\right) +3\kappa \beta _{1}.
\label{kappa wins}
\end{equation}
\end{enumerate}
\end{lemma}

\begin{proof}
Let $\gamma :\left[ 0,2l\right] \longrightarrow B$ be a minimal geodesic
from $b_{1}$ to $b_{2}$ with $\left\vert \gamma ^{\prime }\left( t\right)
\right\vert \equiv 1.$ Let $\tilde{\gamma}_{1}$ be a $p_{1}$--horizontal
lift of $\gamma |_{\left[ 0,l\right] }$ starting at $x_{1},$ and let $\tilde{%
\gamma}_{2}$ be a $p_{2}$--horizontal lift of $\gamma |_{\left[ l,2l\right]
} $ starting at $\tilde{\gamma}_{1}\left( l\right) .$ Set 
\begin{equation*}
\tilde{\gamma}=\tilde{\gamma}_{1}\ast \tilde{\gamma}_{2}.
\end{equation*}%
Then $\tilde{\gamma}\left( 2l\right) \in p_{2}^{-1}\left( b_{2}\right) \cap
\digamma _{2}\left( 5\beta \right) $. Since $\tilde{\gamma}_{i}$ is
everywhere $p_{i}$--horizontal, 
\begin{equation*}
\left\vert \left\vert \tilde{\gamma}^{\prime }\left( t\right) \right\vert
-1\right\vert \leq \varkappa ,
\end{equation*}%
so%
\begin{equation*}
\left\vert \mathrm{len}\left( \tilde{\gamma}\right) -\mathrm{len}\left(
\gamma \right) \right\vert \leq \varkappa \mathrm{len}\left( \gamma \right)
=\varkappa \mathrm{dist}\left( b_{1},b_{2}\right) .
\end{equation*}%
Thus 
\begin{eqnarray*}
\mathrm{dist}\left( x_{1},p_{2}^{-1}\left( b_{2}\right) \cap \digamma
_{2}\left( 5\beta \right) \right) &\leq &\mathrm{len}\left( \tilde{\gamma}%
\right) \\
&\leq &\mathrm{len}\left( \gamma \right) +\varkappa \mathrm{len}\left(
\gamma \right) \\
&=&\left( 1+\varkappa \right) \mathrm{dist}\left( b_{1},b_{2}\right) ,
\end{eqnarray*}%
as desired.

To prove Part 2, let $\tilde{b}_{1}:=I\circ p_{2}\left( \tilde{x}_{1}\right)
.$ By (\ref{JLP}), 
\begin{equation}
\mathrm{dist}\left( b_{1},\tilde{b}_{1}\right) =\mathrm{dist}\left(
p_{1}\left( \tilde{x}_{1}\right) ,I\circ p_{2}\left( \tilde{x}_{1}\right)
\right) <2\kappa \beta _{1}^{+}.  \label{red soxe}
\end{equation}%
Applying Proposition \ref{alm equi dist fibers} to $\tilde{x}_{1}$ and $%
I\circ p_{2}\left( b_{2}\right) $, we get 
\begin{eqnarray*}
\mathrm{dist}\left( \tilde{b}_{1},b_{2}\right) &\leq &\mathrm{dist}\left( 
\tilde{x}_{1},\left( I\circ p_{2}\right) ^{-1}\left( b_{2}\right) \cap
\digamma _{2}\left( 5\beta _{1}\right) \right) +\varkappa \mathrm{dist}%
\left( \tilde{b}_{1},b_{2}\right) \\
&\leq &\mathrm{dist}\left( \tilde{x}_{1},\left( I\circ p_{2}\right)
^{-1}\left( b_{2}\right) \cap \digamma _{2}\left( 5\beta _{1}\right) \right)
+\varkappa \left( \mathrm{dist}\left( b_{1},b_{2}\right) +\mathrm{dist}%
\left( b_{1},\tilde{b}_{1}\right) \right) \\
&\leq &\mathrm{dist}\left( \tilde{x}_{1},\left( I\circ p_{2}\right)
^{-1}\left( b_{2}\right) \cap \digamma _{2}\left( 5\beta _{1}\right) \right)
+\varkappa \left( \mathrm{dist}\left( b_{1},b_{2}\right) +2\kappa \beta
_{1}^{+}\right) ,
\end{eqnarray*}%
where we used (\ref{red soxe}) to get the last term in the last inequality.
The previous two inequalities give us%
\begin{eqnarray*}
\mathrm{dist}\left( b_{1},b_{2}\right) &\leq &\mathrm{dist}\left( b_{1},%
\tilde{b}_{1}\right) +\mathrm{dist}\left( \tilde{b}_{1},b_{2}\right) \\
&\leq &2\kappa \beta _{1}^{+}+\mathrm{dist}\left( \tilde{x}_{1},\left(
I\circ p_{2}\right) ^{-1}\left( b_{2}\right) \cap \digamma _{2}\left( 5\beta
_{1}\right) \right) +\varkappa \mathrm{dist}\left( b_{1},b_{2}\right)
+2\varkappa \kappa \beta _{1}^{+} \\
&\leq &3\kappa \beta _{1}+\mathrm{dist}\left( \tilde{x}_{1},\left( I\circ
p_{2}\right) ^{-1}\left( b_{2}\right) \cap \digamma _{2}\left( 5\beta
_{1}\right) \right) +\varkappa \mathrm{dist}\left( b_{1},b_{2}\right) ,
\end{eqnarray*}

which together with Part 1 gives us Part 2.
\end{proof}

An application of the previous result gives us the following additional tool
to create respectful framed collections.

\begin{proposition}
\label{glued submer prop}Let $\left\{ \left( A_{i}^{+},q_{i}\right) \right\}
_{i=1}^{L}$ be a family of artificially $\varkappa $--framed sets that are
lined up and subordinate to a single $2$--framed set $\left( \digamma \left(
\beta \right) ,p_{\digamma }\right) .$ If each $q_{i}$ is $\kappa $-lined up
with $p_{\digamma }$, then

\begin{enumerate}
\item There is a submersion 
\begin{equation*}
P:\cup _{i}A_{i}^{+}\longrightarrow \mathbb{R}^{2}
\end{equation*}%
so that for all $i$, $P$ is lined up with $q_{i}$.

\item If there is a $c>0$ so that 
\begin{equation}
\min_{i}\left\{ \mathrm{dist}\left( A_{i}^{-},A_{i}^{+}\setminus
A_{i}\right) \right\} \geq c\beta  \label{horiz wid}
\end{equation}%
and $\kappa $ and $\varkappa $ are sufficiently small compared with $c,$
then for all $x\in A_{i}^{-},$ 
\begin{equation*}
q_{i}^{-1}\left( q_{i}\left( x\right) \right) =P^{-1}\left( P\left( x\right)
\right) .
\end{equation*}
\end{enumerate}
\end{proposition}

\begin{proof}
To prove Part 1, we argue by induction on $i.$ For $i=1$ we take $P_1=q_{1}.$
By hypothesis, $P_1$ is lined up with all $q_i$. For $l\leq L,$ assume there
is a submersion 
\begin{equation*}
P_{l}:\cup _{i=1}^{l}A_{i}^{+}\longrightarrow \mathbb{R}^{2}
\end{equation*}%
that is lined up with all $q_{i}.$ Since $P_l$ is lined up with $q_{l+1}$,
there is a diffeomorphism $I_{l+1}:\mathbb{R}^{2}\rightarrow \mathbb{R}^{2}$
such that $P_{l}=I_{l+1}\circ q_{l+1}$ on $A_{l+1}^{+}\cap \left(\cup
_{i=1}^{l}A_{i}^{+}\right)$. Setting $P_{l+1} = P_l$ on $\cup
_{i=1}^{l}A_{i}^{+}$ and $P_{l+1}=I_{l+1}\circ q_{l+1}$ on $A^+_{l+1}$
defines a submersion $P_{l+1}$ on $\cup_{i=1}^{l+1} A_{i}^{+}$ that is lined
up with all $q_i$ and Part 1 is proven.

To prove Part 2, notice that since $P|_{A_{i}^{+}}$ is lined up with $q_{i},$
\begin{equation*}
q_{i}^{-1}\left( q_{i}\left( x\right) \right) =q_{i}^{-1}\left( q_{i}\left(
x\right) \right) \cap A_{i}^{+}=P^{-1}\left( P\left( x\right) \right) \cap
A_{i}^{+}\subset P^{-1}\left( P\left( x\right) \right) .
\end{equation*}%
The opposite containment follows unless for $j\neq k$ there are $x_{j}\in
A_{j}^{-}$ and $x_{k}\in A_{k}^{-}$ so that $P\left( x_{j}\right) =P\left(
x_{k}\right) $. Since $P|_{A_{l}^{+}}$ is lined up with $q_{l},$ this does
not happen if there is an $A_{l}$ with $x_{j},x_{k}\in A_{l}^{+}.$ If there
is no such $A_{l},$ then since $P$ is $\kappa $--lined up with $p_{\digamma
},$ and $\kappa $ and $\varkappa $ are small compared with $c,$ (\ref{kappa
wins}) and (\ref{horiz wid}) give us that $P\left( x_{j}\right) \neq P\left(
x_{k}\right) .$
\end{proof}

\part{Covering the Natural $1$-Stratum\label{Part 1}}

Let 
\begin{equation*}
\mathcal{K}=\mathcal{K}^{0}\cup \mathcal{K}^{1}\cup \mathcal{K}^{2}
\end{equation*}%
be the framed\ cover of Lemma \ref{art Ms Pacman lemma}. In this part, we
replace $\mathcal{K}$ by a framed cover 
\begin{equation*}
\mathcal{A}=\mathcal{A}^{0}\cup \mathcal{A}^{1}\cup \mathcal{A}^{2}
\end{equation*}%
whose collections of $0$- and $1$-framed sets satisfy Part 1 of the
conclusion of Limit Lemma \ref{cover of X lemma}.

To make the $1$-framed sets disjoint, we refine $\mathcal{K}^{1}$ to a
collection that consists of artificially $0$- and $1$-framed sets. To
accomplish this, for each $1$-framed set $(K,p_{K})$ in $\mathcal{K}^{1}$,
we take a finite subset $\{\mathfrak{v}_{1},\ldots ,\mathfrak{v}%
_{n}\}\subset p_{K}(K^{+})\subset \mathbb{R}$ and consider it as the
vertices of a graph. 
\begin{wrapfigure}[9]{r}{0.35\textwidth}\centering 
\vspace{-.2in}
\includegraphics[scale=0.18]{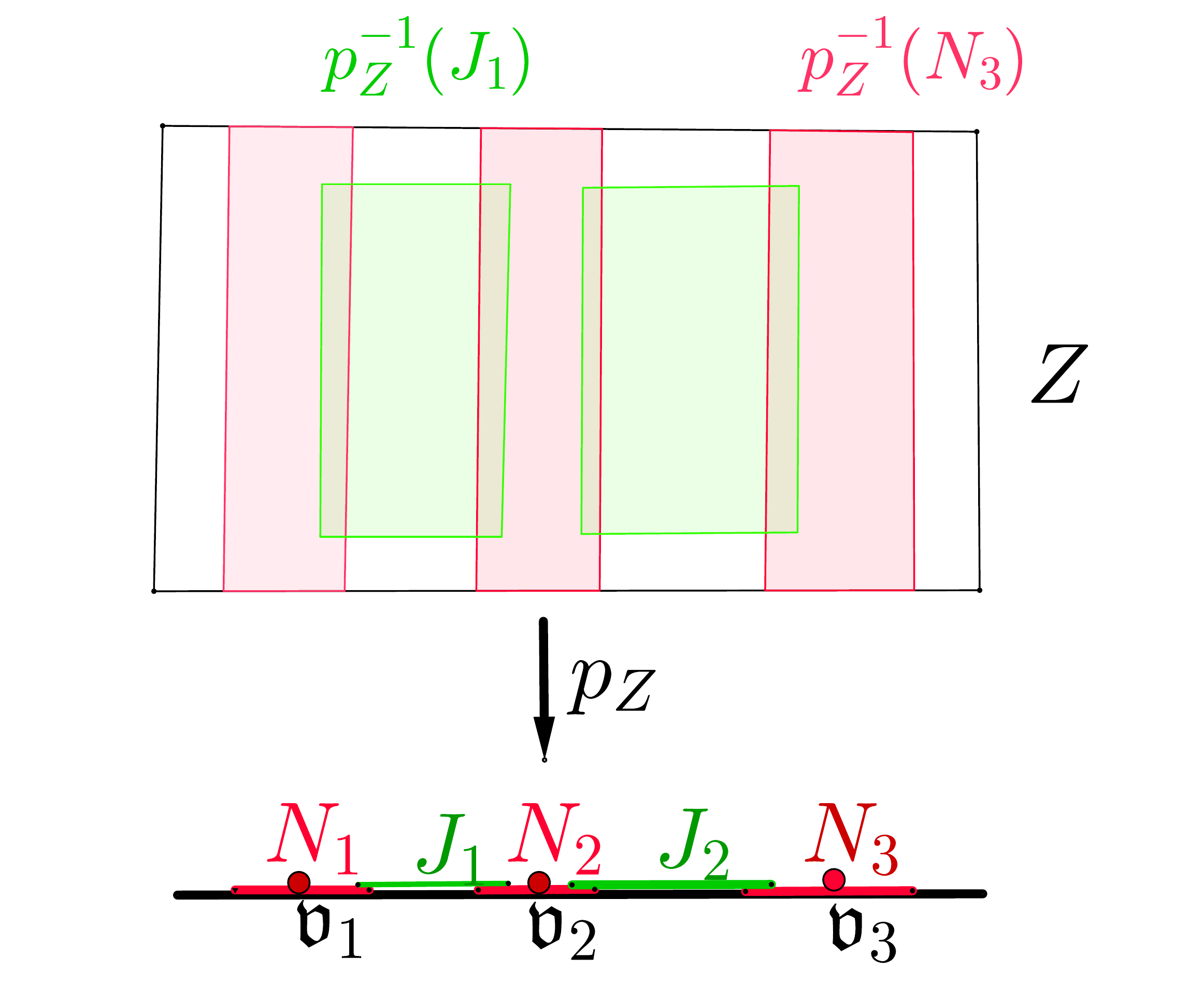}
\vspace{-.185in}
\caption{{\protect\tiny{A graph cover and a local \newline cograph}}}
\label{figure 2.5}
\end{wrapfigure}
We then take disjoint open intervals $N_{i}$ around the $\mathfrak{v}_{i}$
and declare the subsets $p_{K}^{-1}\left( N_{i}\right) \subset K$ to be
artificially $0$--framed. We also choose a family of disjoint intervals $%
\left\{ J_{i}\right\} $ between the $N_{i}$ so that $\left\{ N_{i}\right\}
\cup \left\{ J_{i}\right\} $ covers $p_{K}\left( K\right) $ and declare the
subsets $p_{K}^{-1}\left( J_{i}\right) \subset K$ to be artificially $1$%
-framed. As the open cover $\left\{ p_{K}^{-1}\left( N_{i}\right) \right\}
\cup \left\{ p_{K}^{-1}\left( J_{i}\right) \right\} $ comes from a graph, we
call it a cograph (see Figures \ref{figure 2.5}, \ref{coedge}, and \ref%
{5cograph}). Our declaration that the $p_{K}^{-1}\left( N_{i}\right) $ are $%
0 $--framed then yields a disjoint collection $\mathcal{A}^{0}$ of $0$%
--framed sets and a disjoint collection $\mathcal{A}^{1}$ of $1$--framed
sets. In this part, we exploit these ideas to prove

\begin{covlem}
\label{Cove lem 1}Given $\kappa ,\delta >0,$ there is a $\varkappa _{0}>0$
and for all $\varkappa \in \left( 0,\varkappa _{0}\right) ,$ a respectful, $%
\varkappa $--framed cover $\mathcal{A}=\mathcal{A}^{0}\cup \mathcal{A}%
^{1}\cup \mathcal{A}^{2}$ of $X\setminus X_{\delta }^{4}$ with the following
properties.

\begin{enumerate}
\item For $k=0\ $or $1,$ the closures of the elements of $\mathcal{A}^{k}$
are pairwise disjoint.

\item $\mathcal{A}^{2}$ is $\kappa $--lined up, vertically separated, and
fiber swallowing.
\end{enumerate}
\end{covlem}

This part has three sections. In Section \ref{cocgraph scet}, we define
cographs, explain why the existence of a certain type of cograph implies
that Covering Lemma \ref{Cove lem 1} holds, and assert in Lemma \ref{its
enough lemma} that such a cograph exists. Section \ref{Constr cographs sect}
develops the notion of geometric cograph, and in Section \ref{how to
cographs sect}, we construct the cograph that proves Lemma \ref{its enough
lemma}.

\section{Cographs\label{cocgraph scet}}

In this section, we define cographs and discuss their basic properties. In
Subsection \ref{graph covers subsect}, we define graph covers, which are
collections of open subsets of $\mathbb{R}$ that are precursors of cographs.
In Subsection \ref{What is a cograph subsect}, we explain why the existence
of a certain type of cograph implies that Covering Lemma \ref{Cove lem 1}
holds, and we assert in Lemma \ref{its enough lemma} that such a cograph
exists.

\addtocounter{algorithm}{1}

\subsection{Graph Covers\label{graph covers subsect}}

Recall that if $\mathcal{T}$ is a simplicial complex, then we let $\mathcal{T%
}_{i}$ be the collection of $i$--simplices of $\mathcal{T}$. We adopt the
convention that all graphs $\mathcal{T}$ whose vertex set $\mathcal{T}_{0}$
is a subset of $\mathbb{R}$ have the property that vertices are connected by
an edge only if they are consecutive. Let $\mathcal{T}$ be such a graph. The
starting point for the definition of cograph is a method to create a special
type of open cover $\mathcal{O}\left( \mathcal{T}\right) $ of $\left\vert 
\mathcal{T}\right\vert .$

\begin{definition}
Let $\mathcal{T}\subset \mathbb{R}$ be a graph and let $\left\vert \mathcal{T%
}\right\vert = \cup_{\sigma \in \mathcal{T}} |\sigma|$. A collection $%
\mathcal{O}\left( \mathcal{T}\right) $ of open subsets of $\mathbb{R}$ is
called a \textbf{Graph Cover of }$\mathcal{T}$\textbf{\ }provided:

\begin{enumerate}
\item $\mathcal{O}\left( \mathcal{T}\right) $ covers $\left\vert \mathcal{T}%
\right\vert .$

\item There is a one-to-one correspondence $g_{\mathcal{O},\mathcal{T}}:%
\mathcal{T}\longrightarrow \mathcal{O}\left( \mathcal{T}\right) $ so that
for both $k=0$ and $k=1,$ the collection $g_{\mathcal{O},\mathcal{T}}(%
\mathcal{T}_{k})$ consists of disjoint intervals.
\end{enumerate}

We say that $\mathcal{T}$ \textbf{generates} $\mathcal{O}\left( \mathcal{T}%
\right) .$
\end{definition}

\addtocounter{algorithm}{1}

\subsection{What is a Cograph?\label{What is a cograph subsect}}

In this subsection, we define the concept of a cograph. We start with the
local notion.

\begin{figure}[tbp]
\includegraphics[scale=0.21]{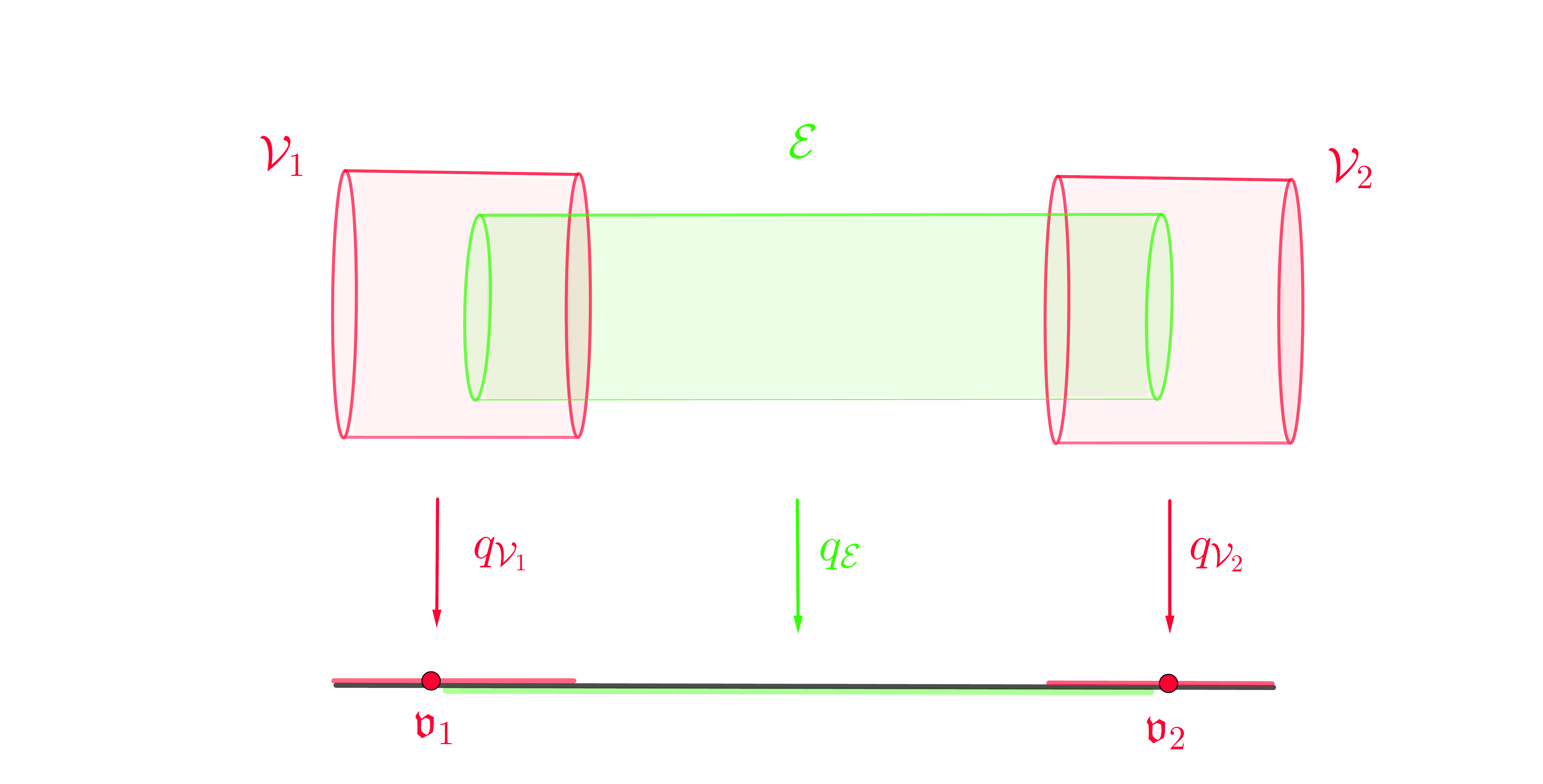}
\caption{{\protect\tiny {A cograph that corresponds to a graph with just two
vertices $\mathfrak{v}_{1}$ and $\mathfrak{v}_{1}$. The green cylinder is
called a coedge. The red cylinders are called covertices. The 3 submersions $%
q_{\mathcal{V}_{1}},q_{\mathcal{E}}$, and $q_{\mathcal{V}_{2}}$ are $\protect%
\kappa $--lined up, so on this scale their fibers appear parallel.}}}
\label{coedge}
\end{figure}

\begin{definition}[see Figures \protect\ref{figure 2.5} and \protect\ref%
{coedge}]
\label{local cograph dfn} Let $\digamma $ be a $1$--framed set and $\mathcal{%
V}$ and $\mathcal{E}$ collections of pairwise disjoint open subsets of $%
\digamma ^{+}.$ Then 
\begin{equation*}
\mathcal{G}:=\mathcal{V}\cup \mathcal{E}
\end{equation*}%
is a \textbf{local cograph} subordinate to $\digamma $ provided

\begin{enumerate}
\item There is a graph $\mathcal{T}$ with $p_{\digamma }\left( \digamma
\right) \subset \left\vert \mathcal{T}\right\vert \subset p_{\digamma
}\left( \digamma ^{+}\right) $ with associated graph cover $O\left( \mathcal{%
T}\right) $.

\item There is a one-to-one correspondence between $\mathcal{G}$ and a
subset $K_{\mathcal{T}}\subset \mathcal{T}$. Each $S\in \mathcal{G}$ is
generated by the corresponding $\mathfrak{s\in }K_{\mathcal{T}}$ as follows:
There is a submersion $q_{S}:$ $\digamma \rightarrow \mathbb{R}$ so that $%
(S,q_{S},U_{S},\mathfrak{s)}$ is an artificially $k$-framed set with $%
\mathfrak{s}\in \mathcal{T}_{k}\cap K_{\mathcal{T}}$, $U_{S}=g_{\mathcal{O},%
\mathcal{T}}\left( \mathfrak{s}\right) $, and 
\begin{equation}
\begin{array}{ll}
\text{if }\mathfrak{s\in }\mathcal{T}_{0}, & \text{then }S=q_{S}^{-1}\left(
U_{S}\right) \cap \digamma (\beta ,\beta ^{+})\in \mathcal{V}, \\ 
\text{if }\mathfrak{s\in }\mathcal{T}_{1}, & \text{then }S=q_{S}^{-1}\left(
U_{S}\right) \cap \digamma (\beta ,\beta )\in \mathcal{E}.%
\end{array}
\label{stop re-writing Curtis}
\end{equation}
\end{enumerate}
\end{definition}

For $\mathcal{G}$, $\mathcal{T},$ and $K_{\mathcal{T}}$ as above we say that 
$K_{\mathcal{T}}$ generates $\mathcal{G}.$ In general, $K_{\mathcal{T}}$
will be a proper subset of $\mathcal{T}.$ This is crucial to our plan to fit
local cographs together with each other and with the natural $0$ stratum
(see e.g. Figure \ref{cogrph inte nat 0 str}).

\begin{wrapfigure}[6]{r}{0.3\textwidth}\centering 
\vspace{-0.1in}
\includegraphics[scale=1.2]{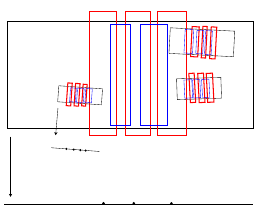}
\vspace{-0.09in}
\caption{A cograph with corresponding local cographs}
\label{5cograph}
\end{wrapfigure}

\begin{wrapfigure}[10]{l}{0.25\textwidth}\centering
\vspace{-12pt}
\includegraphics[scale=0.15]{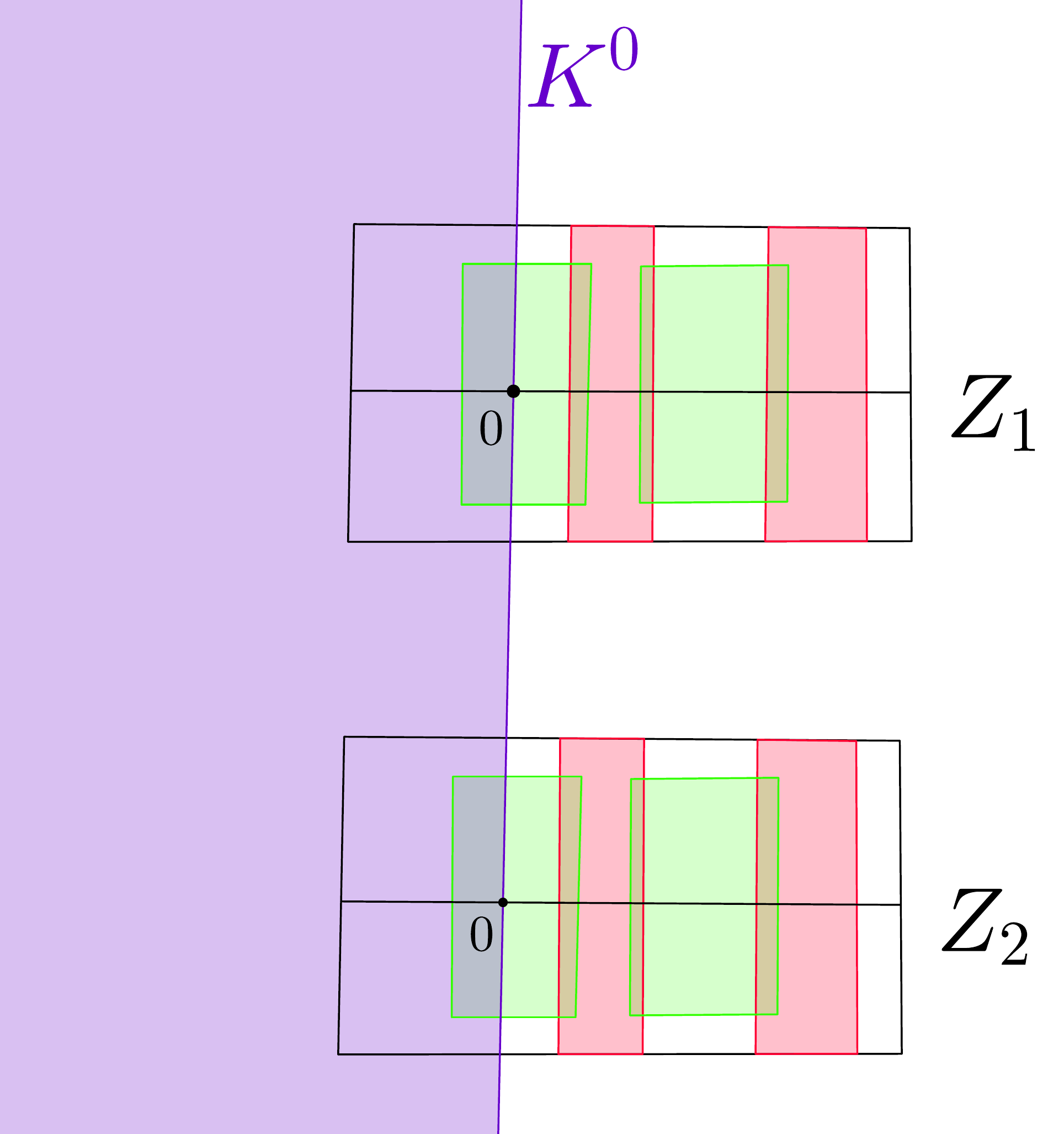}
\vspace{-8pt}
\caption{{\protect\tiny {This depicts a cograph relative to $\Omega_0$ with $\mathcal{V}$ replaced with $\mathcal{V} \backslash \mathcal{V}^0$.  } }}
\label{cogrph inte nat 0 str}
\end{wrapfigure}

\begin{definition}[see Figures \protect\ref{figure 2.5}, \protect\ref{coedge}%
, and \protect\ref{5cograph}]
\label{cograph dfn} Let $\mathcal{U}\equiv \{\digamma _{i}(\beta
_{i})\}_{i=1}^{l}$ be a collection of $\left( 1,\varkappa \right) $--framed
sets that are $\kappa $--lined up, vertically separated, and fiber
swallowing. We say a collection 
\begin{equation*}
\mathcal{G}\equiv \mathcal{V}\cup \mathcal{E}
\end{equation*}
\end{definition}

\noindent \emph{of open subsets of $X$ is a \textbf{cograph} \textbf{%
subordinate to }$\mathcal{U}$ \vspace{.09in}\vspace{1pt} provided \noindent
1. For each $i\in \{1,\ldots ,l\},$ there is a local cograph $\mathcal{G}%
_{i}\equiv \mathcal{V}_{i}\cup \mathcal{E}_{i}$ subordinate to $\digamma
_{i}\ $with $\mathcal{V=}\cup _{i}\mathcal{V}_{i},\ \mathcal{E=}\cup _{i}%
\mathcal{E}_{i},\text{ and }\mathcal{G}=\cup _{i} \mathcal{G}_{i}.$ \vspace{%
6pt} \newline
\noindent 2. Each of the collections $\mathcal{V}\ $and $\mathcal{E}$ is
pairwise \vspace{.07in} disjoint. \vspace{1pt} \newline
\noindent 3. $\mathcal{G}$ covers $\cup \mathcal{U}^{-}.$}

\begin{definition}
A cograph will be called \textbf{respectful} provided it is respectful as a
framed collection (see Definition \ref{dfn resp framed collect}). The
elements of $\mathcal{G},$ $\mathcal{V},$ and $\mathcal{E}$ are called
cosimplices, covertices, and coedges, respectively.
\end{definition}

Next we use the results of Subsection \ref{creating resp} to deform a
cograph to a respectful cograph.

\begin{proposition}
\label{25}Let $\mathcal{G}= \mathcal{V}\cup \mathcal{E}$ be a cograph. There
are fiber exhaustion functions for the elements of $\mathcal{V}$ so that $%
\mathcal{G}$ is a respectful cograph.
\end{proposition}

\begin{proof}
Suppose that $V\in \mathcal{V}$ and $E\in \mathcal{E}$ intersect, and%
\begin{equation*}
q_{V}:V\longrightarrow \mathbb{R}\text{ and }q_{E}:E\longrightarrow \mathbb{R%
}
\end{equation*}%
are the submersions that define $V$ and $E$ as artificially framed sets.
Then $q_{V}$ and $q_{E}$ are $\varepsilon $-lined up on $V\cap E\ $for some
small $\varepsilon >0.$ So after post-composing $q_{E}$ with an isometry $I$
we get an almost Riemannian submersion $I\circ q_{E}:V\cap
E^{+}\longrightarrow \mathbb{R}$ whose gradients are within $\varepsilon $
of those of $q_{V}.$ Using this and Proposition \ref{deform funct}, we get a
deformation $\tilde{q}_{V}$ of $q_{V}$ that is lined up with $q_{E}$ on $%
E\setminus E^{-}.$ It follows from this and Lemma \ref{line up respct lemma}
that there is a fiber exhaustion function for $V$ that makes the
intersection of $E$ and $V$ respectful.
\end{proof}

The gluing method of the previous proof also allows us to extend coedge
submersions so that they respect the submersions of neighboring covertices.

\begin{corollary}
\label{extedn coedge submer}Let $\mathcal{G}=\mathcal{V}\cup \mathcal{E}$ be
a cograph subordinate to $\mathcal{U}$. For $E\in \mathcal{E}$ let $\mathcal{%
V}\left( E\right) $ be the collection of covertices that intersect $E,$ and
let $\digamma _{E}$ be the element of $\mathcal{U}$ so that $E$ is
subordinate to $\digamma _{E}.$ There is a submersion 
\begin{equation*}
q_{E}^{\mathrm{ext}}:E\cup \left( \cup \mathcal{V}\left( E\right) \right)
\cap \digamma _{E}^{+}\longrightarrow \mathbb{R}
\end{equation*}%
so that

\begin{enumerate}
\item $q_{E}^{\mathrm{ext}}$ extends $q_{E}|_{E^{-}}$ and

\item For each $V\in \mathcal{V}\left( E\right) ,$ $q_{E}^{\mathrm{ext}}$ is
lined up with $q_{V}$ on $V\setminus E.$
\end{enumerate}
\end{corollary}

\begin{proof}
For intersecting $V\in \mathcal{V}$ and $E\in \mathcal{E},$ let $I$ and $%
\tilde{q}_{V}$ be as in the proof of Proposition \ref{25}, and set 
\begin{equation*}
q_{E}^{\mathrm{ext}}:=\left\{ 
\begin{array}{ll}
q_{E} & \text{on }E \\ 
I^{-1}\circ \tilde{q}_{V} & \text{on }\left( V\cup \left( E^{+}\setminus
E^{-}\right) \right) \cap \digamma _{E}^{+}.%
\end{array}%
\right.
\end{equation*}
\end{proof}

We will see that the existence of the following type of cograph implies that
Covering Lemma \ref{Cove lem 1} holds.

\begin{definition}
Let $\mathcal{K}=\mathcal{K}^{0}\cup \mathcal{K}^{1}\cup \mathcal{K}^{2}$ be
as in Lemma \ref{art Ms Pacman lemma}. Let $\Omega _{0}:=\cup \mathcal{K}%
^{0} $, and let $\mathcal{G}=\mathcal{V}\cup \mathcal{E}$ be a cograph
subordinate to $\mathcal{K}^{1}$. Then $\mathcal{G}$ is \textbf{relative to }%
$\Omega _{0}$ provided for all $\digamma \left( \beta _{\digamma }\right)
\in \mathcal{K}^{1}$ with $\digamma \cap \Omega _{0}\neq \varnothing ,$

\begin{enumerate}
\item If $(S,q_{S})\in \mathcal{G}$ is an artificially framed set
subordinate to $\digamma ,$ then $q_{S}=p_{\digamma }.$

\item $p_{\digamma }\left( \bar{\digamma}\cap \left( \bar{\Omega}^{0}\right)
^{-}\right) =\{0\}\in \left[ 0,\beta _{\digamma }\right] =p_{\digamma
}\left( \digamma \right) $ generates a covertex of $\mathcal{G}$ that is
subordinate to $\digamma ^{+}.$ Let $\mathcal{V}^{0}$ be the collection of
all such covertices.
\end{enumerate}
\end{definition}

Combining this definition with Part 2 of Lemma \ref{art Ms Pacman lemma}, we
get

\begin{corollary}[cf Figure \protect\ref{cogrph inte nat 0 str}]
\label{resp 0-str prop} Let $\mathcal{G}=\mathcal{V}\cup \mathcal{E}$ be a
cograph subordinate to $\mathcal{K}^{1}$ and relative to $\Omega _{0}.$ Then

\begin{enumerate}
\item $\mathcal{V\setminus V}^{0}$ is disjoint from $\mathcal{K}^{0}.$

\item The collection $\mathcal{E}$ of coedges$\mathcal{\ }$of $\mathcal{G}$
respects $\mathcal{K}^{0}.$
\end{enumerate}
\end{corollary}

%

Combining Proposition \ref{25} and Corollary \ref{resp 0-str prop}, we see
that Covering Lemma \ref{Cove lem 1} holds provided we can construct a
cograph $\mathcal{G}=\mathcal{V}\cup \mathcal{E}$ that is subordinate to $%
\mathcal{K}^{1}$ and relative to $\Omega _{0}.$ Indeed, if $\mathcal{G}=%
\mathcal{V}\cup \mathcal{E} $ is the cograph, then the desired framed cover
is 
\begin{equation}
\hspace*{1.85in} \mathcal{A}^{0}\left( \mathcal{G}\right) :=\mathcal{K}%
^{0}\cup \left( \mathcal{V\setminus V}^{0}\right) ,\text{ }\mathcal{A}%
^{1}\left( \mathcal{G}\right) :=\mathcal{E},\text{ and }\mathcal{A}^{2}=%
\mathcal{K}^{2}.\text{ \label{cov ekm 1 cover}}
\end{equation}

\noindent Thus for the remainder of this part, our goal is to prove

\begin{lemma}
\label{its enough lemma}There is a cograph $\mathcal{G}$ that is subordinate
to $\mathcal{K}^{1}$\ and relative to $\Omega _{0}$.
\end{lemma}

\section{Geometric Cographs\label{Constr cographs sect}}

In this section, we impose geometric restrictions on our cographs. While
Lemma \ref{its enough lemma} can probably be proven without these
constraints, their presence will be crucial to our arguments in Part \ref%
{Part 2} of the paper. The main constraint is called $\eta $--geometric,
where $\eta $ is a positive number that will parameterize lower bounds for
the distances between the vertices of the underlying graph, the horizontal
thickness of the cosimplices, and the distances between the cosimplices of
the same type. We define geometric graph covers and geometric cographs in
Subsections \ref{Geom graph covers sect} and \ref{geom cographs subsect},
respectively.

\addtocounter{algorithm}{1}

\subsection{Geometric Graph Covers\label{Geom graph covers sect}}

In this subsection, we define $\eta $--geometric\ graph covers and discuss
their key properties.

\begin{definition}
Let $\mathcal{T}$ be an abstract graph for which the vertex set $\mathcal{T}%
_{0}$ is a metric space. Given $\eta >0$ and $\delta \in \left[ 0,\frac{1}{%
100}\right] ,$ we say that a graph $\mathcal{T}$ is an \textbf{abstract }$%
\left( \eta ,\delta \eta \right) $\textbf{--geometric graph} provided its
vertices are $\eta \left( 1-\delta \right) $--separated and its edge lengths
are all in $\left[ \eta \left( 1-\delta \right) ,2\eta \left( 1+\delta
\right) \right] .$ \textbf{\ }If $\mathcal{T}_{0}\subset \mathbb{R},$ we
drop the word abstract and simply say that $\mathcal{T}$ is an $\left( \eta
,\delta \eta \right) $--geometric graph. If $\delta =0,$ we say that $%
\mathcal{T}$ is $\eta $--\textbf{geometric.}
\end{definition}

\begin{figure}[tbp]
\includegraphics[scale=0.175]{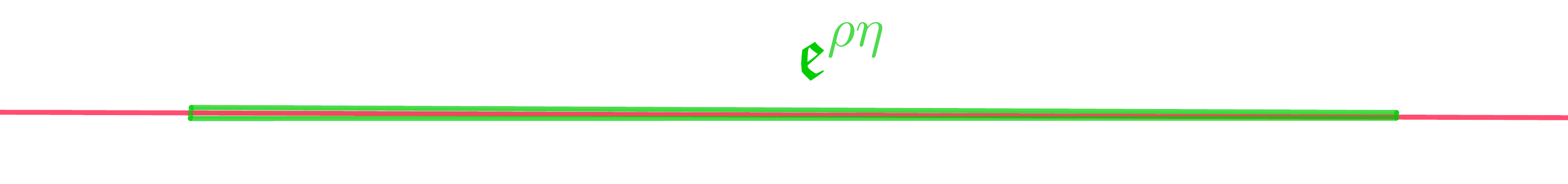}
\caption{{\protect\tiny {The larger interval is $\mathfrak{e}$, and the
smaller is $\mathfrak{e}^{\protect\rho \protect\eta }$. }}}
\label{figure bar in cograph}
\end{figure}

\begin{definition}[see Figure \protect\ref{figure bar in cograph}]
Let $\mathcal{T}$ be an $\left( \eta ,\delta \eta \right) $--geometric
graph. Let $\mathfrak{e}\in \mathcal{T}_{1}$ and let $\mathfrak{v}_1,%
\mathfrak{v}_2 \in \mathcal{T}_0$ be the boundary of $\mathfrak{e}$. For $%
\rho >0$ set 
\begin{equation*}
\mathfrak{e}^{\rho \eta } = \mathfrak{e}\setminus \left( \overline{B(%
\mathfrak{v}_1,2\rho \eta )}\cup \overline{B(\mathfrak{v}_2,2\rho \eta )}%
\right).
\end{equation*}
\end{definition}

\begin{definition}[see Figure \protect\ref{figure 2.5}]
Let $\mathcal{O}\left( \mathcal{T}\right) $ be a graph cover of an $\left(
\eta ,\delta \eta \right) $--geometric graph\textbf{\ }$\mathcal{T}.$ Given $%
\rho ,d>0,$ we say that\textbf{\ }$\mathcal{O}\left( \mathcal{T}\right) $ is 
$\left( \eta ;\rho ,d\right) $--geometric provided

\begin{enumerate}
\item For all $\mathfrak{v\in }\mathcal{T}_{0}$, $g_{\mathcal{O},\mathcal{T}%
}\left( \mathfrak{v}\right) =B\left( \mathfrak{v},5\rho \eta \right) .$

\item For all $\mathfrak{e\in }\mathcal{T}_{1},$ $g_{\mathcal{O},\mathcal{T}%
}\left( \mathfrak{e}\right) =\mathfrak{e}^{\rho \eta }.$

\item Each of the collections 
\begin{equation*}
\left\{ g_{\mathcal{O},\mathcal{T}}\left( \mathfrak{v}\right) \right\} _{%
\mathfrak{v\in \mathcal{T}}_{0}}\text{ and }\left\{ g_{\mathcal{O},\mathcal{T%
}}\left( \mathfrak{e}\right) \right\} _{\mathfrak{e\in \mathcal{T}}_{1}}
\end{equation*}%
is $d\eta $--separated.
\end{enumerate}
\end{definition}

The following is an immediate consequence of this definition.

\begin{proposition}
\label{graph cover prop}There are universal constants $\rho ,d>0$ so that
any $\left( \eta ,\delta \eta \right) $--geometric graph has an $\left( \eta
;\rho ,d\right) $--geometric graph cover $\mathcal{O}\left( \mathcal{T}%
\right) .$ Moreover, if $\mathcal{O}\left( \mathcal{\tilde{T}}\right) $ is
any $\left( \eta ;\rho ,d\right) $--geometric graph cover of a subcomplex $%
\mathcal{\tilde{T}\subset T},$ then $\mathcal{O}\left( \mathcal{\tilde{T}}%
\right) $ extends to an $\left( \eta ;\rho ,d\right) $--geometric graph
cover of $\mathcal{T}$.
\end{proposition}

Since the constants $\rho $ and $d$ are universal, for brevity we will omit
them from the notation and refer to these covers as $\eta $--graph covers.
Thus an $\left( \eta ,\delta \eta \right) $--geometric graph $\mathcal{T}$
determines an $\eta $--graph cover $\mathcal{O}\left( \mathcal{T}\right) .$
If, in addition, $\left\vert \mathcal{T}\right\vert \subset p_{\digamma
}\left( \digamma ^{+}\right) ,$ where $\digamma $ is $1$--framed, then $%
\mathcal{O}\left( \mathcal{T}\right) $ determines a local cograph via the
following construction.

\begin{definition}[see Figure \protect\ref{figure 2.5}]
\label{stand local cographs dfn}Let $\digamma $ be a $1$--framed set, and
let $\mathcal{T}$ be an $\left( \eta ,\delta \eta \right) $--geometric graph 
$\mathcal{T}$ in $p_{\digamma }\left( \digamma ^{+}\right) $ with $\eta $%
--graph cover $\mathcal{O}\left( \mathcal{T}\right) .$ A cosimplex $S$ that
is subordinate to $\digamma $ and generated by a $U\in \mathcal{O}\left( 
\mathcal{T}\right) $ via (\ref{stop re-writing Curtis}) is called a \textbf{%
standard local }$\eta $--\textbf{cosimplex associated to }$\left( \digamma ,%
\mathcal{T}\right) .$ A collection $\mathcal{G}=\mathcal{V}\cup \mathcal{E}$
of \textbf{standard local }$\eta $--\textbf{cosimplices }that are generated
by simplices of $\mathcal{T}$ is called a \textbf{standard local }$\eta $--%
\textbf{cograph}, provided the collection $\mathcal{V}$ of covertices is
disjoint and the collection $\mathcal{E}$ of coedges is disjoint. Otherwise
we call $\mathcal{G}$ a standard \textbf{local precograph.}
\end{definition}

The following are immediate corollaries of this definition and the
definition of $\kappa $--lined up.

\begin{corollary}
\label{cover a one framed set cor}Let $\mathcal{T}$ be an $\eta $--geometric
graph in $p_{\digamma }\left( \digamma ^{+}\right) $, and let 
\begin{equation*}
\left\vert K_{\mathcal{T}}\right\vert :=\cup _{\mathfrak{s}\in K_{\mathcal{T}%
}}\left\vert \mathfrak{s}\right\vert .
\end{equation*}%
There is a respectful, standard local $\eta $--cograph $\mathcal{G}\left( 
\mathcal{T}\right) $ associated to $\left( \digamma ,\mathcal{T}\right) $
which is generated by $K_{\mathcal{T}}$ and whose cosimplices are all
defined using the submersion $p_{\digamma }.$ Moreover, if $p_{\digamma
}\left( \bar{\digamma}\right) \subset \left\vert K_{\mathcal{T}}\right\vert
, $ then $\mathcal{G}\left( \mathcal{T}\right) $ covers $\digamma $.
\end{corollary}

\begin{corollary}
\label{dsij from o stra cor}Let $\mathcal{G=}\mathcal{V}\cup \mathcal{E}$ be
a standard local $\eta $--precograph associated to $\left( \digamma \left(
\beta \right) ,\mathcal{T}\right) .$ If $\eta \geq \frac{\beta }{500}$ and $%
\kappa $ and $\varkappa $ are sufficiently small$,$ then $\mathcal{G}$ is a
cograph.
\end{corollary}

\addtocounter{algorithm}{1}

\subsection{Geometric Cographs\label{geom cographs subsect}}

In this subsection, we define two constraints on cographs which we call $%
\eta $\textbf{--geometric }and \textbf{legal. }

\begin{definition}
\label{dfn eta geom cograph}Let $\{\digamma _{i}(\beta _{i})\}_{i=1}^{L}$ be
a family of $\left( 1,\varkappa \right) $--framed sets that are $\kappa $%
--lined up, vertically separated, and fiber swallowing. For $l\in \left\{
1,\ldots ,L\right\} ,$ let%
\begin{equation*}
\mathcal{G}=\mathcal{V}\cup \mathcal{E}=\bigcup\limits_{i=1}^{l}\mathcal{G}%
_{i}=\bigcup\limits_{i=1}^{l}\left( \mathcal{V}_{i}\cup \mathcal{E}%
_{i}\right)
\end{equation*}%
be a cograph that is subordinate to $\{\digamma _{i}(\beta
_{i})\}_{i=1}^{l}. $ We say that $\mathcal{G}$ is $\eta $--\textbf{geometric}
provided for $i\in \left\{ 1,\ldots ,l\right\} ,$ the local cograph $%
\mathcal{G}_{i}$ is at Gromov-Hausdorff distance $\leq \tau \left( \kappa
\eta \right) $ from a standard local $\eta $--cograph.
\end{definition}

For technical reasons that arise in Part \ref{Part 2} of this paper, we will
need to further constrain some of our cographs as follows.

\begin{definition}
\label{legal dfbn}Let $\mathcal{T}$ be an $\left( \eta ,\delta \eta \right) $%
-geometric graph. If the edges of $\mathcal{T}$ have length strictly less
than $(\sqrt{2}-\frac{1}{10})\eta $, then we say that $\mathcal{T}$ is a 
\textbf{legal} $\eta $--geometric graph.

We say that an $\eta $--geometric cograph is \textbf{legal }if\textbf{\ }all
of its generating graphs are legal$.$
\end{definition}

The only significance of the constant $\frac{1}{10}$ in the definition above
is that it yields a definitive value strictly less than $\sqrt{2}\eta $.

\section{How to Construct Cographs\label{how to cographs sect}}

In this section, we prove Lemma \ref{its enough lemma} and hence Covering
Lemma \ref{Cove lem 1}. In fact, we prove the following more general result,
which will assist in the project of sewing the cograph constructed here with
the framed sets that we construct in Part \ref{Part 2}.

\begin{theorem}
\label{cograph constr them}Let $\mathcal{A}=\left\{ \digamma _{l}(\beta
_{l}^{+})\right\} _{l=1}^{L}$ be a collection of $\left( 1,\varkappa \right) 
$-framed sets that are $\kappa $--lined up, vertically separated, fiber
swallowing, and ordered so that $\beta _{i}\geq \beta _{i+1}.$ There is an $%
\eta _{L}\in \left[ \frac{\beta _{L}}{500},\frac{\beta _{L}}{50}\right] $
and a legal $\eta _{L}$-cograph $\mathcal{G}^{L}=\mathcal{V}^{L}\cup 
\mathcal{E}^{L}$ subordinate to $\mathcal{A}$.

Furthermore, if $\mathcal{A}=\mathcal{K}^{1}$ where $\mathcal{K}^{1}$ is as
in Lemma \ref{art Ms Pacman lemma}, then we can choose $\mathcal{G}^{L}$ to
be relative to $\Omega_{0}$.
\end{theorem}

The proof is by induction on $l$ and involves two key tools---the Extension
and Subdivision Lemmas. These lemmas are stated below as Lemmas \ref{cograph
ext lemma} and \ref{cographsubdiv lemma} and proved in Subsections \ref%
{constrcogrphs subsect} and \ref{subdivided cogrph subsect}, respectively.
The induction statement asserts, among other things, that there is a
sequence of positive numbers $\left\{ \eta _{l}\right\} _{l=1}^{L}$ and a
cograph $\mathcal{G}^{l}$ so that Theorem \ref{cograph constr them} holds
with $L$ replaced by $l.$ The sequence of $\eta _{i}$ comes from the
following proposition.

\begin{proposition}
\label{beta eta sequence}Let $\left\{ \beta _{i}^{+}\right\}
_{i=1}^{M}\subset \left( 0,\infty \right) $ be ordered so that $\beta
_{i}\geq \beta _{i+1}.$ Given any $\eta _{1}\in \left[ \frac{\beta _{1}}{500}%
,\frac{\beta _{1}}{50}\right] $ there is a sequence $\left\{ \eta
_{i}\right\} _{i=1}^{M}$ so that 
\begin{equation*}
\eta _{i}\in \left[ \frac{\beta _{i}}{500},\frac{\beta _{i}}{50}\right]
\end{equation*}%
and 
\begin{equation*}
\eta _{i-1}=10^{k_{i}}\eta _{i}
\end{equation*}%
for some $k_{i}=0,1,2,3,\ldots $.
\end{proposition}

\begin{proof}
For each $i\in \left\{ 2,3,\ldots ,M\right\} ,$ we let $n\left( i\right) $
be the least nonnegative integer so that 
\begin{equation}
\eta _{i}:=\frac{\eta _{1}}{10^{n\left( i\right) }}\leq \frac{\beta _{i}}{50}%
.  \label{n of l dfn}
\end{equation}%
It follows that $\eta _{i}\in \left[ \frac{\beta _{i}}{500},\frac{\beta _{i}%
}{50}\right] .$
\end{proof}

To begin the proof of Theorem \ref{cograph constr them}, we take any $\eta
_{1}\in \left[ \frac{\beta _{1}}{500},\frac{\beta _{1}}{50}\right] $ and let 
$\{\eta _{i}\}_{i=1}^{L}$ be a sequence that satisfies the conclusion of the
proposition above. To anchor the induction, we then take 
\begin{wrapfigure}[17]{r}{0.4\textwidth}\centering
\vspace{-8pt}
\includegraphics[scale=0.14]{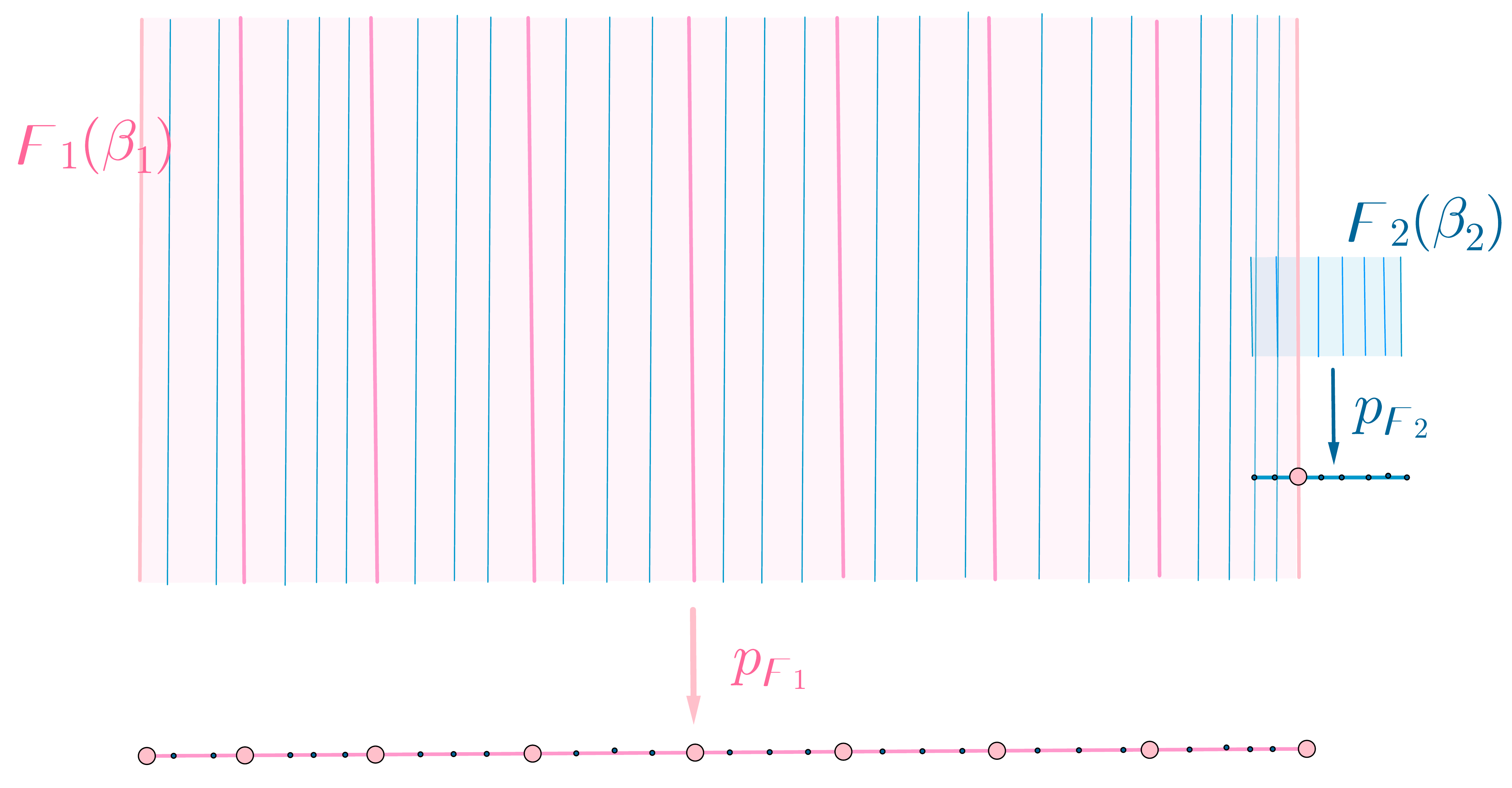}
\caption{ {\protect\tiny {\ The pink and blue rectangles are $2$%
--dimensional depictions of $1$--framed sets. The pink subset of the
interval is the vertex set of an $\frac{\protect\beta _{1}}{8}$--geometric
graph in $p_{\digamma _{1}}(\digamma _{1}(\protect\beta _{1}))$, but only
one of its fibers intersects $\digamma _{2}(\protect\beta _{2})$. Thus we
add the blue points to the pink points to obtain a geometric graph that is
at an appropriate scale for the blue set. If the picture were realistic, the
sides of the pink rectangle would be $10$ rather than $5$ times longer than
the sides of the blue rectangle (cf Figure \protect\ref{etas fig} for the
analogous construction within 2-framed sets in Part \protect\ref{Part 2}).} }
}
\label{candy stripe fig}
\end{wrapfigure}an $\eta _{1}$--geometric graph that covers $p_{\digamma
_{1}}\left( \digamma _{1}\right) $ and apply Corollary \ref{cover a one
framed set cor} to get an associated standard local $\eta _{1}$--cograph $%
\mathcal{G}^{1}$ that covers $\digamma _{1}.$

For an outline of the induction step, assume $\mathcal{G}^{l-1}$ is $\eta
_{l-1}$-geometric. In the event that $\eta _{l}<\eta _{l-1}$, we use the
Subdivision Lemma to replace $\mathcal{G}^{l-1}$ with an $\eta _{l}$%
-geometric cograph $\tilde{\mathcal{G}}^{l-1}$ subordinate to $\left\{
\digamma _{i}(\beta _{i}^{+})\right\} _{i=1}^{l-1}$ (see Figure \ref{candy
stripe fig}). We then use the Extension Lemma to extend $\tilde{\mathcal{G}}%
^{l-1}$ to an $\eta _{l}$-geometric cograph $\mathcal{G}^{l}$ subordinate to 
$\left\{ \digamma _{i}(\beta _{i}^{+})\right\} _{i=1}^{l}$. If ${\eta }%
_{l}=\eta _{l-1},$ the procedure is the same except that we skip the
subdivision step.

Roughly, the Extension Lemma starts with an $\eta _{l}$-geometric cograph $%
\tilde{\mathcal{G}}^{l-1}$ subordinate to $\{\digamma _{i}\}_{i=1}^{l-1}$
and extends it to an $\eta _{l}$-geometric cograph $\tilde{\mathcal{G}}^{l}$
subordinate to $\{\digamma _{i}\}_{i=1}^{l-1}\cup \{\digamma _{l}\}$. To do
this extension, we project $\tilde{\mathcal{G}}^{l-1}$ to $p_{\digamma
_{l}}(\digamma _{l})$ using $p_{\digamma _{l}}$. Since $p_{\digamma _{l}}$
is $\varkappa $-almost Riemannian and the sets in $\mathcal{A}$ are $\kappa $%
-lined up, this yields an $\left( \eta _{l},\delta \eta _{l}\right) $%
--geometric graph $\mathcal{T}^{\mathrm{pr,}l},$ where the error, $\delta $,
depends on $\varkappa $ and $\kappa .$ Since the cardinality $L$ of $%
\mathcal{A}$ has no known relationship with $\varkappa $ and $\kappa ,$ for
our induction argument, we will need a new method to verify that the errors
do not accumulate. This new technique will allow us to extend $\mathcal{T}^{%
\mathrm{pr,}l}$ to a graph $\mathcal{T}^{l}\subset p_{\digamma
_{l}}(\digamma _{l})$ for which the collection of new vertices is error
free. We call this extension an \textbf{error correcting extension}, and in
the Extension Lemma, we make the additional hypothesis and conclusion that $%
\tilde{\mathcal{G}}^{l-1}$ and $\tilde{\mathcal{G}}^{l}$ are \textbf{error
correcting cographs. }This is a term that we define in the next subsection.

In Subsection \ref{subdivided cogrph subsect}, we state and prove the
Subdivision Lemma, and we wrap up the proof of Theorem \ref{cograph constr
them} in Subsection \ref{constr cogrphs subsect}.

\addtocounter{algorithm}{1}

\subsection{Error Correcting Cographs and the Extension Lemma\label%
{constrcogrphs subsect}}

In this subsection, we define the term \textbf{error correcting cograph }and%
\textbf{\ }state and prove the \textbf{Extension Lemma}. We start with error
correcting extensions of geometric graphs.

\begin{definition}
Let $\mathcal{T}$ and $\mathcal{\tilde{T}}$ be $\left( \eta ,\delta \eta
\right) $--geometric graphs with $\mathcal{T}\subset \mathcal{\tilde{T}}$
and $\left\vert \mathcal{T}\right\vert ,|\mathcal{\tilde{T}}|\subset \left[
a,b\right] .$ We say that $\mathcal{\tilde{T}}$ is an \textbf{error
correcting extension} of $\mathcal{T}$ provided $\mathcal{\tilde{T}}_{0}$ is
a subset of $\left[ a,b\right] $ so that%
\begin{equation*}
\left[ a,b\right] \subset B\left( \mathcal{\tilde{T}}_{0},\eta \right) \text{
and}
\end{equation*}%
\begin{equation}
\text{ }\mathrm{dist}\left( \mathfrak{v},\mathfrak{w}\right) \geq \eta
\label{error
cor ext}
\end{equation}%
for all $\mathfrak{v\in }\mathcal{\tilde{T}}_{0}\ $and all $\mathfrak{w\in }%
\mathcal{\tilde{T}}_{0}\setminus \mathcal{T}_{0}.$
\end{definition}

By extending $\mathcal{T}_{0}$ to a maximal subset of $\left[ a,b\right] $
that satisfies (\ref{error cor ext}), we get

\begin{proposition}
\label{peerf ext exists propr}Let $\mathcal{T}$ be an $\left( \eta ,\delta
\eta \right) $--geometric graph with $\left\vert \mathcal{T}\right\vert
\subset \left[ a,b\right] .$ Then there is an error correcting extension $%
\mathcal{\tilde{T}}$ of $\mathcal{T}.$
\end{proposition}

Before proceeding further we detail the process mentioned above wherein we
project a geometric cograph to get an $\left( \eta ,\delta \eta \right) $%
--geometric graph$\mathrm{.}$ To do so, let $\mathcal{G}=\mathcal{V}\cup 
\mathcal{E}$ be an $\eta $--geometric cograph that is subordinate to $%
\left\{ \digamma _{i}(\beta _{i}^{+})\right\} _{i=1}^{l-1}$. Call the set of
vertex fibers $\mathrm{VF},$ that is, 
\begin{equation*}
\mathrm{VF}:=\bigcup\limits_{i=1}^{l-1}\left\{ \left. q_{\mathfrak{v}%
}^{-1}\left( \mathfrak{v}\right) \text{ }\right\vert \mathfrak{v}\in (K_{%
\mathcal{T}^{i}})_{0}\right\} \text{,}
\end{equation*}%
where $q_{\mathfrak{v}}$ is the submersion that defines the covertex that
corresponds to $\mathfrak{v}.$ Let 
\begin{equation*}
\mathcal{T}_{0}^{\mathrm{pr,}l}:=p_{\digamma _{l}}\left( \mathrm{VF}\right) 
\text{ (see, e.g., the pink dot in the blue interval in Figure \ref{candy
stripe fig}).}
\end{equation*}

The elements of $\mathcal{T}_{0}^{\mathrm{pr,}l}$ are not points. However,
since our sets are $\kappa $-lined up, the elements of $\mathcal{T}_{0}^{%
\mathrm{pr,}l}$ all have diameter less than $\tau \left( 2\kappa \beta
_{l}\right) ,$ so to avoid excess technicalities, we will treat the elements
of $\mathcal{T}_{0}^{\mathrm{pr,}l}$ as though they were points.

Let $\mathcal{T}^{\mathrm{pr,}l}$ be the graph whose vertices are $\mathcal{T%
}_{0}^{\mathrm{pr,}l}$ and in which two vertices are connected by an edge if
and only if the corresponding covertices intersect a common coedge. We can
now define the concept of an error correcting cograph.

\begin{definition}
\label{perfe cogrph dfn}Let $\mathcal{A}=\{\digamma _{i}\}_{i=1}^{L}$ be a
collection of $(1,\varkappa )$-framed sets that are $\kappa $--lined up,
vertically separated, and fiber swallowing. For $l\in \left\{ 1,\ldots
,L\right\} $ and $k\geq l,$ let $\mathcal{G}=\mathcal{V}\cup \mathcal{E}$ be
a respectful, $\eta _{k}$--geometric cograph subordinate to $\{\digamma
_{i}\}_{i=1}^{l}$. We say that $\mathcal{G}$ is \textbf{error correcting}
provided

\begin{enumerate}
\item For each local cograph $\mathcal{G}_{j_{0}}$, the associated graph $%
\mathcal{T}^{j_{0}}$ is an error correcting extension of $\mathcal{T}^{%
\mathrm{pr,}j_{0}}.$

\item The generating set $K_{\mathcal{T}^{j_{0}}}$ for $\mathcal{G}_{j_{0}}$
is%
\begin{equation*}
K_{\mathcal{T}^{j_{0}}}=\mathcal{T}^{j_{0}}\setminus \mathcal{T}^{\mathrm{pr,%
}j_{0}}.
\end{equation*}
\end{enumerate}
\end{definition}

With the definition of error correcting cographs in hand, we can state the
Extension Lemma.

\begin{lemma}
(Extension Lemma)\label{cograph ext lemma} Let $\mathcal{A}:=\left\{
\digamma _{j}(\beta _{j}^{+})\right\} _{j=1}^{L}$ be a collection of $%
(1,\varkappa )$--framed sets that are $\kappa $--lined up, vertically
separated, and fiber swallowing. For $l\in \left\{ 1,\ldots ,L\right\} ,$
let $\mathcal{G}^{l-1}=\mathcal{V}^{l-1}\cup \mathcal{E}^{l-1}$ be an error
correcting, $\eta _{l}$-geometric cograph that is subordinate to $\left\{
\digamma _{i}(\beta _{i}^{+})\right\} _{i=1}^{l-1}.$ Then $\mathcal{G}^{l-1}$
extends to an error correcting, $\eta _{l}$--geometric cograph $\mathcal{G}%
^{l}=\mathcal{V}^{l}\cup \mathcal{E}^{l}$ that is subordinate to $\left\{
\digamma _{i}(\beta _{i}^{+})\right\} _{i=1}^{l-1}\cup \{\digamma _{l}(\beta
_{l}^{+})\}.$

Moreover, if $\mathcal{A}=\mathcal{K}^{1}$ and $\mathcal{G}^{l-1}$ is
relative to $\Omega _{0},$ then we can choose $\mathcal{G}^{l}$ so that it
is relative to $\Omega _{0}.$
\end{lemma}

To begin the proof, let $\mathrm{VF}_{i}\left( \nu \right) $ be the
intersection of $\mathrm{VF}$ with $\digamma _{i}\left( \beta _{i},\nu
\right) $. Let $\mathrm{VT}_{i}$ be the abstract graph whose vertex set is $%
\mathrm{VF}_{i}\left( \nu \right) $ and in which two vertices are connected
by an edge if and only if the corresponding covertices intersect a common
coedge.

Next we apply Lemma \ref{Cor Q disj} and observe that for each $j_{0}\in
\left\{ 1,\ldots ,l\right\} ,$ $\mathrm{VT}_{j_{0}}\left( \beta
_{l}^{+}\right) $ is an abstract geometric graph with respect to the
following metric: For any $q_{\mathfrak{v}}^{-1}\left( \mathfrak{v}\right)
,q_{\mathfrak{\tilde{v}}}^{-1}\left( \mathfrak{\tilde{v}}\right) \in \mathrm{%
VT}_{j_{0}}\left( \beta _{l}^{+}\right) ,$ choose any point in $x_{0}\in q_{%
\mathfrak{v}}^{-1}\left( \mathfrak{v}\right) \cap \digamma _{j_{0}}\left(
\beta _{j_{0}},\beta _{l}\right) $ and define 
\begin{equation*}
\mathrm{dist}\left( q_{\mathfrak{v}}^{-1}\left( \mathfrak{v}\right) ,q_{%
\mathfrak{\tilde{v}}}^{-1}\left( \mathfrak{\tilde{v}}\right) \right) :=%
\mathrm{dist}\left( x_{0},q_{\mathfrak{\tilde{v}}}^{-1}\left( \mathfrak{%
\tilde{v}}\right) \cap \digamma _{i}\left( \beta _{i},5\beta _{l}\right)
\right) .
\end{equation*}%
Of course this definition depends on the many possible choices of $x_{0},$
but since the $q_{\mathfrak{v}}$ are $3\kappa $--lined up, $\varkappa $%
--almost Riemannian submersions, the fibers are almost equidistant. So
regardless of these choices, Lemma \ref{Cor Q disj} together with an
induction argument and the fact that $\frac{\beta _{l}}{500}\leq \eta _{l}$
yields

\begin{corollary}
\label{local whole cogrh prop}If $\mathcal{G}$ is an error correcting, $\eta
_{l}$--geometric cograph that is subordinate to a collection $\left\{
\digamma _{i}^{+}\right\} _{i=1}^{l-1}$ of $\kappa $--lined up, vertically
separated, fiber swallowing $(1,\varkappa )$-framed sets, then for all $%
j_{0}\in \left\{ 1,2,3,\ldots ,l-1\right\} ,$ $\mathrm{VT}_{j_{0}}\left(
\beta _{l}^{+}\right) $ is an abstract $\left( \eta _{l},\delta \eta
_{l}\right) $--graph$,$ for $\delta =\varkappa +5,000\kappa .$
\end{corollary}

We omit the proof since Lemma \ref{perf implies abstr cor} in Part \ref{Part
2} is analogous and a complete proof is given. Combining Corollary \ref%
{local whole cogrh prop} with Lemma \ref{Cor Q disj}, we get

\begin{corollary}
\label{proj cograph propr}Let $\mathcal{A}:=\left\{ \digamma _{j}(\beta
_{j}^{+})\right\} _{j=1}^{L}$ be a collection of $(1,\varkappa )$--framed
sets that are $\kappa $--lined up, vertically separated, and fiber
swallowing. For $l\in \left\{ 1,\ldots ,L\right\} ,$ if $\mathcal{G}$ is an
error correcting, $\eta _{l}$--geometric cograph that is subordinate to $%
\left\{ \digamma _{i}^{+}(\beta _{i})\right\} _{i=1}^{l-1}$, then $\mathcal{T%
}^{\mathrm{pr,}l}$ is an $\left( \eta _{l},2\delta \eta _{l}\right) $%
-geometric cograph where $\delta =\varkappa +5,000\kappa .$
\end{corollary}

\begin{proof}[Proof of the Extension Lemma (\protect\ref{cograph ext lemma})]

Since $\mathcal{G}^{l-1}=\mathcal{V}^{l-1}\cup \mathcal{E}^{l-1}$ is an
error correcting, $\eta _{l}$-geometric cograph, Corollary \ref{proj cograph
propr} gives us that $\mathcal{T}^{\mathrm{pr,}l}$ is $\left( \eta
_{l},2\delta \eta _{l}\right) $-geometric, where $\delta =\varkappa
+5,000\kappa $. By Proposition \ref{peerf ext exists propr}, $\mathcal{T}^{%
\mathrm{pr,}l}$ has an error correcting extension to an $\left( \eta
_{l},2\delta \eta _{l}\right) $--geometric graph $\mathcal{T}^{l}$ that
covers $p_{\digamma _{l}}\left( \digamma _{l}\right) .$

Set 
\begin{equation}
K_{\mathcal{T}^{l}}:=\mathcal{T}^{l}\setminus \mathcal{T}^{\mathrm{pr,}l}.
\label{co grph K}
\end{equation}%
By Corollary \ref{cover a one framed set cor}, there is a standard local
cograph $\mathcal{G}_{l}$ whose cosimplices are in one-to-one correspondence
with $K_{\mathcal{T}^{l}}$ and are defined using the submersion $p_{\digamma
_{l}}$ (see Definition \ref{stand local cographs dfn}). Since $\mathcal{T}%
^{l}$ is an error correcting extension of $\mathcal{T}^{\mathcal{T}^{\mathrm{%
pr,}l}}$, Lemma \ref{Cor Q disj} gives us that 
\begin{eqnarray}
&&\text{the intersection of }\mathcal{G}^{l-1}\mathfrak{\cup }\mathcal{G}_{l}%
\text{ with }\digamma _{l}^{+}\text{ is no more than }\tau \left( \kappa
\eta _{l}\right) \text{ from a standard}  \notag \\
&&\text{local cograph associated to the }\eta _{l}\text{--geometric graph }%
\mathcal{T}^{l}.  \label{is stand cogrph}
\end{eqnarray}

To show that 
\begin{equation*}
\mathcal{G}^{l}:=\mathcal{G}^{l-1}\cup \mathcal{G}_{l}=\mathcal{V}^{l-1}\cup 
\mathcal{V}_{l}\cup \mathcal{E}^{l-1}\cup \mathcal{E}_{l}
\end{equation*}%
is a cograph$,$ we combine the fact that $\mathcal{G}^{l-1}$ is a cograph
with the following additional arguments.

\begin{proof}[Proof of Property 1 of Definition \protect\ref{cograph dfn}: ]

Since $\mathcal{G}^{l-1}$ is a cograph and $\mathcal{G}_{l}$ is a local
cograph, Property $1$ holds.
\end{proof}

\begin{proof}[Proof of Property 2 of Definition \protect\ref{cograph dfn}:]
The collections $\mathcal{V}^{l-1}$ and $\mathcal{E}^{l-1}$ are disjoint, by
hypothesis. So it suffices to show that the collections $\mathcal{V}_{l}$
and $\mathcal{E}_{l}$ are disjoint from each other and from $\mathcal{V}%
^{l-1}$ and $\mathcal{E}^{l-1},$ respectively. Since each element of $%
\mathcal{V}_{l}$ and $\mathcal{E}_{l}$ is contained in $\digamma _{l}^{+},$
it is enough to consider the intersection of $\mathcal{G}^{l-1}\cup \mathcal{%
G}_{l}$ with $\digamma _{l}^{+}.$

Since $\eta _{l}\geq \frac{\beta _{l}}{500},$ (\ref{is stand cogrph})
together with Corollary \ref{dsij from o stra cor} gives us that the
intersections of each of $\mathcal{V}^{l-1}\cup \mathcal{V}_{l}$ and $%
\mathcal{E}^{l-1}\cup \mathcal{E}_{l}$ with $\digamma _{l}^{+}$ is a
disjoint collection, and hence each of $\mathcal{V}^{l-1}\cup \mathcal{V}%
_{l} $ and $\mathcal{E}^{l-1}\cup \mathcal{E}_{l}$ is a disjoint collection.
\end{proof}

\begin{proof}[Proof of Property 3 of Definition \protect\ref{cograph dfn}:]
By hypothesis, $\mathcal{G}^{l-1}$ covers $\cup _{i=1}^{l-1}\digamma
_{i}(\beta _{i}^{-}).$ Since $|\mathcal{T}^{l}|$ covers $p_{\digamma
_{l}}\left( \digamma _{l}\right) $, Corollary \ref{cover a one framed set
cor}, (\ref{co grph K}), and (\ref{is stand cogrph}) imply that $\mathcal{G}%
^{l}=\mathcal{G}^{l-1}\cup \mathcal{G}_{l}$ covers $\digamma _{l}(\beta
_{l}^{-})$.
\end{proof}

Thus $\mathcal{G}^{l}$ is a cograph, which by construction, is $\eta _{l}$%
-geometric and error correcting. It remains to prove that if $\mathcal{A}=%
\mathcal{K}^{1}$ and $\mathcal{G}^{l-1}$ is relative to $\Omega _{0},$ then
we can choose $\mathcal{G}^{l}$ to be relative to $\Omega _{0}.$ If $%
\digamma _{l}\cap \Omega _{0}=\emptyset ,$ there is nothing to prove.
Otherwise, by Parts 2 and 3 of Lemma \ref{art Ms Pacman lemma}, $\digamma
_{l}\left( \beta _{l}\right) \cap \bigcup\limits_{i=1}^{l-1}\digamma
_{i}\left( \beta _{i}^{+}\right) =\varnothing $, and $p_{\digamma }\left( 
\bar{\digamma}\cap \bar{\Omega}_{0}^{-}\right) =\left\{ 0\right\} .$ So our
construction gives that $\mathcal{G}^{l}$ is relative to $\Omega _{0}$
provided we choose $\mathcal{T}^{l}$ so that $0\in \mathcal{T}_{0}^{l}.$
This is possible since $\digamma _{l}\left( \beta _{l}\right) \cap
\bigcup\limits_{i=1}^{l-1}\digamma _{i}\left( \beta _{i}^{+}\right)
=\varnothing .$
\end{proof}

\addtocounter{algorithm}{1}

\subsection{Subdividing Cographs\label{subdivided cogrph subsect}}

In this subsection, we prove the cograph subdivision lemma.

\begin{lemma}[Subdivision Lemma]
\label{cographsubdiv lemma} Let $\mathcal{A}:=\left\{ \digamma _{j}(\beta
_{j}^{+})\right\} _{j=1}^{L}$ be a collection of $1$--framed sets that are $%
\kappa $--lined up, vertically separated, and fiber swallowing. Let 
\begin{equation*}
\mathcal{G}:=\mathcal{V}\cup \mathcal{E}=\cup _{j}\mathcal{G}%
_{j}=\bigcup\limits_{j}\mathcal{V}_{j}\cup \mathcal{E}_{j}
\end{equation*}%
be an error correcting, $\eta _{l-1}$--cograph that is subordinate to $%
\left\{ \digamma _{j}(\beta _{j}^{+})\right\} _{j=1}^{l-1}$, where $l\in
\left\{ 1,\ldots ,L\right\} .$ Let $\mathcal{T}_{j}$ be the graph that
generates $\mathcal{G}_{j}$.

There is a subdivision $\mathcal{L}\left( \mathcal{T}_{j}\right) $ of $%
\mathcal{T}_{j}$ that is legal and generates a local cograph $\mathcal{L}%
\left( \mathcal{G}_{j}\right) $ so that $\mathcal{L}\left( \mathcal{G}%
\right) :=\cup _{j}\mathcal{L}\left( \mathcal{G}_{j}\right) $ is a legal,
error correcting $\eta _{l}$--cograph subordinate to $\left\{ \digamma
_{j}(\beta _{j}^{+})\right\} _{j=1}^{l-1}$.

Furthermore, if $\mathcal{A=K}^{1}$ and $\mathcal{G}$ is relative to $\Omega
_{0},$ then we can choose $\mathcal{L}\left( \mathcal{G}\right) $ to be
relative to $\Omega _{0}.$
\end{lemma}

When $\mathcal{G}$ and $\mathcal{L}\left( \mathcal{G}\right) $ are related
as in the Cograph Subdivision Lemma, we will say that $\mathcal{L}\left( 
\mathcal{G}\right) $ is a legal subdivision of $\mathcal{G}.$ The proof of
the Cograph Subdivision Lemma uses the following result from \cite{ProWilh2}.

\begin{lemma}[Lemma 3.6 and Corollary 4.25 of \protect\cite{ProWilh2}]
\label{perfect and subdiv cor} Let $\mathcal{T}$ be an $\left( \eta ,\delta
\eta \right) $--geometric graph with $|\mathcal{T}|\subset \left[ 0,\beta %
\right] .$ There is a legal $\left( \frac{\eta }{19},\delta \frac{\eta }{10}%
\right) $--geometric graph $\mathcal{L}\left( \mathcal{T}\right) $ that
subdivides $\mathcal{T}.$ Moreover, if $\mathcal{\tilde{T}}$ is an error
correcting extension of $\mathcal{T},$ then we can choose legal subdivisions
of $\mathcal{T}$ and $\mathcal{\tilde{T}}$ so that $\mathcal{L}\left( 
\mathcal{\tilde{T}}\right) $ is an error correcting extension of $\mathcal{L}%
\left( \mathcal{T}\right) .$
\end{lemma}

We also make use of the following notation.

\begin{definition}
\label{vert exp}Let $(\digamma \left( \beta ,\nu \right) ,p_{\digamma
},r_{\digamma })$ be $k$--framed. We set 
\begin{equation*}
\digamma ^{v,-}:=\digamma \left( \beta ,\nu ^{-}\right) \text{ and }\digamma
^{v,+}:=\digamma \left( \beta ,\nu ^{+}\right) \text{.}
\end{equation*}
\end{definition}

\begin{proof}[Proof of Lemma \protect\ref{cographsubdiv lemma}]
The existence of the $\mathcal{L}\left( \mathcal{T}_{j}\right) $ follows
from Lemma \ref{perfect and subdiv cor}.

For each $E\in \mathcal{E}_{j}$ let $\mathfrak{e}_{E}$ be the corresponding
edge. Replace $E$ with standard $\frac{\eta }{10}$--geometric coedges and
covertices that are subordinate to $\digamma _{j}$ and that correspond to
the simplices of $\mathcal{L}\left( \mathcal{T}_{j}\right) $ that cover the
interior of $\left\vert \mathfrak{e}_{E}\right\vert .$ The submersion that
defines them is $q_{E}^{\mathrm{ext}},$ where $q_{E}^{\mathrm{ext}}$ is as
in Corollary \ref{extedn coedge submer}.

Given $V\in \mathcal{V}_{j},$ let $\mathfrak{v}_{V}$ be the corresponding
vertex. We replace $V\in \mathcal{V}_{j}$ with a standard $\frac{\eta }{10}$%
--geometric covertex that corresponds to $\mathfrak{v}_{V}$ and is
subordinate to $\digamma _{j}.$ The submersion that defines it is $q_{V}.$

Call the resulting collection $\mathcal{L}\left( \mathcal{G}\right) .$ From
the construction of $\mathcal{L}\left( \mathcal{G}\right) $ and Lemma \ref%
{perfect and subdiv cor}, we see that we can choose the $\mathcal{L}\left( 
\mathcal{T}_{j}\right) $ so that $\mathcal{L}\left( \mathcal{G}\right) $ is
error correcting. Since different fibers of the same submersion are
disjoint, $\mathcal{L}\left( \mathcal{G}\right) $ satisfies the required
disjointness properties.

Notice that $\cup \mathcal{G}^{v,-}\subset \cup \mathcal{L}\left( \mathcal{G}%
\right) .$ So $\mathcal{L}\left( \mathcal{G}\right) $ covers $\cup \mathcal{A%
}^{-}$ because $\cup \mathcal{A}^{-}\subset \cup \mathcal{G}^{v,-}.$ Finally
if $\mathcal{A=K}^{1},$ and $\mathcal{G}$ is relative to $\Omega ^{0},$ then
by construction, $\mathcal{L}\left( \mathcal{G}\right) $ is relative to $%
\Omega _{0}.$
\end{proof}

\addtocounter{algorithm}{1}

\subsection{Proof of\label{constr cogrphs subsect} Theorem \protect\ref%
{cograph constr them}}

\begin{proof}[ Proof of Theorem \protect\ref{cograph constr them}]
For simplicity, we only detail the case when $\mathcal{A=K}^{1}.$ Apply the
Cograph Extension Lemma (\ref{cograph ext lemma}) with $l=1$ and $\mathcal{G}%
^{0}=\emptyset ,$ to get an error correcting, $\eta _{1}$--geometric cograph 
$\mathcal{G}^{1}$ that is subordinate to $\digamma _{1}(\beta _{1}^{+})$ and
relative to $\Omega _{0}$.

For $k\in \{1,\ldots ,l-1\}$, assume there is an error correcting, $\eta
_{k-1}$-geometric cograph $\mathcal{G}^{k-1}$ subordinate to $\left\{
\digamma _{j}(\beta _{j}^{+})\right\} _{j=1}^{k-1}$ which is relative to $%
\Omega _{0}$. Use Lemma \ref{cographsubdiv lemma} to subdivide $\mathcal{G}%
^{k-1}$ to get an error correcting, $\eta _{k}$-geometric cograph $\tilde{%
\mathcal{G}}^{k-1}$ which is relative to $\Omega _{0}.$

By the Cograph Extension Lemma (\ref{cograph ext lemma}), $\tilde{\mathcal{G}%
}^{k-1}$ extends to the desired error correcting geometric cograph $\mathcal{%
G}^{k}.$ By induction, this produces the framed cover that proves Theorem %
\ref{cograph constr them} and thus Covering Lemma \ref{Cove lem 1}.
\end{proof}

\part{Covering the Natural $2$--stratum\label{Part 2}}

In this, part we establish the first conclusion of Limit Lemma \ref{cover of
X lemma} by proving

\begin{covlem}
\label{cov em 2}Given $\delta >0$ there is a $\varkappa >0$ and a
respectful, $\varkappa $--framed cover $\mathcal{B}=\mathcal{B}^{0}\cup 
\mathcal{B}^{1}\cup \mathcal{B}^{2}$ of $X\setminus X_{\delta }^{4}$ so that
for $k=0,1,$ or $2,$ the closures of the elements of $\mathcal{B}^{k}$ are
pairwise disjoint.
\end{covlem}

We establish the second conclusion of Limit Lemma \ref{cover of X lemma} and
finish the proof of Theorem \ref{finiteness} in the appendix. To explain our
strategy to prove Covering Lemma \ref{cov em 2}, let 
\begin{equation*}
\mathcal{A}=\mathcal{A}^{0}\cup \mathcal{A}^{1}\cup \mathcal{A}^{2}
\end{equation*}%
be as in Covering Lemma \ref{Cove lem 1} with $\mathcal{A}^{2}=\left\{
\digamma _{j}\right\} _{j=1}^{\tilde{L}}.$ The goal is to replace $\mathcal{A%
}$ with a respectful framed cover 
\begin{equation*}
\mathcal{B}=\mathcal{B}^{0}\cup \mathcal{B}^{1}\cup \mathcal{B}^{2}
\end{equation*}%
so that each collection $\mathcal{B}^{i}$, $i=0,1,$ or $2$, is a disjoint
collection of $i$-framed sets. Since the collections $\mathcal{A}^{0}$ and $%
\mathcal{A}^{1}$ from Covering Lemma \ref{Cove lem 1} are each disjoint and
respectful, the focus is on separating the elements of $\mathcal{A}^{2}$
using additional $0$- and $1$-framed sets. In analogy with Part 1,%
\begin{wrapfigure}[7]{r}{0.33\textwidth}\centering
\vspace{-16pt}
\includegraphics[scale=.24]{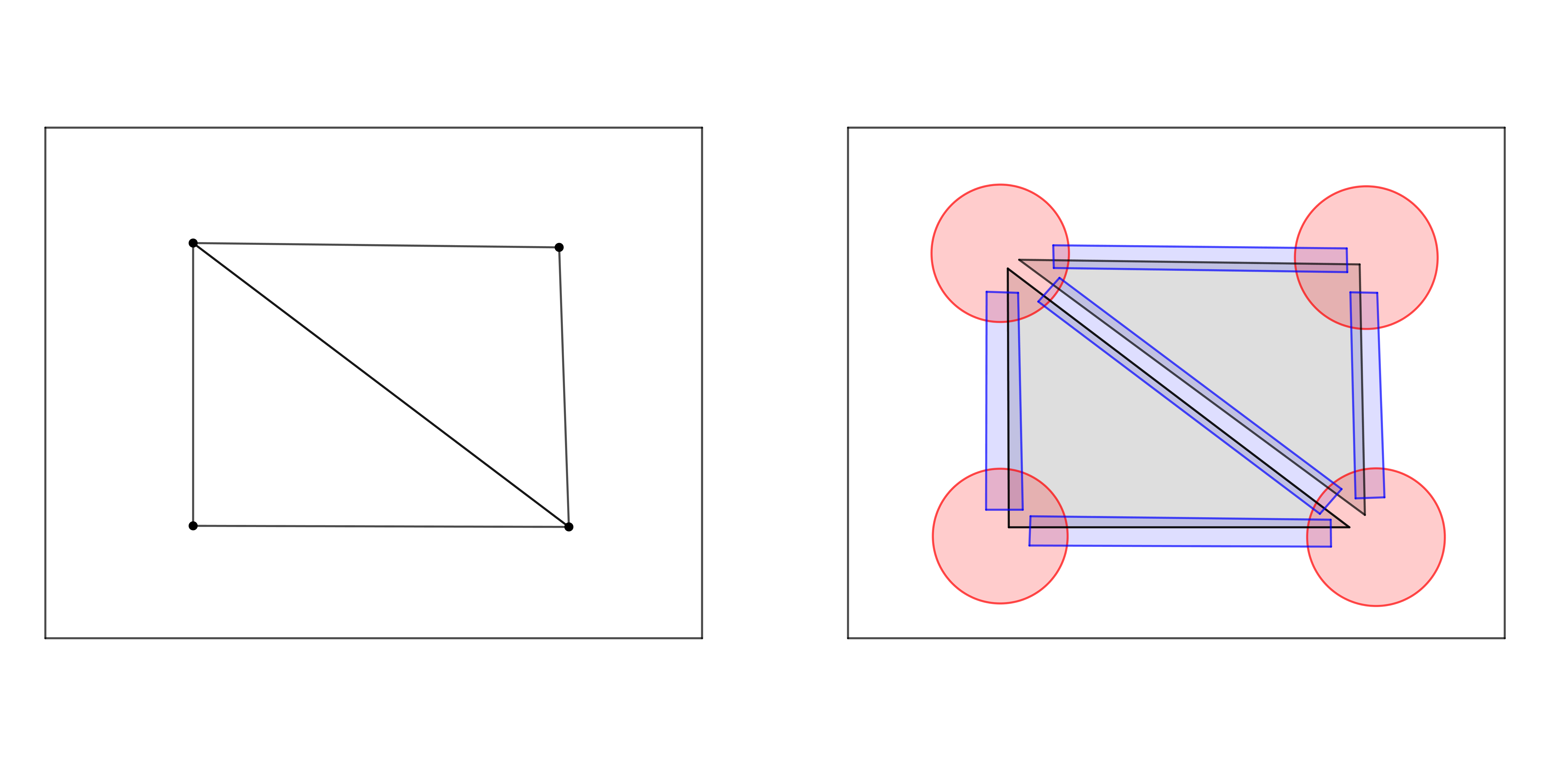}
\vspace{-16pt}
\caption{A triangulation with a corresponding triangulation cover}
\label{tricov}
\end{wrapfigure} these extra $0$- and $1$-framed sets are artificially
framed subsets of the elements of $\mathcal{A}^{2}$ and form a combinatorial
structure that we call a \textbf{cotriangulation.} Roughly speaking,
cotriangulations are analogs of cographs, except the starting point is a
triangulation in $\mathbb{R}^{2}$ rather than a graph in $\mathbb{R}.$

\begin{wrapfigure}[17]{l}{0.4\textwidth}\centering
\vspace{-5pt}
\includegraphics[scale=.26]{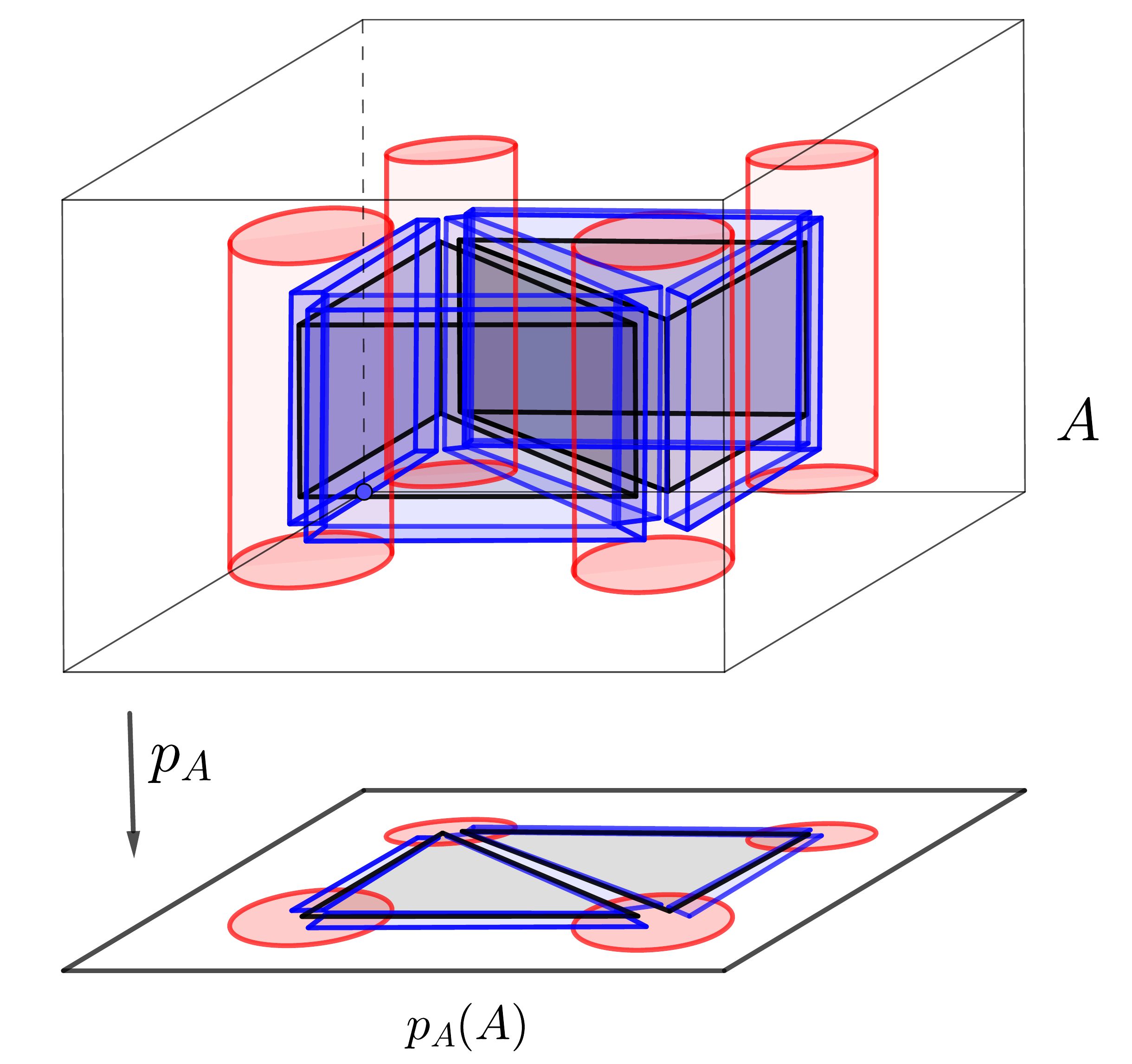}
\vspace{3pt}
\caption{{\protect\tiny {A cotriangulation in a 3--dimensional Alexandrov
space. The red cylinders are the covertices. The coedges are the thin blue
boxes.}}}
\label{cotri figure}
\end{wrapfigure}

More specifically, for each $2$-framed set $(A,p_{A})\in \mathcal{A}^{2},$
we take a triangulation $\mathcal{T}$ in $p_{A}\left( A^{+}\right) \subset 
\mathbb{R}^{2}$. We then construct an open cover of $\left\vert \mathcal{T}%
\right\vert $ which we call a \textbf{triangulation cover} (see Figure \ref%
{tricov}). The open sets of a triangulation cover are disjoint open balls
around each vertex (the red balls in Figure \ref{tricov}), disjoint tubular
neighborhoods of shrunken versions of the edges (the thin blue rectangles in
Figure \ref{tricov}), and disjoint open sets that are shrunken versions of
the faces (the light grey triangles in Figure \ref{tricov}). We then
construct a collection of artificially $0$-, $1$-, and $2$-framed subsets of 
$A$ by taking $p_{A}$--preimages of the open sets in the triangulation
cover. We call this collection a \textbf{cotriangulation}. We call its
elements \textbf{cosimplices}, and we call the $0$-, $1$-, and $2$-framed
sets \textbf{covertices, coedges, }and \textbf{cofaces,} respectively (see
Figure \ref{cotri figure}). Together with an inductive argument, this will
yield a cover of $\cup _{A\in \mathcal{A}^{2}}A^{-}$ with the desired
disjointness properties. 

To create a framed cover that also is respectful, we will apply the
Submersion Deformation Lemma (\ref{const rank def}). Doing this requires
that our collections of cosimplices satisfy Inequality (\ref{Sub def hyp}).
In other words, the quantities $\kappa $ and $\varkappa $ need to be small
compared with the constant $C$ in Inequality (\ref{Sub def hyp}). The
constant $C$ is determined by our cotriangulation, and yet in our
construction we choose $\kappa $ and $\varkappa $ before constructing $%
\mathcal{A}^{2}$. So, to construct a cotriangulation to which we can apply
the Submersion Deformation Lemma, we impose a further constraint on the
triangulations that generate our cotriangulation. The additional constraint
is that our triangulations are what we call \textbf{CDG}s, which is short
for \textbf{Chew-Delaunay-Gromov}. In \cite{ProWilh2}, we prove that CDGs
can be extended, subdivided, and have smallest angle $\geq \pi /6$. The
relevant definitions and results are restated in Subsection \ref{Delaunay
sect}.

In this part, we combine these properties with an inductive argument to
construct a cotriangulation that proves Covering Lemma \ref{cov em 2}.

Throughout Part \ref{Part 2}, there is a strong analogy with Part \ref{Part
1}, as indicated in the following table.%
\begin{equation*}
\begin{tabular}{|l|l|}
\hline
Part 1 concept & Part 2 analog \\ \hline
graph & triangulation \\ \hline
graph cover & triangulation cover, Definition \ref{dfn trig cover} \\ \hline
$\eta $--geometric graph & $\eta $--CDG triangulation, Definition \ref{CDG
dfn} \\ \hline
$\eta $--geometric graph cover & geometric triangulation cover, Definition %
\ref{geom triag cover} \\ \hline
$\eta $--geometric cograph & geometric cotriangulation, Definition \ref{Dfn
Grm Del} \\ \hline
error correcting cograph & error correcting cotriangulation, Definition \ref%
{perfectt dfn} \\ \hline
\end{tabular}%
\end{equation*}

Besides the definition of a geometric cotriangulation, there are three main
ingredients in the proof of Covering Lemma \ref{cov em 2}:%
\begin{equation*}
\begin{tabular}{|l|l|}
\hline
Result & Shows how to \\ \hline
Subdivision Lemma (\ref{subdivide co prop}) & subdivide a geometric
cotriangulation. \\ \hline
Lemma \ref{Curtiss Cor} & deform a geometric cotriangulation to a respectful
one. \\ \hline
Extension Lemma (\ref{rel cotri ext lemma}) & extend a geometric
cotriangulation to an additional $2$--framed set. \\ \hline
\end{tabular}%
\end{equation*}

Part \ref{Part 2} is divided into four sections. In Section \ref%
{cotriangulations sect}, we define cotriangulations and discuss some purely
topological/set theoretic notions that are related to them. In Section \ref%
{Part 1.5}, we show the cover from Covering Lemma \ref{Cove lem 1} can be
constructed in a way that will facilitate the meshing of our cotriangulation
with the natural $0$ and $1$ strata. Section \ref{tools to constr cotri sect}
develops the notion of geometric cotriangulations, and in Section \ref{Sect
Disresp cotri}, we construct the cotriangulation that proves Covering Lemma %
\ref{cov em 2}.

\section{Cotriangulations\label{cotriangulations sect}}

In Subsection \ref{tri covers subsect}, we define triangulation covers,
which are precursors of cotriangulations, and in Subsection \ref{what is a
cotri subsect}, we\ define cotriangulations.

\addtocounter{algorithm}{1}

\subsection{Triangulation Covers\label{tri covers subsect}}

The starting point for the definition of cotriangulation is a method for
creating an open cover $\mathcal{O}\left( \mathcal{T}\right) $ from a
triangulation $\mathcal{T}$ in $\mathbb{R}^{2}$.\ 

\begin{definition}[see Figure \protect\ref{tricov}]
\label{dfn trig cover} Let $\mathcal{T}$ be a triangulation of a subset $%
\left\vert \mathcal{T}\right\vert $ of $\mathbb{R}^{2}$. A collection of
open sets $\mathcal{O}\left( \mathcal{T}\right) $ is called a \textbf{%
triangulation cover }provided

\begin{enumerate}
\item $\mathcal{O}\left( \mathcal{T}\right) $ covers $\left\vert \mathcal{T}%
\right\vert .$

\item There is a one-to-one correspondence $g_{\mathcal{O},\mathcal{T}}:%
\mathcal{T}\longrightarrow \mathcal{O}\left( \mathcal{T}\right) $ so that
for each $k\in \left\{ 0,1,2\right\} ,$ the collection $g_{\mathcal{O},%
\mathcal{T}}\left( \mathcal{T}_{k}\right) $ consists of disjoint open
subsets of $\mathbb{R}^{2}$.

\item For each vertex $\mathfrak{v},$ $g_{\mathcal{O},\mathcal{T}}\left( 
\mathfrak{v}\right) $ is a neighborhood of $\mathfrak{v}$ that is
diffeomorphic to a $2$--ball.

\item For each $1$--simplex $\mathfrak{e},$ $g_{\mathcal{O},\mathcal{T}%
}\left( \mathfrak{e}\right) $ is a neighborhood of a proper open subset $%
\mathfrak{e}^{-}$ of $\mathfrak{e}$ so that $\mathfrak{e}^{-}$ is
diffeomorphic to an interval, and $g_{\mathcal{O},\mathcal{T}}\left( 
\mathfrak{e}\right) $ is a trivial line bundle over $\mathfrak{e}^{-}.$

\item For each $2$--simplex $\mathfrak{f},$ $g_{\mathcal{O},\mathcal{T}%
}\left( \mathfrak{f}\right) $ is a proper subset of $\mathfrak{f}$ that is
diffeomorphic to a $2$--ball.
\end{enumerate}

We say that $\mathcal{O}\left( \mathcal{T}\right) $ is \textbf{associated to}
$\mathcal{T},$ and that $\mathcal{T}$ \textbf{generates} $\mathcal{O}\left( 
\mathcal{T}\right) .$
\end{definition}

\addtocounter{algorithm}{1}

\subsection{What is a Cotriangulation?\label{what is a cotri subsect}}

\begin{figure}[tbp]
\vspace{-10pt} \includegraphics[scale=0.2]{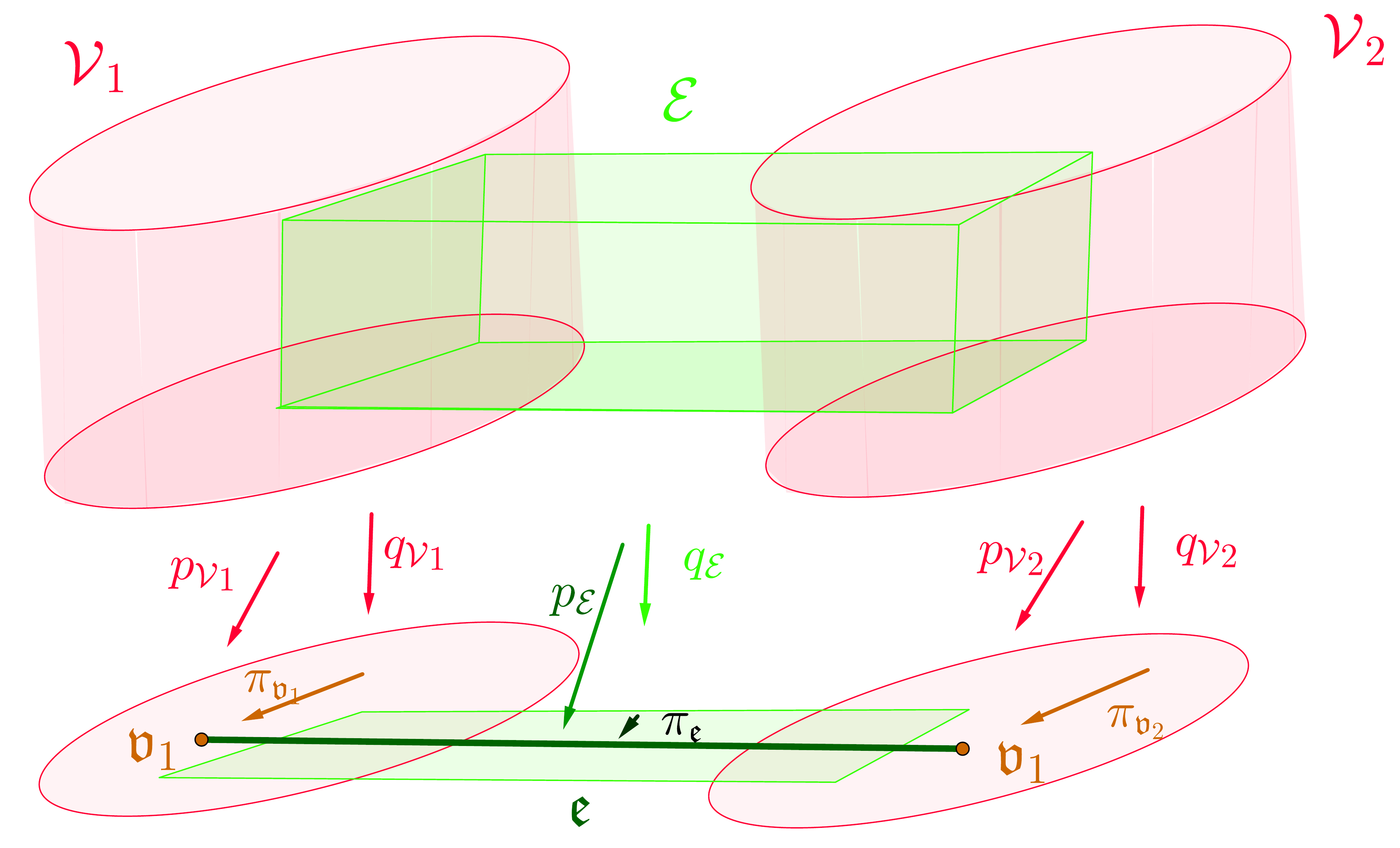} 
\vspace{-10pt}
\caption{{\protect\tiny {The green box $\mathcal{E}$ corresponds to the edge 
$\mathfrak{e}$ and hence is called a coedge. The red cylinders $\mathcal{V}%
_1 $ and $\mathcal{V}_2$ correspond to the vertices $\mathfrak{v}_1$ and $%
\mathfrak{v}_2$, and hence are called covertices. These are two of the three
types of open sets in a cotriangulation. The other type is cofaces (see
Figure \protect\ref{coface}). All three sets are contained in a $2$--framed
set $\digamma$ (not pictured). The three submersions $q_{ \mathcal{V}_1 }$, $%
q_{ \mathcal{E} }$, and $q_{ \mathcal{V}_2 }$ have images in $\mathbb{R}^2$
and are $\protect\kappa$--lined up with $p_{\digamma}$. By composing $q_{ 
\mathcal{V}_1 }$, $q_{ \mathcal{E} }$, and $q_{ \mathcal{V}_2 }$ with the
orthogonal projections $\protect\pi_{ \mathfrak{v}_1 } $, $\protect\pi_{ 
\mathfrak{e} } $, and $\protect\pi_{ \mathfrak{v}_2 } $ onto $\mathfrak{v}_1$%
, $\mathfrak{e} $, and $\mathfrak{v}_2$, we see that $\mathcal{V}_1$, $%
\mathcal{E}$, and $\mathcal{V}_2$ are artificially $0$-, $1$-, and $0$%
--framed, respectively, in the $2$- framed set $\digamma$. }}}
\label{coedge in cotri}
\end{figure}
In this subsection, we define a combinatorial structure associated to a
collection $\mathcal{A}^{2}$ of $2$-framed sets which we call a
cotriangulation. We start with the local notion.

\mbox{}\vspace{1pt}%
\begin{wrapfigure}[11]{r}{0.33\textwidth}\centering
\vspace{-10pt}
\includegraphics[scale=0.13]{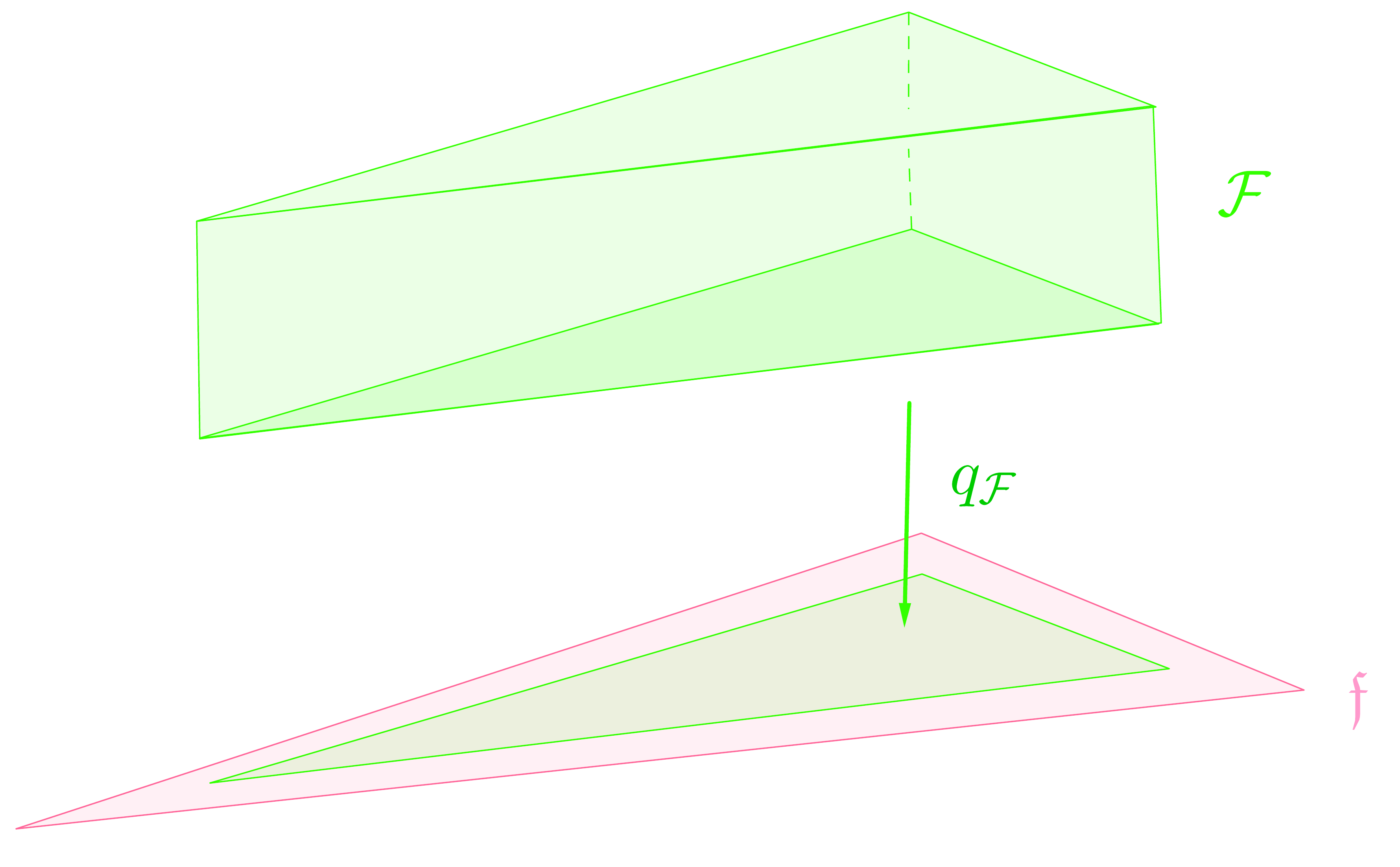}
\vspace{-4pt}
\caption{{\protect\tiny {The green prism, $\mathcal{F}$, corresponds to the
face, $\mathfrak{f} $, and hence is called a coface.  }}}
\label{coface}
\end{wrapfigure}
\vspace{-23pt}

\begin{wrapfigure}[13]{l}{0.47\textwidth}\centering
\vspace{-11pt}
\includegraphics[scale=0.19]{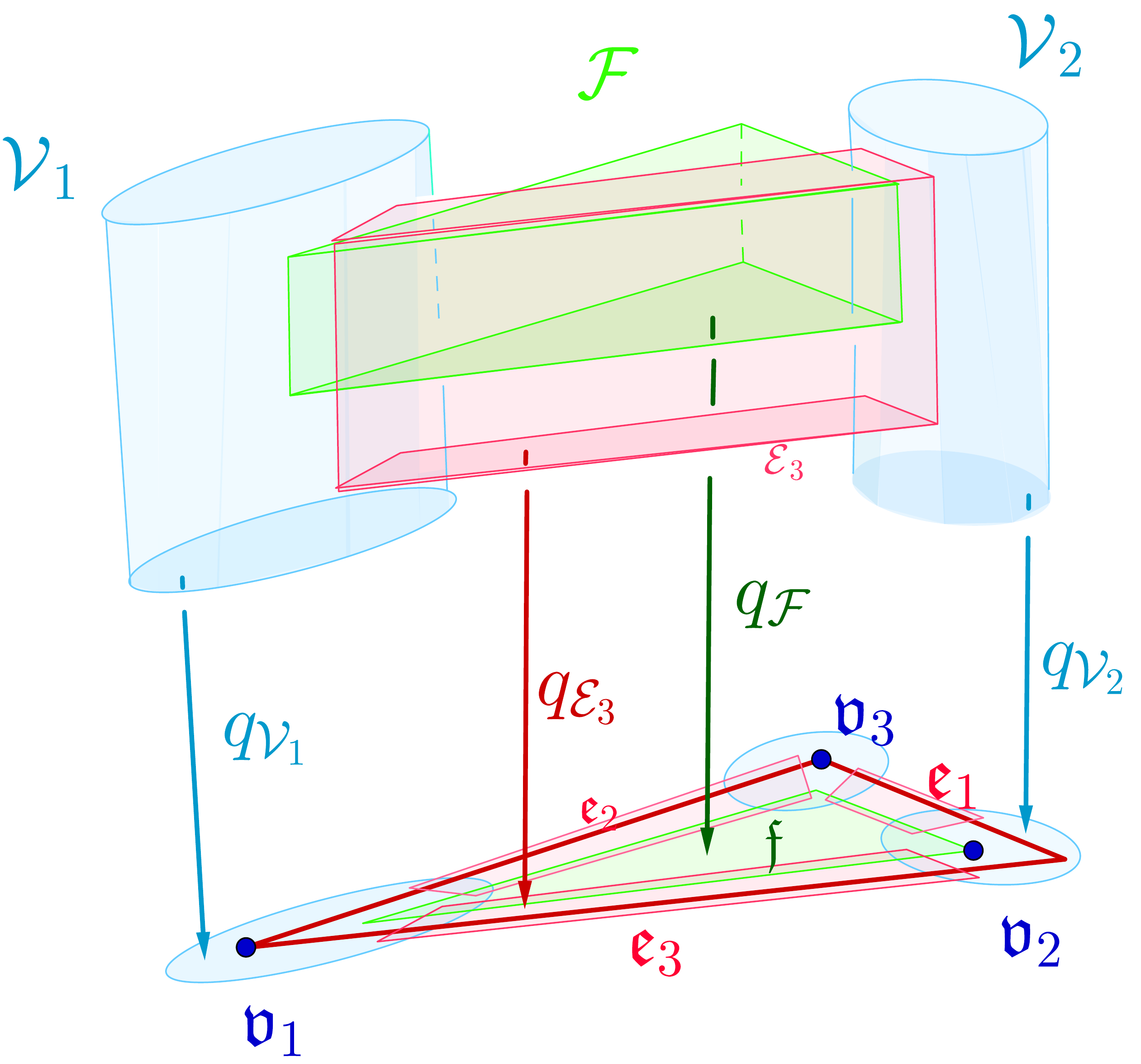}
\vspace{-6pt}
\caption{{\protect\tiny {\ Pictured are $\mathcal{V}_1$, $\mathcal{E }_3$, $%
\mathcal{F}$, and $\mathcal{V}_2$, which are four of the six open sets of a
cotriangulation that correspond to a single $2$--simplex $\mathfrak{f}$.
}}}
\label{cotri one simplex}
\end{wrapfigure}

\begin{definition}[see Figures \protect\ref{cotri figure}, \protect\ref%
{coedge in cotri}, \protect\ref{coface}, and \protect\ref{cotri one simplex}]

\label{local cotri dnf} Let $\digamma $ be a $2$--framed set and $\mathcal{V}%
,$ $\mathcal{E},$ and $\mathcal{F}$ collections of pairwise disjoint open
subsets of $\digamma ^{+}.$ Then 
\begin{equation*}
\hspace*{-2.5in} \mathcal{S}\equiv \mathcal{V}\cup \mathcal{E}\cup \mathcal{F%
}
\end{equation*}%
is a \textbf{local cotriangulation} subordinate to $\digamma $ provided

\vspace{3pt} \noindent 1. There is a triangulation $\mathcal{T}$ with $%
p_{\digamma }\left( \digamma \right) \subset \left\vert \mathcal{T}%
\right\vert \subset p_{\digamma }\left( \digamma ^{+}\right) $ with
associated triangulation cover $O\left( \mathcal{T}\right) $.

\vspace{3pt} \noindent 2. There is a one-to-one correspondence between $%
\mathcal{S}$ and a subset $K_{\mathcal{T}}\subset \mathcal{T}$. Each $S\in 
\mathcal{S}$ is generated by the corresponding $\mathfrak{s\in }K_{\mathcal{T%
}}$ as follows: There is a submersion $q_{S}:$ $\digamma \rightarrow \mathbb{%
R}^{2}$ so that $(S,q_{S},U_{S},\mathfrak{s)}$ is an artificially $k$-framed
set with $\mathfrak{s}\in \mathcal{T}_{k}\cap K_{\mathcal{T}}$, $U_{S}=g_{%
\mathcal{O},\mathcal{T}}\left( \mathfrak{s}\right) $, and 
\begin{equation}
\begin{array}{ll}
\text{if }\mathfrak{s}\in \mathcal{T}_{0}, & \text{then }S=q_{S}^{-1}\left(
U_{S}\right) \cap \digamma (\beta ,\beta ^{+})\in \mathcal{V}, \\ 
\text{if }\mathfrak{s}\in \mathcal{T}_{1}, & \text{then }S=q_{S}^{-1}\left(
U_{S}\right) \cap \digamma (\beta ,\beta )\in \mathcal{E}, \\ 
\text{if }\mathfrak{s}\in \mathcal{T}_{2}, & \text{then }S=q_{S}^{-1}\left(
U_{S}\right) \cap \digamma (\beta ,\beta ^{-})\in \mathcal{F}.%
\end{array}
\label{stand co-simpliex dfn}
\end{equation}
\end{definition}

For $\mathcal{S}$, $\mathcal{T},$ and $K_{\mathcal{T}}$ as above, we say
that $K_{\mathcal{T}}$ generates $\mathcal{S}.$ In general, $K_{\mathcal{T}}$
will be a proper 
subset of $\mathcal{T}$ to allow local cotriangulations to fit together with
each other and with the natural $0$ and $1$ strata (see Figure \ref{cotri
nat str interface figure}).

\begin{definition}[see Figures \protect\ref{cotri figure}, \protect\ref%
{coedge in cotri}, \protect\ref{coface}, and \protect\ref{cotri one simplex}]

\label{Dfn cotri} Let $\mathcal{U}\equiv \{\digamma _{i}(\beta
_{i})\}_{i=1}^{l}$ be a collection of $\left( 2,\varkappa \right) $--framed
sets that are $\kappa $--lined up, vertically separated, and fiber
swallowing. We say a collection 
\begin{equation*}
\hspace*{3in}\mathcal{S}\equiv \mathcal{V}\cup \mathcal{E}\cup \mathcal{F}
\end{equation*}%
of open subsets of $X$ is a \textbf{cotriangulation} \textbf{subordinate to }%
$\mathcal{U}$ provided

\begin{enumerate}
\item For each $i,$ there is a local cotriangulation 
\begin{equation*}
\mathcal{S}_{i}:=\mathcal{V}_{i}\cup \mathcal{E}_{i}\cup \mathcal{F}_{i}
\end{equation*}%
subordinate to $\digamma _{i}$ so that 
\begin{equation*}
\mathcal{V=}\cup _{i}\mathcal{V}_{i},\text{ }\mathcal{E=}\cup _{i}\mathcal{E}%
_{i},\text{ }\mathcal{F=}\cup _{i}\mathcal{F}_{i}\text{, and }\mathcal{S=}%
\cup _{i}\mathcal{S}_{i}.
\end{equation*}

\item Each of the collections $\mathcal{V},$ $\mathcal{E},$ and $\mathcal{F}$
is pairwise disjoint.

\item $\mathcal{S}$ covers $\cup \mathcal{U}^{-}.$
\end{enumerate}
\end{definition}

\begin{definition}
The elements of $\mathcal{V},$ $\mathcal{E}$, $\mathcal{F},$ and $\mathcal{S}
$ are called covertices, coedges, cofaces, and cosimplices, respectively.
\end{definition}

\begin{definition}
We say a cotriangulation is \textbf{respectful} provided it is a respectful
collection of framed sets (see Definition \ref{dfn resp framed collect}).
\end{definition}

At the end of the next section we show that Covering Lemma \ref{cov em 2}
holds provided there is a certain type of cotriangulation. We then state
Lemma \ref{cotr exist lemma}, which asserts that such a cotriangulation
exists.

\section{Interface of the 2-Framed Sets with the Lower Strata\label{Part 1.5}%
}

Before proceeding further with the proof of Covering Lemma \ref{cov em 2},
we constrain the geometry of the $2$--framed sets of Covering Lemma \ref%
{Cove lem 1} and their intersections with the $0$- and $1$-framed sets of
Covering Lemma \ref{Cove lem 1}. These constraints are stated in the next
lemma.

\begin{figure}[h]
\includegraphics[scale=0.21]{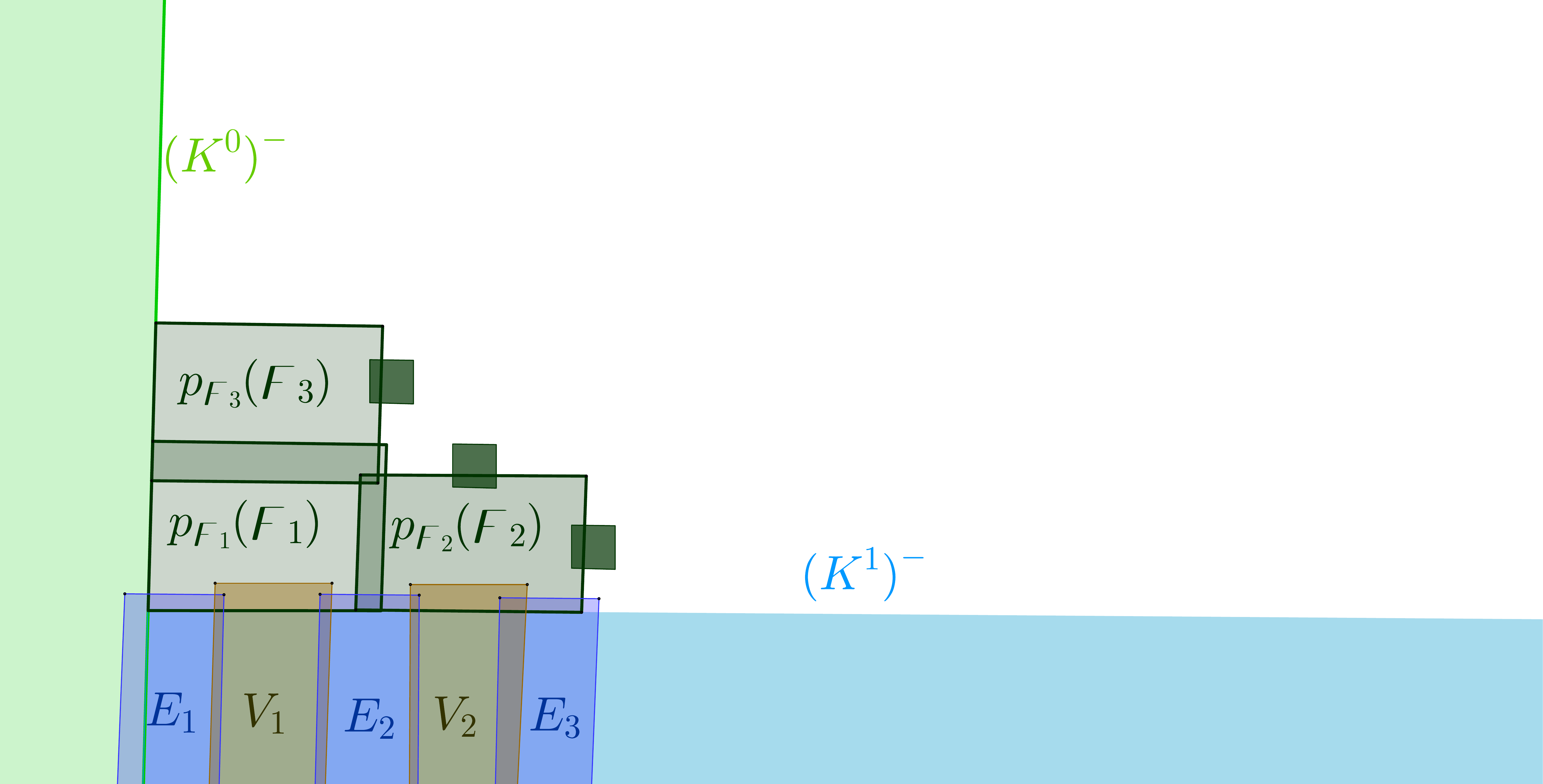}
\caption{{\protect\tiny {This depicts how elements $\digamma _{1}$, $%
\digamma _{2}$, and $\digamma _{3}$ of $\mathcal{A}^{2}$ intersect elements $%
K^{0}$ of $\mathcal{K}^{0}$ and $K^{1}$ of $\mathcal{K}^{1}$, as described
in Lemma \protect\ref{jacked up cov lem 1}. To be more precise, it is a
picture of the image of the configuration under $p_{\digamma_{1}}$, assuming
that $p_{\digamma_{1}} $ is defined in a region this large. Pictured also
are $\{V_{1},V_{2}\} $ and $\{E_{1},E_{2},E_{3}\}$, which are part of a
cograph subordinate to $K^{1}$. The three dark green squares are images of
elements of $\mathcal{A}^{2}$ that intersect neither $\cup \mathcal{K}^{0}$
nor $\cup \mathcal{K}^{1} $. The $\digamma$ pictured thoroughly respect $K^0$%
. They also thoroughly respect the initial cograph described in Lemma 
\protect\ref{jacked up cov lem 1}. They do not, however, thoroughly respect
the cograph depicted here, as it has been subdivided via Proposition \protect
\ref{bottom prop}. }}}
\label{cov lem 1+ fig 1}
\end{figure}

\begin{lemma}[see Figure \protect\ref{cov lem 1+ fig 1}]
\label{jacked up cov lem 1} Let $\mathcal{K}=\mathcal{K}^{0}\cup \mathcal{K}%
^{1}\cup \mathcal{K}^{2}$ be the framed cover from Lemma \ref{art Ms Pacman
lemma}, and let $\mathcal{G}^{L}:\mathfrak{=}\mathcal{V}^{L}\cup \mathcal{E}%
^{L}$ be the $\eta _{L}$-geometric cograph subordinate to $\mathcal{K}^{1}$
that was constructed to prove Theorem \ref{cograph constr them}. There is a
framed cover 
\begin{equation*}
\mathcal{A}:=\mathcal{A}^{0}\cup \mathcal{A}^{1}\cup \mathcal{A}^{2}
\end{equation*}%
that has the following properties:

\begin{enumerate}
\item $\mathcal{A}$ satisfies the conclusion of Covering Lemma \ref{Cove lem
1}.

\item $\mathcal{A}^{0}:=\mathcal{K}^{0}\cup \mathcal{V}^L,$ and $\mathcal{A}%
^{1}:=\mathcal{E}^L.$

\item For each $\left( \digamma \left( \beta _{\digamma }\right)
,p_{\digamma }\right) \in \mathcal{A}^{2},$ $\beta _{\digamma }<\frac{\eta
_{L}}{50}$ and 
\begin{equation}
p_{\digamma }\left( \overline{\digamma (\beta _{\digamma })}\right) =\left[
0,\beta _{\digamma }\right] \times \left[ 0,\beta _{\digamma }\right] .
\label{2 framed images}
\end{equation}

\item $\mathcal{A}^{2}$ thoroughly respects $\mathcal{A}^{0}\cup \mathcal{A}%
^{1}.$

\item If $\digamma \left( \beta _{\digamma }\right) \in \mathcal{A}^{2}$ and 
$\left( K,r_{K}\right) \in \mathcal{K}^{0}$ satisfy $\digamma \left( \beta
_{\digamma }\right) \cap K\neq \emptyset $, then 
\begin{equation*}
p_{\digamma }\left( \overline{\digamma \cap K^{-}}\right) =\left\{ 0\right\}
\times \left[ 0,\beta _{\digamma }\right] ,\text{ }
\end{equation*}%
and $\ r_{K}$ is lined up with $\pi _{1}\circ p_{\digamma },$ where $\pi
_{1}:\mathbb{R}^{2}\longrightarrow \mathbb{R}$ is orthogonal projection to
the first factor.

\item If $\digamma \in \mathcal{A}^{2}$ and $\left( V,r_{V}\right) \in 
\mathcal{V}^{L}$ satisfy $\digamma \left( \beta _{\digamma }\right) \cap
V\neq \emptyset $, then 
\begin{equation}
p_{\digamma }\left( \overline{\digamma \cap V^{-}}\right) =\left[ 0,\beta
_{\digamma }\right] \times \left\{ 0\right\} ,  \label{covert struct}
\end{equation}%
$r_{V}$ is lined up with $\pi _{2}\circ p_{\digamma },$ and $q_{V}$ is lined
up with $\pi _{1}\circ p_{\digamma }$ where $\pi _{i}:\mathbb{R}%
^{2}\longrightarrow \mathbb{R}$ is orthogonal projection to the $i^{th}$%
--factor.

\item If $\digamma \in \mathcal{A}^{2}$ and $\left( E,p_{E},r_{E}\right) \in 
\mathcal{E}^{L}$ satisfy $\digamma \left( \beta _{\digamma }\right) \cap
E\neq \emptyset $, then 
\begin{equation}
p_{\digamma }\left( \overline{\digamma \cap E^{-}}\right) =\left[ 0,\beta
_{\digamma }\right] \times \left\{ 0\right\} ,  \label{coedge
struct}
\end{equation}%
$r_{E}$ is lined up with $\pi _{2}\circ p_{\digamma },$ and $p_{E}$ is lined
up with $\pi _{1}\circ p_{\digamma }.$
\end{enumerate}
\end{lemma}

The proof is a rehash of the ideas of the proof of Lemma \ref{art Ms Pacman
lemma}, so we omit it.

We continue to call the natural $0$-stratum $\Omega _{0}$ and define the
natural $1$-stratum $\Omega _{1}$ to be 
\begin{equation*}
\Omega _{1}:=\cup \mathcal{G}^{L}.
\end{equation*}

\addtocounter{algorithm}{1}

\subsection{The Seam Between the Cotriangulation and the Lower Strata.}

Let $\mathcal{G}^{L}$, $\mathcal{A},$ and $\mathcal{K}$ be the cograph and
framed collections from the lemma above. When a $2$--framed set $\digamma
\in \mathcal{A}^{2}$ intersects $\Omega _{0}\cup \Omega _{1}$, we will need
to understand how its cotriangulation intersects $\Omega _{0}\cup \Omega
_{1}.$ By Part 5 of Lemma \ref{jacked up cov lem 1}, if $\digamma $
intersects $K\in \mathcal{K}^{0},$ then $p_{\digamma }$ takes the
intersection to the left hand edge of $p_{\digamma }\left( \digamma \right) =%
\left[ 0,\beta _{\digamma }\right] \times \left[ 0,\beta _{\digamma }\right]
,$ and by Parts 6, and 7 of Lemma \ref{jacked up cov lem 1} if $\digamma $
intersects $V\in \mathcal{V}^{L}$ or $E\in \mathcal{E}^{L},$ then $%
p_{\digamma }$ takes the intersection to the bottom edge of $p_{\digamma
}\left( \digamma \right) =\left[ 0,\beta _{\digamma }\right] \times \left[
0,\beta _{\digamma }\right] .$ In this subsection, we construct vertex sets
along these edges that later we will extend to triangulations in the $%
p_{\digamma }\left( \digamma \right) =\left[ 0,\beta _{\digamma }\right]
\times \left[ 0,\beta _{\digamma }\right] $ (see Figures \ref{auxillary
vertex sets figure} and \ref{cotri nat str interface figure}).

To fix notation, assume $\mathcal{A}^{2}=\{\digamma _{i}(\beta
_{i})\}_{i=1}^{\tilde{L}}$ is ordered so that $\beta _{i}\geq \beta _{i+1}$.
Choose $\tilde{\eta}_{1}\in \left[ \frac{\beta _{1}}{500},\frac{\beta _{1}}{%
50}\right] $ and let $\{\tilde{\eta}_{i}\}_{i=1}^{\tilde{L}}$ be the
corresponding sequence from the conclusion of Proposition \ref{beta eta
sequence}. We specify the \textquotedblleft bottom\textquotedblright\ and
\textquotedblleft left\textquotedblright\ pieces of $\mathcal{A}^{2}$ to be%
\begin{equation*}
\mathcal{B}:=\left\{ \left. \digamma \in \mathcal{A}^{2}\text{ }\right\vert 
\text{ }\digamma \cap \Omega _{1}\neq \emptyset \right\} \text{ and}
\end{equation*}%
\begin{equation}
\mathcal{L}:=\left\{ \left. \digamma \in \mathcal{A}^{2}\text{ }\right\vert 
\text{ }\digamma \cap \Omega _{0}\neq \emptyset \right\} .
\label{divison
of A2}
\end{equation}%
For $l\in \left\{ 1,2,\ldots ,\tilde{L}\right\} ,$ set 
\begin{equation}
\mathcal{B}_{l}:=\{\digamma _{i}\}_{i=1}^{l}\cap \mathcal{B}\text{ and }%
\mathcal{L}_{l}:=\{\digamma _{i}\}_{i=1}^{l}\cap \mathcal{L}.
\label{B and L dfn}
\end{equation}%
Notice that for $l\in \left\{ 1,2,\ldots ,\tilde{L}\right\} ,$ 
\begin{equation*}
\mathcal{A}_{l}:=\left\{ \left. \left( \digamma ,\pi _{2}\circ p_{\digamma
}\right) \text{ }\right\vert \digamma \in \mathcal{L}_{l}\right\}
\end{equation*}%
is a collection of $1$--framed sets which is $\kappa $--lined up, vertically
separated, and fiber swallowing. So Theorem \ref{cograph constr them} and
its proof give us

\begin{corollary}
\label{left side cor}For each $l\in \left\{ 1,2,\ldots ,\tilde{L}\right\} ,$
there is a legal, error correcting $\tilde{\eta}_{l}$--cograph $\mathcal{LG}%
^{l}$ subordinate to $\mathcal{A}_{l}.$
\end{corollary}

For convenience, for each $(\digamma _{k}(\beta _{k}),\pi _{2}\circ
p_{\digamma _{k}})\in \mathcal{A}_{l}$, we identify 
\begin{equation*}
(\pi _{2}\circ p_{\digamma _{k}})(\digamma _{k})=[0,\beta _{k}]
\end{equation*}%
with $\{0\}\times \lbrack 0,\beta _{k}],$ and we let $\mathcal{L}\mathcal{T}%
^{k}$ be the graph in $\{0\}\times \lbrack 0,\beta _{k}]$ that generates the
local cograph $\mathcal{L}\mathcal{G}_{k}^{l}\subset \mathcal{L}\mathcal{G}%
^{l}$ that is subordinate to $\digamma _{k}.$ Throughout the rest of Part %
\ref{Part 2}, we denote the vertex set $(\mathcal{L}\mathcal{T}^{k})_{0}$ by 
\begin{equation}
L_{k}^{l}:=(\mathcal{L}\mathcal{T}^{k})_{0}\subset \{0\}\times \lbrack
0,\beta _{k}]=(\pi _{2}\circ p_{\digamma _{k}})(\digamma _{k}).
\label{left prop}
\end{equation}
\begin{wrapfigure}[16]{r}{0.4\textwidth}\centering
\vspace{-8pt}
\includegraphics[scale=0.13]{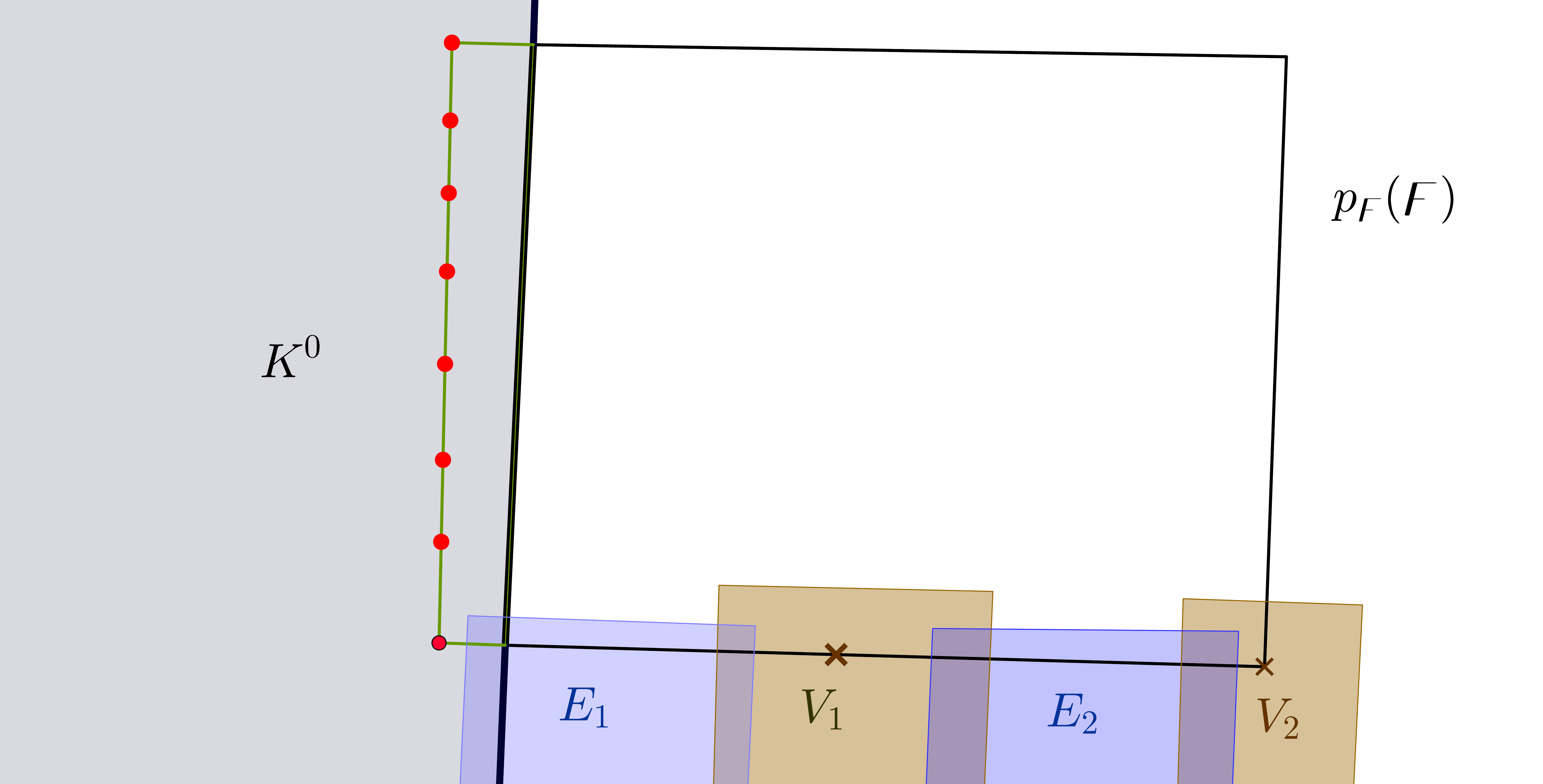}
\caption{{\protect\tiny {Here the $2$--framed set $\digamma $ intersects the
naturally $0$--framed set $K^{0}$ as well as the cograph that covers the
naturally $1$--framed region. The vertex set from (\protect\ref{left prop})
consists of the red points on the left-hand edge. The vertex set from
Proposition \protect\ref{bottom prop} consists of the two brown x's on the
bottom edge.}}}
\label{auxillary vertex sets figure}
\end{wrapfigure}
These vertex sets will play a central role in the interface of our
cotriangulation with the natural $0$--stratum (see Figures \ref{auxillary
vertex sets figure} and \ref{cotri nat str interface figure}).

If $K_{\mathcal{L}\mathcal{T}^{k}}\subset \mathcal{L}\mathcal{T}^{k}$ is the
generating set for $\mathcal{L}\mathcal{G}_{k}^{l}$, and $\mathfrak{v\in }(%
\mathcal{L}\mathcal{T}^{k})_{0}\cap K_{\mathcal{L}\mathcal{T}^{k}}$, we let%
\begin{equation*}
\hspace{-2.5in} q_{\mathfrak{v}}^{l}:\digamma _{k}\longrightarrow \left[
0,\beta _{k}\right]
\end{equation*}%
be the submersion that defines the covertex of $\mathcal{LG}^{l}$ that
corresponds to $\mathfrak{v}$. We then define $q_{\mathfrak{v}}:\digamma
_{k}\longrightarrow \mathbb{R}^{2}$ by 
\begin{equation*}
\hspace{-2.5in} q_{\mathfrak{v}}=\left( \pi _{1}\circ p_{\digamma _{k}},q_{%
\mathfrak{v}}^{l}\right) .
\end{equation*}

The vertex sets $B_{k}^{l}$ in the next result will play an analogous role
in the interface of our cotriangulation with the natural $1$--stratum.

\begin{proposition}[see the two brown $x$'s along the bottom edge in Figure 
\protect\ref{auxillary vertex sets figure}]
\label{bottom prop} Let $\mathcal{A}$, $\mathcal{K}$, and $\mathcal{G}^{L}$
be as in Lemma \ref{jacked up cov lem 1}. If $\mathcal{A}^{2}=\{\digamma
_{i}(\beta _{i})\}_{i=1}^{\tilde{L}}$, there is a choice of $\tilde{\eta}%
_{1}\in \left[ \frac{\beta _{1}}{500},\frac{\beta _{1}}{50}\right] $ so that
if $\{\tilde{\eta}_{i}\}_{i=1}^{\tilde{L}}$ is the corresponding sequence
from Proposition \ref{beta eta sequence}, then for each $l\in \{1,2,\ldots ,%
\tilde{L}\}$, there is a legal, error correcting, $\tilde{\eta}_{l}$%
-geometric cograph 
\begin{equation*}
\tilde{\mathcal{G}}^{l}:=\tilde{\mathcal{V}}^{l}\cup \tilde{\mathcal{E}}^{l}
\end{equation*}%
that has the following properties:

\begin{enumerate}
\item The cograph $\tilde{\mathcal{G}}^{l}$ is a legal subdivision of $%
\mathcal{G}^{L}$.

\item For each $\digamma _{k}\in \mathcal{B}_{l},$ there is a subset $%
B_{k}^{l}$ of $p_{\digamma _{k}}\left( \overline{\digamma _{k}\cap \Omega
_{1}^{-}}\right) =\left[ 0,\beta _{l}\right] \times \left\{ 0\right\} $ with
the property that for each element of $\mathfrak{v\in }B_{k}^{l},$ there is
a $V\in \tilde{\mathcal{V}}^{l}$ so that $\mathfrak{v}$ is the center of the
interval $p_{\digamma _{k}}\left( \partial \overline{V^{-}}\cap \bar{\digamma%
}_{k}\right) .$

\item $B_{k}^{l}$ is the vertex set of a legal $\tilde{\eta}_{l}$--geometric
graph in $\left[ 0,\beta _{k}\right] \times \left\{ 0\right\} .$

\item If $V\in \tilde{\mathcal{V}}^{l}$ and $\digamma \in \mathcal{A}^{2}$
intersect, then $\overline{\digamma \cap V}$ satisfies 
\begin{equation}
p_{\digamma }\left( \overline{\digamma \cap V}\right) \subset I\times \left[
0,\frac{\tilde{\eta}_{l}}{50}\right] ,  \label{covertex vert cut off}
\end{equation}%
where $I\subset \left[ 0,\beta _{\digamma }\right] $ is an interval.

\item If $E\in \tilde{\mathcal{E}}^{l}$ and $\digamma \in \mathcal{A}^{2}$
intersect, then $\overline{\digamma \cap E}$ satisfies 
\begin{equation}
p_{\digamma }\left( \overline{\digamma \cap E}\right) \subset I\times \left[
0,\frac{\tilde{\eta}_{l}}{100}\right] ,  \label{coedge vert cut off}
\end{equation}%
where $I\subset \left[ 0,\beta _{\digamma }\right] $ is an interval.

\item Let%
\begin{equation*}
\mathcal{A}\left( \tilde{\mathcal{G}}^{l}\right) :=\mathcal{A}^{0}\left( 
\tilde{\mathcal{G}}^{l}\right) \cup \mathcal{A}^{1}\left( \tilde{\mathcal{G}}%
^{l}\right) \cup \mathcal{A}^{2},
\end{equation*}%
where 
\begin{equation*}
\mathcal{A}^{0}\left( \tilde{\mathcal{G}}^{l}\right) :=\mathcal{K}^{0}\cup 
\tilde{\mathcal{V}}^{l}\text{ and }\mathcal{A}^{1}\left( \tilde{\mathcal{G}}%
^{l}\right) :=\tilde{\mathcal{E}}^{l}.
\end{equation*}%
Then $\mathcal{A}\left( \tilde{\mathcal{G}}^{l}\right) $ is a framed cover
that satisfies all of the conclusions of Lemma \ref{jacked up cov lem 1},
except for Part 4 and (\ref{covert struct}) and (\ref{coedge struct}).
\end{enumerate}
\end{proposition}

\begin{proof}
Let $\mathcal{A}^{2}=\{\digamma _{i}(\beta _{i})\}_{i=1}^{\tilde{L}},$ and
let $\mathcal{G}^{L}$ be the $\eta _{L}$-geometric cograph from Lemma \ref%
{jacked up cov lem 1}. By Part 3 of Lemma \ref{jacked up cov lem 1}, $\frac{%
\beta _{1}}{50}\leq \eta _{L}$. Let $k\geq 1$ be the smallest integer so
that $\frac{\eta _{L}}{10^{k}}\in \left[ \frac{\beta _{1}}{500},\frac{\beta
_{1}}{50}\right] $. Legally subdivide $\mathcal{G}^{L}$ $k$ times to get a
legal, error correcting, $\tilde{\eta}_{1}$--geometric cograph $\tilde{%
\mathcal{G}}^{1}$ with $\tilde{\eta}_{1}:=\frac{\eta _{L}}{10^{k}}$. Let $\{%
\tilde{\eta}_{i}\}_{i=1}^{\tilde{L}}$ be the corresponding sequence from
Proposition \ref{beta eta sequence}. For each $l\in \left\{ 1,2,\ldots ,%
\tilde{L}\right\} ,$ legally subdivide $\tilde{\mathcal{G}}^{1}$ to get $%
\tilde{\mathcal{G}}^{l}.$ This gives us Part $1.$ Part $2$ follows from the
construction of $\tilde{\mathcal{G}}^{l}$ and Parts 6 and 7 of Lemma \ref%
{jacked up cov lem 1}. Part 3 holds because the submersions that define the
cosimplices of $\tilde{\mathcal{G}}^{l}$ are lined up, almost Riemannian
submersions and because a legal subdivision of a legal cograph is legal.

Parts 4 and 5 follow from Parts 6 and 7 of Lemma \ref{jacked up cov lem 1}
after adjusting the covertices and coedges so that (\ref{covertex vert cut
off}) and (\ref{coedge vert cut off}) hold. To do this, we exploit the fact
that a cosimplex is parameterized by its fiber exhaustion function. So each
cosimplex can be shrunken by replacing it by a super level set of its fiber
exhaustion function.

To prove Part $6,$ we note that all of the properties of Lemma \ref{jacked
up cov lem 1} are preserved under subdivision except of course the three
mentioned.
\end{proof}

\begin{remarknonum}
$\mathcal{A}^{2}$ is not likely to thoroughly respect $\mathcal{A}^{0}\left( 
\tilde{\mathcal{G}}^{l}\right) \cup \mathcal{A}^{1}\left( \tilde{\mathcal{G}}%
^{l}\right) $ as $\mathcal{A}^{2}$ is fixed and we have subdivided our
original cograph $\mathcal{G}^{L}$ to obtain $\tilde{\mathcal{G}}^{l}$ (see
Figure \ref{cov lem 1+ fig 1}).

Similarly, subdividing the cograph changes the width of the cosimplices, so (%
\ref{covert struct}) and (\ref{coedge struct}) are replaced with the
statements of Parts 4 and 5 of the result above.
\end{remarknonum}

\addtocounter{algorithm}{1}

\subsection{Cotriangulations Relative to the Natural $0$ and $1$ Strata}

The main results of this subsection are Proposition \ref{dda} and Lemma \ref%
{cotr exist lemma}. The first result specifies a type of cotriangulation
whose existence implies Covering Lemma \ref{cov em 2}. Lemma \ref{cotr exist
lemma} asserts that this type of cotriangulation exists. The remainder of
Part \ref{Part 2}, after this subsection, is a proof of Lemma \ref{cotr
exist lemma}. Before stating Proposition \ref{dda}, we specify some more
terminology.

Let $\mathcal{B}_{l}$ be as in \ref{B and L dfn}. For $\digamma \in \mathcal{%
B}_{l}$ set 
\begin{equation*}
\digamma _{\mathrm{bot}}:=p_{\digamma }^{-1}\left( \left[ 0,\beta _{\digamma
}\right] \times \left[ 0,\frac{\tilde{\eta}_{l}}{50}\right] \right) ,
\end{equation*}%
and 
\begin{equation*}
\mathcal{A}_{\mathrm{bot}}:=\left\{ \left. \left( \digamma _{\mathrm{bot}%
},\pi _{1}\circ p_{\digamma }\right) \text{ }\right\vert \text{ }\digamma
\in \mathcal{B}_{l}\right\} .
\end{equation*}%
Then $\mathcal{A}_{\mathrm{bot}}$ is a collection of $1$--framed sets that
is $\kappa $--lined up, vertically separated, and fiber 
\begin{wrapfigure}[12]{l}{0.37\textwidth}\centering
\vspace{-8pt}
\includegraphics[scale=0.145]{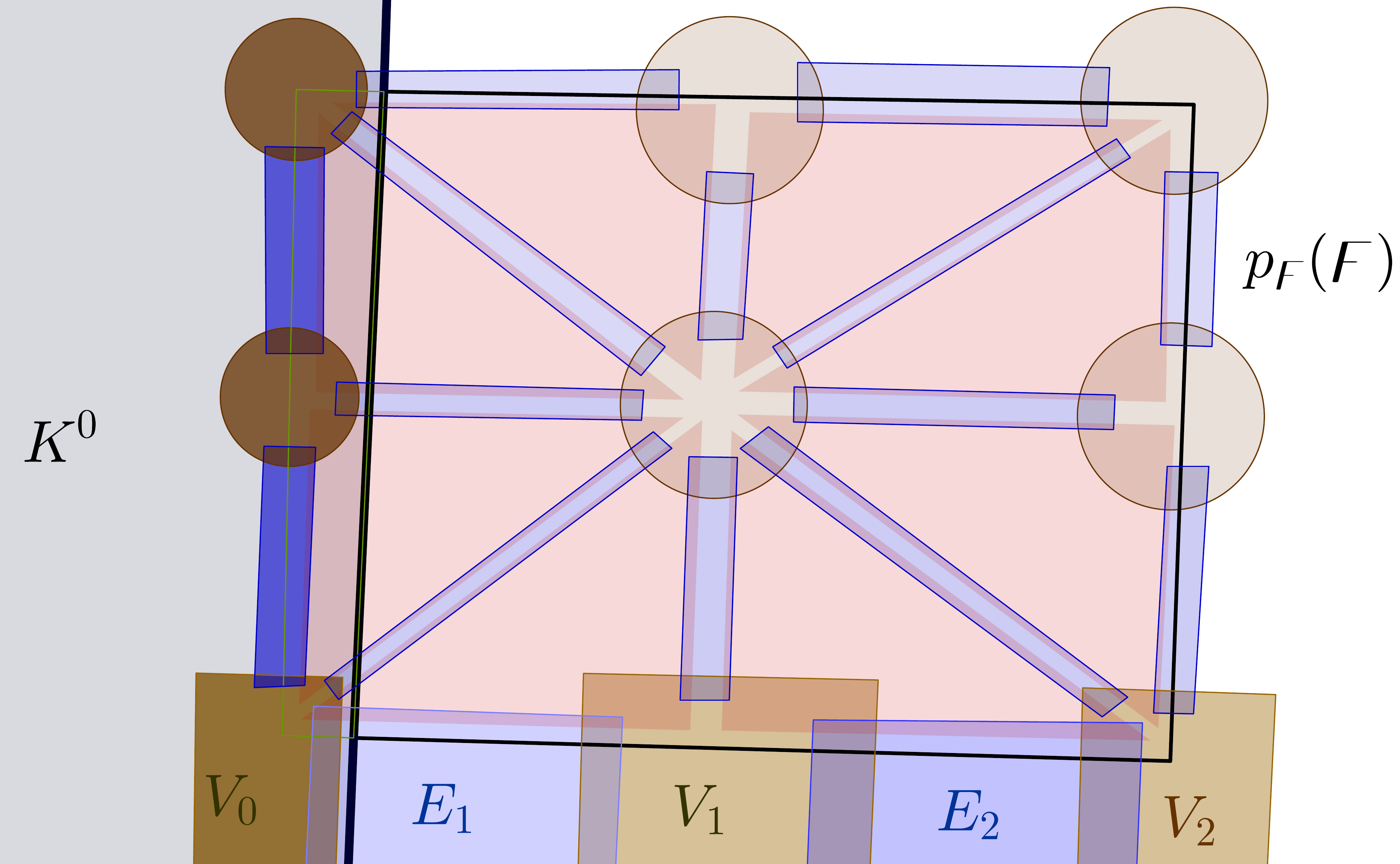}
\caption{ {\protect\tiny {A triangulation cover that generates a cotriangulation relative to   $\mathcal{LG}^{l}\cup \mathcal{B}\tilde{\mathcal{G}}^{l}$. The projection of the cograph ${\mathcal{L}}{%
\mathcal{G}^{l}}$ lies along the left hand edge. The projection of the
cograph ${\mathcal{B}}{\mathcal{\tilde{G}}^{l}}$ lies along the bottom edge. 
} }}
\label{cotri nat str interface figure}%
\end{wrapfigure}
swallowing. For $l\in \left\{ 1,2,\ldots ,\tilde{L}\right\} ,$ let $\tilde{%
\mathcal{G}}^{l}$ be as in Proposition \ref{bottom prop}. Then 
\begin{equation}
\hspace*{2 in} \mathcal{B}\tilde{\mathcal{G}}^{l}:=\left\{ \left. S\cap
(\cup \mathcal{A}_{\mathrm{bot}})\right\vert S\in \tilde{\mathcal{G}}%
^{l}\right\}  \label{bottom
dfn}
\end{equation}%
is a legal, error correcting cograph that is subordinate to $\mathcal{A}_{%
\mathrm{bot}}.$

\begin{definition}
\label{rel to dfn}For $k\leq l,$ let $\mathcal{S}=\mathcal{V}\cup \mathcal{E}%
\cup \mathcal{F}$ be a cotriangulation subordinate to $\{\digamma
_{i}\}_{i=1}^{k}\subset \mathcal{A}^{2}$ which contains the intersection of $%
\mathcal{LG}^{l}\cup \mathcal{B}\tilde{\mathcal{G}}^{l}$ with $\cup
_{i=1}^{k}\digamma _{i}$. Let $\mathcal{V}\left( \mathcal{LG}^{l}\cup 
\mathcal{B}\tilde{\mathcal{G}}^{l}\right) $ and $\mathcal{E}\left( \mathcal{%
LG}^{l}\cup \mathcal{B}\tilde{\mathcal{G}}^{l}\right) $ be the collections
of covertices and coedges of $\mathcal{LG}^{l}\cup \mathcal{B}\tilde{%
\mathcal{G}}^{l},$ respectively. We say that $\mathcal{S}$ is \textbf{%
relative to }$\mathcal{LG}^{l}\cup \mathcal{B}\tilde{\mathcal{G}}^{l}$
provided

\begin{enumerate}
\item The collection of closures of 
\begin{equation*}
\mathcal{V}^{\mathrm{in}}:=\mathcal{V}\setminus \mathcal{V}\left( \mathcal{LG%
}^{l}\cup \mathcal{B}\tilde{\mathcal{G}}^{l}\right)
\end{equation*}%
is disjoint from the collection of closures of $\mathcal{A}^{0}(\tilde{%
\mathcal{G}}^{l}),$ and the collection of closures of 
\begin{equation*}
\mathcal{E}^{\mathrm{in}}:=\mathcal{E}\setminus \mathcal{E}\left( \mathcal{LG%
}^{l}\cup \mathcal{B}\tilde{\mathcal{G}}^{l}\right)
\end{equation*}%
is disjoint from the collection of closures of $\mathcal{A}^{1}(\tilde{%
\mathcal{G}}^{l}).$

\item For each $E$ in $\mathcal{E}^{\mathrm{in}},$ if $E$ intersects an
element of $\mathcal{A}^{0}(\tilde{\mathcal{G}}^{l})\cup \mathcal{A}^{1}(%
\tilde{\mathcal{G}}^{l}),$ then this intersection is contained in $\cup 
\mathcal{V}\left( \mathcal{LG}^{l}\cup \mathcal{B}\tilde{\mathcal{G}}%
^{l}\right) .$

\item For each $F\in \mathcal{F}$, if $F$ intersects $\cup \left( \mathcal{A}%
^{0}(\tilde{\mathcal{G}}^{l})\cup \mathcal{A}^{1}(\tilde{\mathcal{G}}%
^{l})\right) ,$ then this intersection is contained in $\cup \left( \mathcal{%
LG}^{l}\cup \mathcal{B}\tilde{\mathcal{G}}^{l}\right) .$

\item For $F\in \mathcal{F}$ as in Part 3, suppose that $F$ is subordinate
to $\digamma _{F}\in \mathcal{A}^{2}.$ Then along $F\cap \left( \cup \left( 
\mathcal{A}^{0}(\tilde{\mathcal{G}}^{l})\cup \mathcal{A}^{1}(\tilde{\mathcal{%
G}}^{l})\right) \right) ,$ $p_{F}=q_{F}=p_{\digamma _{F}}.$
\end{enumerate}
\end{definition}

\begin{wrapfigure}[17]{l}{0.29\textwidth}\centering
\vspace{0pt}
\includegraphics[scale=0.16]{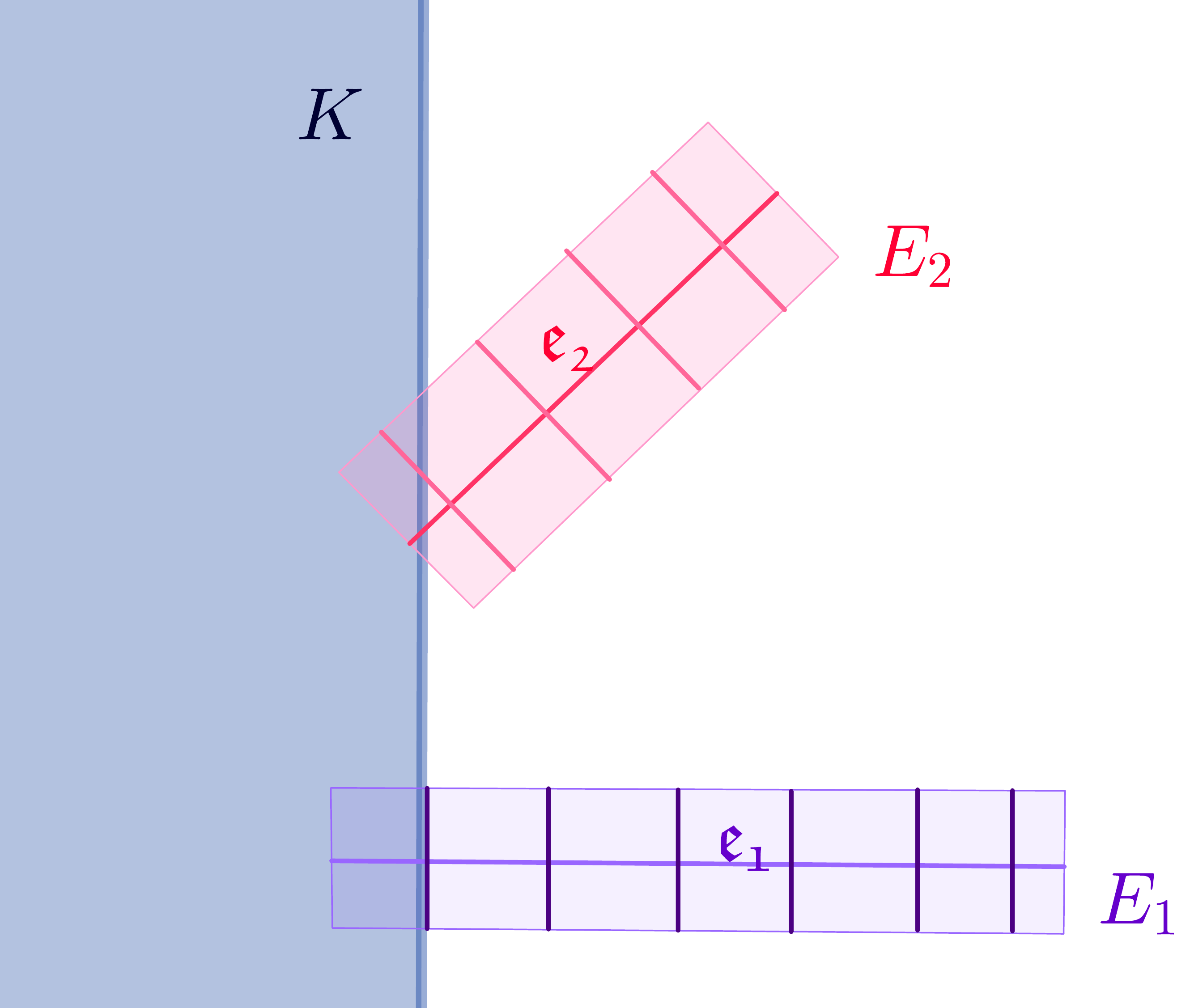}
\caption{ {\protect\tiny {\ Pictured are the images of the coedges $E_{1}$, $%
E_{2}$ and the naturally $0$-framed set $K$ under the submersion $%
p_{\digamma }$ of a naturally $2$-framed set $\digamma $ (not pictured). The
generating edge $\mathfrak{e}_{1}$ of $E_{1}$ is parallel to the $x$--axis,
and hence $E_{1}$ respects $K$. The generating edge $\mathfrak{e}_{2} $ of $%
E_{2}$ is not parallel to the $x$--axis, and hence $E_{2}$ does not respect $%
K$.} }}
\label{lin functs fig}
\end{wrapfigure}
For $k\leq l,$ let $\mathcal{S}=\mathcal{V}\cup \mathcal{E}\cup \mathcal{F}$
be a cotriangulation subordinate to $\{\digamma _{i}\}_{i=1}^{k}\subset 
\mathcal{A}^{2}$ and relative to $\mathcal{LG}^{l}\cup \mathcal{B}\tilde{%
\mathcal{G}}^{l}.$ Let $\mathcal{E}^{\mathrm{left}}$ be the coedges of $%
\mathcal{E}^{\mathrm{in}}$ that intersect $\mathcal{K}^{0}\mathcal{\ }$and
let $\mathcal{E}^{\mathrm{bot}}$ $\ $be the coedges of $\mathcal{E}^{\mathrm{%
in}}$ that intersect $\mathcal{V}( \tilde{\mathcal{G}}^{\tilde{L}}) .$
Suppose that $E^{l}\in \mathcal{E}^{\mathrm{left}}$ and $E^{b}\in \mathcal{E}%
^{\mathrm{bot}}$ are generated by the $1$--simplices $\mathfrak{e}^{l}$ and $%
\mathfrak{e}^{b},$ respectively. It follows from Parts 5 and 6 of Lemma \ref%
{jacked up cov lem 1} that $\mathfrak{e}^{l}$ intersects $\left\{ 0\right\}
\times \mathbb{R},$ and $\mathfrak{e}^{b}$ intersects $\mathbb{R}\times
\left\{ 0\right\} .$ For the next result, we assume that for all $E^{l}\in 
\mathcal{E}^{\mathrm{left}}$ and all $E^{b}\in \mathcal{E}^{\mathrm{bot}},$ 
\begin{equation*}
\hspace{2in} \sphericalangle \left( \mathfrak{e}^{l},\left\{ 0\right\}
\times \mathbb{R}\right) >\left( \frac{\pi }{6}\right) ^{-}\text{ and }
\end{equation*}

\begin{equation}
\hspace{2in}\sphericalangle \left( \mathfrak{e}^{b},\mathbb{R}\times \left\{
0\right\} \right) >\left( \frac{\pi }{6}\right)^{-}.  \label{transverese}
\end{equation}

\begin{wrapfigure}[11]{r}{0.33\textwidth}\centering
\vspace{-20pt}
\includegraphics[scale=0.13]{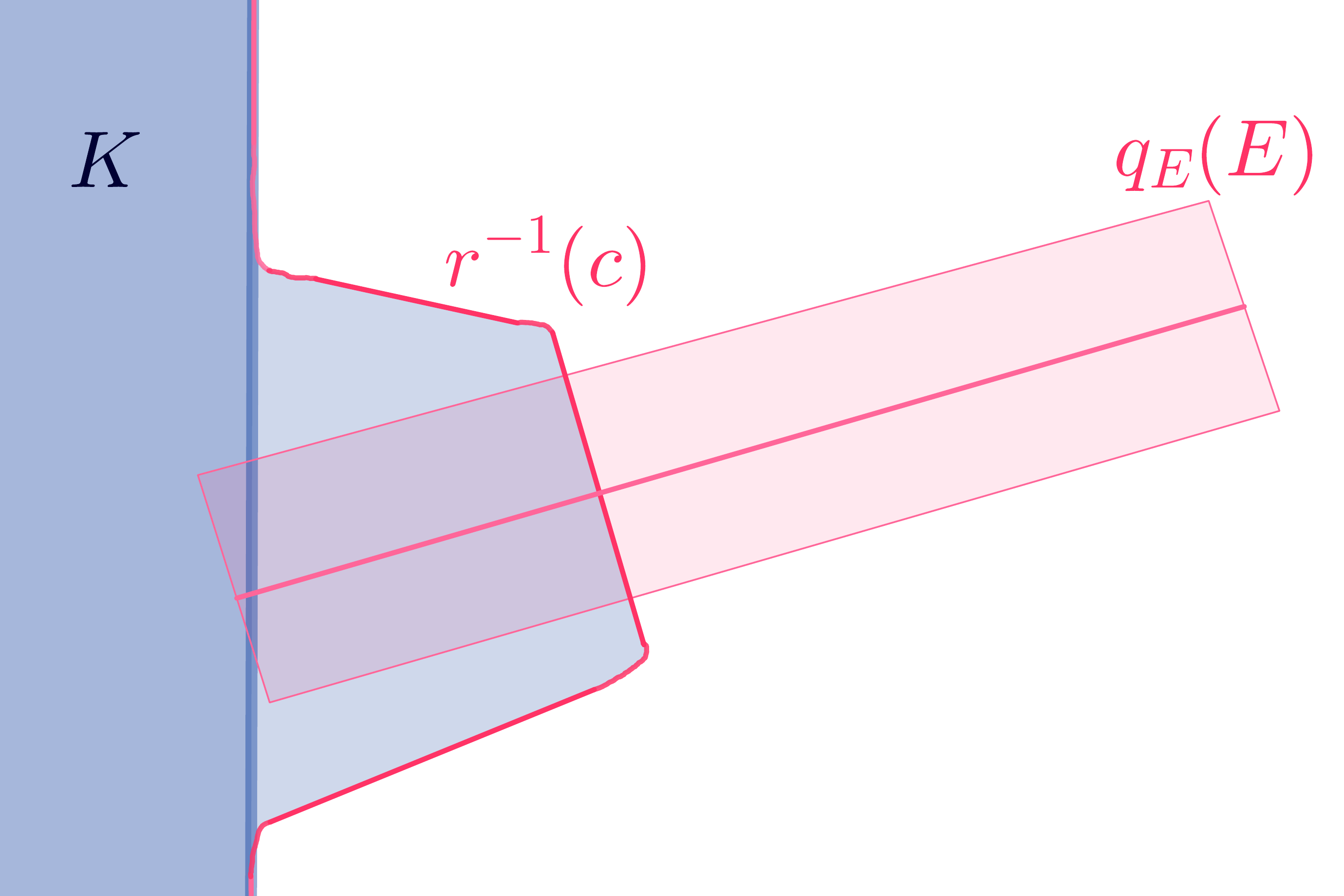}
\caption{ {\protect\tiny {\ This coedge does not intersect the naturally $0$%
--framed set $K $ respectfully, but after modifying the fiber exhaustion
function of $K$ using Proposition \protect\ref{deform funct}, the
intersection becomes respectful. To be more precise, we also replace $K$
with $r^{-1}(-\infty ,c)$. } }}
\label{pimple}
\end{wrapfigure}

\begin{proposition}
\label{dda} For $l=\tilde{L}$, let $\mathcal{A}(\tilde{\mathcal{G}}^{\tilde{L%
}})$ be the framed cover from Part 6 of Proposition \ref{bottom prop}. Let $%
\mathcal{S=V}\cup \mathcal{E\cup F}$ be a respectful cotriangulation that is
subordinate to $\mathcal{A}^{2}$, relative to $\mathcal{LG}^{\tilde{L}}\cup 
\mathcal{B}\tilde{\mathcal{G}}^{\tilde{L}},$ and that satisfies (\ref%
{transverese}).

Then there is a collection $\mathcal{D}^{0}$ of $0$-framed sets which are
deformations of the elements of $\mathcal{A}^{0}\left( \tilde{\mathcal{G}}^{%
\tilde{L}}\right) $ so that Covering Lemma \ref{cov em 2} holds with 
\begin{equation*}
\mathcal{B}^{0}=\mathcal{V}^{\mathrm{in}}\mathcal{\cup D}^{0},\text{ }%
\mathcal{B}^{1}=\mathcal{E}^{\mathrm{in}}\mathcal{\cup A}^{1}\left( \tilde{%
\mathcal{G}}^{\tilde{L}}\right) ,\text{ and }\mathcal{B}^{2}=\mathcal{F}%
\text{ (cf Figure \ref{pimple}).}
\end{equation*}
\end{proposition}

\noindent \textit{Proof. }It follows from Part 1 of Definition \ref{rel to
dfn} that we can replace $\mathcal{A}^{0}(\tilde{\mathcal{G}}^{l})\cup 
\mathcal{A}^{1}(\tilde{\mathcal{G}}^{l})$ with an appropriate choice of $%
\mathcal{A}^{0}(\tilde{\mathcal{G}}^{l})^{+}\cup \mathcal{A}^{1}(\tilde{%
\mathcal{G}}^{l})^{+}$ and 
\begin{equation}
\mathcal{S}\text{ will still be relative to }\mathcal{LG}^{l}\cup \mathcal{B}%
\tilde{\mathcal{G}}^{l}.  \label{topologists geometry}
\end{equation}%
We must construct $\mathcal{D}^{0}$ and show that the resulting collection $%
\mathcal{B}:=\mathcal{B}^{0}\mathcal{\cup B}^{1}\mathcal{\cup B}^{2}$ is
respectful, covers $X\setminus X_{\delta }^{4},$ and has the required
disjointness properties.

\noindent \textbf{Respect: }Since $\mathcal{S\ }$is respectful, we only need
concern ourselves with how $\mathcal{E}^{\mathrm{in}}\mathcal{\cup F}$
intersect the natural $0$- and $1$-strata. It follows from Parts $6$ and $7$
of Lemma \ref{jacked up cov lem 1}, Parts 3 and 4 of Definition \ref{rel to
dfn}, and the construction of $\mathcal{B}\tilde{\mathcal{G}}^{l}$ that $%
\mathcal{F}$ respects $\mathcal{A}^{1}(\tilde{\mathcal{G}}^{\tilde{L}}).$
The main project is to construct a deformation $\mathcal{D}^{0}$ of $%
\mathcal{A}^{0}( \tilde{\mathcal{G}}^{\tilde{L}}) $ so that $\mathcal{E}^{%
\mathrm{in}}$ respects $\mathcal{D}^{0}.$ To accomplish this, we deform $%
\mathcal{K}^{0}$ and $\mathcal{V}^{0}( \tilde{\mathcal{G}}^{\tilde{L}}) $
separately. Except for notation, the two arguments are the same, so we only
give the details for $\mathcal{K}^{0}.$

Suppose that $K\in \mathcal{K}^{0}$ and $E\in \mathcal{E}^{\mathrm{left}%
}\subset \mathcal{E}^{\mathrm{in}}$ intersect. Let $r_{K}:K\rightarrow 
\mathbb{R}$ be the fiber exhaustion function for $K,$ and let $%
q_{E}:E\rightarrow \left[ 0,\beta \right] \times \left[ 0,\beta \right]
\subset \mathbb{R}^{2}$ and $p_{E}=\pi _{\mathfrak{e}}\circ
q_{E}:E\rightarrow \mathfrak{e}\subset \mathbb{R}^{2}$ be the submersions
that define $E$ as a coedge.

Let $\tilde{p}_{E}=\pi _{1}\circ p_{E}:E^{l}\longrightarrow \mathbb{R},$ and
suppose that $E$ is subordinate to $\digamma \in \mathcal{A}^{2}.$ By Part 5
of Lemma \ref{jacked up cov lem 1}, $r_{K}$ is lined up with $\pi _{1}\circ
p_{\digamma },$ so for simplicity we assume 
\begin{equation*}
r_{K}=\pi _{1}\circ p_{\digamma }.
\end{equation*}%
This, together with (\ref{transverese}), gives us that the angles between
the gradients of $\tilde{p}_{E}$ and $r_{K}=\pi _{1}\circ p_{\digamma }$ are 
$\leq \left( \frac{5\pi }{6}\right) ^{+}.$ Thus, after choosing the
appropriate bump function $\tilde{\lambda}:$ $\mathbb{R}\longrightarrow 
\mathbb{R}$ and setting $\lambda :=\tilde{\lambda}\circ \pi _{1}\circ
p_{\digamma },$ we see that $r_{K},$ $\tilde{p}_{E},$ and $\lambda $ satisfy
the hypotheses of Proposition \ref{deform funct}. Thus by Proposition \ref%
{deform funct}, there is a deformation $\tilde{r}_{K}$ of $r_{K}$ with $%
\tilde{r}_{K}=r_{K}$ on ${K\setminus E^{+}}$ and $\tilde{r}_{K}=\tilde{p}%
_{E} $ on $E.$ By redefining $K$ to be the appropriate sublevel of $\tilde{r}%
_{K}, $ we have that $E$ respects $K$ (cf Figure \ref{pimple}).

It remains to show that $\mathcal{F}$ respects $\mathcal{D}^{0}.$ As before,
we only give the details for our deformation of $\mathcal{K}^{0}.$ Except
for notation, the argument for the deformation of $\mathcal{V}^{0}\left( 
\tilde{\mathcal{G}}^{\tilde{L}}\right) $ is the same. From Parts $3$ and $4$
of Definition \ref{rel to dfn}, we have that if $F\in \mathcal{F}$
intersects $K\in \mathcal{K}^{0},$ then along this intersection, 
\begin{equation*}
p_{F}=p_{\digamma },
\end{equation*}%
where $\digamma $ is the element of $\mathcal{A}^{2}$ to which $F$ is
subordinate. Since $r_{K}=\pi _{1}\circ p_{\digamma },$ it follows that the
fibers of $p_{F}$ are contained in levels of $r_{K},$ so $F$ respects $K,$
and it remains to see that this property is not destroyed when we replace $K$
with $\tilde{K}.$

Suppose that $E$ is a coedge that intersects $K$ and $F.$ Since $\mathcal{S\ 
}$is respectful the fibers of $p_{F}$ are contained in levels of $p_{E}$ and
hence also in the levels of $\tilde{p}_{E}=\pi _{1}\circ p_{E}.$ Since our
bump function $\lambda $ is $\tilde{\lambda}\circ \pi _{1}\circ p_{\digamma
},$ we also have that the fibers of $p_{F}$ are contained in levels of $%
\lambda .$ Thus the fibers of $p_{F}$ are contained in levels of 
\begin{equation*}
\tilde{r}_{K}=\left( 1-\lambda \right) \tilde{p}_{E}+\lambda r_{K},
\end{equation*}%
and $\mathcal{B}$ is respectful.

To assist in the remainder of the proof, we note that the deformation
outlined above can be performed while not affecting the existing properties
of the $\mathcal{A}^{0}\left( \tilde{\mathcal{G}}^{\tilde{L}}\right) .$
Indeed, by (\ref{topologists geometry}) there is a choice of $K^{+}$ with $%
K\Subset K^{+}$ so that $\mathcal{S}$ is also relative to $\mathcal{LG}%
^{l}\cup \mathcal{B}\tilde{\mathcal{G}}^{l}$ even with $K$ replaced by $%
K^{+}.$ Proposition \ref{deform funct} and the argument above allow us to
construct $\tilde{K}$ so that 
\begin{equation}
K\subset \tilde{K}\subset K^{+}\text{.\label{samll deformation}}
\end{equation}

\noindent \textbf{Covers: }Since $\mathcal{S}$ covers $\left( \mathcal{A}%
^{2}\right) ^{-}$ and $\mathcal{A}^{0}\left( \tilde{\mathcal{G}}^{\tilde{L}%
}\right) \cup \mathcal{A}^{1}\left( \tilde{\mathcal{G}}^{\tilde{L}}\right) $
covers $\Omega _{0}\cup \Omega _{1},$ 
\begin{equation*}
\mathcal{A}^{0}\left( \tilde{\mathcal{G}}^{\tilde{L}}\right) \cup \mathcal{A}%
^{1}\left( \tilde{\mathcal{G}}^{\tilde{L}}\right) \cup \mathcal{S}
\end{equation*}%
covers $\Omega _{0}\cup \Omega _{1}\cup \left( \mathcal{A}^{2}\right) ^{-}.$
Since $\cup \left( \mathcal{V}\setminus \mathcal{V}^{\mathrm{in}}\right) $
and $\cup \left( \mathcal{E}\setminus \mathcal{E}^{\mathrm{in}}\right) $ are
contained in $\Omega _{0}\cup \Omega _{1},$ 
\begin{equation*}
\mathcal{A}^{0}\left( \tilde{\mathcal{G}}^{\tilde{L}}\right) \cup \mathcal{A}%
^{1}\left( \tilde{\mathcal{G}}^{\tilde{L}}\right) \cup \mathcal{V}^{\mathrm{%
in}}\cup \mathcal{E}^{\mathrm{in}}\cup \mathcal{F}
\end{equation*}%
covers $\Omega _{0}\cup \Omega _{1}\cup \left( \mathcal{A}^{2}\right) ^{-},$
and it follows from (\ref{samll deformation}) that $\mathcal{B}$ covers $%
X\setminus X_{\delta }^{4}.$

\noindent \textbf{Disjointness: } Since $\mathcal{S}$ is a cotriangulation, $%
\mathcal{V}^{\mathrm{in}},$ $\mathcal{E}^{\mathrm{in}}$, and $\mathcal{B}%
^{2}=\mathcal{F}$ are disjoint collections. By Part 6 of Proposition \ref%
{bottom prop}, $\mathcal{A}^{0}(\mathcal{G}^{\tilde{L}})$ and $\mathcal{A}%
^{1}(\mathcal{G}^{\tilde{L}})$ are each disjoint collections. By Part 1 of
Definition \ref{rel to dfn}, $\mathcal{E}^{\mathrm{in}}$ is disjoint from $%
\mathcal{A}^{1}(\mathcal{G}^{\tilde{L}})$ and $\mathcal{V}^{\mathrm{in}}$ is
disjoint from $\mathcal{A}^{0}(\mathcal{G}^{\tilde{L}}).$ So $\mathcal{B}%
^{1} $ is a disjoint collection, and using also (\ref{samll deformation}) we
have that $\mathcal{B}^{0}$ is a disjoint collection. \mbox{}\hfil $\square$ 
\vspace{5pt}

So our goal for the remainder of this part of the paper is to prove

\begin{lemma}
\label{cotr exist lemma}There is a respectful cotriangulation $\mathcal{S}=%
\mathcal{V}\cup \mathcal{E\cup F}$ that is subordinate to $\mathcal{A}^{2},$
is relative to $\mathcal{LG}^{\tilde{L}}\cup \mathcal{B}\tilde{\mathcal{G}}^{%
\tilde{L}},$ and satisfies (\ref{transverese}).
\end{lemma}



\section{Geometric Cotriangulations\label{tools to constr cotri sect}}

In spite of the fact that the concept of a cotriangulation is rather
topological, we do not know how to prove Lemma \ref{cotr exist lemma}
without imposing some further geometric restrictions on our
cotriangulations. To do this, we will exploit a special fact about $\mathbb{R%
}^{2}$, namely Chew's Theorem (\ref{Gromov-Del Prop}, below), which says
that all angles of a special type of simplicial complex in $\mathbb{R}^{2}$
are $\geq \frac{\pi }{6}.$ We call these special complexes CDGs.

The only significance of the number $\frac{\pi }{6}$ is that it is a
universal lower bound for all angles of all CDGs. In the end, this will
allow us to verify that our cotriangulations satisfy the geometric
inequality in (\ref{Sub def hyp}), and thus can be deformed to respectful
cotriangulations.

In Subsection \ref{Delaunay sect}, we review the relevant results from \cite%
{ProWilh2} on Chew's Theorem and its consequences. In Subsection \ref%
{geomcotriansubsect}, we define the concept of a geometric triangulation
cover and in Subsection \ref{Del Grm subsect} we define geometric
cotriangulations. Roughly speaking, geometric cotriangulations are
cotriangulations whose underlying triangulations are CDGs.

\addtocounter{algorithm}{1}

\subsection{Chew's Theorem and Its Consequences\label{Delaunay sect}}

A CDG is a special type of triangulation in $\mathbb{R}^{2}.$ Roughly
speaking it is a Delaunay triangulation of a Gromov set. In this subsection,
we define CDGs and review their properties that were proven in \cite%
{ProWilh2}. We start with the definition of Gromov sets.

\begin{definition}
\label{dfn eta grmv}Given $\eta >0,$ we say a subset $A$ of a metric space $%
X,$ is an $\eta $--\textbf{Gromov set} if and only if 
\begin{equation*}
B\left( A,\eta \right) =X
\end{equation*}%
and for distinct $a,b\in A,$ 
\begin{equation*}
\mathrm{dist}\left( a,b\right) \geq \eta .
\end{equation*}%
Equivalently, $A$ is $\eta $--separated and not properly contained in an $%
\eta $--separated subset of $X.$
\end{definition}

The definitions of Delaunay triangulation and CDGs are:

\begin{definition}
\label{dfn dela}(see e.g. \cite{GartHoff}, Chapter 6 and Proposition 1.5 in 
\cite{ProWilh2}) A triangulation $\mathcal{T}$ of a discrete point set $S$
of $\mathbb{R}^{2}$ is called Delaunay if and only if for every edge $e$
that bounds two faces, the sum of the opposite angles is $\leq \pi .$
Equivalently, the circumdisk of each $2$--simplex contains no points of $S$
in its interior.
\end{definition}

\begin{definition}
\label{CDG dfn}We call a simplicial complex $\mathcal{T}$ in $\mathbb{R}^{2}$
an $\eta $\textbf{--}CDG, provided there is a Delaunay triangulation $%
\mathcal{\hat{T}}$ of an $\eta $\textbf{--}Gromov subset of $\mathbb{R}^{2},$
and $\mathcal{T}$ is a subcomplex of $\mathcal{\hat{T}}.$
\end{definition}

The algorithm on Pages 7-8 of \cite{Chew} combined with the theorem on Page
9 and the corollary on Page 10 of \cite{Chew} give the following.

\begin{theorem}
\label{Gromov-Del Prop}(Chew's Angle Theorem) Let $\mathcal{T}$ be an $\eta $%
--CDG. Then all angles of $\mathcal{T}$ are $\geq \frac{\pi }{6}$ and the
lengths of all edges of $\mathcal{T}$ are in $\left[ \eta ,2\eta \right] .$
\end{theorem}

\noindent (See \cite{ProWilh2} for an alternative exposition.) Next we
discuss deforming the concept of a $\mathrm{CDG}$ to what we call almost $%
\mathrm{CDGs}.$ This begins with

\begin{definition}
Given $\eta ,\delta >0,$ we say a subset $A$ of a metric space $X,$ is an $%
\left( \eta ,\delta \eta \right) $--\textbf{Gromov set} if and only if 
\begin{equation*}
B\left( A,\eta \left( 1+\delta \right) \right) =X
\end{equation*}%
and for distinct $a,b\in A,$ 
\begin{equation*}
\mathrm{dist}\left( a,b\right) \geq \eta \left( 1-\delta \right) .
\end{equation*}
\end{definition}

\begin{figure}[tbp]
\includegraphics[scale=.13]{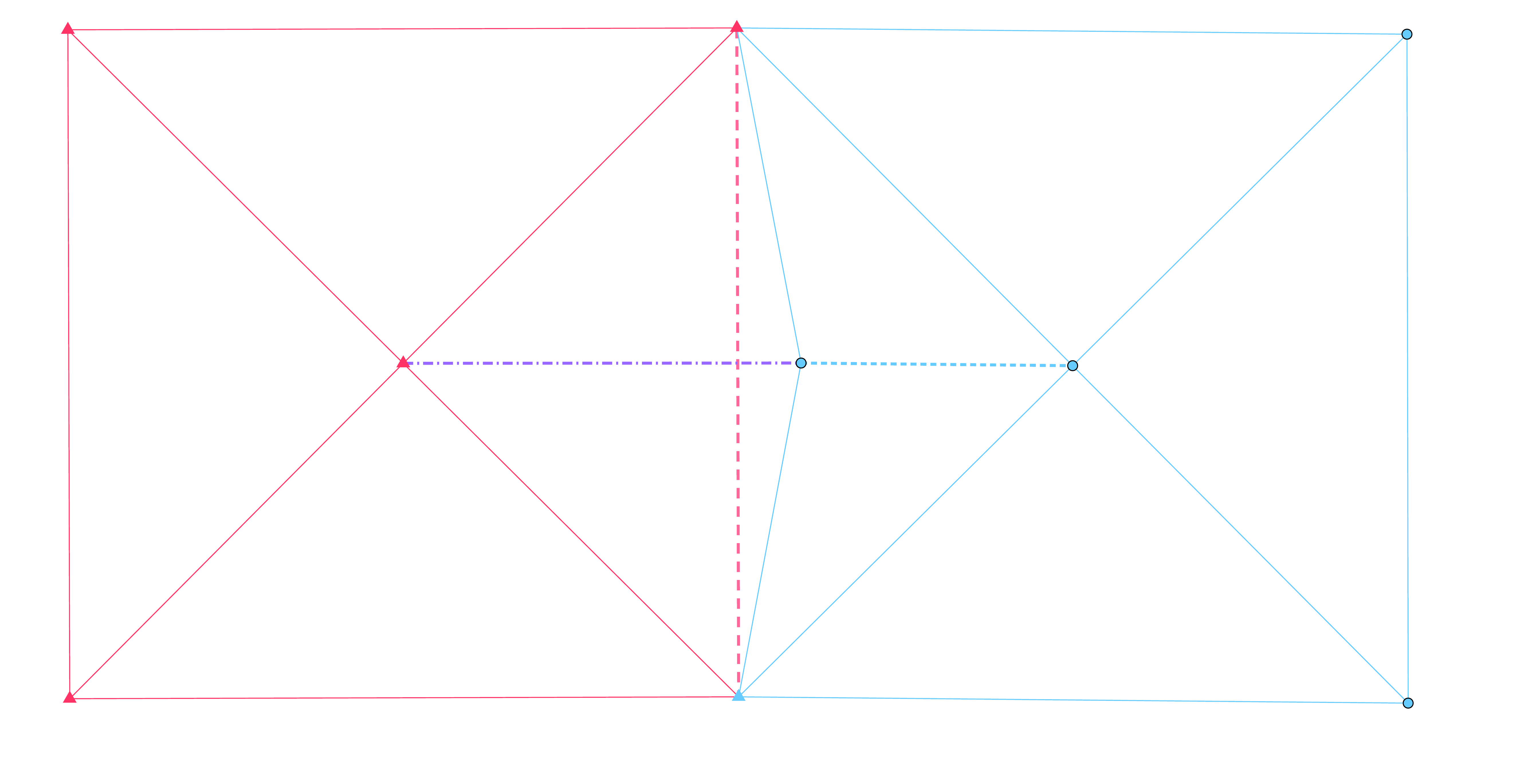}
\caption{{\protect\tiny {Assume the side lengths of these squares is just
under $2$. Then the red vertices are a $1$--Gromov subset of the red square,
and the red edges are a Delaunay triangulation of the red vertices. The blue
vertices are an extension of the red ones to a $1$--Gromov subset of the
union of the two squares, but the dashed red edge fails the angle flip test
and hence is not Delaunay for the extended configuration. The Delaunay
triangulation of the extended set consists of all of the edges pictured
except for the dashed red one.}}}
\label{sqrt3}
\end{figure}

To avoid technical difficulties near the boundaries of Delaunay
triangulations of Gromov sets the definition of almost \textrm{CDGs }is
phrased in terms of the following technical concept (cf Figure \ref{sqrt3}).

\begin{definitionnonum}
\label{ext grm dfn}Let $X$ be a metric space with $U,U_{\mathrm{x}}\subset X$
open in $X.$ Given $\eta >0$ and $\delta \geq 0,$ let $A$ and $A_{\mathrm{x}%
} $ be $\left( \eta ,\delta \eta \right) $--Gromov subsets of the closures
of $U$ and $U_{\mathrm{x}}$ respectively. We say that $A_{\mathrm{x}}$ is a 
\textbf{buffer }of $A$ provided:

\begin{enumerate}
\item The closed ball $D(U,6\eta )\subset U_{\mathrm{x}}$.

\item $A=A_{\mathrm{x}}\cap \mathrm{closure}\left( U\right) .$
\end{enumerate}
\end{definitionnonum}

A $\mathrm{CDG}$ passes the following test with $\varepsilon =0.$

\begin{definition}
Let $\mathcal{T}$ be a simplicial complex in $\mathbb{R}^{2}.$ Let $e$ be an
edge of $\mathcal{T}$ which bounds two faces. We say that $e$ passes the $%
\varepsilon $--angle Flip Test if and only if the sum of the angles opposite 
$e$ is $\leq \pi +\varepsilon .$
\end{definition}

We are ready for the definition of almost $\mathrm{CDGs.}$

\begin{definition}
\label{alm CDG dfn}Given $\eta ,\delta ,\varepsilon >0,$ let $\mathcal{T}%
_{0} $ be an $\left( \eta ,\delta \eta \right) $--Gromov subset of the
closure of an open set $U\subset \mathbb{R}^{2}$ with buffer $\left( 
\mathcal{T}_{0}\right) _{\mathrm{x}}.$ A triangulation $\mathcal{T}$ of $%
\mathcal{T}_{0} $ is called an $\left( \eta ,\delta \eta ,\varepsilon
\right) $--CDG, provided $\mathcal{T}$ is a subcomplex of a triangulation $%
\mathcal{T}_{\mathrm{x}}$ of $\left( \mathcal{T}_{0}\right) _{\mathrm{x}}$,
and each edge of $\mathcal{T}_{\mathrm{x}}$ passes the $\varepsilon $--angle
flip test.
\end{definition}

We will say $\mathcal{T}_{\mathrm{x}}$ is a \textbf{buffer} of $\mathcal{T}.$
In \cite{ProWilh2} we observe

\begin{proposition}
(Proposition 4.15 of \cite{ProWilh2}) \label{alm CDG Chew}Given $d,\theta >0$%
, there are $\varepsilon ,\delta >0$ so that for every $\left( \eta ,\delta
\eta ,\varepsilon \right) $-CDG,

\begin{enumerate}
\item Every edge length of $\mathcal{T}$ is in the interval $\left[ \eta
-d,2\eta +d\right] .$

\item All angles of $\mathcal{T}$ are $\geq \frac{\pi }{6}-\theta .$
\end{enumerate}
\end{proposition}

A given $S\subset \mathbb{R}^{2}$ can have more than one Delaunay
triangulation, and hence sets $S_{i}$ arbitrarily close to $S$ need not have
Delaunay triangulations that are close to an initially prescribed
triangulation of $S.$ The following three definitions are part of our
strategy to deal with these issues.

\begin{definition}
\label{del grph dfn}\emph{(\cite{GartHoff}, page 74) }Given a discrete $%
S\subset \mathbb{R}^{2},$ a segment $e$ between two points of $S$ is in the
Delaunay graph of $S$ if and only if $e$ is an edge of every Delaunay
triangulation of $S.$
\end{definition}

\begin{definition}
Given a discrete $S\subset \mathbb{R}^{2},$ $\xi >0$, and a segment $ab$
between two points $a$ and $b$ of $S$ we say that $ab$ is $\xi $--stable
provided the following holds. For any embedding%
\begin{eqnarray*}
\iota &:&S\longrightarrow \mathbb{R}^{2}\text{ so that} \\
\left\vert \iota -\mathrm{id}_{S}\right\vert &<&\xi ,
\end{eqnarray*}%
the segment $\iota \left( a\right) \iota \left( b\right) $ is in the
Delaunay graph of $\iota \left( S\right) .$
\end{definition}

\begin{definition}
We say that an $\left( \eta ,\delta \eta ,\varepsilon \right) $--CDG is $\xi
\eta $--stable provided each of its boundary edges $\mathcal{T}$ is $\xi
\eta $--stable$.$
\end{definition}

The definition of legal graph (Definition \ref{legal dfbn}) is motivated by
the following result.

\begin{proposition}
\label{Gromov Guarantee Lemma}(\cite{ProWilh2}, Proposition 3.3) For $\eta
>0,$ let $S$ be an $\eta $--Gromov subset of $\mathbb{R}^{2}$. If $a,b\in S$
satisfy dist$(a,b)<\sqrt{2}\eta ,$ then the segment $ab$ between $a$ and $b$
is in Delaunay graph of $S.$
\end{proposition}

We can now state the main result of the section, which is the following
extension theorem for almost $\mathrm{CDGs.}$ It has the added feature that
the extension is \textbf{error correcting }in the sense that the error
estimates, $\varepsilon $ and $\delta ,$ for the new simplices is $0$ (cf
Definition \ref{error correcting dfn}, below).

\begin{theorem}[Theorem 4.17, \protect\cite{ProWilh2}]
\label{real chew thm} There are $\varepsilon ,\xi ,\delta >0$ with the
following property. For $\beta \geq 50\eta ,$ let $\mathcal{T}$ be a $\xi
\eta $--stable, $\left( \eta ,\delta \eta ,\varepsilon \right) $--CDG in $%
\left[ 0,\beta \right] \times \left[ 0,\beta \right] \subset \mathbb{R}^{2}.$
There\ is a $\xi \eta $--stable, $\left( \eta ,\delta \eta ,\varepsilon
\right) $--CDG, $\mathcal{\tilde{T}},$\ which extends $\mathcal{T}$\ and has
the following properties.

\begin{enumerate}
\item $\left[ 0,\beta \right] \times \left[ 0,\beta \right] \subset B\left( 
\mathcal{\tilde{T}}_{0},\eta \right) .$

\item The buffer of $\mathcal{\tilde{T}}$ extends the buffer of $\mathcal{T}$
that is, 
\begin{equation*}
\left( \mathcal{T}_{\mathrm{x}}\right) _{0}\subset \left( \mathcal{\tilde{T}}%
_{\mathrm{x}}\right) _{0}.
\end{equation*}

\item $\left( \mathcal{\tilde{T}}_{\mathrm{x}}\right) _{0}$ is a subset of $%
\left[ 0,\beta +6\eta \right] \times \left[ 0,\beta +6\eta \right] $ so that%
\begin{equation*}
\mathrm{dist}\left( \left( \mathcal{\tilde{T}}_{\mathrm{x}}\right)
_{0}\setminus \left( \mathcal{T}_{\mathrm{x}}\right) _{0},\mathcal{T}%
_{0}\right) \geq 3\eta ,
\end{equation*}%
\begin{equation*}
\left[ 0,\beta +6\eta \right] \times \left[ 0,\beta +6\eta \right] \subset
B\left( \left( \mathcal{\tilde{T}}_{\mathrm{x}}\right) _{0},\eta \right) ,%
\text{ and}
\end{equation*}%
\begin{equation*}
\mathrm{dist}\left( v,w\right) \geq \eta ,
\end{equation*}

for $v\in \left( \mathcal{\tilde{T}}_{\mathrm{x}}\right) _{0}\setminus
\left( \mathcal{T}_{\mathrm{x}}\right) _{0}$ and $w\in \left( \mathcal{%
\tilde{T}}_{\mathrm{x}}\right) _{0}.$

\item $\mathcal{\tilde{T}}_{\mathrm{x}}$ contains all legal subgraphs of $%
\mathcal{T}_{\mathrm{x}}.$

\item Every edge of $\mathcal{\tilde{T}}$ that has a vertex in $\left( 
\mathcal{\tilde{T}}_{0}\right) _{\mathrm{x}}\setminus \mathcal{T}_{0}$ is
Delaunay, in the sense that it passes the $\varepsilon $--angle Flip Test
with $\varepsilon =0.$

\item Every edge with a vertex in $\left( \mathcal{\tilde{T}}_{0}\right) _{%
\mathrm{x}}\setminus \left( \mathcal{T}_{0}\right) _{\mathrm{x}}$ has length 
$\leq 2\eta .$
\end{enumerate}
\end{theorem}

\addtocounter{algorithm}{1}

\subsection{From Triangulations to Covers in Dimension 2\label%
{geomcotriansubsect}}

In this subsection, we define geometric triangulation covers and discuss
their key properties (see Figure \ref{tricov}). We assume $\left\{ \tilde{%
\eta}_{l}\right\} _{l=1}^{\tilde{L}}$ and $\mathcal{A}^{2}=\left\{ \digamma
_{l}\right\} _{l=1}^{\tilde{L}}$ are as in Proposition \ref{bottom prop}.

\begin{definition}[see Figure \protect\ref{cotri nat str interface figure}]
For $\digamma _{k}\in \mathcal{A}^{2}$ and $k\leq l,$ let $\mathcal{T}$ be
an $\left( \tilde{\eta}_{l},\delta \tilde{\eta}_{l}\right) $--CDG in $%
p_{\digamma _{k}}\left( \digamma _{k}^{+}\right) ,$ and let $\mathcal{B}_{l}$
and $\mathcal{L}_{l}$ be as in (\ref{B and L dfn}). We say that $\mathcal{T}$
is an $\tilde{\eta}_{l}$--\textbf{geometric triangulation} provided:

\begin{enumerate}
\item If $\digamma _{k}\in \mathcal{L}_{l},$ then $\mathcal{T}_{0}\cap
\left\{ \left\{ 0\right\} \times \left[ 0,\beta _{l}\right] \right\} $ is
the set $L_{k}^{l}$ from (\ref{left prop}).

\item If $\digamma _{k}\in \mathcal{B}_{l},$ then $\mathcal{T}_{0}\cap
\left\{ \left[ 0,\beta _{l}\right] \times \left\{ 0\right\} \right\} $ is
the set $B_{k}^{l}$ from Part 2 of Proposition \ref{bottom prop}.

\item The boundary edges of $\mathcal{T}$ are $\xi \tilde{\eta}_{l}$--stable
where $\xi $ is as in Theorem \ref{real chew thm}. 
\end{enumerate}
\end{definition}

\mbox{}\vspace{-35pt}%
\begin{wrapfigure}[6]{r}{0.33\textwidth}\centering
\includegraphics[scale=.15]{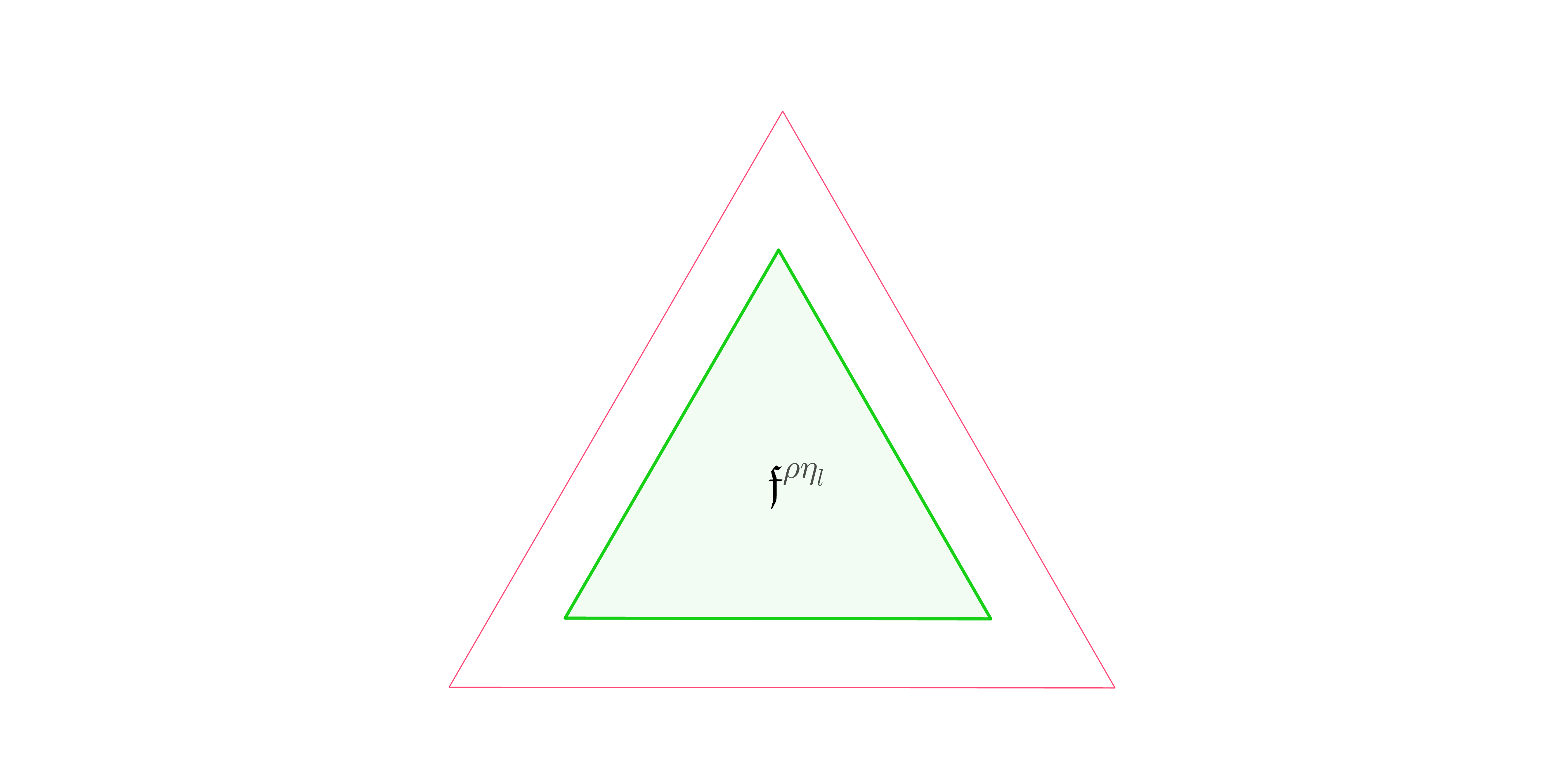}

\end{wrapfigure}

\begin{definition}
\label{tiles and bars dfn} Let $\mathcal{T}$ be an $\tilde{\eta}_{l}$%
--geometric triangulation. For $\mathfrak{e}\in \left( \mathcal{T}_{\mathrm{x%
}}\right) _{1}$ and $\rho >0$ set 
\begin{eqnarray*}
\hspace*{-4in}\mathfrak{e}^{\mathrm{Shr}} &:= &\mathfrak{e}\setminus \left( 
\overline{B(\partial \mathfrak{e},40\rho \tilde{\eta}_{l})}\right) \text{ and%
} \\
\mathfrak{e}^{\rho \tilde{\eta}_{l}} &:= &\left\{ \left. p+v\text{ }%
\right\vert \text{ }p\in \mathfrak{e}^{\mathrm{Shr}},\text{ }v\perp 
\mathfrak{e},\text{ and }\left\vert v\right\vert <2\rho \tilde{\eta}%
_{l}\right\} .
\end{eqnarray*}
For $\mathfrak{f}\in \left( \mathcal{T}_{\mathrm{x}}\right) _{2},$ set 
\begin{equation*}
\hspace{-1.18in} \mathfrak{f}^{\rho \tilde{\eta}_{l}}:= \left\{ x\in 
\mathfrak{f}\mid \mathrm{dist}(x,\partial \mathfrak{f})>\rho \tilde{\eta}%
_{l}\right\} \text{.}
\end{equation*}
\end{definition}

\begin{definition}[see Figure \protect\ref{tricov}]
\label{geom triag cover} Let $\mathcal{O}\left( \mathcal{T}\right) $ be a
triangulation cover of an $\tilde{\eta}_{l}$--geometric triangulation $%
\mathcal{T}.$ Given $\rho ,d>0$ we say that\textbf{\ }$\mathcal{O}\left( 
\mathcal{T}\right) $ is $\left( \tilde{\eta}_{l};\rho ,d\right) $--geometric
provided:

\begin{enumerate}
\item For all $\mathfrak{v\in }\mathcal{T}_{0},g_{\mathcal{O},\mathcal{T}%
}\left( \mathfrak{v}\right) =B\left( \mathfrak{v},50\rho \tilde{\eta}%
_{l}\right) .$

\item For all $\mathfrak{e\in }\mathcal{T}_{1},$ $g_{\mathcal{O},\mathcal{T}%
}\left( \mathfrak{e}\right) =\mathfrak{e}^{\rho \tilde{\eta}_{l}}.$

\item For all $\mathfrak{f\in }\mathcal{T}_{2},g_{\mathcal{O},\mathcal{T}%
}\left( \mathfrak{f}\right) =\mathfrak{f}^{\rho \tilde{\eta}_{l}}.$

\item Each of the three collections, 
\begin{equation*}
\left\{ g_{\mathcal{O},\mathcal{T}}\left( \mathfrak{v}\right) \right\} _{%
\mathfrak{v\in }\mathcal{T}_{0}},\text{ }\left\{ g_{\mathcal{O},\mathcal{T}%
}\left( \mathfrak{e}\right) \right\} _{\mathfrak{e\in }\mathcal{T}_{1}},%
\text{ and }\left\{ g_{\mathcal{O},\mathcal{T}}\left( \mathfrak{f}\right)
\right\} _{\mathfrak{f\in }\mathcal{T}_{2}}
\end{equation*}%
is $d\tilde{\eta}_{l}$--separated.
\end{enumerate}
\end{definition}

Combining Proposition \ref{alm CDG Chew} with this definition gives us

\begin{corollary}
\label{cover triang lemma}There exist universal constants $\rho ,d>0$ so
that for any $\tilde{\eta}_{l}$--geometric triangulation $\mathcal{T}$ there
is an $\left( \tilde{\eta}_{l};\rho ,d\right) $--geometric triangulation
cover $\mathcal{O}\left( \mathcal{T}\right) $ of $\mathcal{T}.$ Moreover, if 
$\mathcal{O}\left( \mathcal{\tilde{T}}\right) $ is any $\left( \tilde{\eta}%
_{l};\rho ,d\right) $--geometric triangulation cover of a subcomplex $%
\mathcal{\tilde{T}\subset T},$ then $\mathcal{O}\left( \mathcal{\tilde{T}}%
\right) $ extends to an $\left( \tilde{\eta}_{l};\rho ,d\right) $--geometric
triangulation cover of $\mathcal{T}.$
\end{corollary}

Since the constants $\rho $ and $d$ are universal, for brevity we will omit
them from the notation and refer to these covers as $\tilde{\eta}_{l}$%
--geometric triangulation covers. Next we record the following method of
getting a cotriangulation from a triangulation, which is canonical up to a
Gromov-Hausdorff error.

\begin{definition}[see Figures \protect\ref{cotri figure} and \protect\ref%
{cotri nat str interface figure}]
\label{stand local cotri dfn}

For $k\leq l,$ let $\mathcal{T}^{k}$ be an $\tilde{\eta}_{l}$--geometric
triangulation in $p_{\digamma _{k}}\left( \digamma _{k}^{+}\right) \ $with $%
\tilde{\eta}_{l}$--geometric triangulation cover $\mathcal{O}\left( \mathcal{%
T}\right) .$ A cosimplex $S$ that is subordinate to $\digamma ^{+}$ and
generated by a $U\in \mathcal{O}\left( \mathcal{T}\right) $ via (\ref{stand
co-simpliex dfn}) is called a \textbf{standard local }$\tilde{\eta}_{l}$--%
\textbf{cosimplex associated to }$\left( \digamma ,\mathcal{T}\right) .$ A
collection $\mathcal{S}=\mathcal{V}\cup \mathcal{E}\cup $ $\mathcal{F}$ of 
\textbf{standard local }$\tilde{\eta}_{l}$--\textbf{cosimplices }is called a 
\textbf{standard local }$\tilde{\eta}_{l}$--\textbf{cotriangulation},
provided the collection $\mathcal{V}$ of covertices is disjoint, the
collection $\mathcal{E}$ of coedges is disjoint, and the collection $%
\mathcal{F}$ of cofaces is disjoint. Otherwise, we call $\mathcal{S}$ a
standard \textbf{local precotriangulation.}
\end{definition}

\begin{definition}
\emph{(cf Definition \ref{alm CDG dfn})} Let $\mathcal{S}_{k}$ be a standard
local cotriangulation associated to\textbf{\ }$\left( \digamma _{k},\mathcal{%
T}^{k}\right) .$ An extension $\left( \mathcal{S}_{k}\right) _{\mathrm{x}}$
of $\mathcal{S}_{k}$ is called a \textbf{buffer }of $\mathcal{S}_{k}$\textbf{%
\ }provided $\left( \mathcal{S}_{k}\right) _{\mathrm{x}}$ is a standard
local cotriangulation associated to\textbf{\ }$\left( \digamma
_{k}^{+},\left( \mathcal{T}^{k}\right) _{\mathrm{x}}\right) $.
\end{definition}

The following corollaries follow by combining these definitions with the
fact that $\tilde{\eta}_{l}\geq \frac{\beta _{l}}{500},$ and our global
hypothesis that $\kappa $ and $\varkappa $ are sufficiently small.

\begin{corollary}
\label{cover a convex set cor}Let $\mathcal{T}$ be an $\tilde{\eta}_{l}$%
--geometric triangulation in $p_{\digamma _{l}}\left( \digamma
_{l}^{+}\right) $. Given any $K_{\mathcal{T}}\subset \mathcal{T}$, there is
a respectful, standard local $\tilde{\eta}_{l}$--cotriangulation $\mathcal{S}%
\left( \mathcal{T}\right) $ associated to $\left( \digamma _{l},\mathcal{T}%
\right) \ $whose generating set is $K_{\mathcal{T}}\ $and whose cosimplices
are defined using the submersion $p_{\digamma _{l}}.$ Moreover, setting $%
\left\vert K_{\mathcal{T}}\right\vert :=\cup _{\mathfrak{s}\in K_{\mathcal{T}%
}}\left\vert \mathfrak{s}\right\vert ,$ if 
\begin{equation*}
p_{\digamma _{l}}\left( \digamma _{l}\right) \subset \left\vert K_{\mathcal{T%
}}\right\vert ,
\end{equation*}%
then $\mathcal{S}\left( \mathcal{T}\right) $ covers $\digamma _{l}^{-}.$
\end{corollary}

\begin{corollary}
\label{dijs with natur cor}Let $\mathcal{S=}\mathcal{V}\cup \mathcal{E}\cup 
\mathcal{F}$ be a standard local $\tilde{\eta}_{l}$--precotriangulation
associated to $\left( \digamma _{l}\left( \beta _{l}\right) ,\mathcal{T}%
\right) .$ Then $\mathcal{S}$ is a cotriangulation, provided $\kappa $ and $%
\varkappa $ are sufficiently small. Moreover, if the cofaces $F\in \mathcal{F%
}$ that intersect $\cup \left( \mathcal{A}^{0}(\tilde{\mathcal{G}}^{l})\cup 
\mathcal{A}^{1}(\tilde{\mathcal{G}}^{l})\right) $ satisfy Part 4 of
Definition \ref{rel to dfn}, then $\mathcal{S}$ is relative to $\mathcal{LG}%
^{l}\cup \mathcal{B}\tilde{\mathcal{G}}^{l}.$
\end{corollary}

\addtocounter{algorithm}{1}

\subsection{Geometric Cotriangulations\label{Del Grm subsect}}

In this subsection, we define geometric cotriangulations and establish some
of their key properties. The most important is that we can use Lemmas \ref%
{const rank def} and \ref{line up respct lemma} to deform many of them to
respectful cotriangulations (see Lemma \ref{Curtiss Cor} below).

\begin{definition}[Geometric Cotriangulations]
\label{Dfn Grm Del}Let%
\begin{equation*}
\mathcal{S}=\mathcal{V}\cup \mathcal{E}\cup \mathcal{F}=\bigcup%
\limits_{i=1}^{l}\mathcal{S}_{i}=\bigcup\limits_{i=1}^{l}\left( \mathcal{V}%
_{i}\cup \mathcal{E}_{i}\cup \mathcal{F}_{i}\right)
\end{equation*}%
be a cotriangulation that is subordinate to $\{\digamma _{i}(\beta
_{i})\}_{i=1}^{l}.$ For $l\leq m,$ we say that $\mathcal{S}$ is $\tilde{\eta}%
_{m}$--\textbf{geometric} provided for $k\in \left\{ 1,\ldots ,l\right\} :$

\begin{enumerate}
\item There is a standard local cotriangulation whose Gromov-Hausdorff
distance from $\mathcal{S}_{k}$ is no more than $\tau \left( \kappa \tilde{%
\eta}_{m}\right) $.

\item Let the triangulation $\mathcal{T}^{k}$ that generates $\mathcal{S}%
_{k} $ have buffer $\mathcal{T}_{\mathrm{x}}^{k}$. Then $\mathcal{T}_{%
\mathrm{x}}^{k}$ generates a buffer $\left( \mathcal{S}_{k}\right) _{\mathrm{%
x}}$ of $\mathcal{S}_{k}$ which also has Gromov-Hausdorff distance at most $%
\tau \left( \kappa \tilde{\eta}_{m}\right) $--from a standard local
cotriangulation.
\end{enumerate}
\end{definition}

In the remainder of this subsection we let $\mathcal{S}=\mathcal{V}\cup 
\mathcal{E}\cup \mathcal{F}$ be an $\tilde{\eta}_{l}$--geometric
cotriangulation subordinate to $\{\digamma _{i}(\beta _{i})\}_{i=1}^{l}$ and
explain why, with a few additional hypotheses, $\mathcal{S}$ can be deformed
to be respectful. To do this, we will deform the submersions of neighboring
cosimplices so that they are lined up and then use the following corollary
to Lemma \ref{line up respct lemma} to obtain a respectful collection.

\begin{corollary}
\label{linup to respect}Suppose all pairs $\left( L,q_{L}\right) $, $\left(
H,q_{H}\right) $ of $l$-framed and $h$-framed sets with $l<h\ $and $L\cap
H\neq \emptyset $ in $\mathcal{S}$ are lined up on $L\cap \left\{ H\setminus
H^{-}\right\} .$ Then we can choose fiber exhaustion functions for the
covertices and coedges of $\mathcal{S}$ so that $\mathcal{S}$ is a
respectful framed collection.
\end{corollary}

We can now state the main result of this subsection, which provides a
criterion for when a geometric cotriangulation can be deformed to a
respectful geometric cotriangulation.

\begin{lemma}
\textbf{\label{Curtiss Cor}}There is a $C>1$ with the following property.
Let $\mathcal{S}=\cup_{i = 1}^l \mathcal{S}_i$ be an $\tilde{\eta}_{l}$%
--geometric cotriangulation that is subordinate to $\left\{ \digamma
_{i}^{+}\right\} _{i=1}^{l}$. Set 
\begin{equation*}
\mathcal{I}:=\left\{ \left. S\in \bigcup\limits_{i=1}^{l-1}\mathcal{S}_{i}%
\text{ }\right\vert \text{ }S\cap \left( \cup \mathcal{S}_{l}\right) \neq
\emptyset \right\} .
\end{equation*}%
Suppose that $\bigcup\limits_{i=1}^{l-1}\mathcal{S}_{i}$ is respectful and
every cosimplex of $\mathcal{S}_{l}$ is defined via the submersion $%
p_{\digamma _{l}}.$

If $\varkappa $ and $\kappa $ are sufficiently small, then there is a
respectful geometric cotriangulation $\mathcal{\tilde{S}}=\cup _{i=1}^{l}%
\tilde{\mathcal{S}}_{i}$ subordinate to $\left\{ \digamma _{i}(\beta
_{i}^{+})\right\} _{i=1}^{l}$ and a bijection 
\begin{eqnarray*}
\mathcal{S} &\rightarrow &\mathcal{\tilde{S}} \\
\left( S,q_{S}\right)  &\longmapsto &(\tilde{S},q_{\tilde{S}})
\end{eqnarray*}%
that have the following properties:

\begin{description}
\item[a] If $S\in \mathcal{S\setminus I},$ then $(\tilde{S},q_{\tilde{S}%
})=\left( S,q_{S}\right) .$

\item[b] For $S\in \mathcal{I},$ $q_{\tilde{S}}=q_{S}$ on $\left\{
S\setminus \cup \mathcal{S}_{l}\right\} $.

\item[c] For $S\in \mathcal{I}$ and any open $U\subset S,$ if $q_{S}|_{U}$
is $\varkappa $--almost Riemannian and $\kappa $--lined up with $p_{\digamma
_{l}},$ then $q_{\tilde{S}}|_{U}$ is $\tilde{\varkappa}$--almost Riemannian
and $\tilde{\kappa}$--lined with $p_{\digamma _{l}},$ where 
\begin{equation}
\tilde{\varkappa}:=\varkappa +\kappa \left( 1+2C\right) \text{ and }\tilde{%
\kappa}:=\kappa \left( 2+2C\right) .  \label{tilde kappa}
\end{equation}
\end{description}
\end{lemma}

Set 
\begin{equation*}
\mathcal{I}^{l}:=\left\{ \left. S\cap \digamma _{l}^{+}\text{ }\right\vert 
\text{ }S\in \mathcal{I}\right\} .
\end{equation*}

Since the cosimplices of $\mathcal{I}^{l}$ are lined up and subordinate to $%
\digamma _{l}^{+},$ it follows from Proposition \ref{glued submer prop} that
there is a submersion 
\begin{equation}
P:\cup \mathcal{I}^{l}\longrightarrow \mathbb{R}^{2}  \label{P}
\end{equation}%
that is lined up with each of the submersions $q_{S}$ that define the
cosimplices of $\mathcal{I}^{l}.$

The definition of geometric cotriangulation and the hypothesis that $%
\{\digamma _{i}(\beta _{i})\}_{i=1}^{l}$ is vertically separated imply that
there is a universal constant $C>1$ so that (\ref{Sub def hyp}) holds with $%
W=\cup \mathcal{I}^{l}\cap \left( \cup \mathcal{S}_{l}^{-}\right) ,$ $U=\cup 
\mathcal{I}^{l,+}\cap \left( \cup \mathcal{S}_{l}\right) ,$ and $\Omega
:=\cup \mathcal{I}^{l,+}.$ Applying the Submersion Deformation Lemma we then
get

\begin{lemma}
\label{52}If $\varkappa $ and $\kappa $ are sufficiently small then, there
is a Riemannian submersion $\tilde{P}:\cup \mathcal{I}^{l,+}\longrightarrow 
\mathbb{R}^{2}$ and a diffeomorphism $T_{l}:\mathbb{R}^{2}\longrightarrow 
\mathbb{R}^{2}$ so that 
\begin{equation}
\tilde{P}=\left\{ 
\begin{array}{cc}
T_{l}\circ p_{\digamma _{l}} & \text{on }\cup \mathcal{I}^{l}\cap \left(
\cup \mathcal{S}_{l}^{-}\right)  \\ 
P & \text{on }\cup \mathcal{I}^{l,+}\setminus \cup \mathcal{S}_{l}.%
\end{array}%
\right.   \label{85}
\end{equation}%
Moreover, for any open $U\subset \cup \mathcal{I}^{l,+},$ if $P|_{U}$ is $%
\varkappa $--almost Riemannian and $\kappa $--lined with $p_{\digamma _{l}},$
then $\tilde{P}|_{U}$ is $\tilde{\varkappa}$--almost Riemannian and $\tilde{%
\kappa}$--lined with $p_{\digamma _{l}},$ where $\tilde{\varkappa}$ and $%
\tilde{\kappa}$ are as in (\ref{tilde kappa}).
\end{lemma}

\begin{proof}[Proof of Lemma \protect\ref{Curtiss Cor}]
Given $S\in \mathcal{I},$ there is a diffeomorphism $T_{S}:\mathbb{R}%
^{2}\longrightarrow \mathbb{R}^{2}$ so that on $S\cap \digamma _{l}^{+},$ $%
q_{S}=T_{S}\circ P,$ where $P$ is as in (\ref{P}). Set 
\begin{equation}
\tilde{q}_{S}=\left\{ 
\begin{array}{cc}
T_{S}\circ \tilde{P} & \text{on }\digamma _{l}^{+}\cap S \\ 
q_{S} & \text{on }S\setminus \cup \mathcal{S}_{l},%
\end{array}%
\right.  \label{deformation to line up}
\end{equation}%
where $\tilde{P}$ is as in Lemma \ref{52}. Since the elements of $\mathcal{S}%
_{l}$ are undeformed, it follows from (\ref{85}) and (\ref{deformation to
line up}) that $\tilde{q}_{S}$ is lined up with each element of $\mathcal{S}%
_{l}$ that intersects $S.$ It also follows from (\ref{deformation to line up}%
) that the process of replacing each $q_{S}$ with $\tilde{q}_{S}$ does not
alter the fact that the collection $\bigcup\limits_{i=1}^{l-1}\mathcal{S}%
_{i} $ is lined up.

Conclusion a holds since we only deform the cosimplices of $\mathcal{I}.$
Conclusion $b$ holds since $\tilde{q}_{S}=q_{S}$ on $\left\{ S\setminus \cup 
\mathcal{S}_{l}\right\} .$ Conclusion $c$ follows from the \textquotedblleft
moreover\textquotedblright\ conclusion of Lemma \ref{52}.
\end{proof}

\section{How to Construct Cotriangulations\label{Sect Disresp cotri}}

In this section, we complete the proof of Covering Lemma \ref{cov em 2}. By
Proposition \ref{dda} it suffices to prove Lemma \ref{cotr exist lemma}.

Recall that $\mathcal{LG}^{l}$ and $\mathcal{B}\tilde{\mathcal{G}}^{l}$ are
the cographs from Corollary \ref{left side cor} and Equation (\ref{bottom
dfn}). By Proposition \ref{alm CDG Chew} a geometric cotriangulation that is
relative to $\mathcal{LG}^{\tilde{L}}\cup \mathcal{B}\tilde{\mathcal{G}}^{%
\tilde{L}}$ automatically satisfies (\ref{transverese}). Thus it suffices to
prove

\begin{lemma}
\label{ind st}For each $l=1,2,\ldots ,\tilde{L},$ there is a respectful $%
\tilde{\eta}_{l}$--geometric cotriangulation $\mathcal{S}^{l}=\mathcal{V}%
^{l}\cup \mathcal{E}^{l}\cup \mathcal{F}^{l}$ that is subordinate $\left\{
\digamma _{i}(\beta _{i}^{+})\right\} _{i=1}^{l}$\ and relative to $\mathcal{%
LG}^{l}\cup \mathcal{B}\tilde{\mathcal{G}}^{l}$.
\end{lemma}

The proof is by induction. As in Part \ref{Part 1}, the main ingredients of
the proof are extension and subdivision lemmas. They are stated below as
Lemmas \ref{rel cotri ext lemma} and \ref{subdivide co prop}.

The induction is anchored by taking an $\tilde{\eta}_{1}$--\textrm{CDG }that
covers $p_{\digamma _{1}}\left( \digamma _{1}\right) $ and contains $%
\mathcal{LG}^{1}$ and $\mathcal{B}\tilde{\mathcal{G}}^{1}$ and then applying
Corollary \ref{cover a convex set cor} to get an associated standard local $%
\tilde{\eta}_{1}$--cotriangulation $\mathcal{S}^{1}$ that covers $\digamma
_{1}^{-}.$

For an outline of the induction step, assume that $\mathcal{S}^{l-1}$ is a
respectful $\tilde{\eta}_{l-1}$--geometric cotriangulation subordinate to $%
\left\{ \digamma _{i}(\beta _{i}^{+})\right\} _{i=1}^{l-1}.$ In the event
that $\tilde{\eta}_{l}<\tilde{\eta}_{l-1},$ we use the Subdivision Lemma to
replace $\mathcal{S}^{l-1}$ with an $\tilde{\eta}_{l}$-geometric
cotriangulation $\tilde{\mathcal{S}}^{l-1}$ subordinate to $\left\{ \digamma
_{i}(\beta _{i}^{+})\right\} _{i=1}^{l-1}$ 
\begin{wrapfigure}[13]{r}{0.39\textwidth}\centering
\vspace{-4pt}
\includegraphics[scale=0.135]{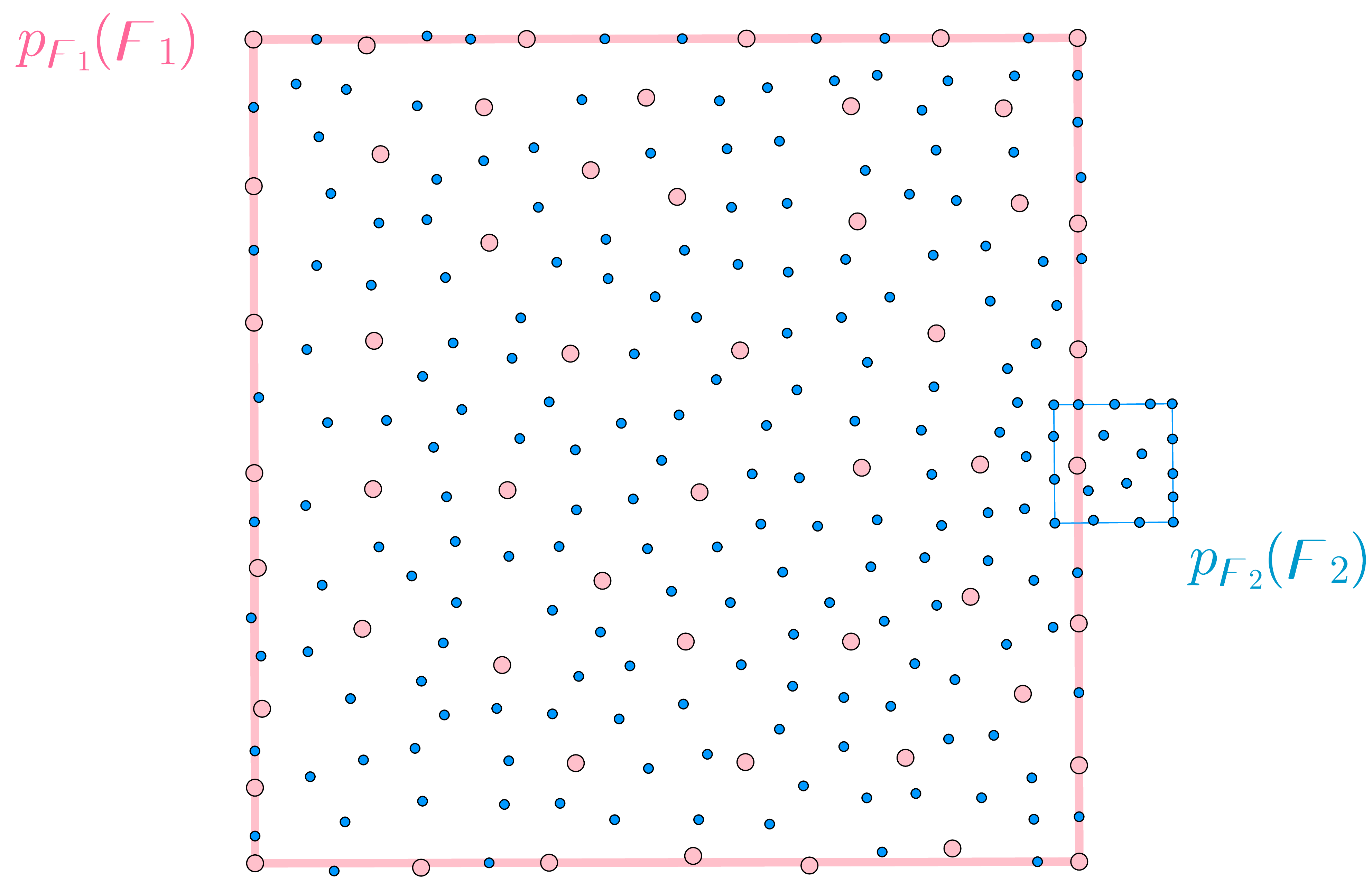}
\vspace{-5pt}
\caption{ {\protect\tiny {\ The pink and blue rectangles are the images of
two $2$--framed sets $\digamma _{1}(\protect\beta _{1})$ and $\digamma _{2}(%
\protect\beta _{1})$. The pink set is roughly a $\frac{\protect\beta _{1}}{8}
$--Gromov subset of $p_{\digamma _{1}}(\digamma _{1}(\protect\beta _{1}))$,
but only it intersects $p_{\digamma _{2}}(\digamma _{2}(\protect\beta _{2}))$
in one point. Thus we add the blue points to the pink points to obtain a
Gromov subset that is at an appropriate scale for the blue set. } }}
\label{etas fig}
\end{wrapfigure}
(see Figure \ref{etas fig}). We then use the Extension Lemma to extend $%
\tilde{\mathcal{S}}^{l-1}$ to an $\tilde{\eta}_{l}$-geometric
cotriangulation $\tilde{\mathcal{S}}^{l}$ subordinate to $\left\{ \digamma
_{i}(\beta _{i}^{+})\right\} _{i=1}^{l}$. If $\tilde{\eta}_{l}=\tilde{\eta}%
_{l-1},$ the procedure is the same except that we skip the subdivision step.

The first step in the proof of the Extension Lemma is to project $\mathcal{%
\tilde{S}}^{l-1}$ to $\digamma _{l}$ using $p_{\digamma _{l}}$. This yields
an $\left( \tilde{\eta}_{l},\delta \tilde{\eta}_{l},\varepsilon \right) $--$%
\mathrm{CDG,}$ $\mathcal{T}^{\mathrm{pr,}l},$ where the errors, $\delta $
and $\varepsilon $ depend on $\varkappa $ and $\kappa .$ Since the
cardinality $\tilde{L}$ of $\mathcal{A}^{2}$ has no known relationship with $%
\varkappa $ and $\kappa ,$ we need a way to see that these errors do no
accumulate to something significant. This is the role played by the as yet
undefined term, \textbf{error correcting cotriangulation.}\emph{\textbf{\ }}%
We define this term in the next subsection where we also state and prove the
Extension Lemma. In Subsection \ref{sub per cotri}, we state and prove the
Subdivision Lemma, and in Subsection \ref{constrcuting cotrig subsect}, we
complete the proof of Lemma \ref{ind st}.

\addtocounter{algorithm}{1}

\subsection{Error Correcting Cotriangulations and the Extension Lemma}

\label{error correcting and Ext Lem subsect}

In this subsection, we state and prove the Extension Lemma. Before doing
this we define the term \textbf{error correcting cotriangulation. }This also
requires several preliminary notions.

To start, notice that Theorem \ref{real chew thm} accepts a $\xi $--stable $%
\left( \eta ,\delta \eta ,\varepsilon \right) $--$\mathrm{CDG}$ as its input
and extends it to a $\xi $--stable $\left( \eta ,\delta \eta ,\varepsilon
\right) $--$\mathrm{CDG}$ for which the new simplices are error free. This
motivates

\begin{definition}
\label{error correcting dfn}If $\mathcal{T}$ and $\mathcal{\tilde{T}}$ are
related as in Theorem \ref{real chew thm}, then we will say that $\mathcal{%
\tilde{T}}$ is an \textbf{error correcting extension }of $\mathcal{T}.$
\end{definition}

\begin{wrapfigure}[18]{l}{0.37\textwidth}\centering
\vspace{-13pt}
\includegraphics[scale=0.11]{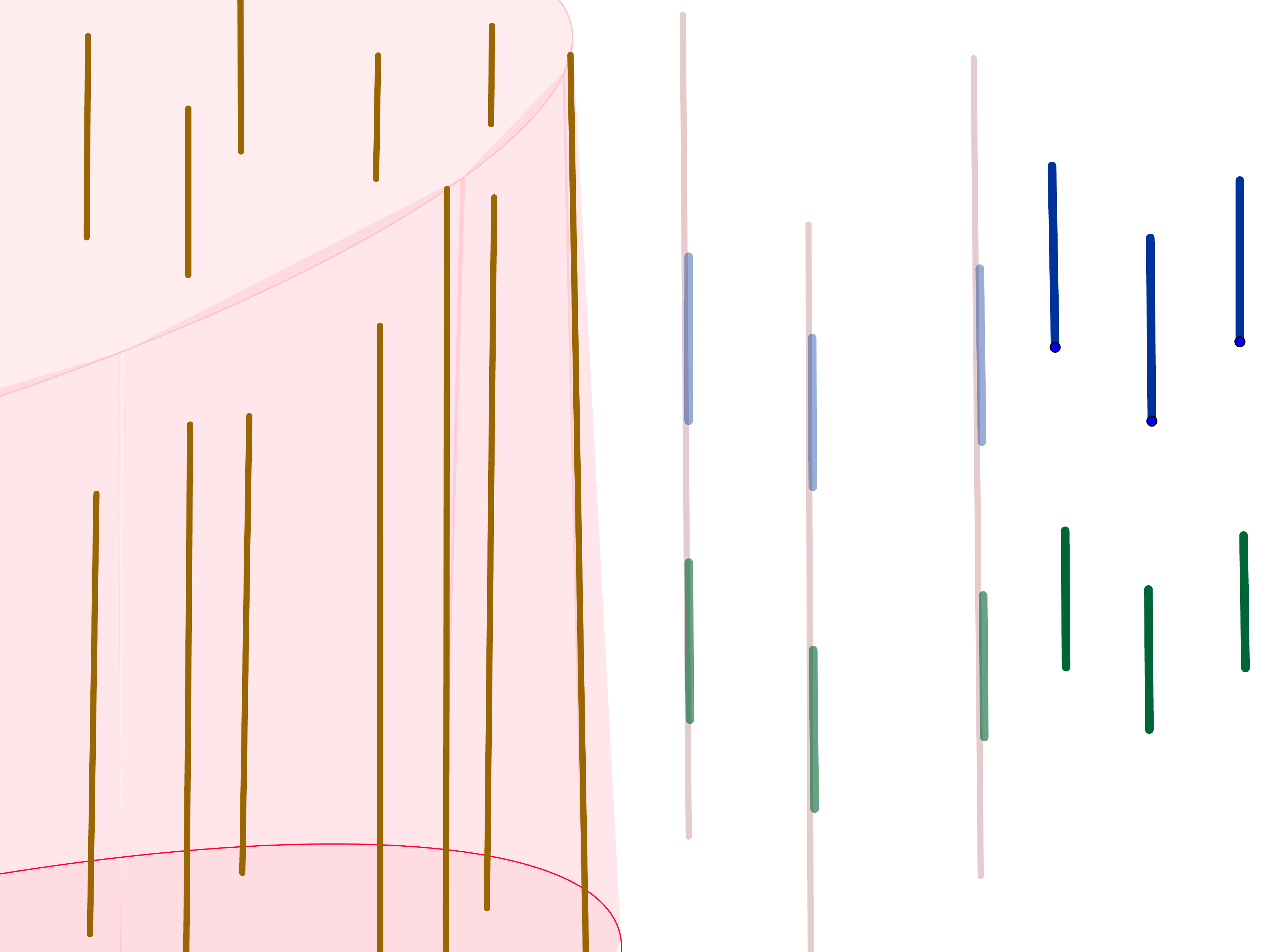}
\caption{{\protect\tiny {\ \textbf{The Extension Process:} The long brown
segments are vertex fibers of the existing cotriangulation. The short
segments are extensions of the brown vertex fibers to two different $2$%
--framed sets. The vertical splitting depicted here is caused by a
bifurcation of the $2$--strained points (cf Figure \protect\ref{bifurc}).
The long pink segments are part of the buffer of the existing, pink-brown,
cotriangulation. Its presence is necessitated by technicalities involved
with extending Delaunay triangulations of Gromov sets (cf Figure \protect\ref%
{sqrt3}). For a more detailed, blown up version of this region see Figure 
\protect\ref{buffer blow up}.}}}
\label{buffers}
\end{wrapfigure}
Before proceeding further we detail the process mentioned above wherein we
project a geometric cotriangulation to get an almost $\mathrm{CDG.}$ To do
so, for $k\geq l-1,$ let $\mathcal{S}=\cup _{i=1}^{l-1}\mathcal{S}_{i}$ be
an $\tilde{\eta}_{k}$--geometric cotriangulation that is subordinate to $%
\left\{ \digamma _{i}(\beta _{i}^{+})\right\} _{i=1}^{l-1}\subset \mathcal{A}%
^{2}$. For $\mathfrak{v}\in (K_{\mathcal{T}^{i}})_{0},$ let $q_{\mathfrak{v}%
} $ be the submersion that defines the covertex that corresponds to $%
\mathfrak{v}.$ We call the set of vertex fibers 
\begin{equation}
\hspace{3in} \mathrm{VF}:=\bigcup\limits_{i=1}^{l-1}\left\{ \left. q_{%
\mathfrak{v}}^{-1}\left( \mathfrak{v}\right) \text{ }\right\vert \mathfrak{v}%
\in (K_{\mathcal{T}^{i}})_{0}\right\} ,\text{\label{north star dfn}}
\end{equation}%
and define 
\begin{equation*}
\hspace{2.5in} \mathcal{T}_{0}^{\mathrm{pr,}l}:=p_{\digamma _{l}}\left( 
\mathrm{VF}\right) .
\end{equation*}

Strictly speaking, the elements of $\mathcal{T}_{0}^{\mathrm{pr},l}$ are not
points. So a triangulation of $\mathcal{T}_{0}^{\mathrm{pr},l}$ is,
technically, an abstract simplicial complex rather than a complex in $%
\mathbb{R}^{2}.$ However, since our $2$--framed sets are $\kappa $-lined up,
the elements of \newpage \noindent $\mathcal{T}_{0}^{\mathrm{pr},l}$ all
have diameter less than $\tau \left( 2\kappa \beta _{l}\right) .$ To avoid
excess technicalities we will treat the elements of $\mathcal{T}_{0}^{%
\mathrm{pr},l}$ as though they were points, with no further comment.

Let $\mathcal{T}^{\mathrm{pr,}l}$ be the triangulation of $\mathcal{T}_{0}^{%
\mathrm{pr},l}$ in which two vertices are connected by an edge if and only
if the corresponding covertices intersect a common coedge. Define buffers $%
\left( K_{\mathcal{T}^{l}}\right) _{\mathrm{x}},$ $\mathrm{VF}_{\mathrm{x}}$%
, and $\mathcal{T}_{\mathrm{x}}^{\mathrm{pr,}l}$ analogously. Before stating
the definition of error correcting cotriangulation, we need two more notions.

\begin{definition}
Let $S\in \mathcal{V\cup E}$ be generated by $\mathfrak{s}\in \mathcal{K}_{%
\mathcal{T}^{i}}\subset \mathcal{T}^{i}.$ If $\mathfrak{s}\in \partial 
\mathcal{T}^{i},$ then we say that $S$ is a boundary cosimplex of $\mathcal{S%
}$, and we let $\partial \mathcal{V},$ $\mathcal{\partial E},$ and $\partial 
\mathcal{S}$ be the collection of all boundary covertices, coedges, and
cosimplices, respectively.
\end{definition}

\begin{definition}
Let $S\in \mathcal{S}_{i}\subset \mathcal{S}.$ If $q_{S}$ is $\tilde{%
\varkappa}$--almost Riemannian and $\tilde{\kappa}$--lined up with $%
p_{\digamma _{i}},$ then we say that $S$ is $\left( \tilde{\varkappa},\tilde{%
\kappa}\right) $--deformed. For $U\subset S$, we say that $q_{S}$ is
undeformed on $U$ if $q_{S}|_{U}$ and $p_{\digamma _{i}}|_{U}$ are lined up.
If $q_{S}$ is undeformed on $S,$ then we say that $S$ is undeformed.
\end{definition}

\begin{definition}
\label{perfectt dfn}(\textbf{Error Correcting Cotriangulations) }For $k\geq
l,$ let 
\begin{equation*}
\mathcal{S}=\mathcal{V}\cup \mathcal{E}\cup \mathcal{F}=\bigcup%
\limits_{i=1}^{l}\mathcal{S}_{i}
\end{equation*}%
be a respectful $\tilde{\eta}_{k}$--geometric cotriangulation that is
subordinate to $\{\digamma _{i}(\beta _{i})\}_{i=1}^{l}.$ We say that $%
\mathcal{S}$ is \textbf{error correcting} provided for each $j_{0}\in $ $%
\left\{ 1,2,\ldots ,l\right\} ,$

\begin{enumerate}
\item The triangulation $\mathcal{T}^{j_{0}}$ associated to $\mathcal{S}%
_{j_{0}}$ is an error correcting extension of $\mathcal{T}^{\mathrm{pr,}%
j_{0}}$.

\item The generating set $K_{\mathcal{T}^{j_{0}}}$ for $\mathcal{S}_{j_{0}}$
is%
\begin{equation*}
K_{\mathcal{T}^{j_{0}}}=\mathcal{T}^{j_{0}}\setminus \mathcal{T}^{\mathrm{pr,%
}j_{0}}.
\end{equation*}

\item Every $S\in \mathcal{S}\setminus \partial \mathcal{S}$ is undeformed.

\item Every element of $E\in \partial \mathcal{E}$ is undeformed on $\left(
E\setminus \cup \mathfrak{F}\left( E\right) ^{+}\right) $, where $\mathfrak{F%
}\left( E\right) $ is the set of the cofaces $F\in \mathcal{F}$ that
intersect $E$.

\item For every covertex of $V\in \partial \left( \mathcal{V}\right) ,$ $%
q_{V}$ is undeformed on $V\setminus \left( \cup \mathfrak{E}\left( V\right)
^{+}\bigcup \cup \mathfrak{F}\left( V\right) ^{+}\right) $, where $\mathfrak{%
E}\left( V\right) $ and $\mathfrak{F}\left( V\right) $ are the sets of
coedges and cofaces that intersect $V.$

\item For every covertex of $V\in \partial \left( \mathcal{V}\right) ,$ $%
q_{V}$ is at most $\left( \varkappa _{2},\kappa _{2}\right) $--deformed
everywhere and at most $\left( \varkappa _{1},\kappa _{1}\right) $--deformed
on $V\setminus \left( \cup \mathfrak{F}\left( V\right) ^{+}\right) ,$ where 
\begin{eqnarray*}
\varkappa _{1} &:&=\varkappa +\kappa \left( 1+2C\right) \text{, }\kappa
_{1}:=\kappa \left( 2+2C\right) , \\
\varkappa _{2} &:&=\varkappa _{1}+\kappa _{1}\left( 1+2C\right) \text{, }%
\kappa _{2}:=\kappa _{1}\left( 2+2C\right)
\end{eqnarray*}%
and $C$ is the constant from Lemma \ref{Curtiss Cor}.
\end{enumerate}
\end{definition}

We can now state the Extension Lemma.

\begin{lemma}
\label{rel cotri ext lemma}(Extension Lemma) Let $\mathcal{S}=\mathcal{V}%
\cup \mathcal{E}\cup \mathcal{F}$ be an error correcting, $\tilde{\eta}_{l}$%
--geometric cotriangulation which is subordinate to $\left\{ \digamma
_{i}(\beta _{i}^{+})\right\} _{i=1}^{l-1}\ $and relative to $\mathcal{LG}%
^{l}\cup \mathcal{B}\tilde{\mathcal{G}}^{l}.$

Then $\mathcal{S}$ extends to an error correcting, $\tilde{\eta}_{l}$%
--geometric cotriangulation $\tilde{\mathcal{S}}=\tilde{\mathcal{V}}\cup 
\tilde{\mathcal{E}}\cup \tilde{\mathcal{F}}$ which is subordinate to $%
\left\{ \digamma _{i}(\beta _{i}^{+})\right\} _{i=1}^{l-1}\cup \left\{
\digamma _{l}(\beta _{l}^{+})\right\} $ and relative to $\mathcal{LG}%
^{l}\cup \mathcal{B}\tilde{\mathcal{G}}^{l}.$
\end{lemma}

To begin the proof, define 
\begin{equation*}
\mathrm{VF}_{i}\left( \nu \right) :=\left\{ \left. q_{\mathfrak{v}%
}^{-1}\left( \mathfrak{v}\right) \cap \digamma _{i}\left( \beta _{i},\nu
\right) \text{ }\right\vert \text{ }q_{\mathfrak{v}}^{-1}\left( \mathfrak{v}%
\right) \in \mathrm{VF}\right\} ,
\end{equation*}%
and let $\mathrm{VT}_{i}\left( \nu \right) $ be the simplicial complex whose
vertex set is $\mathrm{VF}_{i}\left( \nu \right) $ and whose edge set is
pairs of fibers in $\mathrm{VF}_{i}\left( \nu \right) $ whose corresponding
covertices intersect a common coedge. Define buffers $\mathrm{VF}_{i}\left(
\nu \right) _{\mathrm{x}}$ and $\mathrm{VT}_{i}\left( \nu \right) _{\mathrm{x%
}}$ analogously. The error correcting hypothesis implies that $\mathrm{VT}%
_{i}\left( \nu \right) $ is geometrically similar to a $\mathrm{CDG.}$ To
formulate this precisely, we propose

\begin{definition}
A simplicial complex is \textbf{metric} provided the vertex set $\mathcal{T}%
_{0}$ is a subset of a metric space $X.$
\end{definition}

\begin{definition}
We say that a $2$--dimensional simplicial complex is \textbf{surface--like}
provided each edge $e$ is contained in at most two $2$--simplices.
\end{definition}

\begin{definition}
Let $\mathcal{T}$ be $2$--dimensional simplicial complex which is metric and
surface--like$.$ Let $e$ be an edge of $\mathcal{T}$ which bounds two faces.
We say that $e$ passes the $\varepsilon $--\textbf{angle flip test} if and
only if the sum of the comparison angles opposite $e$ is $\leq \pi
+\varepsilon ,$ where our comparison triangles are in $\mathbb{R}^{2}.$
\end{definition}

\begin{definition}
Let $\mathcal{T}$ be $2$--dimensional simplicial complex which is metric and
surface--like$.$ We say that $\mathcal{T}$ is an \textbf{abstract }$\left(
\eta ,\delta \eta ,\varepsilon \right) $\textbf{--CDG} provided $\mathcal{T}%
_{0}$ is an $\left( \eta ,\delta \eta \right) $--Gromov subset of $X$ and
each edge passes the $\varepsilon $--angle flip test$.$
\end{definition}

Next we prove that for each $j_{0}\in \left\{ 1,\ldots ,l\right\} ,$ $%
\mathrm{VT}_{j_{0}}\left( \beta _{l}^{+}\right) $ is an abstract \textrm{CDG}
with respect to the following metric: For any $q_{\mathfrak{v}}^{-1}\left( 
\mathfrak{v}\right) ,q_{\mathfrak{\tilde{v}}}^{-1}\left( \mathfrak{\tilde{v}}%
\right) \in \mathrm{VT}_{j_{0}}\left( \beta _{l}^{+}\right) ,$ choose any
point in $x_{0}\in q_{\mathfrak{v}}^{-1}\left( \mathfrak{v}\right) \cap
\digamma _{j_{0}}\left( \beta _{j_{0}},\beta _{l}\right) $ and define 
\begin{equation}
\mathrm{dist}\left( q_{\mathfrak{v}}^{-1}\left( \mathfrak{v}\right) ,q_{%
\mathfrak{\tilde{v}}}^{-1}\left( \mathfrak{\tilde{v}}\right) \right) :=%
\mathrm{dist}\left( x_{0},q_{\mathfrak{\tilde{v}}}^{-1}\left( \mathfrak{%
\tilde{v}}\right) \cap \digamma _{j_{0}}\left( \beta _{j_{0}},5\beta
_{l}\right) \right) .  \label{the metric}
\end{equation}%
Of course this definition depends on the many possible choices of $x_{0},$
but since the $q_{\mathfrak{v}}$ are $3\kappa _{2}$--lined up $\varkappa
_{2} $--almost Riemannian submersions, the fibers are almost equidistant. So
regardless of these choices we have

\begin{lemma}
\label{perf implies abstr cor}If $\mathcal{S}$ is an error correcting $%
\tilde{\eta}_{l}$--geometric cotriangulation that is subordinate to $\left\{
\digamma _{i}^{+}\right\} _{i=1}^{l-1}$, then for all $j_{0}\in \left\{
1,2,\ldots ,l-1\right\} ,$ $\mathrm{VT}_{j_{0}}\left( \beta _{l}^{+}\right) $
is an abstract $\left( \tilde{\eta}_{l},\tilde{\eta}_{l}\delta ,\tau \left(
\delta \right) \right) $--\textrm{CDG}$,$ for $\delta =1-\varkappa
_{2}-5,000\kappa _{2}.$
\end{lemma}

\begin{remarknonum}
The quantity $\tau \left( \delta \right) $ of the previous lemma is
universal in the sense that it only depends on $\delta .$
\end{remarknonum}

\begin{proof}
It follows from the definition of error correcting cotriangulations and (\ref%
{fibers close}) that $\mathrm{VF}_{j_{0}}\left( \beta _{l}^{+}\right) $ is $%
\tilde{\eta}_{l}\left( 1+\varkappa _{2}\right) $ dense in $\digamma
_{j_{0}}\left( \beta _{j_{0}},\beta _{l}^{+}\right) .$ We will show that the
vertex fibers $q_{\mathfrak{v}_{i}}^{-1}\left( \mathfrak{v}_{i}\right) $, $%
q_{\mathfrak{v}_{j}}^{-1}\left( \mathfrak{v}_{j}\right) $ of the buffer $%
\left( \mathrm{VT}_{j_{0}}\left( \beta _{l}^{+}\right) \right) _{\mathrm{x}}$
of $\mathrm{VT}_{j_{0}}\left( \beta _{l}^{+}\right) $ are $\left( \tilde{\eta%
}_{l}\left( 1-\varkappa _{2}\right) -9\kappa _{2}\beta _{l}^{+}\right) $%
-separated, and if a pair is connected by an edge, then the distance between
them lies in the interval%
\begin{equation*}
\left[ \tilde{\eta}_{l}\left( 1-\varkappa _{2}\right) -9\kappa _{2}\beta
_{l}^{+},2\tilde{\eta}_{l}\left( 1+\varkappa _{2}\right) \right] .
\end{equation*}%
Since $\beta _{l}^{+}\leq 501\tilde{\eta}_{l},$ it will then follow that our
vertex fibers are $\delta \tilde{\eta}_{l}$-separated, where $\delta
=1-\varkappa _{2}-5,000\kappa _{2},$ and if a pair is connected by an edge,
then the distance between them lies in the interval%
\begin{equation*}
\left[ \tilde{\eta}_{l}\left( 1-\delta \right) ,2\tilde{\eta}_{l}\left(
1+\varkappa _{2}\right) \right] .
\end{equation*}

By combining the definition of error correcting cotriangulation with
Proposition \ref{alm equi dist fibers} we see that both estimates hold when $%
j_{0}=1.$ Assume that both statements hold for $j_{0}\leq k-1\leq l-2.$

Let $\mathcal{S}=\cup _{i=1}^{l-1}\mathcal{S}_{i}$ be the decomposition of $%
\mathcal{S}$ into local cotriangulations. Suppose 
\begin{equation}
\mathfrak{v}_{i}\in \left( K_{\mathcal{T}^{i}}\right) _{\mathrm{x}},%
\mathfrak{v}_{j}\in \left( K_{\mathcal{T}^{j}}\right) _{\mathrm{x}}\text{
and }i\leq j\leq k.  \label{typ vertices}
\end{equation}%
If $j\leq k-1,$ then by induction the required estimates hold.

Otherwise, $\mathfrak{v}_{j}\in \left( K_{\mathcal{T}^{k}}\right) _{\mathrm{x%
}}\setminus \left( \mathcal{T}_{\mathrm{x}}^{\mathrm{pr,}k}\right) _{0}$ (cf
Figure \ref{buffers})$.$ Set $\mathfrak{\tilde{v}}_{i}=p_{\digamma
_{k}}\left( q_{\mathfrak{v}_{i}}^{-1}(\mathfrak{v}_{i})\right) .$ The
definition of error correcting cotriangulation then gives us that $\mathrm{%
dist}\left( \mathfrak{\tilde{v}}_{i},\mathfrak{v}_{j}\right) \geq \tilde{\eta%
}_{l}.$ Together with (\ref{kappa wins}), this gives us that $q_{\mathfrak{v}%
_{i}}^{-1}\left( \mathfrak{v}_{i}\right) $ and $q_{\mathfrak{v}%
_{j}}^{-1}\left( \mathfrak{v}_{j}\right) $ are $\left( \tilde{\eta}%
_{l}\left( 1-\varkappa _{2}\right) -9\kappa _{2}\beta _{l}^{+}\right) $%
-separated with respect to the metric of (\ref{the metric}). If, in
addition, $q_{\mathfrak{v}_{i}}^{-1}\left( \mathfrak{v}_{i}\right) $ and $q_{%
\mathfrak{v}_{j}}^{-1}\left( \mathfrak{v}_{j}\right) $ are connected by an
edge, then $\mathfrak{\tilde{v}}_{i}$ and $\mathfrak{v}_{j}$ are connected
by an edge of $\mathcal{T}^{k},$ which by Part 6 of Theorem \ref{real chew
thm} has length $\leq 2\tilde{\eta}_{l}.$ So, as before, (\ref{fibers close}%
) gives us that $\mathrm{dist}\left( q_{\mathfrak{v}_{i}}^{-1}\left( 
\mathfrak{v}_{i}\right) ,q_{\mathfrak{v}_{j}}^{-1}\left( \mathfrak{v}%
_{j}\right) \right) \leq 2\tilde{\eta}_{l}\left( 1+\varkappa _{2}\right) .$

It remains to prove that the edges of $\mathrm{VT}_{j_{0}}\left( \beta
_{l}^{+}\right) $ pass the $\tau \left( \delta \right) $--comparison angle
flip test. The proof is also by induction and the case $j_{0}=1$ is a
consequence of Proposition \ref{alm equi dist fibers}. We take 
\begin{equation*}
\mathfrak{v}_{i}\in K_{\mathcal{T}^{i}},\mathfrak{v}_{j}\in K_{\mathcal{T}%
^{j}}\text{ with }i\leq j\leq k
\end{equation*}%
and assume that $q_{\mathfrak{v}_{i}}^{-1}\left( \mathfrak{v}_{i}\right) $
and $q_{\mathfrak{v}_{j}}^{-1}\left( \mathfrak{v}_{j}\right) $ are connected
by an edge of $\mathrm{VT}_{k}\left( \beta _{l}^{+}\right) .$ By induction
we may assume that $j=k$ and hence that $\mathfrak{v}_{j}\in K_{\mathcal{T}%
^{k}}:=\mathcal{T}^{k}\setminus \mathcal{T}^{\mathrm{pr,}k},$ so by Part 5
of Theorem \ref{real chew thm}, the edge $\mathfrak{\tilde{v}}_{i}\mathfrak{v%
}_{j}$ is Delaunay. This, together with Lemma \ref{Cor Q disj}, gives us
that the edge between $q_{\mathfrak{v}_{i}}^{-1}\left( \mathfrak{v}%
_{i}\right) $ and $q_{\mathfrak{v}_{j}}^{-1}\left( \mathfrak{v}_{j}\right) $
passes the $\tau \left( \delta \right) $--comparison angle flip test. So $%
\mathrm{VT}_{j_{0}}\left( \beta _{l}^{+}\right) $ is an abstract $\left( 
\tilde{\eta}_{l},\delta \tilde{\eta}_{l},\tau \left( \delta \right) \right) $%
--\textrm{CDG.}
\end{proof}

\begin{figure}[tbp]
\includegraphics[scale=0.16]{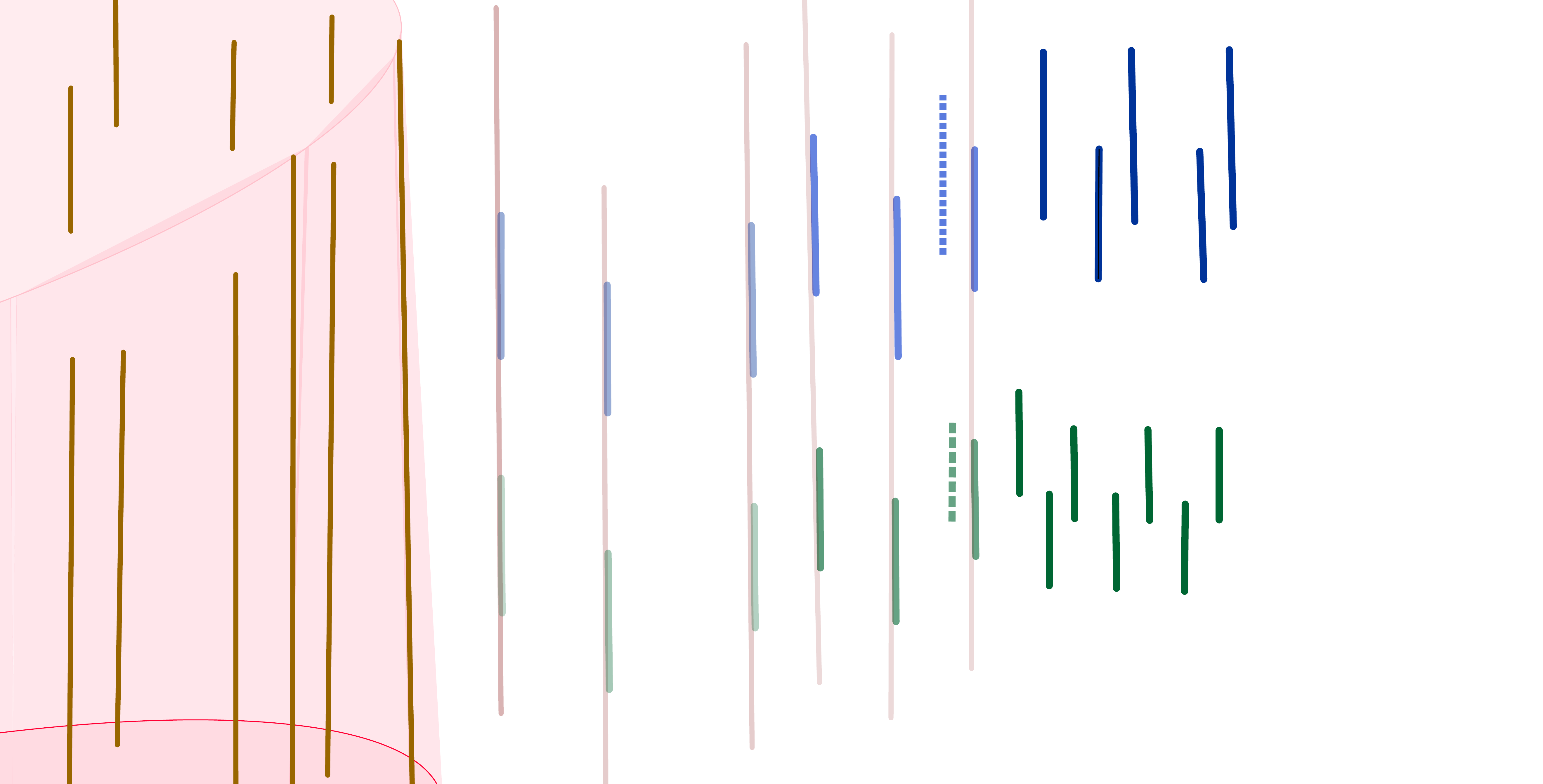}
\caption{{\protect\tiny {\ \textbf{Buffer Blown Up:} It follows from Part $3$
of Theorem \protect\ref{real chew thm} that no new vertex fibers are added
to the part of the buffer (light blue and light grey-green) that is within $3%
\tilde{\protect\eta}_{l}$ of the existing cotriangulation. However, the
buffer extends out to the $6\tilde{\protect\eta}_{l}$ neighborhood. New
vertex fibers (the broken segments) can be added that are further than $3%
\tilde{\protect\eta}_{l}$ from the old cotriangulation, but are still in the
buffer region. This would only happen if an old edge was originally almost
length $2\tilde{\protect\eta}_{l}$ and was lengthened by the lifting and
projecting process. The codege set can also change in the buffer region. }}}
\label{buffer blow up}
\end{figure}

Combining the previous lemma with Lemma \ref{Cor Q disj} yields

\begin{corollary}
\label{proj CDG prop}If $\mathcal{S}$ is an error correcting $\tilde{\eta}%
_{k}$--geometric cotriangulation that is subordinate to $\left\{ \digamma
_{i}^{+}\right\} _{i=1}^{l-1}$, then\textrm{\ }$\mathcal{T}^{\mathrm{pr,}l}$
is an $\left( \tilde{\eta}_{k},2\delta \tilde{\eta}_{k},2\tau \left( \delta
\right) \right) $--\textrm{CDG, }where $\delta $ and\textrm{\ }$\tau $ are
as in the previous lemma.
\end{corollary}

\begin{proof}[Proof of the Extension Lemma]
For simplicity, we write $\eta $ for $\tilde{\eta}_{l}.$ Since $\mathcal{S}=%
\mathcal{V}\cup \mathcal{E}\cup \mathcal{F}$ is error correcting, Corollary %
\ref{proj CDG prop} gives us that $\mathcal{T}^{\mathrm{pr,}l}$ is an $%
\left( \eta ,2\delta \eta ,2\tau \left( \delta \right) \right) $--\textrm{%
CDG, }where $\delta $ and $\tau $ are as in Corollary \ref{proj CDG prop}. 
\textrm{\ }By Theorem \ref{real chew thm}, $\mathcal{T}^{\mathrm{pr,}l}$ has
an error correcting extension to an $\left( \eta ,2\delta \eta ,2\tau \left(
\delta \right) \right) $--\textrm{CDG }$\mathcal{T}^{l}$ which covers $\left[
0,\beta _{l}\right] \times \left[ 0,\beta _{l}\right] $.

Set 
\begin{equation}
K_{\mathcal{T}^{l}}=\mathcal{T}^{l}\setminus \mathcal{T}^{\mathrm{pr,}l},
\label{genrating set}
\end{equation}%
and let $\mathcal{S}_{l}=\mathcal{V}_{l}\cup \mathcal{E}_{l}\cup \mathcal{F}%
_{l}$ be the standard local cotriangulation whose cosimplices are in
one-to-one correspondence with $K_{\mathcal{T}^{l}}$ and are defined using
the submersion $p_{\digamma _{l}}$ (see Definition \ref{stand local cotri
dfn}). Since $\mathcal{T}^{l}$ is an error correcting extension of $\mathcal{%
T}^{\mathrm{pr}}$, Lemma \ref{Cor Q disj} gives us that 
\begin{eqnarray}
&&\text{the intersection of }\mathcal{S}\mathfrak{\cup }\mathcal{S}_{l}\text{
with }\digamma _{l}^{+}\text{ is no more than }\tau \left( \kappa \eta
\right) \text{ from a standard}  \notag \\
&&\text{local cotriangulation associated to the }\eta \text{--geometric
triangulation }\mathcal{T}^{l}.  \label{is stand}
\end{eqnarray}

To show that 
\begin{equation*}
\widehat{\mathcal{S}}:=\mathcal{S}\cup \mathcal{S}_{l}=\mathcal{V}\cup 
\mathcal{E}\cup \mathcal{F}\cup \mathcal{V}_{l}\cup \mathcal{E}_{l}\cup 
\mathcal{F}_{l}
\end{equation*}%
is a cotriangulation, we combine the fact that $\mathcal{S}$ is a
cotriangulation with the following additional arguments.

\noindent \textbf{Proof of Property 1 of Definition \ref{Dfn cotri}. }Since $%
\mathcal{S}$ is a cotriangulation and $\mathcal{S}_{l}\ $is a local
cotriangulation, Property $1$ holds $\square $

\noindent \textbf{Proof of Property 2 of Definition \ref{Dfn cotri}. }The
collections $\mathcal{V},$ $\mathcal{E},$ and $\mathcal{F}$ are disjoint, by
hypothesis, so it suffices to show that the collections $\mathcal{V}_{l}%
\mathfrak{,}$ $\mathcal{E}_{l}\mathfrak{,}$ and $\mathcal{F}_{l}$ of the
newly defined covertices, coedges, and cofaces are disjoint from each other
and from $\mathcal{V},$ $\mathcal{E},$ and $\mathcal{F},$ respectively.
Since each element of $\mathcal{V}_{l}\mathfrak{,}$ $\mathcal{E}_{l}%
\mathfrak{,}$ and $\mathcal{F}_{l}$ is contained in $\digamma _{l}^{+},$ it
is enough to consider the intersection of $\mathcal{S}\mathfrak{\cup }%
\mathcal{S}_{l}$ with $\digamma _{l}^{+}.$

By combining (\ref{is stand}) with Corollary \ref{dijs with natur cor}, we
see that the intersections of each of $\mathcal{V}\mathfrak{\cup }\mathcal{V}%
_{l},$ $\mathcal{E}\mathfrak{\cup }\mathcal{E}_{l},$ and $\mathcal{F}%
\mathfrak{\cup }\mathcal{F}_{l}$ with $\digamma _{l}^{+}$ is a disjoint
collection, and hence each of $\mathcal{V}\mathfrak{\cup }\mathcal{V}_{l},$ $%
\mathcal{E}\mathfrak{\cup }\mathcal{E}_{l},$ and $\mathcal{F}\mathfrak{\cup }%
\mathcal{F}_{l}$ is a disjoint collection. $\square $

\noindent \textbf{Proof of Property 3 of Definition \ref{Dfn cotri}. }By
hypothesis, $\mathcal{S}$ covers $\cup _{i=1}^{l-1}\digamma _{i}(\beta
_{i}^{-}).$ Since $|\mathcal{T}^{l}|$ covers $p_{\digamma _{l}}\left(
\digamma _{l}(\beta _{l})\right) $, Corollary \ref{cover a convex set cor}, (%
\ref{genrating set}), and (\ref{is stand}) together imply that $\widehat{%
\mathcal{S}}$ covers $\digamma _{l}(\beta _{l}^{-})$. $\square $

Thus $\widehat{\mathcal{S}}$ is a cotriangulation. By construction, $%
\widehat{\mathcal{S}}$ is $\tilde{\eta}_{l}$--geometric. Our original
cotriangulation, $\mathcal{S},$ is relative to $\mathcal{LG}^{l}\cup 
\mathcal{B}\tilde{\mathcal{G}}^{l}$, so it contains the intersection of $%
\mathcal{LG}^{l}\cup \mathcal{B}\tilde{\mathcal{G}}^{l}$ with $\cup
_{i=1}^{l-1}\digamma _{i}(\beta _{i}^{+})$. Since Part 4 of Theorem \ref%
{real chew thm} guarantees that $\mathcal{T}^{l}$ contains all legal
subgraphs of $\mathcal{T}^{\mathrm{pr,}l},$ we can choose $\mathcal{T}^{l}$
to contain the graphs that generate $\mathcal{LG}^{l}\cup \mathcal{B}\tilde{%
\mathcal{G}}^{l}.$ This, together with (\ref{is stand}), Corollary \ref{dijs
with natur cor}, and the fact that the cosimplices of $\mathcal{S}_{l}$ are
defined using the submersion $p_{\digamma _{l}}$ gives us that $\widehat{%
\mathcal{S}}$ is relative to $\mathcal{LG}^{l}\cup \mathcal{B}\tilde{%
\mathcal{G}}^{l}.$

Our original cotriangulation, $\mathcal{S},$ is respectful and every
cosimplex of $\mathcal{S}_{l}$ is defined via the submersion $p_{\digamma
_{l}}.$ Together with the fact that $\mathcal{S}$ is $\eta $--geometric and
error correcting, this gives us that $\widehat{\mathcal{S}}$ satisfies the
hypotheses of Lemma \ref{Curtiss Cor}. So there is a deformation $\mathcal{%
\tilde{S}}$ of $\widehat{\mathcal{S}}$ which is a respectful $\eta $\textbf{%
--}geometric cotriangulation that is subordinate to $\left\{ \digamma
_{i}(\beta _{i}^{+})\right\} _{i=1}^{l}$.

Let $\mathcal{I}\subset $ $\widehat{\mathcal{S}}$ be as in Lemma \ref%
{Curtiss Cor}. Because $\mathcal{S}$ is error correcting, $\mathcal{I}%
\subset \partial \widehat{\mathcal{S}}.$ Combining this with the conclusion
of Lemma \ref{Curtiss Cor}, the fact that $\mathcal{S}$ is error correcting,
the construction of $\mathcal{\tilde{S}}$, and the fact that $\mathcal{%
\tilde{S}}$ is respectful and $\eta $\textbf{--}geometric, we see that $%
\mathcal{\tilde{S}}$ is error correcting.
\end{proof}

\addtocounter{algorithm}{1}

\subsection{Subdividing Error Correcting Cotriangulations\label{sub per
cotri}}

In this subsection, we explain how to subdivide error correcting
cotriangulations. We start with the analogous result for $\mathrm{CDGs}$
from \cite{ProWilh2}, which shows how to subdivide an $\left( \eta ,\delta
\eta ,\varepsilon \right) $--CDG in a manner that preserves the error
correcting property.

\begin{theorem}
\label{perfect ext thm}\emph{(Theorem 4.23, \cite{ProWilh2})} Let $\mathcal{T%
}$ be an $\left( \eta ,\delta \eta ,\varepsilon \right) $--CDG in $\left[
0,\beta \right] \times \left[ 0,\beta \right] \subset \mathbb{R}^{2}.$ If $%
\mathcal{\tilde{T}}$ is an error correcting extension of $\mathcal{T},$ then
there are subdivisions $\mathcal{L}\left( \mathcal{T}\right) ,$ $\mathcal{L}%
\left( \mathcal{\tilde{T}}\right) $ of $\mathcal{T}$ and $\mathcal{\tilde{T}}
$, respectively, so that

\begin{enumerate}
\item $\mathcal{L}\left( \mathcal{T}\right) $ and $\mathcal{L}\left( 
\mathcal{\tilde{T}}\right) $ are $\left( \frac{\eta }{10},\delta \frac{\eta 
}{10}\right) $--\textrm{CDGs}.

\item $\mathcal{L}\left( \mathcal{\tilde{T}}\right) $ is an error correcting
extension of $\mathcal{L}\left( \mathcal{T}\right) .$

\item If $\tilde{G}\subset \mathcal{\tilde{T}}_{1}$ is a legal subgraph of
the $1$--skeleton of $\mathcal{\tilde{T}}$ and $\mathcal{L}\left( \tilde{G}%
\right) $ is a legal $\left( \frac{\eta }{10},\delta \frac{\eta }{10}\right) 
$\textbf{--}geometric\textbf{\ }subdivision of $\tilde{G},$ then we can
choose $\mathcal{L}\left( \mathcal{\tilde{T}}\right) $ so that it contains $%
\mathcal{L}\left( \tilde{G}\right) .$
\end{enumerate}
\end{theorem}

Let 
\begin{equation*}
\mathcal{S}=\cup _{i=1}^{l-1}\mathcal{S}_{i}=\mathcal{V}\cup \mathcal{E}\cup 
\mathcal{F}
\end{equation*}%
be an error correcting, $\tilde{\eta}_{l-1}$--geometric cotriangulation
which is subordinate to $\left\{ \digamma _{i}(\beta _{i}^{+})\right\}
_{i=1}^{l-1}\ $ and relative to $\mathcal{LG}^{l-1}\cup \mathcal{B}\tilde{%
\mathcal{G}}^{l-1}$. We will explain how to subdivide $\mathcal{S}$ and
obtain an error correcting, $\tilde{\eta}_{l}$--cotriangulation which is
subordinate to $\left\{ \digamma _{i}(\beta _{i}^{+})\right\} _{i=1}^{l-1}\ $%
and relative to $\mathcal{LG}^{l}\cup \mathcal{B}\tilde{\mathcal{G}}^{l}.$
We make use of the following corollary of Proposition \ref{glued submer prop}%
.

\begin{corollary}
Let $S\in \mathcal{E}_{i}\cup \mathcal{F}_{i},$ and let $\mathcal{L}_{S}$ be
the collection of lower framed cosimplicies that intersect $S.$ There is a
submersion 
\begin{equation*}
q_{S}^{\mathrm{ext}}:\left( S\cup \cup \mathcal{L}_{S}\right) \cap \digamma
_{i}\longrightarrow \mathbb{R}^{2}
\end{equation*}%
so that

\begin{enumerate}
\item $q_{S}^{\mathrm{ext}}$ extends $q_{S}|_{S^{-}}$ and

\item For each $L\in \mathcal{L}_{S},$ $q_{S}^{\mathrm{ext}}$ is lined up
with $q_{L}$ on $L\setminus S.$
\end{enumerate}
\end{corollary}

\noindent \textbf{The Procedure:}

Since $\mathcal{S}$ is error correcting, for each $k,$ $\mathcal{T}^{k}$ is
an error correcting extension of $\mathcal{T}^{\mathrm{pr,}k}.$ Thus by
Theorem \ref{perfect ext thm}, there are subdivisions $\mathcal{L}\left( 
\mathcal{T}^{k}\right) $ and $\mathcal{L}\left( \mathcal{T}^{\mathrm{pr,}%
k}\right) $ of $\mathcal{T}^{k}$ and $\mathcal{T}^{\mathrm{pr,}k}$,
respectively, so that $\mathcal{L}\left( \mathcal{T}^{k}\right) $ is an
error correcting extension of $\mathcal{L}\left( \mathcal{T}^{\mathrm{pr,}%
k}\right) $ and each of $\mathcal{L}\left( \mathcal{T}^{k}\right) $ and $%
\mathcal{L}\left( \mathcal{T}^{\mathrm{pr,}k}\right) $ is $\frac{\tilde{\eta}%
_{l-1}}{10}$--geometric. Set%
\begin{equation*}
K_{\mathcal{L}\left( \mathcal{T}^{k}\right) }=\mathcal{L}\left( \mathcal{T}%
^{k}\right) \setminus \mathcal{L}\left( \mathcal{T}^{\mathrm{pr,}k}\right) .
\end{equation*}

Each $\mathfrak{s}\in K_{\mathcal{L}\left( \mathcal{T}^{k}\right) }$ is
contained in the interior of a unique $\mathfrak{t}\left( \mathfrak{s}%
\right) \in K_{\mathcal{T}^{k}}.$ Let $S\left( \mathfrak{s}\right) $ be the
standard geometric $\frac{\tilde{\eta}_{l-1}}{10}$--cosimplex that
corresponds to $\mathfrak{s}$ and is defined using the submersion $q_{%
\mathfrak{t}\left( \mathfrak{s}\right) }^{\mathrm{ext}}\ $(cf Definition \ref%
{stand local cotri dfn}).

\noindent \textbf{The Outcome: }

\begin{lemma}
\label{subdivide co prop}(Subdivision Lemma) Let $\mathcal{\tilde{S}=\tilde{V%
}}\cup \mathcal{\tilde{E}}\cup \mathcal{\tilde{F}}$ be the collection
obtained by performing the above procedure $n\left( l\right) -n\left(
l-1\right) $ times, where $n\left( l\right) $ is as in (\ref{n of l dfn}).
Then $\mathcal{\tilde{S}}$ is an error correcting, $\tilde{\eta}_{l}$%
--geometric cotriangulation that is relative to $\mathcal{LG}^{l}\cup 
\mathcal{B}\tilde{\mathcal{G}}^{l}$.
\end{lemma}

The definition of error correcting cotriangulation together with Lemma \ref%
{Cor Q disj} gives us the following result.

\begin{proposition}
\label{vert shrk prop}For $k\geq l-1$, let $\mathcal{S}$ be an error
correcting $\eta _{k}$--geometric cotriangulation that is subordinate to $%
\{\digamma _{i}(\beta _{i})\}_{i=1}^{l-1}.$ For $S,\hat{S}\in \mathcal{S},$
if $S\cap \hat{S}=\emptyset $, then $S^{v,+}\cap \hat{S}^{v,+}=\emptyset $.
(See (\ref{vert exp}) for the definition of $S^{v,+}$.)
\end{proposition}

\begin{proof}[Proof of Lemma \protect\ref{subdivide co prop}]
Notice that $\mathcal{\cup }_{S\in \mathcal{S}}S^{v,-}\mathcal{\subset }\cup 
\mathcal{\tilde{S}},$ so the fact that $\mathcal{S}$ covers $\{\digamma
_{i}(\beta _{i}^{-})\}_{i=1}^{l-1}$ gives us that $\mathcal{\tilde{S}}$
covers $\{\digamma _{i}(\beta _{i}^{-})\}_{i=1}^{l-1}.$

Write $\mathfrak{\tilde{S}}$ for 
\begin{equation*}
\mathcal{\tilde{V}},\text{ }\mathcal{\tilde{E}},\text{ or }\mathcal{\tilde{F}%
}.
\end{equation*}%
Notice that for $\tilde{S}\in \mathfrak{\tilde{S}},$ there is a unique $S\in 
\mathcal{S}$ so that $\left\vert \mathfrak{s}_{\tilde{S}}\right\vert \subset
\left\vert \mathrm{interior}\left( \mathfrak{s}_{S}\right) \right\vert ,$ 
\begin{equation}
q_{\tilde{S}}=q_{S}^{\mathrm{ext}},  \label{same subm}
\end{equation}%
and $\tilde{S}\subset S^{v,+}.$ For distinct $\tilde{S}_{1},\tilde{S}_{2}\in 
\mathfrak{\tilde{S}},$ if $S_{1}\cap S_{2}=\emptyset ,$ then by Proposition %
\ref{vert shrk prop}, $S_{1}^{v,+}\cap S_{2}^{v,+}=\emptyset .$ So 
\begin{equation}
\tilde{S}_{1}\cap \tilde{S}_{2}\subset S_{1}^{v,+}\cap S_{2}^{v,+}=\emptyset
.  \label{disj with vertc shrink}
\end{equation}%
Now suppose that for $\tilde{S}_{1},\tilde{S}_{2}\in \mathfrak{\tilde{S}}$, $%
S_{1}\cap S_{2}\neq \emptyset .$ Then $q_{S_{1}}^{\mathrm{ext}}$ and $%
q_{S_{2}}^{\mathrm{ext}}$ are lined up and $\mathfrak{s}_{\tilde{S}_{1}}\cap 
\mathfrak{s}_{\tilde{S}_{2}}=\emptyset .$ This together with $q_{\tilde{S}%
_{1}}=q_{S_{1}}^{\mathrm{ext}}$ and$\ q_{\tilde{S}_{2}}=q_{S_{2}}^{\mathrm{%
ext}}$ gives us that $\tilde{S}_{1}\cap \tilde{S}_{2}=\emptyset .$ Hence $%
\mathfrak{\tilde{S}}$ is pairwise disjoint.

Since we defined the cosimplices of $\mathcal{\tilde{S}}$ to be
geometrically standard using (\ref{same subm}),\textrm{\ }$\mathcal{\tilde{S}%
}$ is an $\tilde{\eta}_{l}$--geometric cotriangulation. Similarly, using the
equation $q_{\tilde{S}}=q_{S}^{\mathrm{ext}}$ and Lemma \ref{line up respct
lemma}, we can choose fiber exhaustion functions for the cosimplices of $%
\mathcal{\tilde{S}}$ that make $\mathcal{\tilde{S}}$ respectful. Since $%
\mathcal{S}$ is error correcting and $\mathcal{\tilde{S}}$ is respectful, it
follows from Theorem \ref{perfect ext thm} and the construction of $\mathcal{%
\tilde{S}}$ that $\mathcal{\tilde{S}}$ is error correcting. Our original
cotriangulation $\mathcal{S}$ is relative to $\mathcal{LG}^{l-1}\cup 
\mathcal{B}\tilde{\mathcal{G}}^{l-1}.$ It follows from Part 3 of Theorem \ref%
{perfect ext thm} that we can choose $\mathcal{\tilde{S}}$ to contain $%
\mathcal{LG}^{l}\cup \mathcal{B}\tilde{\mathcal{G}}^{l}.$ This, together
with the equation $q_{\tilde{S}}=q_{S}^{\mathrm{ext}}$ and the fact that $%
\mathcal{\tilde{S}}$ is $\tilde{\eta}_{l}$--geometric, gives us that $%
\mathcal{\tilde{S}}$ is relative to $\mathcal{LG}^{l}\cup \mathcal{B}\tilde{%
\mathcal{G}}^{l}.$
\end{proof}

\addtocounter{algorithm}{1}

\subsection{Constructing the cotriangulation\label{constrcuting cotrig
subsect}}

In this subsection, we show how to combine the extension and subdivision
lemmas to prove Lemma \ref{ind st} and hence Covering Lemma \ref{cov em 2}.

\begin{proof}[Proof of Lemma \protect\ref{ind st}]
Apply the Extension Lemma with $l=1$ and $\mathcal{S}=\emptyset $ to get an
error correcting, $\tilde{\eta}_{1}$--geometric cotriangulation $\mathcal{S}%
^{1}$ which is subordinate to $\digamma _{1}(\beta _{1}^{+})$ and relative
to $\mathcal{LG}^{1}\cup \mathcal{B}\tilde{\mathcal{G}}^{1}$.

Assume, by induction that there is an error correcting, $\tilde{\eta}_{l-1}$%
--geometric cotriangulation $\mathcal{S}^{l-1}$ which is subordinate to $%
\left\{ \digamma _{j}(\beta _{j}^{+})\right\} _{j=1}^{l-1}$ and relative to $%
\mathcal{LG}^{l-1}\cup \mathcal{B}\tilde{\mathcal{G}}^{l-1}$. Use the
Subdivision Lemma (\ref{subdivide co prop}) to subdivide $\mathcal{S}^{l-1}$
to get an error correcting, $\tilde{\eta}_{l}$--geometric cotriangulation $%
\mathcal{\tilde{S}}^{l-1}.$ By the Extension Lemma (\ref{rel cotri ext lemma}%
), $\mathcal{\tilde{S}}^{l-1}$ extends to a respectful, error correcting,
geometric cotriangulation $\mathcal{S}^{l}.$
\end{proof}

\section{Appendix: The Diffeomorphism of the Top Stratum}

Part 2 of Limit Lemma \ref{cover of X lemma} follows from the next result,
which we prove in this section. To do so we make liberal use of the notion
of stable submersions and their associated concepts as detailed in
Definitions \ref{stable subm dfn} and \ref{C1 close S-subms dfn}.

\begin{theorem}
\label{append thm}Let $\mathcal{B}$ be as in Covering Lemma \ref{cov em 2}.
There is an $\Omega \subset X_{\delta }^{4}$ so that $\mathcal{B}$ respects $%
\Omega $ and is $\Omega $--proper.
\end{theorem}

We construct $\Omega $ as the union of the following four open sets: 
\begin{eqnarray*}
O_{0} &:&=\bigcup\limits_{B\in \mathcal{B}^{0}}B\setminus
\bigcup\limits_{B\in \mathcal{B}^{0}\cup \mathcal{B}^{1}\cup \mathcal{B}^{2}}%
\bar{B}^{-}, \\
O_{1} &:&=\bigcup\limits_{B\in \mathcal{B}^{1}}B\setminus
\bigcup\limits_{B\in \mathcal{B}^{1}\cup \mathcal{B}^{2}}\bar{B}^{-},\text{ }
\\
O_{2} &:&=\bigcup\limits_{B\in \mathcal{B}2}B\setminus \bar{B}^{-},\text{ and%
} \\
O_{4} &:&=X\setminus \bigcup\limits_{B\in \mathcal{B}^{0}\cup \mathcal{B}%
^{1}\cup \mathcal{B}^{2}}\bar{B}^{-}.
\end{eqnarray*}%
It follows from the Quantitatively Proper Corollary (\ref{real grad pres}), (%
\ref{stop re-writing Curtis}), (\ref{stand co-simpliex dfn}), and Covering
Lemma \ref{cov em 2} that for $i=0,1,$ $2,$ or $4,$ $O_{i}\Subset X_{\delta
}^{4}.$ Hence for 
\begin{equation*}
\Omega :=O_{0}\cup O_{1}\cup O_{2}\cup O_{4},
\end{equation*}%
we have $\Omega \Subset X_{\delta }^{4}.$ It follows from the construction
of $\Omega $ that $\mathcal{B}$ is $\Omega $--proper. It remains to show

\begin{theorem}
There are open sets $\Omega _{\alpha }\subset M_{\alpha }$ that converge to $%
\Omega \subset X$, and for all but finitely many $\alpha ,\beta ,$
diffeomorphisms 
\begin{equation*}
\Phi _{\beta ,\alpha }:\Omega _{\alpha }\longrightarrow \Omega _{\beta },
\end{equation*}%
that have the following property.

For $k\in \left\{ 0,1,2\right\} ,$ all $B\in \mathcal{B}^{k},$ and all but
finitely many $\alpha ,\beta ,$ there are approximations $B_{\alpha
},B_{\beta }\rightarrow B$ so that 
\begin{equation}
\Phi _{\beta ,\alpha }\left( B_{\alpha }\cap \Omega _{\alpha }\right)
=B_{\beta }\cap \Omega _{\beta }.  \label{intersect preserv}
\end{equation}

Moreover, if $p:B\longrightarrow \mathbb{R}^{k}$ and $r:B\longrightarrow 
\mathbb{R}$ are the projection and fiber exhaustion function of $B$. Then
there are approximations $p_{\alpha }:B_{\alpha }\longrightarrow \mathbb{R}%
^{k}$ and $p_{\beta }:B_{\beta }\longrightarrow \mathbb{R}^{k}$ of $p$ and $%
r_{\alpha }:B_{\alpha }\longrightarrow \mathbb{R}$ and $r_{\beta }:B_{\beta
}\longrightarrow \mathbb{R}$ of $r$ so that 
\begin{equation}
p_{\beta }\circ \Phi _{\beta ,\alpha }|_{B_{\alpha }\cap \Omega _{\alpha
}}=p_{\alpha }|_{B_{\alpha }\cap \Omega _{\alpha }}\text{ and }r_{\beta
}\circ \Phi _{\beta ,\alpha }|_{B_{\alpha }\cap \Omega _{\alpha }}=r_{\alpha
}|_{B_{\alpha }\cap \Omega _{\alpha }}.  \label{Four str resp}
\end{equation}
\end{theorem}

The existence of the $\Omega _{\alpha }$ and the diffeomorphisms $\Phi
_{\beta ,\alpha }:\Omega _{\alpha }\longrightarrow \Omega _{\beta }$ follows
from Theorem 6.1 in \cite{KMS}. It only remains to show that the $\Phi
_{\beta ,\alpha }$ can be chosen to satisfy (\ref{intersect preserv}) and (%
\ref{Four str resp}). The argument is nearly identical to the proof of
Proposition 6.4 in \cite{ProWilh1}, so we will only outline it.

A key ingredient is Gromov's Covering Lemma (\cite{ProWilh1}, \cite{Zhu}),
whose statement makes use of the notion of the first order of a cover, which
we review here.

\begin{definitionnonum}
We say that a collection of sets $\mathcal{C}$ has first order $\leq 
\mathfrak{o}$ if and only if each $C\in \mathcal{C}$ intersects no more than 
$\mathfrak{o}-1$ other members of $\mathcal{C}$.
\end{definitionnonum}

\begin{lemma}
\label{Gromov Pack}(Gromov's Covering Lemma) Let $X$ be an $n$--dimensional
Alexandrov space with curvature $\geq k$ for some $k\in \mathbb{R}.$ There
are positive constants $\mathfrak{o}\left( n,k\right) $ and $r_{0}\left(
n,k\right) $ that satisfy the following property.

For all $r\in \left( 0,r_{0}\right) ,$ any compact subset of $A\subset X$
contains a finite subset $\left\{ a_{i}\right\} _{i\in I}$ so that

\begin{itemize}
\item $A\subset \cup _{i}B\left( a_{i},r\right) ,$ and

\item the first order of the cover $\left\{ B\left( a_{i},3r\right) \right\}
_{i}$ is $\leq \mathfrak{o}.$ {\huge \ }
\end{itemize}
\end{lemma}

In the Riemannian case, this follows from relative volume comparison, so one
only needs the corresponding lower bound on Ricci curvature. Since relative
volume comparison holds for rough volume in Alexandrov spaces (\cite{BGP}),
the proof in \cite{Zhu} yields, with minor modifications, Lemma \ref{Gromov
Pack}.

Since $\Omega $ is precompact and every point of $\bar{\Omega}$ is $\left(
4,\delta \right) $--strained, there is an $r>0$ and a finite open cover $%
\widetilde{\mathcal{U}}^{4}=\left\{ \tilde{U}\right\} $ of $\Omega $ by $%
\left( 4,\delta ,r\right) $--strained open sets $\tilde{U}.$ Since $\mathcal{%
B}^{0},$ $\mathcal{B}^{1}$ and $\mathcal{B}^{2}$ each consist of sets whose
closures are pairwise disjoint, we choose the sets of $\widetilde{\mathcal{U}%
}^{4}$ so that

\begin{itemize}
\item For $k\in \left\{ 0,1,2\right\} ,$ if $U\in $ $\widetilde{\mathcal{U}}%
^{4}$, then there is at most one $B\in \mathcal{B}^{k}$ so that $U\cap
B^{-}\neq \emptyset ,$ and in this event, 
\begin{equation}
U\subset B.  \label{then contains}
\end{equation}

\item For $k\in \left\{ 0,1,2\right\} ,$ if $U\in $ $\widetilde{\mathcal{U}}%
^{4},$ $\left( B,p_{B},r_{B}\right) \in \mathcal{B}^{k},$ and $U\cap
B^{-}\neq \emptyset ,$ let $\left\{ \left( a_{i},b_{i}\right) \right\}
_{i=1}^{k+1}$ be the first $\left( k+1\right) $-strainers of $U,$ and define 
$p_{U}:U\longrightarrow \mathbb{R}^{k+1}$ by $p_{U}\left( x\right) =\left( 
\mathrm{dist}_{a_{1}}\left( x\right) ,\ldots ,\mathrm{dist}_{a_{k+1}}\left(
x\right) \right) .$ Then using (\ref{k+1 submer dfn}), we choose $\left\{
\left( a_{i},b_{i}\right) \right\} _{i=1}^{k+1}$ so that 
\begin{equation}
\left\vert dp_{u}-d\left( p_{B},r_{B}\right) \right\vert <\tau \left( \delta
\right) .  \label{ext str}
\end{equation}
\end{itemize}

Let $3\rho $ be a Lebesgue number for $\widetilde{\mathcal{U}}^{4}$. Apply
Lemma \ref{Gromov Pack} to get a refinement $\mathcal{U}^{4}\left( 3\right) $
of $\widetilde{\mathcal{U}}^{4}$ that consists of $3\rho $--balls and
satisfies

\begin{itemize}
\item the corresponding collection $\mathcal{U}^{4}\left( 1\right) $ of $%
\rho $--balls is also a cover of $\Omega ,$ and

\item the first order of $\mathcal{U}^{4}\left( 3\right) $ is $\leq 
\mathfrak{o},$ where $\mathfrak{o}$ is as in Lemma \ref{Gromov Pack}.
\end{itemize}

Next we exploit (\ref{ext str}) to get local embeddings of each element $%
3U\in \mathcal{U}\left( 3\right) $ that respect the projections of our
handles. The result is stated in terms of the following definition.

\begin{definitionnonum}
Given $l,k\in \mathbb{N}$ with $k\leq l,$ let $\pi _{k,l}:$ $\mathbb{R}%
^{l}\longrightarrow \mathbb{R}^{k}$ be orthogonal projection. Given two
stable submersions 
\begin{equation*}
q:U\longrightarrow \mathbb{R}^{l}\text{ and }p:U\longrightarrow \mathbb{R}%
^{k}
\end{equation*}%
we write $p\overset{s}{=}\pi _{k,l}\circ q$ if and only if for all but
finitely many $\alpha ,$ $q_{\alpha }$ and $p_{\alpha }$ can be chosen so
that%
\begin{equation*}
p_{\alpha }=\pi _{k,l}\circ q_{\alpha }.
\end{equation*}
\end{definitionnonum}

\begin{proposition}
For each $3U\in \mathcal{U}^{4}\left( 3\right) ,$ there is a stable $\tau
\left( \delta \right) $--almost Riemannian embedding 
\begin{equation}
\mu :3U\longrightarrow \mathbb{R}^{4}  \label{embeeding}
\end{equation}%
with the following property. For $k\in \left\{ 0,1,2\right\} ,$ if $3U\cap
B^{-}\neq \emptyset $ for some $\left( B,p_{B},r_{B}\right) \in \mathcal{B}%
^{k}$, then%
\begin{equation}
3U\subset B,\text{ and\label{contains 2}}
\end{equation}%
\begin{equation}
\pi _{k+1}\circ \mu \overset{s}{=}\left( p_{B},r_{B}\right) ,
\label{respt top}
\end{equation}%
where $p_{B}:B\longrightarrow \mathbb{R}^{k}$ is the projection of $B,$ $%
r_{B}:B\longrightarrow \mathbb{R}$ is the fiber exhaustion function of $B$
and $\pi _{k+1}:\mathbb{R}^{4}\longrightarrow \mathbb{R}^{k+1}$ is
projection to the first $k+1$ factors.
\end{proposition}

\begin{proof}[Outline of Proof]
The $k=4$ version of Lemma \ref{framedsetprop} gives us the desired stable
embeddings. The fact that the embeddings can be chosen to satisfy (\ref%
{respt top}) follows from (\ref{ext str}) and the Submersion Deformation
Lemma (\ref{const rank def}).
\end{proof}

Let $\mathcal{U}^{4,\alpha }\left( 3\right) =\left\{ 3U_{i}^{\alpha
}\right\} _{i=1}^{L}$ and $\mathcal{U}^{4,\beta }\left( 3\right) =\left\{
3U_{i}^{\beta }\right\} _{i=1}^{L}$ be collections of metric balls in $%
M_{\alpha }$ and $M_{\beta }$ that approximate $\mathcal{U}^{4}\left(
3\right) $. Let 
\begin{equation*}
\mu _{i}^{\alpha }:3U_{i}^{\alpha }\longrightarrow \mathbb{R}^{4}
\end{equation*}%
be an embedding that approximates the embedding 
\begin{equation*}
\mu _{i}:3U_{i}\longrightarrow \mathbb{R}^{4}
\end{equation*}%
from (\ref{embeeding}). For $j\in \left\{ 1,2,\ldots ,L\right\} ,$ we get
embeddings 
\begin{equation*}
\nu _{j}:3U_{j}^{\alpha }\longrightarrow 3U_{j}^{\beta }
\end{equation*}%
by setting 
\begin{equation*}
\nu _{j}=\left( \mu _{j}^{\beta }\right) ^{-1}\circ \mu _{j}^{\alpha }.
\end{equation*}

It remains to glue the $\nu _{j}$ together in a manner that respects the
handles. This is accomplished via Theorem 5.3 of \cite{ProWilh1}. The idea,
which also appears in the proof of Perelman's Stability Theorem (\cite{Kap2}%
, \cite{Perel1}), is to use the Submersion Deformation Lemma (\ref{const
rank def}) to inductively glue the $\nu _{j}$. Since the first order of our
cover is bounded from above by the a priori constant $\mathfrak{o},$ Part 2
of the Submersion Deformation Lemma (\ref{const rank def}) ensures that our
glued embeddings stay close to all of our original embeddings. Since our
original embeddings respect $\mathcal{B}$, Part 4 of Lemma \ref{const rank
def} ensures that the glued embeddings respect $\mathcal{B}$.

\end{document}